\documentclass[dvips,12pt]{article}
\usepackage{amsmath,amssymb,amsfonts,amscd,rotating,pstricks,pst-node}
\def\simarrow{\mathrel{\raise -0.5mm\hbox{$\sim$}}\hspace{-1.8mm}{\rightarrow} } 

\def\bsimarrow{\leftarrow\hspace{-0.7mm}\mathrel{\raise -0.5mm\hbox{$\backsim$}} }

\def\bt{\begin{tabular}}
\def\te{\end{tabular}}

\def\lettrine#1#2#3{\noindent\hangindent#1\hangafter-#2
\hskip-#1\smash{\hbox to #1{#3\hfill}}\ignorespaces}

\newcommand{\To}[1]{\mathop{\to}\limits_{#1}}

\def\BM{\begin{pmatrix}}
\def\EM{\end{pmatrix}}

\def\pt{\mathrel{\scriptstyle \bullet}}

\def\txt{\textstyle}
\def\ds{\displaystyle}

\def\d=f{\buildrel\hbox{\scriptsize d\'{e}f}\over \Longleftrightarrow}

\def\cit{\text{\it I\hskip -6ptC\/}}

\def\square{\hfill\hbox{\vrule height .9ex width .8ex depth -.1ex}}
\def\rit{\text{\it I\hskip -2pt  R}}

\def\nit{\text{\it I\hskip -2pt  N}}

\def\rl {\rit^{\hskip 1pt\ell}}

\def\Bd{{\text B}}

\def\Ed{{\text E}}

\def\Fs{{\cal F}}

\def\Hs{{\cal H}}

\def\Is{{\cal I}}

\def\Ns{{\cal N}}

\def\Os{{\cal O}}

\def\Us{{\cal U}}

\def\ung{\hbox{1\hskip -4.2pt \rm 1}}

\def\be{\begin{equation}}
\def\ee{\end{equation}}
\def\beqn{\begin{eqnarray}}
\def\eeqn{\end{eqnarray}}
\def\nobeqn{\begin{eqnarray*}}
\def\noeeqn{\end{eqnarray*}}
\def\ba{\left(\begin{array}}
\def\ea{\end{array} \right) }

\def\bpr{\paragraph{Proof.}}

\def\epr{\square\vskip 6pt}
\def\eop{\hbox{\vrule height .9ex width .8ex depth -.1ex}}

\def\u{\underline}
\def\o{\overline}
\def\and{\; \mbox{and} \;}

\newcommand{\half}{\frac{1}{2}}

\def\hfl#1#2{\smash{\mathop{\hbox to 12mm{\rightarrowfill}}
\limits^{\scriptstyle #1}_{\scriptstyle #2}}}

\def\Ker{\mathop{\rm Ker}\nolimits}

\def\mod{\mathop{\rm mod}\nolimits}

\def\Be{\begin{enumerate}}
\def\Ee{\end{enumerate}}

\def\Bena{\begin{enumerate}
\def\labelenumi{\theenumi)}
\def\theenumi{\arabic{enumi}}
\def\labelenumii{\theenumii)}
\def\theenumii{\alph{enumii}}}

\def\Bean{\begin{enumerate}
\def\labelenumii{\theenumii)}
\def\theenumii{\arabic{enumii}}
\def\labelenumi{\theenumi)}
\def\theenumi{\alph{enumi}}}

\def\Bero{\begin{enumerate}
\def\labelenumii{\theenumii)}
\def\theenumii{\arabic{enumii}}
\def\labelenumi{(\theenumi)}
\def\theenumi{\roman{enumi}}}

\def\BeRo{\begin{enumerate}
\def\labelenumii{\theenumii)}
\def\theenumii{\arabic{enumii}}
\def\labelenumi{(\theenumi)}
\def\theenumi{\Roman{enumi}}}

\def\Bi{\vskip 11pt\begin{itemize}\itemsep=18pt}
\def\Ei{\end{itemize}}

\def\Bd{\begin{description}}
\def\Ed{\end{description}}

\def\R{\right}
\def\L{\left}
\def\F{\frac}

\usepackage[latin1]{inputenc}

\def\bigoplus{\mathop{\oplus}\limits}
\def\bigotimes{\mathop{\otimes}\limits}

\def\prod{\mathop{\Pi}\limits}
\def\sum{\mathop{\Sigma}\limits}

\def\bbf{\boldmath\bf}
\def\o{\overline}
\def\wt{\widetilde}

\def\Bi{\begin{itemize}}
\def\Ei{\end{itemize}}

\newcommand{\Aa}{\mathbb{A}\,}
\newcommand{\NN}{\mathbb{N}\,}
\newcommand{\ZZ}{\mathbb{Z}\,}
\def\hZ{\hbox{$Z$\hskip -10 pt\raisebox{-3pt}{$\bar{\hspace{4mm}}$}}}
\def\pZ{{\rm{Z}}\,}

\newcommand{\CC}{\mathbb{C}\,}
\newcommand{\RR}{\mathbb{R}\,}
\newcommand{\QQ}{\mathbb{Q}\,}
\newcommand{\FF}{\mathbb{F}\,}

\newcommand{\WW}{\mathbb{W}\,}
\renewcommand{\SS}{\mathbb{S}\,}

\newcommand{\HH}{{\rm{H}}\,}

\def\det{\operatorname{det}}
\def\Gal{\operatorname{Gal}}
\def\Frob{\operatorname{Frob}}
\def\End{\operatorname{End}}
\def\mod{\operatorname{mod}}
\def\hol{\operatorname{hol}}

\def\Mot{\operatorname{Mot}}

\def\Hom{\operatorname{Hom}}
\def\Corr{\operatorname{Corr}}

\def\cusp{\operatorname{cusp}}
\def\Irr{\operatorname{Irr}}
\def\Lie{\operatorname{Lie}}

\def\Rep{\operatorname{Rep}}
\def\Reps{\operatorname{Repsp}}
\def\FReps{\operatorname{FRepsp}}
\def\Repsp{\operatorname{Repsp}}
\def\FREPSP{\operatorname{FREPSP}}

\def\Red{\operatorname{Red}}

\def\Int{\operatorname{Int}}
\def\Aut{\operatorname{Aut}}
\def\out{\operatorname{Out}}
\def\Eis{\operatorname{EIS}}
\def\ELLIP{\operatorname{ELLIP}}
\def\CY {\operatorname{CY}}

\def\CEL {\operatorname{CEL}}
\def\CH{\operatorname{CH}}

\def\GL{\operatorname{GL}}
\def\Fr{\operatorname{Fr}}

\def\Pic{\operatorname{Pic}}
\def\alg{{\rm alg}}
\def\aut{{\rm aut}}
\def\op{{\rm op}}
\def\Tr{\operatorname{tr}}

\def\bM{\begin{matrix}}
\def\eM{\end{matrix}}

\def\pt{\bullet}
\def\lr{left (resp. right) }
\def\rl{right (resp. left) }

\def\bpr{\noindent{\bf{Proof\/}}:\;\;}
\def\To{\begin{CD} @>>>\end{CD}}

\def\SRL{_{S_{R\times L}}}
\def\FRL{_{F_{R\times L}}}
\def\RL{_{R\times L}}
\def\LR{_{L\times R}}
\def\Rl{_{R,L}}
\def\Lr{_{L,R}}

\begin{document}

\setcounter{page}{0}
{\pagestyle{empty}
\null\vfill
\begin{center}
{\LARGE $n$-dimensional global correspondences of Langlands}
\vfill
{\sc C. Pierre\/}
\vskip 11pt

Institut de Mathématique pure et appliquée\\
Université de Louvain\\
Chemin du Cyclotron, 2\\
B-1348 Louvain-la-Neuve,  Belgium\\
pierre@math.ucl.ac.be

\vfill

\begin{abstract}
The program of Langlands is studied here on the basis of:
\Bi
\item new concepts of global class field theory related to the explicit construction of global class fields and of reciprocity laws;
\item the representations of the reductive algebraic groups $\GL(n)$ constituting the $n$-dimensional representations of the associated global Weil-Deligne groups;
\item a toroidal compactification of the conjugacy classes of these reductive algebraic groups whose analytic representations constitute the cuspidal representations of these groups $\GL(n)$ in the context of harmonic analysis.
\Ei

This leads us to build two types of $n$-dimensional global bilinear correspondences of Langlands by taking into account the irreducibility or the reducibility of the representations of the considered bilinear algebraic semigroups.  

The major outcome of this global approach is the generation of general algebraic symmetric structures, consisting of double symmetric towers of conjugacy class representatives of algebraic groups, so that the analytic toroidal representations of these conjugacy classes are the cuspidal conjugacy class representations of these algebraic groups.

\end{abstract}
\vfill
\eject
\end{center}

\tableofcontents
\vfill\eject}
\setcounter{page}{1}
\def\thepage{\arabic{page}}
{\parindent=0pt 
\setcounter{section}{0}
\section*{Introduction}
\addcontentsline{toc}{section}{Introduction}

After the pioneer works of G. Laumon, M. Rapoport, U. Stuhler and others, recent significant advances were realized these past few years in the Langlands program \cite{Lan1}, \cite{Kna} by M. Harris, R. Taylor \cite{H-T} and G. Henniart \cite{Hen} with respect to the proof of the Langlands correspondences in the case of $p$-adic number fields and by L. Lafforgue \cite{Laf}, \cite{Lau} in the case of function fields.
\vskip 11pt

However, the approach of the Langlands program from a global point of view was relatively left sideways until now: it is the aim of this paper to try to fill up partially this gap by studying the Langlands correspondences over global number fields.
\vskip 11pt

The global methods used in this context are based on:
\Bi
\item new concepts in global class field theory related to the explicit construction of global class fields and of reciprocity laws;
\item the representations of the reductive algebraic groups $\GL(n)$ constituting the $n$-dimensional representations of the associated global Weil-Deligne groups;
\item a toroidal compactification of the conjugacy classes of these algebraic groups whose analytic representations constitute the cuspidal representations of these groups $\GL(n)$ in the context of harmonic analysis.  It will be shown that the global approach of the Langlands program consists in studying the algebraic groups and their representations to bridge the gap between analytic and algebraic problems as described in the diagram:
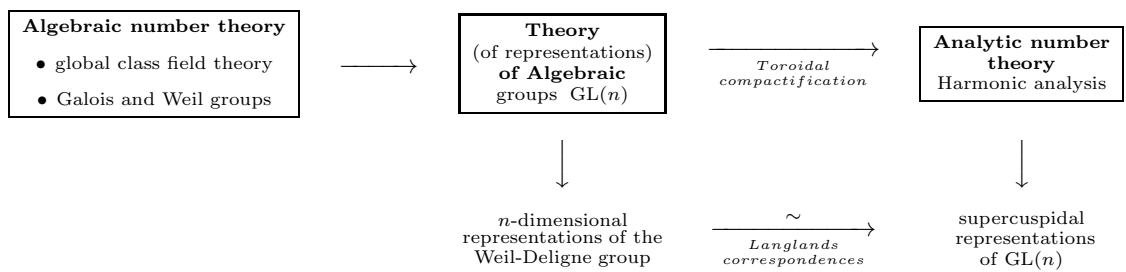
\Ei
\begin{footnotesize}
\[\bM
\boxed{\begin{smallmatrix}
\text{\bf Algebraic number theory}\\[6pt]
\bullet\  \text{global class field theory}\\[6pt] \bullet\  \text{Galois and Weil groups}\end{smallmatrix}}
& \To &
\boxed{\begin{smallmatrix} 
\text{\bf Theory}\\ \text{(of representations)}\\ \text{\bf of Algebraic}\\ \text{groups}\ \GL(n)\end{smallmatrix}}
& \begin{CD}
@>>{{Toroidal}\atop{compactification}}> \end{CD} &
\boxed{\begin{smallmatrix} 
\text{\bf Analytic number}\\ \text{\bf theory}\\ \text{Harmonic analysis}\end{smallmatrix}} \\
&& \begin{CD} @VVV \end{CD} && \begin{CD} @VVV \end{CD} \\
&& \begin{smallmatrix} \text{$n$-dimensional}\\ \text{representations of the}\\ \text{Weil-Deligne group}\end{smallmatrix}
& \begin{CD} @>{\sim}>{{Langlands}\atop{correspondences}}> \end{CD}
&   \begin{smallmatrix} \text{supercuspidal}\\\text{ representations}\\ \text{of $\GL(n)$} \end{smallmatrix} \eM\]\end{footnotesize}

\def\rop{\rm op}

The developments considered here will essentially concern symmetric objects in such a way that the envisaged mathematical objects will be cut into two symmetric semiobjects $\Os _R$ and $\Os _L$ localized respectively in (or referring to) the lower and in the upper half space.
\vskip 11pt

The right semiobject $\Os _R$ is then the dual of the left semiobject $\Os _L$ and the interest of considering an object, decomposed into two dual symmetric semiobjects $\Os _R$ and $\Os _L$~, is that the informations concerning the internal mathematical structure of the object ``~$\Os $~'' can be obtained from the product  $\Os _R\times \Os _L$ of the semiobjects $\Os _R$ and $\Os _L$~.  Indeed, every endomorphism $E$ of  ``~$\Os $~'' can be decomposed into the product $E_R\times E_L$ of a right endomorphism $E_R$ acting on the right semiobject $\Os _R$ by the ``opposite'' left endomorphism $E_L$ acting on the corresponding left semiobject $\Os _L$ in such a way that $E_R\simeq E_L^{-1}$~.
\vskip 11pt

The general existence of symmetric objects can be shown by the following considerations on function fields.  Let $k$ denote a global number field of characteristic zero.  Let $k[x_1,\cdots,x_n]$ be a polynomial ring over $k$ 
and let  $R=R_R\cup R_L$ be a symmetric 
(closed) extension of $k$ as introduced in the following.

Let $k[x_1,\cdots,x_n]$ be a polynomial ring at $n$ indeterminates over $k$~.

Let $I_L=\{P_\mu (x_1,\cdots,x_n)\mid P_\mu (V_L)=0\}$ be the ideal of $k[x_n,\cdots,x_n]$ in such a way that:
\Bean
\item $P_\mu (V_L)$ be the polynomial function in $k[V_L]$ represented by $P_\mu (x_1,\cdots,x_n)$~.
\item
$V_L\subset R_L$ be an affine semispace restricted to the upper half space.
\Ee
Let $I_R=\{P_\mu (-x_1,\cdots,-x_n)\mid P_\mu (V_R)=0\}$ be the symmetrical ideal of $I_L$ obtained by the involution:
\[ \tau : \quad P_\mu (x_1,\cdots,x_n) \To P_\mu (-x_1,\cdots,-x_n)\]
in such a way that
\Bean
\item $P_\mu (V_R)$ be the polynomial function in $k[V_R]$ represented by $P_\mu (-x_1,\cdots,-x_n)$~.
\item
$V_R\subset R_R$ be an affine semispace
\Bi
\item of dimension $n$~.
\item restricted to the lower half space.
\item symmetric of $V_L$~.
\item disjoint of $V_L$ or possibly connected to $V_L$ on the symmetric axis, plane, \ldots
\Ei
\Ee
The quotient ring obtained modulo the ideal $I_L$ (resp. $I_R$~) is $Q_L=k[x_1,\cdots,x_n]\big/ I_L$ (resp. $Q_R=k[x_1,\cdots,x_n]\big/ I_R$~) \cite{Wat}.

$Q_L$ and $Q_R$ are quotient algebras characterized by the
corresponding homomorphisms:
\[ \phi _L : \quad Q_L\To R_L\;, \qquad\qquad  \phi _R: \quad Q_R\To R_R\;,\]
where $R_L$ and $R_R$ are commutative (division) semirings localized respectively in (or referring to) the upper and the lower half space.  (Commutative (division) semirings are recalled (or introduced) in the appendix).
\vskip 11pt

So, the pair of homomorphisms $\phi _L$ and $\phi _R$ sends the ``general'' solution to a pair of symmetric solutions respectively in $R_L$ and $R_R$~.

On the other hand, let $\widetilde F=\widetilde F_R\cup \widetilde F_L$ be a symmetric ``algebraic'' finite extension of $k$~.

Let $T_n(\widetilde F_L)$ (resp. $T^t_n(\widetilde F_R)$~) $\subset \GL_n(\centerdot)$ be the (semi)group of matrices of dimension $n$ over $\widetilde F_L$ (resp. $\widetilde F_R$~) viewed as an operator sending $\widetilde F_L$ (resp. $\widetilde F_R$~) into the affine semispace
$T^{(n)}(\widetilde F_L)$ (resp. $T^{(n)}(\widetilde F_R)$~) of dimension $n$~:
\begin{align*}
T_n(\centerdot) : \quad \widetilde F_L &\To T^{(n)}(\widetilde F_L)\\
\text{(resp.} \quad
T^t_n(\centerdot) : \quad \widetilde F_R &\To T^{(n)}(\widetilde F_R)\ )\end{align*}
in such a way that to the indeterminates $(x_1,\cdots,x_\ell ,\cdots,x_n,\cdots, x_{1\ell },\cdots, x_{\ell n},\cdots)$ of $Q_L$ (resp. $(-x_1,\cdots,-x_\ell ,\cdots,-x_n,\cdots, -x_{1\ell },\cdots, -x_{\ell n},\cdots)$ of $Q_R$~), $\forall\ x_{\ell n}=x_\ell \cdot x_n$~, corresponds the homomorphism \cite{Bor2}:
\begin{align*}
\phi' _L: \quad Q_L &\To  \widetilde F_L  (x_1\to e_{11},\cdots,x_\ell\to e_{\ell\ell} ,\cdots, x_{1\ell }\to e_{1\ell},\cdots, x_{n\ell }\to e_{n\ell})\\
\text{(resp.} \quad
\phi' _R: \quad Q_R &\To  \widetilde F_R  (-x_1\to -e_{11},\cdots,-x_\ell\to -e_{\ell\ell} ,\cdots,\\
& \hspace{6cm} x_{\ell 1}\to e_{\ell 1},\cdots, x_{n\ell }\to e_{n\ell})\ )\end{align*}
where:
\begin{align*}
T_n(\widetilde F_L) &= \{e_L=(e_{\ell n})\in T_n(\widetilde F_L)\mid P_{T_\mu }(e_{\ell n})=0\}\\
\text{(resp.} \quad
T^t_n(\widetilde F_R) &= \{e_R=(e_{\ell n})\in T^t_n(\widetilde F_R)\mid P_{T_\mu }(e_{\ell n})=0\}\ )\end{align*}
with the polynomials $P_{T_\mu }(e_{\ell n})\in k[x]$~.

On the other hand, let $X_L$ (resp. $X_R$~) be  a functor from the quotient ring $Q_L$ (resp. $Q_R$~) to the affine semispace $V_L$ (resp. $V_R$~) in such a way that the diagram:
\begin{align*}
&\begin{CD} Q_L @>{\phi' _L}>> \widetilde F_L(e_{11},\cdots,e_{nn},\cdots, e_{\ell n})\\
@VV{X_L}V    @VV{T_n(\centerdot)}V\\
V_L @>>{\psi _L}> T^{(n)}(\widetilde F_L)
\end{CD}\\[15pt]
\mbox{\Huge(} \text{resp} \quad
&\begin{CD} Q_R @>{\phi' _R}>> \widetilde F_R(-e_{11},\cdots,-e_{nn},\cdots, e_{\ell n})\\
@VV{X_R}V    @VV{T^t_n(\centerdot)}V\\
V_R @>>{\psi _R}> T^{(n)}(\widetilde F_R)
\end{CD}\qquad \mbox{\Huge)}\end{align*}
commutes, leading to:
\[ X_L=\psi _L^{-1}\circ T_n\circ \phi' _L \qquad \text{(resp.} \quad
X_R=\psi _L^{-1}\circ T^t_n\circ \phi _R\ ).\]
Thus, if the elements $v_L\in V_L$ (resp. $v_R\in V_R$~), resulting from $X_L\circ \phi _L^{'-1}(\widetilde F_L)$ (resp. $X_R\circ \phi _R^{'-1}(\widetilde F_R)$~), correspond to solutions in $\widetilde F_L$ (resp. $\widetilde F_R$~) of some family of equations, there is a $k$-algebra $Q_L$ (resp. $Q_R$~) and a natural correspondence between $X_L\circ \phi _L^{'-1}\widetilde F_L$ (resp. $X_R\circ \phi _R^{'-1}\widetilde F_R$~) and $\Hom_{k}(Q_L,\widetilde F_L)$ (resp. $\Hom_{k}(Q_R,\widetilde F_R)$~), taking into account the algebraic (semi)group of matrices $T_n(\widetilde F_L)$ (resp. $T^t_n(\widetilde F_R)$~) at the condition that $\psi _L:V_L\to T^{(n)}(\widetilde F_L)$ (resp. $\psi _R:V_R\to T^{(n)}(\widetilde F_R)$~) be a homeomorphism.

\vskip 11pt

Such $X_L$ (resp. $X_R)$ is then called representable in the sense that a \lr affine semigroup scheme over $k$ is a  representable functor from the $k$-algebra $Q_L$ (resp. $Q_R$~) to the affine semispace $V_L$ (resp. $V_R$~), homeomorphic to $T^{(n)}(\widetilde F_L)$
(resp. $T^{(n)}(\widetilde F_R)$~).
\vskip 11pt

The left and right affine semigroup schemes $X_L$ and $X_R$ are said to be symmetric if every element $v_L\in V_L$~, localized in the upper half space, is symmetric to every element $v_R\in V_R$~, localized in the lower half space with respect to an element or a set of elements of symmetry.
\vskip 11pt 

Thus, the consideration of a sufficiently big polynomial ring as $k[x_1,\cdots,x_n]$ allows to consider objects as being generally symmetric and able to be cut into two symmetric semiobjects in one-to-one correspondence.
\vskip 11pt

In this respect, as the endomorphisms $\End_{\widetilde F}(B)$ of a division $\widetilde F$-algebra $B$ can be handled throughout its enveloping algebra $B^e=B\otimes_{ \widetilde F}B^{\rop}$~, where $B^{\rop}$ denotes the opposite algebra of $B$~, because $B^e\simeq \End_{\widetilde F}(B)$ and as fundamental algebras, such as the algebra of automorphic forms, are essentially defined on half spaces, we shall work in a bilinear mathematical framework.  So, enveloping semialgebras as well as bisemialgebras will be considered, which has the following advantages:
\Be 
\item As the representation of automorphic forms is intrinsically defined in the upper half space, the bialgebra of automorphic biforms, associated to the algebra of automorphic forms, will be envisaged in such a way that the coalgebra will refer to dual automorphic forms localized in the lower half space.  One of the interests of considering a bialgebra of automorphic biforms is that the endomorphism of this bialgebra leads to take into account the Hecke bialgebra of tensor products of Hecke bioperators having a nice matricial representation \cite{Pie3}.
\item The representation of a reducible general bilinear (semi)group of order $2n$ decomposes diagonally following the direct sum of irreducible bilinear representations and off diagonally following the direct sum of complementary irreducible bilinear representations (see section 0.3 of the introduction and chapter 4).\pagebreak

So, the nonorthogonal reducibility of the representation of a bilinear general Lie semigroup generating off diagonal  complementary irreducible bilinear representations could lead to a new approach of:
\Bi
\item the endoscopy problem considered in trace formulas;
\item the functoriality envisaged in automorphic representations of general ``linear'' groups \cite{Lan2}, \cite{C-K-P-S}.
\Ei

\item Some conjectures of algebraic number theory proceed probably from the fact that the techniques of algebra endomorphisms could be more worked out, as for example, by considering judicious enveloping (semi)algebras \cite{Pie3}.
\Ee
In this context, let us point out that the appendix of this paper introduces the used semistructures and bisemistructures \cite{Pie5}.
\vskip 11pt

In a few words, let us say that the Langlands global program developed here is related to the generation of a reductive algebraic group with entries in a double tower of extensions of $k$ which are characterized by increasing Galois extension degrees.
\vskip 11pt

The conjugacy class representatives of this reductive algebraic group then appear in symmetric pairs in such a way that a double tower of conjugacy class representatives, characterized by increasing ranks, is generated symmetrically respectively in the upper and in the lower half space and constitutes the $n$-dimensional representation of a global Weil-Deligne group.
\vskip 11pt

The Langlands global correspondences then consist in finding the toroidal analytic representations of these conjugacy class representatives by means of a suitable toroidal compactification of these in such a way that the sum of these toroidal analytical representatives of conjugacy classes constitutes the searched supercuspidal representation of the reductive algebraic group, i.e. the Fourier development of the $n$-dimensional automorphic form on this algebraic group.
\vskip 11pt

So, the Langlands global program is directly based on the generation of a general algebraic symmetric structure being in one-to-one correspondence with its automorphic analytical counterpart.


\subsection{Nonabelian global class field concepts}

More concretely, the new concepts introduced in global class field theory are inspired by the factorization of a prime
 $p$ into primes in the ring of integers $\Os_E$ of a finite extension $E$ of $\QQ$ and by the Artin's reciprocity law which 
asserts that a Dirichlet character 
$\chi _\sigma :(\ZZ\big/N\ \ZZ)^*\to \CC^*$ exists such that $\sigma (\Fr p)=\chi _\sigma (p)$ where $\sigma :\Gal(E\big/\QQ)\to \CC^*$~.  This leads us to introduce pseudo-ramified completions, covered by the symmetric algebraic extensions of a global number field $k$ of characteristic zero, at the infinite places such that their degrees or ranks are integers modulo $N$ where $N$ is the order of the global inertia subgroups.  The subsets of  pseudo-ramified completions are thus covered by finite subsets of Galois extensions and generate a double symmetric tower of one-dimensional $k$-semimodules characterized by increasing ranks.
  A character $\chi _{\omega _j}$ is associated to each basic completion $F_{\omega _j}$ in an infinite or global place $\omega _j$ and the Dirichlet character $\chi _{\omega _{j,m_j}}$~, corresponding to an equivalent completion $F_{\omega _{j,m_j}}$ in $\omega _j$~, can be obtained by the action of the decomposition group element $D_{\omega _{j,m_j}}$ on $\chi _{\omega _j}$~.
As we are concerned with pairs of symmetric (semi)structures, left algebraic extension semifields
$\wt F_\omega =\{\wt F_{\omega _1},\dots,\wt F_{\omega _{j,m_j}},\dots,\wt F_{\omega _r}\}$
 referring to the upper half space and right algebraic extension semifields 
$\wt F_{\o\omega} =\{\wt F_{\o\omega _1},\dots,\wt F_{\o\omega _{j,m_j}},\dots,\wt F_{\o\omega _r}\}$
referring to the lower half space are considered.  

In correspondence with these left and right algebraic extension semifields, we have a left and a right tower $F_\omega =\{F_{\omega _1},\cdots,F_{\omega _{j,m_j}},\cdots,F_{\omega _r}\}$ and $F_{\o\omega} =\{F_{\o\omega _1},\cdots,F_{\o\omega _{j,m_j}},\cdots,F_{\o\omega _r}\}$ of packets of equivalent completions characterized by increasing ranks.

From these, we can consider:
\Bean \item their direct sums
\[ F_{\omega _\oplus}=\bigoplus_j\bigoplus_{m_j}F_{\omega _{j,mj}} \quad \and\quad
F_{\o\omega _\oplus}=\bigoplus_j\bigoplus_{m_j}F_{\o\omega _{j,mj}} \;,\]
\item and the products
\[ 
\Aa^\infty _{F_\omega }
=\prod_{j_p}F_{\omega _{j_p}}
\quad \and\quad
\Aa^\infty _{F_{\o\omega }}=\prod_{j_p}F_{\o\omega _{j_p}}
\]
\Ee
of packets of primary  completions 
$F_{\omega _{j_p}}$ and $F_{\o\omega _{j_p}}$
where $\Aa^\infty _{F_\omega }$ and $\Aa^\infty _{F_{\o\omega }}$ are infinite adele semirings.

This allows to introduce the Galois subgroups of these left and right algebraic extensions and the global Weil groups $W^{ab}_{F_L}$ and $W^{ab}_{F_R}$ referring respectively to (the sums of) the automorphisms of left and right algebraic extensions explicited subsequently.  

The set of \lr pseudo-ramified extensions considered until now constitute a \lr affine semigroup  $\SS^1_L$ (resp. $\SS^1_R$~).
\vskip 11pt 

All that constitutes elementary concepts of abelian global class field theory.  To enter into the $n$-dimensional global Langlands program, nonabelian global class field concepts have to be set up \cite{Lan2}.  The challenge then consists in building up the $n$-dimensional equivalent of the affine semigroup  $\SS^1_L$ (resp. $\SS^1_R$~), and from a bilinear point of view, the $n$-dimensional equivalent of the bilinear affine semigroup $\SS^1_R\times \SS^1_L$~,
which needs the introduction of the injective morphisms \cite{Gel}:
\[ \bM 
\sigma _L : & \quad & W^{ab}_{F_L} & \To & T_n(\wt F_{\omega_\oplus} )\;, \\[-6pt]
\sigma _R : & \quad & W^{ab}_{F_R} & \To & T^t_n(\wt  F_{\o\omega_\oplus} )\;, \\[-6pt]
\sigma \RL : & \quad & W^{ab}_{F_R}\times W^{ab}_{F_L} & \To & 
\GL_n(\wt F_{\o\omega_\oplus} \times \wt F_{\omega_\oplus} ) \equiv T^t_n( \wt F_{\o\omega_\oplus} )\times T_n(\wt F_{\omega_\oplus} )
 \;.\eM\]
The affine bilinear algebraic semigroup $G^{(2n)}(\wt F_{\o\omega} \times \wt F_{\omega} )\equiv T^{(2n)}(\wt F_{\o\omega })\times T^{(2n)}(\wt F_\omega )$ results from the action 
 of the bilinear 
algebraic semigroup of matrices $\GL_n( \wt F_{\o\omega} \times \wt F_{\omega} )\equiv T^t_n(\wt F_{\o\omega} )
\times T_n(\wt F_{\omega} )$ where $T^t_n(\wt F_{\o\omega} )$ is the group of lower triangular matrices with entries in 
${\wt F_{\o\omega} }$ while $T_n({\wt F_{\omega} })$ is the group of upper triangular matrices with entries in ${\wt F_\omega }$~.  The algebraic bilinear affine semigroup  $G^{(2n)}( {\wt F_{\o\omega }}\times {\wt F_\omega } )$ over 
$ {\wt F_{\o\omega }}\times {\wt F_\omega } $ is a $\GL_n( {\wt F_{\o\omega }}\times {\wt F_\omega } )$-bisemimodule $\wt M_R\otimes \wt M_L$~.  And, the algebraic representation of $\GL_n( {\wt F_{\o\omega }}\times {\wt F_\omega } )$ in $(\wt M_R\otimes\wt M_L)$ refers to an algebraic morphism from $\GL_n( {\wt F_{\o\omega }}\times {\wt F_\omega } )$ into $\GL(\wt M_R\otimes \wt M_L)$ which denotes the group of automorphisms of $(\wt M_R\otimes \wt M_L)$~.

%
The reasons for considering bilinear algebraic (semi)groups $\GL_n({\wt F_{\o\omega }}\times {\wt F_\omega })$ instead of linear algebraic groups $\GL_n({\wt F_{\o\omega -\omega }})$ result from the facts that:
\Be
\item a bilinear algebraic semigroup covers the corresponding linear algebraic group (as it is proved in chapter 2) in the sense that the representation space of this linear algebraic group is a $n^2$-complex dimensional vector space $W$ corresponding to the $n^2$-complex dimensional representation space $(\wt M_R\otimes \wt M_L)$ of the bilinear algebraic (semi)group $\GL_n({\wt F_{\o \omega }}\times {\wt F_\omega })$~.  So, we have: \quad $W\simeq \wt M_R\otimes \wt M_L$~;
\item a bilinear algebraic semigroup is directly connected to the enveloping algebras allowing to generate endomorphisms as it was explained above.
\Ee

The bilinear semigroup $G^{(2n)}(F_{\o \omega }\times F_\omega )$ over the product, right by left, of pseudo-ramified completions $F_{\o\omega }$ and $F_\omega $, is a complete ``algebraic'' bilinear semigroup  or an abstract bisemivariety since the sets of completions $F_{\o\omega }$ and $F_\omega $ are covered by the corresponding sets of algebraic extensions $\wt F_{\o\omega }$ and $\wt F_\omega $.  By concern of abbreviation, this complete ``algebraic'' bilinear semigroup $G^{(2n)}(F_{\o\omega }\times F_\omega )$ will sometimes be simply called ``algebraic'' bilinear semigroup, taking into account a universal property between these bilinear semigroups.

Let $\wt M_{R_\oplus}\otimes \wt M_{L_\oplus}$ denote the representation space of $\GL_n(\wt F_{\o\omega _\oplus}\times \wt F_{\omega _\oplus})$ in such a way that 
$\wt M_{R_\oplus}\otimes \wt M_{L_\oplus}$ decomposes into the direct sum of subbisemimodules $\wt M_{\o\omega _{j,m_j}}\otimes \wt M_{\omega _{j,m_j}}$ being the conjugacy class representatives of $\GL_n(\wt F_{\o\omega }\times \wt F_\omega )$~.

Then, $\GL(\wt M_{R_\oplus}\otimes \wt M_{L_\oplus})$ constitutes the $n$-dimensional complex equivalent of the product 
$W^{ab}_{F_R}\times W^{ab}_{F_L}$ of the global Weil groups.  As $\GL(\wt M_{R_\oplus}\otimes \wt M_{L_{\oplus}})$ is isomorphic to 
$\GL_n( {\wt F_{\o\omega_\oplus }}\times {\wt F_{\omega_\oplus} } )$~, the bilinear algebraic semigroup 
$G^{(2n)}( {\wt F_{\o\omega_\oplus }}\times {\wt F_{\omega_\oplus} } )$ becomes the $2n$-dimensional (irreducible) complex representation (space) 
$(\Irr)\Rep^{(2n)}_{W_{F\RL}}(W^{ab}_{F_R}\times W^{ab}_{F_L})$ of $W^{ab}_{F_R}\times W^{ab}_{F_L}$~:
\[ (\Irr)\Rep^{(2n)}_{W_{F\RL}}(W^{ab}_{F_R}\times W^{ab}_{F_L}): \quad
\GL(\wt M_{R_\oplus }\otimes \wt M_{L_\oplus }) \To G^{(2n)}( {\wt F_{\o\omega_\oplus }}\times {\wt F_{\omega _\oplus}} )\;.\]

As the bilinear algebraic semigroup $G^{(2n)}({\wt F_{\o\omega }}\times {\wt F_\omega } )$ is built over $( {\wt F_{\o\omega }}\times {\wt F_\omega } )$~, it is composed of $r$ conjugacy classes, $1\le j\le r\le \infty$~, having multiplicities $m^{(r)}$ and associated with the $r$ places of ${F_\omega }$ or of ${F_{\o\omega }}$~.  Consequently,
$(\Irr)\Rep^{(2n)}_{W_{F\RL}}(W^{ab}_{F_R}\times W^{ab}_{F_L})$ also decomposes into $r$ conjugacy classes with multiplicities $m^{(r)}$~.

This decomposition of $G^{(2n)}( {\wt F_{\o\omega }}\times {\wt F_\omega } )$ into conjugacy classes can also be obtained by considering the cutting of bilattices $\Lambda _{\o\omega }\times \Lambda _\omega $ into subbilattices in $(\wt M_R\otimes \wt M_L)$ under the action of the product $T_R(n;r)\otimes T_L(n;r)$ of Hecke operators having as representation $\GL_{2n}((\ZZ\big/ N\ \ZZ)^2)$~, $N\in\NN$~.

This leads to a one-to-one correspondence between the pseudo-ramified extensions $\wt F_{\omega _j}$ (and $\wt F_{\omega _{j,m_j}}$~), the diagonal (and off-diagonal) conjugacy classes of the bilinear algebraic semigroup $G^{(2n)}({\wt F_{\o\omega }}\times {\wt F_\omega } )$ and the subbilattices of Hecke in the $\GL_n({\wt F_{\o\omega }}\times {\wt F_\omega } )$-bisemimodule $\wt M_R\otimes \wt M_L$~.
\vskip 11pt 

Let $G^{(2n)}( {\wt F^+_{\o v}} \times  {\wt F^+_{v}} )$ denote the bilinear algebraic semigroup over the product of 
the symmetric sets $\wt F^+_{\o v}=\{\wt F^+_{\o v_{1_\delta }},\cdots,\wt F^+_{\o v_{j_\delta ,m_{j_\delta }}},\cdots,\wt F^+_{\o v_{r_\delta }}\}$ and
$\wt F^+_{v}=\{\wt F^+_{v_{1_\delta }},\cdots,\wt F^+_{v_{j_\delta ,m_{j_\delta }}},\cdots,\wt F^+_{v_{r_\delta }}\}$
of real extensions
 and let 
$G^{(2n)}( F^+_{\o v} \times  F^+_{v} )$ denote its compact equivalent.  
The compactification of $G^{(2n)}( \widetilde F^+_{\o v} \times  
\widetilde F^+_{v} )$ into 
$G^{(2n)}( {F^+_{\o v}} \times  {F^+_{v}} )$ is realized by means of the Fulton-McPherson compactification by a 
set of successive blowups \cite{F-M} as developed in chapter 3.
\vskip 11pt

If the conjugacy class representatives $g^{(2n)}\RL[j ,m_{j }]$ of 
$G^{(2n)}( {F^+_{\o \omega }} \times  {F^+_{\omega }} )$ are glued together, the bifunction 
$(f_{\o \omega }(z^*)\otimes f_\omega (z))$ on $G^{(2n)}( {F^+_{\o \omega }} \times  {F^+_{\omega }} )$ 
has a holomorphic development in the multiple power series
\[ f_{\o \omega }(z^*)\otimes f_\omega (z) = \sum^\infty_{j =1}\sum_{m_{j }} c^*_{j ,m_{j }}\ c_{j ,m_{j }}\ 
(z^*_1\ z_1-z^*_{01}\ z_{01})^{j } \ \cdots\ (z^*_n\ z_n-z^*_{0n}\ z'_{0n})^{j }\]
at the bipoint $(z^*_0\times z_0)$~, $z_0\in \cit^n$~, such that each term of $f_{\o \omega }(z^*)\otimes f_\omega (z)$ corresponds to a 
conjugacy class representative $g^{(2n)}\RL[j ,m_{j }]$~. And, the holomorphic bifunction 
$f_{\o \omega }(z^*)\otimes f_\omega (z)$ constitutes an irreducible holomorphic representation 
$\Irr\hol( \GL_n({F_{\o \omega }}\times {F_\omega } ))$ of 
$\GL_n({F_{\o \omega  }}\times {F_\omega } )\approx 
G^{(2n)}({F_{\o \omega }}\times {F^+_\omega } )$~.\vskip 11pt

\subsection{Langlands irreducible global bilinear correspondences}

Let $\widetilde F_L$ (resp. $\widetilde F_R$~) denote a complex symmetric splitting semifield.

Then, the pseudo-ramified lattice bisemispace
\[ X\SRL =\GL_n({\widetilde F_{R }}\times {\widetilde F_L } )\Big/\GL_n((\ZZ\big/ N\ \ZZ)^2)\approx G^{(2n)}(F_{\o\omega }\times F_\omega )\]
was introduced such that its cosets are isomorphic to the conjugacy classes of $G^{(2n)}({\widetilde F_{\o\omega  }}\times{\widetilde F_{\omega } } )$~.

Let $\{\wt F_{\omega^1} \}$ (resp. $\{\wt F_{\o\omega^1 }\}$~) denote  the set of irreducible extensions  $\wt F_{\omega^1_j} $ (resp. ${\wt F_{\o\omega ^1_j}}$~) having a rank $N$~. Then, the smallest pseudo-ramified normal bilinear algebraic subsemigroup $P^{(2n)}({\wt F_{\o\omega^1} }\times {\wt F_{\omega^1} })$ of 
$G^{(2n)}({\wt F_\omega }\times {\wt F_\omega })$ is introduced as being 
the complex bilinear equivalent of the minimal (not standard) 
parabolic subgroup: it corresponds to the irreducible $n$-dimensional complex
representation of the product $I_{F_R}\times I_{F_L}$ of global inertia 
subgroups having an order equal to $(N\cdot m^{(j)})$ where $m^{(j)}=m_j+1$ is the multiplicity of the $j$-th complex place.
\vskip 11pt 

Then, a toroidal compactification of $X\SRL$ is envisaged by mapping $X\SRL$ into the corresponding toroidal compactified lattice bisemispace
\[ \o X\SRL=\GL_n( {F^T_{R }}\times {F^T_L } )\Big/ \GL_n((\ZZ\big/ N\ \ZZ)^2)\approx G^{(2n)}(F^T_{\o\omega }\times F^T_\omega )\]
in such a way that $X\SRL$ may be viewed as the interior of $\o X\SRL$ in the sense of the Borel-Serre compactification.  
${F^T_\omega } $ and ${F^T_{\o\omega }}$ are the sets of toroidal completions.

A double coset decomposition of the bilinear complete algebraic semigroup $\GL_n({F^T_{R }}\times {F^T_L } )$ gives rise to the compactified bisemispace
\[
\o S^{P_n}_{K_n}
= P_n({F^{T}_{\o\omega^1 }}\times {F^{T}_{\omega^1 }} )
\setminus \GL_n({F^T_{R}}\times {F^T_L } )\Big/ 
\GL_n((\ZZ\big/ N\ \ZZ)^2)\;. \]
A general bilinear cohomology  is introduced in section 3.2 as a contravariant functor $H^*$ from smooth abstract (resp. algebraic) bisemivarieties together with a natural transformation $n_{H^*\to H^{[*,*]}}$ from $H^*$ to the associated  the de Rham bilinear cohomology $H^{[*,*]}$.

This bilinear cohomology is:
\Bi
\item of general type in the sense that it is a motivic (or Weil) bilinear cohomology directly related to the singular, de Rham of Betti cohomologies,
\item characterized by Hodge (bisemi)-cycles, a Künneth standard conjecture and a Künneth biisomorphism.
\Ei
It results then that:
\[
H^{2i}(\o S^{P_n}_{K_n},M^{2i}_{T_{R_\oplus}}\otimes M^{2i}_{T_{L_\oplus}})
= \Reps (\GL_{i}({F^T_{\o\omega_\oplus }}\times {F^T_{\omega_\oplus} } ))\]
where 
\Bi
\item $M^{2i}_{T_{R_\oplus}}\otimes M^{2i}_{T_{L_\oplus}}$ is an irreducible $\GL_{i}({F^T_{\o\omega_\oplus }}\times {F^T_{\omega_\oplus} } )$-subbisemimodule of \\
 $M^{2n}_{T_{R_\oplus }}\otimes M^{2n}_{T_{L_\oplus }}$~, $i\le n$ is a complex dimension; $2i$ is the real corresponding dimension~;

\item $\Reps (\GL_{i}({F^T_{\o\omega_\oplus }}\times {F^T_{\omega_\oplus} } ))$ is the  representation space of 
the complete bilinear semigroup of matrices $\GL_{i}({F^T_{\o\omega_\oplus }}\times {F^T_{\omega_\oplus} } )$~.
\Ei
\vskip 11pt

As the toroidal compactification of $X\SRL$ is an isomorphism, we have that:
\begin{align*}
\Irr \Rep W^{(2i)}\FRL (W^{ab}_{F_R}\times W^{ab}_{F_L})
&\simeq G^{(2i)}({F^T_{\o\omega_\oplus }}\times {F^T_{\omega_\oplus} } )\\
&\simeq\Reps (\GL_{i}({F^T_{\o\omega_\oplus }}\times {F^T_{\omega_\oplus} } ))\;.\end{align*}

On the other hand, the bilinear cohomology $H^{2i}(\o S^{P_n}_{K_n},M^{2i}_{ T_{R_\oplus} }\otimes M^{2i}_{ T_{L_\oplus} })$ has a decomposition 
following the equivalent representatives $g^{2i}_{T\RL}[j,m_j]$ of the conjugacy classes of the complete bilinear semigroup 
$G^{2i}({F^T_{\o\omega }}\times {F^T_\omega } )$ according to:
\[
H^{2i}(\o S^{P_n}_{K_n},M^{2i}_{T_{R_\oplus} }\otimes M^{2i}_{T_{L_\oplus} })
= \bigoplus^r_{j=1}\ \bigoplus_{m_j}\ g^{2i}_{T\RL}[j,m_j]\;.\]

As a consequence, the bilinear cohomology $H^{2i}(\o S^{P_n}_{K_n},M^{2i}_{T_{R_\oplus} }\otimes M^{2i}_{ T_{L_\oplus} })$ has an 
analytic development consisting in the product $\Eis_R(2i,j,m_j)\otimes \Eis_L(2i,j,m_j)$ of the (truncated) Fourier development 
of a cusp form of weight 2 by its left equivalent.  
The $n$-dimensional complex cusp biform $\Eis_R(2i,j,m_j)\otimes \Eis_L(2i,j,m_j)$~, 
in one-to-one correspondence with the bilinear  cohomology $H^{2i}(\o S^{P_n}_{K_n},
M^{2i}_{T_{R_\oplus} } \otimes M^{2i}_{T_{L_\oplus} })$~, is constructed according to a solvable way and corresponds to 
an eigenbifunction of the product $T_R(i;r)\otimes T_L(i;r)$, right by left, of the 
Hecke operators; so, $\Eis_R(2i,j,m_j)\otimes \Eis_L(2i,j,m_j)$ constitutes an irreducible 
supercuspidal representation $\Irr\cusp(\GL_i({F^T_{\o\omega }}\times {F^T_\omega }))$~ of the general bilinear semigroup $\GL_i({F^T_{\o\omega }}\times {F^T_\omega })$~.
\vskip 11pt

This leads us to a Langlands irreducible global correspondence:  

\[\begin{CD}
\Irr\Rep^{(2i)}_{W_{F_{R\times L}}}(W^{ab}_{F_{R}}\times W^{ab}_{F_{L}})
@>>>
\Irr\cusp(\GL_i({F^T_{\o\omega }}\times {F^T_\omega }) \\
@|     @| \\
G^{(2i)}({\wt F_{\o\omega_\oplus }}\times {\wt F_{\omega_\oplus} }) @.
 \Eis\RL(2i,j,m_j) \\
@VV{\wr}V @AAA  \\
G^{(2i)}(F^T_{\o\omega _\oplus} \times F^T_{\omega _\oplus} )
@>>>
\widehat G^{(2i)}({F^T_{\o\omega }}\times {F^T_\omega }) 
\\
 @| \\
H^{2i}(\o S^{P_n}_{K_n},M^{2i}_{T_{R_\oplus}}\otimes M^{2i}_{T_{L_\oplus}})\end{CD}\]
\Bi
\item from the sum of the products, right by left, of the conjugacy classes (associated with the places 
of the semifields $F^T_\omega $ and $F^T_{\o\omega }$~) of the irreducible $i$-dimensional 
representation $\Irr\Rep^{(2i)}_{W_{F_{R\times L}}}(W^{ab}_{F_R}\times W^{ab}_{F_L})$ of the bilinear Weil-Deligne group given by 
$G^{(2i)}({F_{\o\omega_\oplus }}
\times {F_{\omega_\oplus} }))$~;
\item to the sum of the products, right by left, of the conjugacy classes of the irreducible 
cuspidal representation $\Irr\cusp(\GL_i({F^T_{\o\omega_\oplus }}\times {F^T_{\omega_\oplus} }))$ of $
\GL_i({F^T_{\o\omega }}\times {F^T_\omega })$ given by $\Eis_R(2i,j,m_j)\otimes \Eis_L(2i,j,m_j)$
\Ei
where $\widehat G^{(2i)}({F^T_{\o\omega }}\times {F^T_\omega }) 
$ is a (bisemi)sheaf over the complete bilinear semigroup 
$G^{(2i)}({F^T_{\o\omega }}\times {F^T_\omega }) $
\vskip 11pt

A second kind of Langlands irreducible global correspondence can be reached by considering the real equivalent of the preceding correspondence: this can be realized by envisaging the boundary $\partial \o X_{S_{R\times L}}= \GL_n({F^{+,T}_{R }}\times {F^{+,T}_L })\big/ \GL_n((\ZZ\big/N\ \ZZ)^2)$ of the compactified lattice bisemispace $\o X_{S_{R\times L}}$~, where 
${F^{+,T} _{L}}$ (resp. ${F^{+,T}_{R}}$~) is a toroidal compact semifield.  This inclusion morphism $
\gamma ^\delta _{R\times L}: \o X_{S_{R\times L}}\to \partial \o X_{S_{R\times L}}$ sends bijectively the (diagonal) ``complex'' class representatives  into the corresponding ``real'' conjugacy class representatives covering them.
This gives rise to a double coset decomposition
\[ \partial \o S^{P_n}_{K_n}=P_n(  {F^{+,T}_{\o v^1}}\times {F^{+,T}_{v^1}})
\setminus \GL_n({F^{+,T}_{R}}\times {F^{+,T}_{L}})\Big/\GL_n((\ZZ\big/N\ \ZZ)^2)\]
of the general  bilinear semigroup $ \GL_n({F^{+,T}_{R}}\times {F^{+,T}_{L}})$~.

$\partial \o S^{P_n}_{K_n}$ is then the equivalent of a Shimura (bisemi)variety whose cohomology:
\[H^{2i}(\partial \o S^{P_n}_{K_n},M^{2i}_{T_{\o v_{R_\oplus}}}\otimes M^{2i}_{T_{v_{_\oplus}L}})=
\Reps ( \GL_{2i}({F^{+,T}_{\o  v_\oplus}}\times {F^{+,T}_{v_\oplus}}))\] is the irreducible Eisenstein bilinear cohomology which has an analytic representation given by the product
\[\ELLIP_R(2i,j^\delta ,m_{j^\delta }) \otimes \ELLIP_L(2i,j^\delta ,m_{j^\delta })\] of a right $2i$-dimensional global elliptic semimodule $\ELLIP_R(2i,j^\delta ,m_{j^\delta })$ by its left equivalent.
$\ELLIP_R(2i,j^\delta ,m_{j^\delta })\otimes \ELLIP_L(i,j^\delta ,m_{j^\delta })$ is also a ``solvable'' (bi)series and is an eigenbifunction of the product of Hecke operators.

On the irreducible bisemivariey $\partial\o S^{P_n}_{K_n}$~, there is a Langlands irreducible global correspondence:
\[\begin{CD}
\Irr\Rep^{(2i)}_{W_{F^+_{R\times L}}}(W^{ab}_{F^{+}_{\o v }}\times W^{ab}_{F^{+}_{v }})
@>>>
\Irr \ELLIP(\GL_{2i}({F^{+,T}_{\o v }}\times {F^{+,T}_v }))\\
@| @|\\
G^{(2i)}({\wt F^+_{\o v_\oplus}}\times {\wt F^+_{v_\oplus}}) 
@>>>
 \ELLIP\RL(2i,j_\delta ,m_{j_\delta })\\
@VV{\wr}V @AAA \\
G^{(2i)}(F^{+,T}_{\o v_\oplus} \times F^{+,T}_{v_\oplus} )
@>>>
\widehat G^{(2i)}({F^{+,T}_{\o v}}\times {F^{+,T}_v}) 
\\
@| @.\\
H^{2i}(\o S^{P_n}_{K_n},M^{2i}_{T_{v_{R_\oplus}}}\otimes M^{2i}_{T_{v_{L_\oplus}}})\end{CD}\]
\Bi
\item from the sum of the products, right by left, of the conjugacy classes of the irreducible 
$n$-dimensional representation $\Irr\Rep^{(2i)}_{W_{F^+_{R\times L}}}(W^{ab}_{F^{+}_{\o v }}
\times W^{ab}_{F^{+}_v })$ of the bilinear global Weil group given by the complete bilinear 
real semigroup $G^{(2i)}({F^+_{\o v_\oplus }}\times {F^+_{v_\oplus }})$~;
\item to the sum of the products, right by left, of the conjugacy classes of the irreducible cuspidal (and elliptic) representation $\Irr \ELLIP(\GL_{2i}({F^{+,T}_{\o v }}\times {F^{+,T}_{v} }))$ of $
\GL_{2i}({F^{+T}_{\o v }}\times {F^{+T}_v })$ given by the ``solvable'' global elliptic bisemimodule 
$\ELLIP_R(2i,j^\delta ,m_{j^\delta })\linebreak\otimes \ELLIP_L(2i,j^\delta ,m_{j^\delta })$~.
\Ei

Two kinds of trace formulas have been  considered.  The first type of trace formula is an adaptation of the Arthur-Selberg 
trace formula to the operator associated with the bilinear parabolic semigroup $P^{(2n)} ( {F_{\o\omega^1 }}\times {F_{\omega^1 }})$ envisaged in the present context as the unitary representation of the complete bilinear semigroup $\GL_n
( {F_{\o\omega }}\times {F_\omega })$~.  This operator acts by convolution on the bialgebra
$L^{1-1}\RL (G^{(2n)}( F^{nr}_{\o\omega }\times F^{nr}_\omega ))$ of smooth continuous bifunctions 
on the pseudo-unramified bilinear complete semigroup $G^{(2n)}( F^{nr}_{\o\omega }\times F^{nr}_\omega )$ and 
decomposes according to the unitary conjugacy classes of the pseudo-ramified bilinear complete semigroup 
$G^{(2n)}( {F_{\o\omega }}\times {F_\omega })$~. The resulting trace formula occurs in the bialgebra 
$L^{1-1}\RL (G^{(2n)}( {F_{\o\omega }}\times {F_\omega }))$ of bifunctions on 
$G^{(2n)}( {F_{\o\omega }}\times {F_\omega })$ and relies on Lefschetz trace formula. Remark that the complete trace formula, referring 
to the set of irreducible representations of $G^{(2n)}( {F_{\o\omega }}\times {F_\omega })$ must be envisaged in the frame of chapter 4.
\vskip 11pt 

The second kind of trace formula, occurring directly in the bialgebra 
$L^{1-1}\RL(G^{(2n)}( {F_{\o\omega }}\times {F_\omega }))$ of bifunctions on the pseudo-ramified algebraic 
complete semigroup $G^{(2n)}( {F_{\o\omega }}\times {F_\omega }))$~, leads to the Plancherel formula 
and corresponds to the first type of trace formula.
\vskip 11pt


\subsection{Langlands reducible global bilinear correspondences}
While the bilinear global correspondences of Langlands deal essentially with the irreducible representations of algebraic general bilinear semigroups, reducible global bilinear correspondences will be constructed with respect to the reducibility of the representations of the considered general bilinear semigroups.  Three kinds of reducibility will be envisaged:
\Be
\item the representation of the complete bilinear semigroup $\GL_n({F_{\o\omega }}\times {F_\omega })$ will be said to be {\bf partially reducible} if it decomposes according to the direct sum of irreducible bilinear representations $\Rep( \GL_{n_\ell }({F_{\o\omega }}\times{F_\omega }))$ taking into account the partition\linebreak  $n=n_1+\cdots + n_\ell +\cdots+n_s$ of $n$~.

\item the representation of $\GL_{2n}({F_{\o\omega }}\times {F_\omega })$ will be said to be {\bf orthogonally completely reducible} if it decomposes ``diagonally'' following the direct sum of irreducible bilinear representations $\Rep(\GL_{2_\ell }({F_{\o\omega }}\times {F_\omega }))$~, $1\le \ell\le n$~.

\item the representation of $\GL_{2n}({F_{\o\omega }}\times {F_\omega })$  will be said to be {\bf nonorthogonally completely reducible} if it decomposes diagonally following the direct sum of irreducible bilinear representations $\Rep(\GL_{2_\ell }({F_{\o\omega }}\times {F_\omega }))$  and off diagonally following the direct sum of irreducible bilinear representations $\Rep(T^t_{2_{k_R}}({F_{\o\omega }})\times T_{2_{\ell_L}}({F_\omega }))$~.
\Ee
Let then $\o X^{n=n_1+\cdots+n_s}_{S_{R\times L}}$~,  $\o X^{2n=2_1+\cdots+2_n}_{S_{R\times L}}$
and $\o X^{2n_R\times 2n_L}_{S_{R\times L}}$ be the three kinds of reducible compactified lattice bisemispaces and let
$\o S^{P_{n=n_1+\cdots+n_s}}_{K_{n=n_1+\cdots+n_s}}$~, $\o S^{P_{2n=2_1+\cdots+2_n}}_{K_{2n=2_1+\cdots+2_n}}$ and $\o S^{P_{2n_R\times 2n_L}}_{K_{2n_R\times 2n_L}}$ be their corresponding double coset decompositions.

Their bilinear cohomologies then decompose according to:
\Bi
\item $\begin{array}[t]{ll}
H^*(\o S^{P_{n=n_1+\cdots+n_s}}_{K_{n=n_1+\cdots+n_s}},\widehat M^{2n}_{T_{R_\oplus}}\otimes \widehat M^{2n}_{T_{L_\oplus}}) &=\txt\txt\bigoplus\limits_{n_\ell} H^{2n_\ell}(\o S^{P_{n}}_{K_{n}},\widehat M^{2n_\ell}_{T_{R_\oplus}}\otimes \widehat M^{2n_\ell}_{T_{L_\oplus}}) \\
&=\txt\txt\bigoplus\limits_{n_\ell} \bigl(\Eis _R(2n_\ell,j_R,m_{j_R}) \otimes  \Eis _L(2n_\ell,j_L,m_{j_L}) \bigr)
\end{array}$
\vskip 11pt

\item $\begin{array}[t]{ll}
H^{2n}(\o S^{P_{2n=2_1+\cdots+2_n}}_{K_{2n=2_1+\cdots+2_n}},
\widehat M^{2n^*}_{T_{R_\oplus}}\otimes \widehat M^{2n^*}_{T_{L_\oplus}}) 
&=\ \txt\bigoplus\limits_{\ell} 
H^{2_\ell}(\o S^{P_{2_n}}_{K_{2_n}},\widehat M^{2_\ell}_{T_{R_\oplus}}\otimes 
\widehat M^{2_\ell}_{T_{L_\oplus}}) \\
& =\ \txt\bigoplus\limits_{\ell} \bigl(\Eis _R(2_\ell,j_R,m_{j_R}) 
\otimes  \Eis _R(2_\ell,j_L,m_{j_L})\bigr)\end{array}$
\vskip 11pt

\item $\begin{array}[t]{ll}
H^{2n}(\o S^{P_{2n_R\times 2n_L}}_{K_{2n_R\times 2n_L}},\widehat M^{2n_R}_{T_{R_\oplus}}\otimes 
\widehat M^{2n_L}_{T_{L_\oplus}}) 
 & =\ \txt\bigoplus\limits_{\ell} H^{2_\ell}(\o S^{P_{2n_R\times 2n_L}}_{K_{2n_R\times 2n_L}},
\widehat M^{2_\ell}_{T_{R_\oplus}}\otimes 
\widehat M^{2_\ell}_{T_{L_\oplus}})\\ 
&\quad \txt\bigoplus\limits_{k_R\neq \ell_L} H^{2_{k_R,\ell_L}}(\o S^{P_{2n_R\times 2n_L}}_{K_{2n_R\times 2n_L}},
\widehat M^{2_{k_R}}_{T_{R_\oplus}} \otimes 
\widehat M^{2_{\ell_L}}_{T_{L_\oplus}})\\
& =\ \txt\bigoplus\limits_{\ell} \bigl(\Eis _R(2_\ell,j_R,m_{j_R}) \otimes  \Eis _L(2_{\ell},j_L,m_{j_L})\bigr)\\
&\quad \txt\bigoplus\limits_{k_R\neq \ell_L} \bigl(\Eis _R(2_{k_R},j_R,m_{j_R}) \otimes  
\Eis _L(2_{\ell_L},j_L,m_{j_L})\bigr)\end{array}$
\Ei
where $\widehat M^{2n_\ell}_{T_{R_\oplus}}\equiv \widehat G^{(2n_\ell)}(F^T_{\o\omega _\oplus})$
(resp. $\widehat M^{2n_\ell}_{T_{L_\oplus}}\equiv \widehat G^{(2n_\ell)}(F^T_{\omega _\oplus})$) is a \rl semisheaf over the \rl linear semigroup
$ G^{(2n_\ell)}(F^T_{\o\omega _\oplus})$
(resp. $ G^{(2n_\ell)}(F^T_{\omega _\oplus})$~)
,\\
and are in one-to-one correspondence with the sums of the products of the (truncated) Fourier developments of the cusp forms of weight two.
This leads to evident reducible global bilinear correspondences of which the partially reducible case is explicitly developed here:
\[\begin{array}{rcl}
\Red(\Rep^{(2n)}_{W_{F_{R\times L}}}(W^{2n=2n_1+\cdots+2n_s}_{F_{\o\omega}}\times 
W^{2n=2n_1+\cdots+2n_s}_{F_{\omega}}))
& \!\!\rightarrow\!\! & 
\Red\cusp (\GL_{n=n_1+\cdots+n_s}({F^T_{\o\omega}}\times {F^T_{\omega}}))\\
\searrow && \nearrow\end{array}\]
\[ \mbox{\qquad} \CY ^{2n=2n_1+\cdots+2n_s}_T(\o X_{R}) \times \CY ^{2n=2n_1+\cdots+2n_s}_T(\o X_{L})\]
where:
\Bi
\item $\Red(\Rep^{(2n)}_{W_{F_{R\times L}}}(W^{2n=2n_1+\cdots+2n_s}_{F_{\o\omega}}\times 
W^{2n=2n_1+\cdots+2n_s}_{F_{\omega}}))$ is given by  
$ \widehat G^{(2n=2n_1+\cdots+2n_s)}({F_{\o\omega_\oplus} }\times{F_{\omega_\oplus} })$~;
\vskip 11pt
\item $\Red\cusp (\GL_{n=n_1+\cdots+n_s}({F^T_{\o\omega}}\times {F^T_{\omega}}))$ is given by 
$\bigoplus_{n_\ell}H^{2n_\ell}(\o S^{P_{n=n_1+\cdots+n_s}}_{K_{n_1+\cdots+n_s}},\widehat M^{2n_\ell}_{T_{R_\oplus}}\otimes \widehat M^{2n_\ell}_{T_{L_\oplus}})=
\bigoplus_{n_\ell}(\Eis_R(2n_\ell,j_R,m_{j_R})\otimes \Eis_L(2n_\ell,j_L,m_{j_L}))$~;
\vskip 11pt
\item $\CY ^{2n=2n_1+\cdots+2n_s}_T(\o X_{L})= \txt\bigoplus\limits_{2n_\ell}\CY ^{2n_\ell}_T(\o X_{L})$ with $\CY ^{2n_\ell}_T(\o X_{L})$ a $2n_\ell$-dimensional toroidal compactified semi ``cycle'' over the (compactified) semischeme $\o X_{L}$~.
\Ei
As in the irreducible case, a second kind of reducible bilinear global correspondences of Langlands has been developed in section 4.2 on the equivalents of the  reducible Shimura bisemivarieties.

This last version of the paper was undertaken in order to precise and define the bilinear cohomology theory in section 3.2.
\vskip 11pt

Finally, I would like to thank J. Arthur, P. Cartier, H. Carayol, P. Deligne, M. Harris, G. Henniart, H. Kim, L. Lafforgue, R. Langlands, G. Laumon and many others for helpful advices and conversations.

\section{Global class field concepts, motivic chain bicomplexes and pure Chow bimotives}

\subsection{Global class field concepts}

\subsubsection{Classical background}

The two main challenges of class field theory are the explicit construction of class fields and of reciprocity laws.  Historically, this program proceeds essentially from the fundamental problem in algebraic number theory consisting in describing how an ordinary prime $p$ factorizes into ``primes'' in the ring of integers $\Os_E$ of finite extension $E$ of $\QQ$~.  If we set that $\Os_E=\ZZ(i)$~, then we obtain the following well known:
\vskip 11pt 

\paragraph{Theorem:}{\em
Suppose $p$ is an odd prime.  Then $p$ can be written following $p=n^2+m^2=(n+im)(n-im)$~, $n,m\in\ZZ$ if and only if $p\equiv 1(4)$~.
}
\vskip 11pt

This is an elementary example of the Frobenius automorphisms $\Fr_p$ in the Galois group $G=\Gal(E/\QQ)$~.  Remark that $\Fr_p=id$ when $p$ splits completely into $E$~, i.e. when the ideal it generates in $\Os_E$ factors into distinct prime ideals of $\Os_E$~.

A major objective consists in proving that the splitting properties of $p$ in $E$ depend only  on its residue modulo some fixed modulus $N$~.  This can be achieved if $G=\Gal(E/\QQ)$ is abelian and if $\sigma:G\to\CC^*$ is a homomorphism.  Then, there exists an integer $N>0$ and a Dirichlet character $\chi _\sigma:(\ZZ/N\ \ZZ)^*\to \CC^*$ such that $\sigma(\Fr_p)=\chi _\sigma(p)$ for all primes $p$ unramified in $E$~: this is E. Artin's fundamental reciprocity law of abelian class field theory.  All that was developed very clearly by S. Gelbart in \cite{Gel}.
\vskip 11pt

\subsubsection{How do global algebraic number fields proceed from this classical background?}

An interesting step of this paper 
in relation with section 1.1.1, consists in replacing the finite adele ring by
the infinite places of an algebraic extension (i.e. a splitting (semi)field) $\widetilde F^+$ of a global number (semi)field $k$ of characteristic zero such that:
\Bean
\item the degrees or ranks of the generated ``pseudo-ramified'' completions of $\widetilde F^+$ are integers modulo $N$ where $N$ is the order of the associated global inertia subgroups.

The completions of $\widetilde F_+$  are constructed from  irreducible closed algebraic subsets in such a way that they are isomorphically covered by the corresponding algebraic extensions.

These completions result from an isomorphism of compactifications of the corresponding algebraic extensions.

\item a character $\chi _{v_{j_\delta}}:(\ZZ/N\ \ZZ)^*\to\CC^*$ is associated to every pseudo-ramified completion 
$F^+_{v_{j_\delta}}$ at the $v_{j_\delta}$-th place of $F^+$ such that its decomposition group $D_{v_{j_\delta}}$ 
generates the Dirichlet character 
$\chi _{v_{j_{\delta,m_{j_\delta}}}}$ of the equivalent pseudo-ramified completion $F^+_{v_{j_{\delta,m_{j_\delta}}}}$ following:
\[ D_{v_{j_\delta}}(\chi _{v_{j_\delta}})= \chi _{v_{j_{\delta,m_{j_\delta}}}}\;, \qquad m_{j_\delta}\in\NN\;.\]
\Ee

The point a) deals with the definition of a set of increasing ``pseudo-ramified'' completions 
$F^+_{v_{j_\delta}}$ which  have increasing ranks equal to the associated extension degrees\linebreak 
$[\wt F^+_{v_{j_\delta}}:k]\simeq j_\delta\cdot N$ such that every (pseudo-)ramified completion
 $F^+_{v_{j_\delta}}$ is generated from an irreducible ``central'' completion $F^+_{v^1_{j_\delta}}$ having a rank equal to 
$N$~.

And the point b) refers to the construction  of equivalent pseudo-ramified completions 
$F^+_{v_{j_{\delta,m_{j_\delta}}}}$~.

As a consequence, the residue subfield of $F^+_{v_{j_{\delta}}}$ is defined by 
$F^{+,nr}_{v_{j_{\delta}}}= F^+_{v_{j_\delta}}\big/ F^+_{v^1_{j_\delta}}$~: this implies that this residue subfield $F^{+,nr}_{v_{j_\delta}}$ is a pseudo-unramified completion of $F^+$  whose rank is given by the extension degree 
$[\wt F^{+,nr}_{v_{j_\delta}}:k]\simeq j_\delta$ of the corresponding subfield $\widetilde F^{+,nr}_{j_\delta}$ of $\wt F^+$~, also called a global residue degree.
\vskip 11pt 

The generation of global algebraic extension (semi)fields thus originates from classical algebraic number theory as it will be developed in the following sections and presents some analogy with the construction of local $p$-adic number fields as envisaged in the following section.
\vskip 11pt 

\subsubsection{Classical notions of local number fields}

The classical notions concerning local number fields can be summarized as follows.
\vskip 11pt 

Let $\Os_K$ denote the ring of integers of a finite extension $K$ of $\QQ_p$~.  Its residue field is $k(v_K)=k(\wp_K)=\Os_K/\wp_K$ where $\wp_K$ is the unique maximal ideal of $\Os_K$~.  Let $v_K:K^*\to\ZZ$ be the unique valuation so that the absolute value on $K$ is given by
\[ |\cdot|_K=|\cdot|_{v_K}\qquad \text{with}\quad |x|_K=(\#k(v_K))^{-v_K(x)}\quad \text{for}\quad x\in K^*\;.\]
The number of elements in $k(\wp_K)$ is $q=p^f$ where $f_{v_K}=[k(v_K):\FF_p]$ is the residue degree over $\QQ_p$~.

The ideal $\wp_K\ \Os_K$ of $\Os_K$ has the form $\wp_K^{e_{v_K}}=\widetilde \omega_K^{e_{v_K}}$ where
\Bi
\item $\widetilde\omega_K$ is a uniformiser in $\Os_K$~;
\item $e_{v_K}$ is the ramification degree of $K$ over $\QQ_p$~.
\Ei

Then, we have $[K:\QQ_p]=e_{v_K}\cdot f_{v_K}$ such that $e_{v_K}=[K:\QQ_p]\big/f_{v_K}$~.
\vskip 11pt 

If $K^{nr}$ denotes the maximal unramified extension of $K$~, the inertia subgroup $I_K$ can be defined by:
\[
\begin{CD} \Gal(K^{ac}\big/K)\Big/I_K 
@>\sim>> \Gal(K^{nr}\big/ K) @>\sim>>\Gal(k(v_K)^{ac}\big/k(v_K))\;.\end{CD}\]
\vskip 11pt 

\subsubsection{Infinite places of a global number field of characteristic 0}

Let $k$ be a global number field of characteristic zero and let $\wt F$ denote a finite (algebraically closed) extension of 
$k$ such that $ \widetilde F=\widetilde F_R\cup \widetilde F_L$ is a symmetric splitting field composed of a right and a left 
algebraic extension semifields $\widetilde F_R$ and $\widetilde F_L$ in one-to-one correspondence.  In the complex case, 
$\widetilde F_L$ (resp. 	$\widetilde F_R$~) is the set of complex (resp. conjugate complex) simple roots of a polynomial ring $k[x]$ 
over $k$~.  In the real case, the symmetric splitting field is noted $\widetilde F^+=\widetilde F^+_R\cup \widetilde F^+_L$ where $\widetilde F^+_L$ (resp. $\widetilde F^+_R$~) is the algebraic extension semifield composed of the set of positive (resp. symmetric negative) simple real roots as developed in the appendix.  Remark that $k$ can also be written following $k=k_R\cup k_L$~.

The left and right equivalence classes of the local completions of $\widetilde F^{(+)}_L$ and $\widetilde F^{(+)}_R$ are the left and right real (resp. complex) infinite places of $F^{(+)}_L$ and $F^{(+)}_R$~:\par they are noted $v=\{v_{1_\delta },\cdots,v_{j_\delta},\cdots,v_{r_\delta }\}$ and 
$\o v=\{\o v_{1_\delta },\cdots,\o v_{j_\delta},\cdots,\o v_{r_\delta }\}$ in the real case and
$\omega =\{\omega _1,\cdots,\omega _{j},\cdots,\omega _r\}$ and 
$\o \omega =\{\o \omega _1,\cdots,\o \omega _{j},\cdots,\o \omega _r\}$ in the complex case and are constructed in such a way that each complex place is covered by its real equivalent.
\vskip 11pt 

Recall that the completions of $\widetilde F^{(+)}_L$ and $\widetilde F^{(+)}_R$ at infinite places are defined for the topology by their archimedean absolute values such that the Cauchy sequences in these completions converge (i.e. have a limit).
\vskip 11pt

\subsubsection{Completions of a global number field}

\Bi
\item Let $F_{\omega _j}$ (resp. $F_{\o \omega _j}$~) denote a \lr complex pseudo-ramified completion of $\widetilde F_L$ (resp. $\widetilde F_R$~) in $\omega _j$ (resp. $\o \omega _j$~) and let
$F^+_{v_{j_\delta}}$ (resp. $F^+_{\o v_{j_\delta}}$~) 
 denote a \lr real pseudo-ramified completion of $\widetilde F^+_L$ (resp. $\widetilde F^+_R$~) in 
$v_{j_\delta }$ (resp. $\o v_{j_\delta }$~): 
they  will be assumed to be generated from an ``irreducible'' central completion $F_{\omega ^1_j} $ (resp.  $F_{\o\omega ^1_j} $~) of rank $N\cdot m^{(j)}$ in the complex case and from a irreducible completion $F_{v^1_{j_\delta}} $ (resp.  $F_{\o v^1_{j_\delta}}$~) of rank $N$ in the real case where $m^{(j)}=\sup(m_j)+1$ is the multiplicity of the $j$-th real completion covering its complex equivalent.

So, the ranks (i.e. the Galois extension degrees of the associated extensions) of the complex and real completions 
will be given by the integers modulo $N$~, $N\in\NN$~, where the integer $N$ corresponds to the Artin conductor.
\vskip 11pt

\item In the complex case, the ranks of the complex pseudo-ramified completions will be given by:
\begin{alignat*}{4}
[ F_{\omega _j}:k]
&= (*+j\cdot N)\ m^{(j)}\qquad  \qquad 
& [F_{\o\omega _j}:k]
&= (*+j\cdot N)\ m^{(j)} \;, \end{alignat*}
where $*$ is an integer inferior to $N$~,
in such a way that $[F_{\omega _j}:k]=[ F_{\o \omega _j}:k]=(j\cdot N)\ m^{(j)}$ if $[F_{\omega _j}:k]=[F_{\o\omega _j}:k]=0\mod N$~, $1\le j\le r\le \infty $~.

The ``irreducible'' complex completions have ranks:
\[ [F_{\omega ^1_j}:k]=N\cdot m^{(j)}\;, \qquad \qquad 
 [F_{\o\omega ^1_j}:k]=N\cdot m^{(j)}\;.\]
\vskip 11pt

\item In the real case, the ranks of the pseudo-ramified completions will be:
\begin{alignat*}{4}
[ F^+_{v_{j_\delta }}:k]
&= *+j\cdot N \qquad \qquad 
&[F^+_{\o v_{j_\delta }}:k]
&= *+j\cdot N \;, \end{alignat*}
in such a way that
$
[ F^+_{v_{j_\delta}}:k]
= [ F^+_{\o v_{j_\delta}}:k]
=j\cdot N$
 if $ [ F^+_{v_{j_\delta}}:k]
= [ F^+_{\o v_{j_\delta}}:k]
=0\mod N$~, $1\le j\le r\le \infty $~.

The irreducible complex completions have ranks:
\[ [F^+_{v^1_{j_\delta}}:k]=N\;, \qquad \qquad 
 [F^+_{\o v^1_{j_\delta}}:k]=N\;.\]
\vskip 11pt 

\item The introduction of ``irreducible'' central completions $F_{\omega ^1_j}$ (resp. $F_{\o\omega ^1_j}$~) in the complex case and of irreducible central completions $F^+_{v^1_{j_\delta }}$ (resp.  $F^+_{\o v^1_{j_\delta }}$~) in the real case leads to consider that:
\Bi
\item the complex pseudo-ramified completion $F_{\omega _j}$ (resp. $F_{\o\omega _j}$~), $1\le j\le r$~, can be cut into a set 
of $j$ irreducible equivalent complex completions $F_{\omega _j^{j'}}$~, $1\le j'\le j$ 
(resp.   $F_{\o\omega _j^{j'}}$~) of rank $N\cdot m^{(j)}$~;

\item the real pseudo-ramified completions $F^+_{v_{j_\delta }}$ (resp. $F^+_{\o v_{j_\delta }}$~), $1\le j_\delta \le r_\delta $~, can be cut into a set of $j_\delta $ irreducible equivalent real completions $F^+_{v_{j_\delta }^{j'_\delta }}$~, $1\le j'_\delta \le j_\delta $ (resp.  $F^+_{\o v_{j_\delta }^{j'_\delta }}$~) of rank $N$~.
\Ei
\vskip 11pt

However, remark that the extension degrees characterize more exactly the splitting subfields in bijection with the corresponding completions of $ \widetilde F^{(+)}_L$ and $\widetilde F^{(+)}_R$~.  
\vskip 11pt

\item The introduction of ``irreducible'' central complex left completions 
$F_{\omega ^1_j}$ (resp. right  
$F_{\o\omega ^1_j}$~) of rank $N\cdot m^{(j)}$ at the complex places $\omega _j$ (resp. $\o\omega _j$~), $1\le j\le r$~, and of irreducible real left completions $F^+_{v^1_{j_\delta}}$ (resp. right 
  $F^+_{\o v^1_{j_\delta}}$~) of rank $N$ at the real places $v_{j_\delta}$ (resp. $\o v_{j_\delta}$~), $1\le j_\delta \le r$~, allow to define the global residue subfields of these completions as quotient subfields given by:
\[ F^{nr}_{\omega _j}=F_{\omega _j}\Big/ F_{\omega ^1_j}\qquad \text{(resp.} \quad
 F^{nr}_{\o\omega _j}=F_{\o\omega _j}\Big/ F_{\o\omega ^1_j}\ )\]
in the complex case and by
\[ F^{+,nr}_{v_{j_\delta}}=F^+_{v_{j_\delta}}\Big/ F^+_{v^1_{j_\delta}} \qquad \text{(resp.} \quad
 F^{+,nr}_{\o v_{j_\delta}}=F^+_{\o v_{j_\delta}}\Big/ F^+_{\o v^1_{j_\delta}} \ )\]
in the real case: they are called global residue completions or pseudo-unramified completions of $F^{(+)}_L$ (resp. $F^{(+)}_R$~) characterized by their ranks, or global residue degrees,  given respectively by
\[ [F^{nr}_{\omega _j} :k]= j \cdot m^{(j)}\qquad \text{(resp.} \quad
[ F^{nr}_{\o\omega _j} :k]=j\cdot m[(j) \ )\]
in the complex case and by
\[[ F^{+,nr}_{v_{j_\delta}}:k] =j \qquad \text{(resp.}\quad
 [F^{+,nr}_{\o v_{j}} :k]= {j} \ )\]
in the real case.
\vskip 11pt

\item As a place is an equivalence class of completions, we have to consider a set of complex completions 
$\{F_{\omega _{j,m_j}}\}_{m_j}$~, $m_j\in\NN$~, equivalent to $F_{\omega _j}$~, for all $1\le j\le r$ and characterized by the same rank 	as $F_{\omega _j}$~.

Similarly, a real place $v_{j_\delta}$ will be given by the basic real completion $F_{v_{j_\delta}}$ and by the set of real equivalent completions $\{F^+_{v_{j_{\delta,m_{j_\delta}}}}\}_{m_{j_\delta}}$~, $m_{j_\delta}\in\NN$~, characterized by the same rank as $F_{v_{j_\delta}}$~.

A complex (resp. real) equivalent completion $F_{\omega _{j,m_j}}$ (resp. $F^+_{v_{j,m_{j_\delta}}}$~) is generated from the 
respective basic completion $F_{\omega _j}$ (resp. $F^+_{v_{j_\delta}}$~) by the action of the nilpotent group element
$u_{\omega _{j,m_{j}}}$ (resp.  $u_{v_{j_\delta ,m_{j_\delta}}}$~) following:
\[ u_{\omega _{j,m_j}}F_{\omega _j}=F_{\omega _{j,m_j}}\quad
\text{(resp.}\quad
u_{v_{j_{\delta,m_{j_\delta}}}}F^+_{v_{j_\delta}}=F^+_{v_{j_{\delta,m_{j_\delta}}}}\ )\]

Since only bilinear cases are relevant, it is more exactly the nilpotent group element $u_{j^2;m^2_j}$ acting on $(F_{\o \omega _j}\times F_{\omega _j})$ which must be taken into account.  It generates the equivalent completion 
$(F_{\o \omega _{j;m_j}}\times F_{\omega _{j;m_j}})$ in the complex case according to:
\[ u_{j^2;m^2_j}(F_{\o\omega_j }\times F_{\omega _j})=F_{\o\omega _{j;m_j}}\times F_{\omega _{j;m_j}}\;;\]
in the real case, we should have:
\[u_{j^2_\delta;m_{j^2_\delta}}(F^+_{\o v_{j_\delta}}\times F^+_{v_{j_\delta}})=
F^+_{\o v_{j_\delta};m_{j_\delta}}\times F^+_{v_{j_\delta};m_{j_\delta}}\]
as it will be developed in the next chapter.
\Ei
\vskip 11pt 

\subsubsection{Infinite adele  semirings and semigroups $F_{\omega _\oplus}$~, 
$F^{nr}_{\omega _\oplus}$~, $F^+_{v_\oplus}$ and $F^{+,nr}_{v_\oplus}$}

\Bi
\item So, a \lr infinite ``pseudo-ramified'' adele semiring $\Aa^\infty _{F_\omega }$ 
(resp. $\Aa^\infty _{F_{\o \omega }}$~) can then be introduced by considering 
the product over Archimedean prime  places $j_p$ of $\widetilde F_L$ 
(resp. $\widetilde F_R$~) of the basic completions 
$F_{\omega _{j_p}}$ (resp. $F_{\o\omega _{j_p}}$~) according to
\begin{align*}
\Aa^\infty _{F_\omega } &= \prod_{j_p} F_{\omega _{j_p}}\;, \quad 1\le {j_p}\le r\le\infty\;,\qquad
&\text{(resp.} \quad 
\Aa^\infty _{F_{\o\omega }} &= \prod_{j_p} F_{\o\omega _{j_p}}\ ).\\[11pt]\end{align*}
In the same way, a \lr ``pseudo-unramified'' adele semiring
$\Aa^{nr,\infty }_{F_\omega }$ (resp. $\Aa^{nr,\infty }_{F_{\o \omega }}$~) 
can then be defined from the product of pseudo-unramified ``prime'' completions following:
\begin{align*}
\Aa^{nr,\infty }_{F_\omega } &= \prod_{j_p} F^{nr}_{\omega _{j_p}}\qquad
&\text{(resp.} \quad 
\Aa^{nr,\infty }_{F_{\o\omega }} &= \prod_{j_p} F^{nr}_{\o\omega _{j_p}}\ ).\end{align*}

And, in the real case, left and right pseudo-ramified and pseudo-unramified adele semirings are introduced similarly by:
\begin{align*}
\Aa^\infty _{F^+_v} &= \prod_{j_{\delta _p}} F^+_{v_{j_p}}\;, 
\quad 1\le j_{\delta _p}\le r_p\le\infty\;,\qquad
&\text{(resp.} \quad 
\Aa^\infty _{F^+_{\o v}} &= \prod_{j_{\delta _p}} F^+_{\o v_{j_{\delta _p}}}\ )
\end{align*}
{and by}
\begin{align*}
\Aa^{nr,\infty }_{F^+_v} &= 
\prod_{j_{\delta _p}} F^{+,nr}_{v_{j_{\delta _p}}}
\qquad
&\text{(resp.} \quad 
\Aa^{nr,\infty }_{F^+_{\o v}} 
&= \prod_{j_{\delta _p}} F^{+,nr}_{\o v_{j_{\delta _p}}}\ ).\end{align*}
\vskip 11pt

\item Let 
\begin{align*}
F_\omega  &= \{F_{\omega _1},\cdots,F_{\omega _j},F_{\omega _{j,m_j}},\cdots,F_{\omega _r}\}\\
\text{(resp.} \quad
F_{\o\omega}  &= \{F_{\o\omega _1},\cdots,F_{\o\omega _j},F_{\o\omega _{j,m_j}},\cdots,F_{\o\omega _r}\}\ )\end{align*}
denote the set of pseudo-ramified completions at all Archimedean complex places $\omega _j$ (resp. $\o\omega _j$~).  Then, the direct sum of these pseudo-ramified completions is given by:
\begin{align*}
F_{\omega_\oplus}  &= \bigoplus_j\bigoplus_{m_j}
F_{\omega _{j,m_j}}\;, \qquad 1\le j\le r\le\infty\;, \\
\text{(resp.} \quad
F_{\o\omega_\oplus}  &= \bigoplus_j\bigoplus_{m_j}
F_{\o\omega _{j,m_j}} \ ).\end{align*}
Similarly, if 
\begin{align*}
F^{nr}_\omega  &= \{F^{nr}_{\omega _1},\cdots,F^{nr}_{\omega _j},F^{nr}_{\omega _{j,m_j}},\cdots,F^{nr}_{\omega _r}\}\\
\text{(resp.} \quad
F^{nr}_{\o\omega}  &= \{F^{nr}_{\o\omega _1},\cdots,F^{nr}_{\o\omega _j},F^{nr}_{\o\omega _{j,m_j}},\cdots,F^{nr}_{\o\omega _r}\}\ )\end{align*}
denote the set of corresponding pseudo-unramified completions, their direct sum is given by:
\begin{align*}
F^{nr}_{\omega_\oplus}  &= \bigoplus_j\bigoplus_{m_j}
F^{nr}_{\omega _{j,m_j}}\;, \qquad 1\le j\le r\le\infty\;, \\
\text{(resp.} \quad
F^{nr}_{\o\omega_\oplus}  &= \bigoplus_j\bigoplus_{m_j}
F^{nr}_{\o\omega _{j,m_j}} \ ).\end{align*}
In the real case, if
\begin{align*}
F^{+}_v&= \{F^{+}_{v_{1_\delta }},\cdots,
F^{+}_{v_{j_\delta ,m_{j_\delta }}},\cdots,F^{+}_{v_{r_\delta }}\} \\
\text{(resp.} \quad
F^{+}_{\o v}&= \{F^{+}_{\o v_{1_\delta }},\cdots,
F^{+}_{\o v_{j_\delta ,m_{j_\delta }}},\cdots,F^{+}_{\o v_{r_\delta }}\} \ )\end{align*}
denotes the set of pseudo-ramified completions at all Armimedean real places $v_{j_\delta }$ (resp. $\o v_{j_\delta }$~), 
 the direct sum of these gives:
\begin{align*}
F^{+}_{v_\oplus}  &= \bigoplus_{j_\delta }\bigoplus_{m_{j_\delta }}
F^{+}_{v_{j_\delta ,m_{j_\delta }}} \;, \qquad 1\le j_\delta \le r_\delta \le\infty\;, \\
\text{(resp.} \quad
F^{+}_{\o v_\oplus}  &= \bigoplus_{j_\delta }\bigoplus_{m_{j_\delta }}
F^{+}_{\o v_{j_\delta ,m_{j_\delta }}}   \ ).\end{align*}

Similarly, the direct sum of all pseudo-unramified real completions gives:
\begin{align*}
F^{+,nr}_{v_\oplus}  &= \bigoplus_{j_\delta }\bigoplus_{m_{j_\delta }}
F^{+,nr}_{v_{j_\delta ,m_{j_\delta }}} \;, \qquad 1\le j_\delta \le r_\delta \le\infty\;, \\
\text{(resp.} \quad
F^{+,nr}_{\o v_\oplus}  &= \bigoplus_{j_\delta }\bigoplus_{m_{j_\delta }}
F^{+,nr}_{\o v_{j_\delta ,m_{j_\delta }}}  \ ).\end{align*}
\Ei
\vskip 11pt

Galois groups and Weil groups will now be introduced in the complex case by taking into account that the real case can be handled similarly.  Furthermore, the considered Galois subgroups refer to Galois subgroups of splitting subfields $\widetilde F_{\omega_j} $ in bijection with the corresponding completions $F_{\omega _j}$~.
\vskip 11pt 

\subsubsection{Galois and inertia subgroups}

Let $\Gal^D(\widetilde F_{\omega _j}\big/k)$ (resp. $\Gal^D(\widetilde F_{\o\omega _j}\big/k)$~),
$\Gal(\widetilde F_{\omega _{j,m_j}}\big/k)$ (resp. $\Gal(\widetilde F_{\o\omega _{j,m_j}}\big/k)$~) and\linebreak
$\Gal(\widetilde F_{\omega _{j}}\big/k)$ (resp. $\Gal(\widetilde F_{\o\omega _{j}}\big/k)$~) denote respectively the  Galois subgroups of the basic (non-compact) extension $\widetilde F_{\omega _j}$ (resp. $\widetilde F_{\o\omega _j}$~), of the $m_j$-th equivalent extension $\widetilde F_{\omega _{j,m_j}}$ (resp. $\widetilde F_{\o\omega _{j,m_j}}$~) and of the set of extensions
$\{\widetilde F_{\omega _{j,m_j}}\}$ (resp. $\{\widetilde F_{\o\omega _{j,m_j}}\}$~) including $\widetilde F_{\omega _j}$ (resp. $\widetilde F_{\o\omega _j}$~).

So, we have that:
\[
\Gal(\widetilde F_{\omega _j}\big/ k)
=  \Gal^D(\widetilde F_{\omega _j}\big/k) \bigoplus_{m_j}\Gal (\widetilde F_{\omega _{j,m_j}}\big/k)\;.\]

For the respective pseudo-unramified extension, we should have:
\[ \Gal(\widetilde F^{nr}_{\omega _j}\big/ k)
=  \Gal^D(\widetilde F^{nr}_{\omega _j}\big/k)
\bigoplus_{m_j}\Gal (\widetilde F^{nr}_{\omega _{j,m_j}}\big/k)\;.\]
The Galois subgroup of the ``irreducible'' central extension 
$\widetilde F_{\omega ^1_j}$ having a rank equal to $N\cdot m^{(j)}$ is obviously the global inertia subgroup 
$I^D_{F_{\omega _j}}$ of $\Gal^D(\widetilde F_{\omega _j}\big/k)$ because it can be defined by:
\[ I^D_{F_{\omega _j}}= \Gal^D(\widetilde F_{\omega _j}\big/k) \Big/\Gal^D(\widetilde F^{nr}_{\omega _j}\big/k)\;.\]
Similarly, the global inertia subgroup $I_{F_{\omega _{j,m_j}}}$ of  $\Gal(\widetilde F_{\omega _{j,m_j}}\big/k)$ is defined by:
\[ I_{F_{\omega _{j,m_j}}}= \Gal(\widetilde F_{\omega _{j,m_j}}\big/k) \Big/\Gal(\widetilde F^{nr}_{\omega _{j,m_j}}\big/k)\;.\]
As the global inertia subgroups are of Galois type, they are all isomorphic:
\[\begin{matrix}
I_{F^D_{\omega _1}} &\simeq & \cdots&\simeq&I_{F^D_{\omega _j}} &\simeq& \cdots& \simeq& I_{F^D_{\omega _r}}\;,\\
I_{F^D_{\omega _1,m_1}} &\simeq & \cdots & \simeq & I_{F^D_{\omega _{j,m_j}}}& \simeq & \cdots & \simeq&  I_{F^D_{\omega _r,m_r}}\; ,\end{matrix} \begin{matrix} \qquad 1\le j\le r\;, \\ \mbox{}\end{matrix}\]
which has for consequence that the kernel of the map:
\[\Gal_j: \quad \Gal(\widetilde F_{\omega _j}\big/k) \To \Gal (\widetilde F^{nr}_{\omega _j}\big/k)\]
is a general global inertia subgroup $I_{F_{\omega _j}}$ verifying
\[ I_{F_{\omega _j}}=\Gal(\widetilde F_{\omega _j}\big/k)\Big/\Gal(\widetilde F^{nr}_{\omega _j}\big/k)\;.\]
\vskip 11pt

\subsubsection{Proposition}  

{\em There exists an injective nilpotent morphism
\[ \Ns_{\Gal_{\omega _j}}: \quad \Gal^D(\widetilde F_{\omega _j}\big/k)\To \Gal (\widetilde F_{\omega _j}\big/k)\]
from the Galois subgroup $\Gal^D(\widetilde F_{\omega _j}\big/k)$ of the basic extension $\widetilde F_{\omega _j}$ to the Galois subgroup $\Gal(\widetilde F_{\omega _j}\big/k)$ of the set of extensions $\{ \widetilde F_{\omega _{j,m_j}}\}$~.}
\vskip 11pt 

\bpr According to section 1.1.5, the extension $\widetilde F_{\omega _{j,m_j}}$ is generated from the basic extension $\widetilde F_{\omega _j}$ by the action of the nilpotent group element $u_{\omega _{j,m_j}}$ following:
\begin{equation}
 u_{\omega _{j,m_j}}\widetilde F_{\omega _j}=\widetilde F_{\omega _{j,m_j}}\;.\tag*{\eop}\end{equation}
\vskip 11pt 

\subsubsection{Definition: global Weil groups}  

Referring to section 1.1.7, we have that:
\Bi
\item $\ds\Gal(\wt F^{ac}_L/k)=\txt\bigoplus\limits^r_{j=1}\Gal(\widetilde F_{\omega _j}/k)$~;\vskip 11pt

\item $\ds\Gal(\wt F^{ac}_R/k)=\txt\bigoplus\limits^r_{j=1} \Gal(\widetilde F_{\o\omega _j}/k)$~;\vskip 11pt

\item $\ds\Gal^D(\wt F^{ac}_L/k)=\txt\bigoplus\limits^r_{j=1} \Gal^D(\widetilde F_{\omega _j}/k)$~;\vskip 11pt

\item $\ds\Gal^D(\wt F^{ac}_R/k)=\txt\bigoplus\limits^r_{j=1} \Gal^D(\widetilde F_{\o\omega _j}/k)$~;
\Ei
where $\wt F^{ac}_L$ is the union of all finite abelian extensions $\widetilde F_{\omega_j}$ of $F_L$ in $\o F_L$ 
leading to:
\Bi
\item $\ds\Gal(\wt F^{ac}_R/k) \times \Gal(\wt F^{ac}_L/k) =\txt\bigoplus\limits^r_{j=1} (\Gal(\widetilde F_{\o\omega _j}/k) \times \Gal( \widetilde F_{\omega _j}/k)) $~;\vskip 11pt

\item $\ds\Gal^D(\wt F^{ac}_R/k) \times \Gal^D(\wt F^{ac}_L/k) =\txt\bigoplus\limits^r_{j=1}(\Gal^D(\widetilde F_{\o\omega _j}/k) \times \Gal^D(\widetilde F_{\omega _j}/k)) $~.
\Ei
Similarly, in the pseudo-unramified case, we get:
\Bi
\item $\ds\Gal(\wt F^{nr}_R/k) \times \Gal(\wt F^{nr}_L/k) =\txt\bigoplus\limits^r_{j=1} (\Gal(\widetilde F^{nr}_{\o\omega _j}/k) \times \Gal(\widetilde F^{nr}_{\omega _j}/k)) $~;\vskip 11pt

\item $\ds\Gal^D(\wt F^{nr}_R/k) \times \Gal^D(\wt F^{nr}_L/k) =\txt\bigoplus\limits^r_{j=1}(\Gal^D(\widetilde F^{nr}_{\o\omega _j}/k) \times \Gal^D(\widetilde F^{nr}_{\omega _j}/k)) $~.
\Ei

As in the $p$-adic case, the Weil group \cite{Tat}, \cite{H-T} is the Galois subgroup of the elements inducing on the residue field an integer power of a Frobenius element, we shall assume that, in the characteristic zero case, the Weil group will be the Galois subgroup of the pseudo-ramified extensions characterized by extension degrees $d=0\mod N\cdot m^{(j)}$~.

In this respect, let $\dot{\widetilde F}_{\omega _j}$ (resp. $\dot{\widetilde F}_{\o\omega _j}$~) denote a Galois extension characterized by a degree
\[ [\dot{\widetilde F}_{\omega _j}:k]=(j\cdot N)\   m^{(j)} \qquad \mbox{(resp.}\quad
 [\dot{\widetilde F}_{\o\omega _j}:k]=(j\cdot N )\   m^{(j)} \ ),\]
and let $\dot{\widetilde F}^{nr}_{\omega _j}$ (resp. $\dot{\widetilde F}^{nr}_{\o\omega _j}$~) denote the respective pseudo-unramified extension characterized by the global residue degree:
\[ [\dot{\widetilde F}^{nr}_{\omega _j}:k]=j \cdot m^{(j)} \qquad \mbox{(resp.}\quad
 [\dot{\widetilde F}^{nr}_{\o\omega _j}:k]=j \cdot m^{(j)} \ ).\]
The sum of the Galois subgroups of such extensions is then given by:
\begin{align*}
\Gal  (\dot{\wt F}^{ac}_{L}/k)&= \txt\bigoplus\limits^r_{j=1}\Gal( \dot{\widetilde F}_{\omega _j}/k)=W^{ab}_{F_L}\\[11pt]
\mbox{(resp.}\quad
\Gal ( \dot{\wt F}^{ac}_{R}/k)&= \txt\bigoplus\limits^r_{j=1}\Gal( \dot{\widetilde F}_{\o\omega _j}/k)=W^{ab}_{F_R}\ )
\end{align*}
and corresponds to the Weil global group $W^{ab}_{F_L}$ (resp. $W^{ab}_{F_R}$~).

The product, right by left, of these Weil groups thus is:
\[ W^{ab}_{F_R} \times W^{ab}_{F_L}
= \Gal  (\dot{\wt F}^{ac}_{R}/k)\times \Gal  (\dot{\wt F}^{ac}_{L}/k)
\subset \Gal  ({\wt F}^{ac}_{R}/k)\times \Gal  ({\wt F}^{ac}_{L}/k)\;.\]
\vskip 11pt

\subsubsection{Remark concerning the pseudo-ramified extensions taken into account the Weil global groups}

As it was mentioned in section 1.1.9, the Weil global groups are the Galois subgroups of the pseudo-ramified extensions characterized by extension degrees $d=0\mod N$~: they were noted 
$\dot{\widetilde F}_{\omega _j}$ and $\dot{\widetilde F}_{\o\omega _j}$~, $1\le j\le r\le \infty $~.

Since the Langlands program concerns the $n$-dimensional representations of Weil global groups in bijection with the corresponding automorphic representations, we shall be essentially interested by  the completions corresponding to these pseudo-ramified extensions $\dot{\widetilde F}_{\omega _{j,m_j}}$ and $\dot{\widetilde F}_{\o\omega _{j,m_j}}$~.  In order to simplify the notations, we shall consider, until the end of chapter 4, these restricted completions with the notations introduced in section 1.1.6 for the general case.
\vskip 11pt

Let us now introduce the Suslin-Voevodsky motivic bicomplexes and the Chow bimotives which will be especially used in chapter 4.

{\bf The proposed treatment shows how it is possible to construct enough algebraic cycles\/}.
\vskip 11pt

\subsection{Suslin-Voevodsky motivic bicomplexes}

Let $F_\omega $ (resp. $F_{\o\omega }$~) denote the set of the finite completions of $F_L$ (resp. $F_R$~) at the set of infinite places $\omega $ (resp. $\o\omega $~).
\vskip 11pt 

\subsubsection{Definition: CW semicomplexes}

Let $T_L$ (resp. $T_R$~) be a \lr topological semispace over $F_\omega$  (resp. $F_{\o\omega }$~) restricted respectively to the upper (resp. lower) half space.  $T_L$ (resp. $T_R$~) is identified with a \lr  CW complex having a partition in \lr closed cells noted $\CEL_L(F_{\omega _j})$ (resp. $\CEL_R(F_{\o\omega _j})$~) such that every \lr closed cell $\CEL_L(F_{\omega _j})$ (resp. $\CEL_R(F_{\o\omega _j})$~) be defined over a \lr place $\omega _j$ (resp. $\o\omega _j$~).
Furthermore, it will be assumed that every \lr reducible closed cell
 $\CEL_L(F_{\omega _j})$ (resp. $\CEL_R(F_{\o\omega _j})$~) decomposes into a set of one-dimensional irreducible subcells $C_L(F_{\omega _j})$ (resp. $C_R(F_{\o\omega _j})$~).
Each \lr CW complex $T_L$ (resp. $T_R$~) is isomorphic to an object of the category $Sm_L(k)$ (resp. $Sm_R(k)$~) of the \lr smooth quasi-projective semivarieties over $k$~.
The Eilenberg-MacLane topological semispace generated by the \lr  CW complex $T_L$ (resp. $T_R$~) is the free topological abelian semigroup generated by $T_L$ (resp. $T_R$~), i.e. the semigroup completion $L^{\text{top}}[T_L]$ (resp. $L^{\text{top}}[T_R]$~) of the topological commutative monoid $\txt\bigsqcup\limits_{i_\ell \ge 0}SP^{i_\ell }(T_L)$ (resp. $\txt\bigsqcup\limits_{i_\ell \ge 0}SP^{i_\ell }(T_R)$~), where $SP^{i_\ell }(T_L)$ (resp. $SP^{i_\ell }(T_R)$~) denotes the $i_\ell $-th symmetric product of $T_L$ (resp. $T_R$~).
\vskip 11pt

\subsubsection{Definition: Suslin-Voevodsky motivic presheaf}  

Let $\Delta _L^{2n_\ell }$ (resp. $\Delta _R^{2n_\ell }$~) 
denote a \lr real topological $2n_\ell $-simplex, $2n_\ell\in\NN$~, 
i.e. a closed subscheme in the affine $(2n_\ell +2)$-semispace $\Aa_{k }^{2n_\ell +2}$ 
(resp. $\Aa_{k}^{2n_\ell +2}$~) and let  $\Delta ^\pt_L$ (resp. $\Delta ^\pt_R$~) 
be a cosimplicial object in $Sm_R(k )$ (resp. $Sm_L(k)$~) from the collection of the $\Delta _R^{2n_\ell }$ (resp. 
$\Delta _L^{2n_\ell }$~).
A presheaf of complexes on $Sm_L(k )$ (resp. $Sm_R({k })$~) is a functor from $Sm_L(k )$ (resp. $Sm_R{k }$~) to the \lr chain complexes of abelian semigroups.
On the other hand, let $X_L^{sv}$ (resp. $X_R^{sv}$~) denote a \lr Suslin-Voevodsky smooth semischeme of real dimension $2\ell $ 
isomorphic to the \lr topological semispace $T_L$ (resp. $T_R$~) over $F_\omega $ (resp. $F_{\o\omega }$~) decomposing into closed cells $\CEL_L(F_{\omega _j})$ (resp. $\CEL_R(F_{\o\omega _j})$~) which are also isomorphic to their algebraic analogues $\CEL_L(\wt F_{\omega _j})$ (resp. $\CEL_R(\wt F_{\o\omega _j})$~).
Fix $2n_\ell=i_\ell\times 2\ell$~.
Then, a Suslin-Voevodsky motivic \lr presheaf of $X_L^{sv}$ (resp. $X_R^{sv}$~) on $Sm_L(k )$ (resp. $Sm_R(k)$~), denoted $\u C_*(X_L^{sv})$ (resp. $\u C_*(X_R^{sv})$~) and called in short a Suslin-Voevodsky \lr motive, is any functor from $(X_L^{sv})$ (resp. $(X_R^{sv})$~) to the \lr chain complex associated to the \lr abelian semigroup $\txt\bigsqcup\limits_{i_\ell }\Hom_{Sm_{L/k}}(\Delta ^\pt_L,SP^{i_\ell }(X_L^{sv}))$ (resp. $\txt\bigsqcup\limits_{i_\ell }\Hom_{Sm_R/k}(\Delta ^\pt_R,SP^{i_\ell }(X_R^{sv}))$~).
\vskip 11pt

\subsubsection{Definition: reducibility of Suslin-Voevodsky presheaf}  

If we take into account the 
decomposition of the \lr smooth semischeme $X_L^{sv}$ (resp. $X_R^{sv}$~) in closed subsemischemes isomorphic to algebraic cells 
$\CEL_L( F_{\omega _j})$ (resp. $\CEL_R(F_{\o\omega _j})$~, then the Suslin-Voevodsky \lr motive $\u C_*(X_L^{sv})$ (resp. $\u C_*(X_R^{sv})$~) is the functor from $X_L^{sv}$ (resp. $X_R^{sv}$~) to 
\[\txt\bigsqcup\limits_{i_\ell }  \txt\bigsqcup\limits_{j}\Hom_{Sm_L/k}(\Delta ^\pt_L,SP^{i_\ell }_L(X_L^{sv}[j]))\qquad \text{(resp.} \quad \txt\bigsqcup\limits_{i_\ell }  \txt\bigsqcup\limits_{j}\Hom_{Sm_R/k}(\Delta ^\pt_L,SP^{i_\ell }_L(X_R^{sv}[j]))\ )\]
where $SP^{i_\ell }(X_L^{sv}[j])$ (resp. $SP^{i_\ell }(X_R^{sv}[j])$~) denotes the $j$-th equivalence class of the $i_\ell $-th symmetric product of $X_L^{sv}$ (resp. $X_R^{sv}$~) associated with the \lr place $\omega _j$ (resp. $\o\omega_j$~) of the completion of the  semifield $\widetilde F_L$ (resp. $\widetilde F_R$~).
On the other hand, if we consider the decomposition of the closed cells into one-dimensional 
irreducible affine curves $C_L(j)$ (resp. $C_R(j)$~), then the Suslin-Voevodsky 
\lr motive $\u C_*(X_L^{sv})$ (resp. $\u C_*(X_R^{sv})$~) will be represented by
\begin{align*}
&\txt\bigsqcup\limits_{i_\ell }  \txt\bigsqcup\limits_j\txt\bigsqcup\limits_{m =1}\Hom_{Sm_L/k}
(\Delta ^\pt_L,C_m(SP^{i_\ell }(X_L^{sv}[j])))\\
\text{(resp.} \qquad & \txt\bigsqcup\limits_{i_\ell }  
\txt\bigsqcup\limits_j\txt\bigsqcup\limits_{m=1}\Hom_{Sm_R/k}(\Delta ^\pt_R, C_m(SP^{i_\ell }(X_R^{sv}
[j])))\ )\end{align*}
where $\{C_m (SP^{i_\ell }(X_L^{sv}[j]))\}^{2n_\ell }_{m=1}$ (resp. $\{C_m (SP^{i_\ell }(X_R^{sv}[j]))\}^{2n_\ell }_{m=1}$~) is a set of   affine curves, having multiplicities $m_{j_\ell }$~.  
Finally, if $2n_\ell=i_\ell\times2 \ell$~, let $L_{2n_\ell ;L}$ (resp. $L_{2n_\ell;R}$~) denote the restriction of $\u C_*(X_L^{sv})$ (resp. $\u C_*(X_R^{sv})$~) to its $2n_\ell $-th element and let $Z_L(2n_\ell )$ (resp. $Z_R(2n_\ell )$~) denote the $2n_\ell $-th \lr Suslin-Voevodsky submotive given by $\u C_*(L_{2n_\ell ;L})[-2n_\ell ]$ (resp. $\u C_*(L_{2n_\ell ;R})[-2n_\ell ]$~) which corresponds to the $2n_\ell $-th desuspension of $\u C_*(L_{2n_\ell ;L})$ (resp. $\u C_*(L_{2n_\ell ;R})$~) (the desuspension sends objects of homological degree $(2n_\ell +2)$ to degree $2n_\ell $~) \cite{Mor}.
\vskip 11pt

\subsubsection{Presheaf with transfers}  

$\u C_*(X_L^{sv})$ (resp. $\u C_*(X_R^{sv})$~) has the property to be a presheaf with transfer, i.e. has the correspondence property of motives.  A \lr correspondence between the smooth \lr semischemes $X_L^{sv}$ (resp. $X_R^{sv}$~) and $Y_L^{sv}$ (resp. $Y_R^{sv}$~), noted $\Corr(X_L^{sv},Y_L^{sv})$ (resp. $\Corr(X_R^{sv},Y_R^{sv})$~), is the free abelian semigroup on the set of closed irreducible subsemischemes $XY_L^{sv}$ (resp. $XY_R^{sv}$~) of $X_L^{sv}\times_{k} Y_L^{sv}$ (resp. $X_R^{sv}\times_{k} Y_R^{sv}$~) \cite{Voe} endowed with the projection from $XY_L^{sv}$ (resp. $XY_R^{sv}$~) on one of the irreducible components of $X_L^{sv}$ or of $Y_L^{sv}$ (resp. $X_R^{sv}$ or of $Y_R^{sv}$~).
\vskip 11pt

\subsubsection{Bilinear correspondence on Suslin-Voevodsky semischemes}  

As left and right semischemes $X_L^{sv}$ and $X_R^{sv}$ are taken into account, we have to introduce a new kind of correspondence called bilinear left-right (resp. right-left) correspondence between the \lr smooth semischeme $X_L^{sv}$ (resp. $X_R^{sv}$~) and the \rl smooth semischeme $Y_R^{sv}$ (resp. $Y_L^{sv}$~).  Then, a bilinear left-right (resp. right-left) correspondence will be given by $\Corr(Y_R^{sv},X_L^{sv})$ (resp. $\Corr(Y_L^{sv},X_R^{sv})$~) and defined as the free abelian bilinear semigroup on the set of closed irreducible subsemischemes $XY_{R\times L}^{sv}$ (resp. $XY_{L\times R}^{sv}$~) of $Y_{R}^{sv}\times X_L^{sv}$ (resp. $Y_{L}^{sv}\times X_R^{sv}$~).
This leads to introduce a \lr Suslin-Voevodsky motivic bilinear presheaf $\u C_*(Y_R^{sv}\times X_L^{sv})$ (resp. $\u C_*(Y_L^{sv}\times X_R^{sv})$~) from $Y_R^{sv}\times X_L^{sv}$ (resp. $Y_L^{sv}\times X_R^{sv}$~) to the  right-left (resp.  left-right) chain complex associated to the bilinear semigroup
\begin{align*}
&\txt\bigsqcup\limits_{i_\ell }\Hom_{Sm_R\otimes Sm_L}(\Delta ^\pt\RL,SP^{i_\ell }
(\Corr(Y_R^{sv},X_L^{sv}))\\
\text{(resp.}\qquad & \txt\bigsqcup\limits_{i_\ell }\Hom_{Sm_R\otimes Sm_L}(\Delta ^\pt\LR,SP^{i_\ell }(\Corr(Y_L^{sv},X_R^{sv}))\ ).\end{align*}
\vskip 11pt

\subsection{Chow motives}
\subsubsection{Equivalence classes of algebraic cycles}  

Let $X_L$ (resp. $X_R$~) denote a \lr smooth semischeme of complex dimension $n$ isomorphic to an algebraic linear semigroup of the same dimension over the extensions of a field $k$ of characteristic $0$ such that
\Bi
\item $2n\ge2 \ell $ (~$2\ell $ is the real dimension of the Suslin-Voevodsky semischeme $X_L^{sv}$ (resp. $X_R^{sv}$~);
\item $2n=\sum\limits_\ell i_\ell \times 2\ell =\sum\limits_\ell 2n_\ell$~;
\item $i$ is an  integer corresponding to the dimension $2n_\ell$ used in section 1.2;

so, we have that: $i\equiv 2n_\ell$~.
\Ei

Let $\hZ ^i(X_L)$ (resp. $\hZ ^i(X_R)$~) be the semigroup of algebraic semicycles 
$\CY^i(X_L)$ (resp. $\CY^i(X_R)$~) of codimension $i$
on $X_L$ (resp. $X_R$~).
As the \lr semischeme $X_L$ (resp. $X_R$~) is isomorphic to an algebraic linear semigroup of the same dimension decomposing into conjugacy classes labeled by the integer $j$ (see the introduction and section 2.4),
 it can be assumed that the semigroup of algebraic semicycles $\hZ^i(X_L)$ (resp. $\hZ^i(X_R)$~) is partitioned into equivalence classes corresponding to the set of \lr places $\omega _1,\cdots,\omega _r$ (resp. $\o\omega _1,\cdots,\o\omega _r$~).  So, let $\hZ^i(X_L[j])$ (resp. $\hZ^i(X_R[j])$~) denote the $j$-th equivalence class of the semigroup of algebraic semicycles 
$\CY^i(X_L[j])$ (resp. $\CY^i(X_R[j])$~) 
associated to the place $\omega _j$ (resp. $\o\omega _j$~) of $ F_L$ (resp. $ F_R$~).
\vskip 11pt

\subsubsection{Definition: rational equivalence}  

On the other hand, the \lr semicycle $\CY^i(X_L[j])$ (resp. $\CY^i(X_R[j])$~) must be rationally equivalent to zero, i.e. there must exist rational functions $f_{\CY^i_L}$ (resp. $f_{\CY^i_R}$~) on irreducible subschemes of $X_L$ (resp. $X_R$~) such that $\CY^i(X_L[j])$ (resp. $\CY^i(X_R[j])$~) can be decomposed into a formal sum of (Weil) divisors on $f_{\CY^i_L}$ (resp. $f_{\CY^i_R}$~) \cite{Mur}.
Let $\hZ^i_{\text{rat}}(X_L[j])$ (resp. $\hZ^i_{\text{rat}}(X_R[j])$~) denote the $j$-th equivalence class of the semigroup of algebraic semicycles of codimension $i$ rationally equivalent to zero, i.e. $\hZ_{\text{rat}}^i(X_L[j])=\{\CY^i(X_L[j])\in\hZ^i(X_L[j])\mid \CY^i(X_L[j])$ be rationally equivalent to zero$\}$ (idem for\linebreak $\hZ^i_{\text{rat}}(X_R[j])$~).  Thus, we have that 
\[\CH^i(X_L)=\F{\hZ^i(X_L)}{\hZ^i_{\text{rat}}(X_L)}\qquad  \text{(resp.} \quad \CH^i(X_R)=\F{\hZ^i(X_R)}{\hZ^i_{\text{rat}}(X_R)}\;)\]
 is the $i$-th \lr Chow semigroup of $X_L$ (resp. $X_R$~) and that
\[\CH^i(X_L[j])=\F{\hZ^i(X_L[j])}{\hZ^i_{\text{rat}}(X_L[j])}\qquad \text{(resp.} \quad  \CH^i(X_R[j])=\F{\hZ^i(X_R[j])}{\hZ^i_{\text{rat}}(X_R[j])}\;)\]
 is the $j$-th equivalence class of the $i$-th \lr Chow semigroup.
\vskip 11pt

\subsubsection{Lemma}  

{\em If a numerical equivalence is considered for the $i$-th \lr Chow semigroup 
$\CH^i(X_L)$ (resp. $\CH^i(X_R)$~), then it can be partitioned according to: 
\begin{align*}
\CH^i(X_L)&=
\txt\bigoplus\limits_{m =1}^i\txt\bigoplus\limits_j \CH^1_m (X_L[j])\\
 \text{(resp.} \quad \CH^i(X_R)&=\txt\bigoplus\limits_{m =1}^i\txt\bigoplus\limits_j \CH^1_m (X_R[j])\;)\end{align*}
 where $\CH^1_m (X_L[j])$ (resp. $\CH^1_m (X_R[j])$~) is the $j$-th equivalence class of the one-dimensional simple \lr Chow semigroup $\CH^1_m (X_L)$ (resp. $\CH^1_m (X_R)$~).}
\vskip 11pt

\bpr an algebraic semicycle $\CY^i(X_L)\in\hZ^i(X_L)$ (resp. $\CY^i(X_R)\in\hZ^i(X_L)$~) is said to 
be numerically equivalent to zero if, for all algebraic semicycles $\CY^{2n-i}(X_L)\in\hZ^{2n-i}(X_L)$ 
(resp. $\CY^{2n-i}(X_R)\in\hZ^{2n-i}(X_R)$~), the intersection number $\#(\CY^i(X_L)\cdot \CY^{2n-i}(X_L))=0$ (idem for right part).\newline
According to U. Jannsen \cite{Jan2}, the category of motives is a semisimple abelian category if there is a numerical equivalence on the cycles of the Chow semigroup $\CH^i(X_L)$\linebreak (resp. $\CH^i(X_R)$~) of the motives.  So, if there is a numerical equivalence on the algebraic semicycles $\CY^i(X_L)$ (resp. $\CH^i(X_R)$~), the $i$-th Chow semigroup can decompose as a direct sum of simple one-dimensional semigroups $\CH^1_m (X_L)$ (resp. $\CH^1_m (X_R)$~).
So, the numerical equivalence on algebraic cycles is stronger than rational equivalence and implies it.\epr
\vskip 11pt

\subsubsection{Definition: reducibility of algebraic semicycles}  

A \lr algebraic semicycle thus decomposes following: \begin{align*}
&\CY^i(X_L)=  \txt\bigoplus\limits_{j}\txt\bigoplus\limits_{m=1}^i \txt\bigoplus\limits_{m_{j_m}}\CY^1_m(X_L[j])\\ 
\text{(resp.}\qquad & \CY^i(X_R)=  \txt\bigoplus\limits_{j}\txt\bigoplus\limits_{m=1}^i \txt\bigoplus\limits_{m_{j_m}}\CY^1_m(X_R[j])\ )
\end{align*} where $m_{j_m}$ denotes the equivalent representatives of the $m$-th one-dimensional semicycle\linebreak 
$\CY_m^1(X_L[j])$ (resp. $\CY_m^1(X_R[j])$~).
\vskip 11pt

\subsubsection{Proposition}  

{\em Let $\CH(X_L)=\txt\bigoplus\limits_{i=1}^{2n} \CH^i(X_L)$ (resp. $\CH(X_R)=\txt\bigoplus\limits_{i=1}^{2n} \CH^i(X_R)$~) be the \lr Chow semiring.\newline
If numerical equivalence is considered for the algebraic semicycles of all codimensions, the Chow semiring $\CH(X_L)$ (resp. $\CH(X_R)$~) develops according to:
\begin{align*}
&\CH(X_L)= \txt\bigoplus\limits_{i=1}^{2n} \txt\bigoplus\limits_{j}\txt\bigoplus\limits_{m=1}^i \CH^1_{i;m}(X_L[j])\\ 
\text{(resp.}\qquad 
& \CH(X_R)= \txt\bigoplus\limits_{i=1}^{2n} \txt\bigoplus\limits_{j}\txt\bigoplus\limits_{m=1}^i \CH^1_{i;m}(X_R[j])\ )\end{align*} 
and has for elements:
\begin{align*}
&\CY(X_L)= \txt\bigoplus\limits_{i=1}^{2n}  \txt\bigoplus\limits_{j}\txt\bigoplus\limits_{m=1}^i\txt\bigoplus\limits_{m_{j_m}} \CY^1_{i;m}(X_L[j])\\ 
\text{(resp.}\qquad 
&\CY(X_R)= \txt\bigoplus\limits_{i=1}^{2n}  \txt\bigoplus\limits_{j}\txt\bigoplus\limits_{m=1}^i\txt\bigoplus\limits_{m_{j_m}} \CY^1_{i;m}(X_R[j])\ ).\end{align*}
$\CH(X_L)$ (resp. $\CH(X_R)$~) is thus a semisimple semiring (see  {\em \cite{D-M}, \cite{Ram}\/} for a recent literature on this subject).}\vskip 11pt

\subsubsection{Definition}  

A {\em pure Chow motive\/} restricted to the left or right case considered here is commonly defined as being a pair $(X_L,\Corr^0(X_L,X_L))$ (resp. $(X_R,\Corr^0(X_R,X_R))$~) 
consisting in a \lr smooth semischeme $X_L$ (resp. $X_R$~) having a decomposition in the \lr Chow semiring $\CH(X_L)$ (resp. $\CH(X_R)$~) as developed in Proposition 1.3.5 and a \lr correspondence \cite{Mur} $\Corr^0(X_L,X_L)=\CH^0(X_L,X_L)$ (resp. $\Corr^0(X_R,X_R)=\CH^0(X_R,X_R)$~) which is the set of \lr projectors of $X_L$ (resp. $X_R$~).
\vskip 11pt

\subsubsection{Proposition} 

{\em The category of \lr Chow semimotives is additive, semisimple and graded by weights.}
\vskip 11pt

\bpr results directly from Proposition 1.3.5.
\vskip 11pt

\subsubsection{Definition: product of Chow semirings}  

As in 1.2.5, a bilinear correspondence $\Corr(X_R,X_L)$ can be introduced and defined by\linebreak $\Corr(X_R,X_L)=\CH(X_R\times X_L)\simeq \CH(X_R)\times \CH(X_L)$~.  This leads to introduce the product between a right and a left Chow semiring which can be developed according to:
\begin{align*}
\CH(X_R)\times \CH(X_L)
= &\left( \txt\bigoplus\limits_i \txt\bigoplus\limits_j \txt\bigoplus\limits_{m_R=1}^i \CH^1_{i;m_R}(X_R[j])\right)
\times \left( \txt\bigoplus\limits_i \txt\bigoplus\limits_j \txt\bigoplus\limits_{i;m_L}^i \CH^1_{i;m_L}(X_L[j])\right)\\
=  &\txt\bigoplus\limits_i \txt\bigoplus\limits_j \txt\bigoplus\limits_{m_R=m_L=m} \CH^1_{i;m}(X_R[j]) \times \CH^1_{i;m}(X_L[j])\\
&\txt\bigoplus\limits_i \txt\bigoplus\limits_j \txt\bigoplus\limits_{m_R\neq m_L} \CH^1_{i;m_R}(X_R[j]) \times \CH^1_{i;m_L}(X_L[j])\end{align*}
where the first direct sum gives rise to a bilinear intersection pairing:
\[\CH^1_{i;m}(X_R[j]) \times \CH^1_{i;m}(X_L[j])\to \CH^1_{i;m}(X_{R\times L}[j]) \] while the second direct sum yields a bilinear mixed intersection pairing:
\[\CH^1_{i;m_R}(X_R[j]) \times \CH^1_{i;m_L}(X_L[j])\to \CH^0_{i;m_R;m_L}(X_{R\times L}[j])\;.\]
\vskip 11pt

\subsubsection{Proposition} 

{\em Let $\u C_*(X_L^{sv})$ (resp. $\u C_*(X_R^{sv})$~) be the \lr Suslin-Voevodsky reducible motivic presheaf represented by \begin{align*}
 & \txt\bigsqcup\limits_{i_\ell }\txt\bigsqcup\limits_{j}\txt\bigsqcup\limits_{m}\Hom_{Sm_L/k}(\Delta ^\pt_L,C_m(SP^{i_\ell}(X_L^{sv}[j])))\\
 \text{(resp.}\quad &\txt\bigsqcup\limits_{i_\ell }\txt\bigsqcup\limits_{j}\txt\bigsqcup\limits_{m}\Hom_{Sm_R/k}(\Delta ^\pt_R,C_m(SP^{i_\ell}(X_R^{sv}[j])))\ )\end{align*}
 and let 
$\CH(X_L)=\txt\bigoplus\limits_i \txt\bigoplus\limits_j \txt\bigoplus\limits_{m} \CH^1_{i;m}(X_L[j]) $ (resp. $\CH(X_R)=\txt\bigoplus\limits_i \txt\bigoplus\limits_j \txt\bigoplus\limits_{m} \CH^1_{i;m}(X_R[j]) $~) be the \lr reducible Chow semiring of $X_L$ (resp. $X_R$~).\newline
Then, we can envisage the maps:
\begin{align*}
&\Mot_L:\quad \u C_*(X_L^{sv})\overset{\sim}{\To} (X_L,\Corr^0(X_L,X_L))\\
 \text{(resp.} \quad &Mot_R:\quad \u C_*(X_R^{sv})\overset{\sim}{\To}(X_R,\Corr^0(X_R,X_R))\;)\end{align*}
which are bijective between the pure \lr Suslin-Voevodsky semimotive and the corresponding pure \lr Chow semimotive $(X_L,\Corr^0(X_L,X_L))$ (resp. $(X_R,\Corr^0(X_R,X_R))$~).}
\vskip 11pt

\bpr this is a direct consequence of the developments of this chapter connecting:
\Bi
\item the dimension $2n$ of $X_L$ (resp. $X_R$~) with the dimension $2\ell $ of $X_L^{sv}$ 
(resp. $X_R^{sv}$~) by the relation $2n=\sum\limits_\ell i_\ell \times2 \ell $~;
\item the reducibility of $\u C_* (X_L^{sv})$ (resp. $\u C_* (X_R^{sv})$~) to the reducibility 
of $\CH(X_L)$ (resp. $\CH(X_R)$~).\epr
\Ei\vskip 11pt

\section{Nonabelian global class field concepts based on the representation of an algebraic bilinear semigroup}

\subsection{The general bilinear algebraic semigroup}

\subsubsection{From the abelian global class field theory towards a nonabelian global class field theory: a summary}

In section 1.1, packets of left and right  equivalent complex pseudo-ramified completions associated respectively with the left and right algebraic extension semifields $\widetilde F_L$ and $\widetilde F_R$ were introduced so that:
\Bean
\item the $j$-th packet of the \lr complex pseudo-ramified completions is composed of the basic completion $F_{\omega _j}$ (resp. $F_{\o\omega _j}$~) and of the equivalent completions $F_{\omega _{j,m_j}}$ (resp. $F_{\o\omega _{j,m_j}}$~): these completions are  characterized by a rank given by the extension degree of the associated extension:
\begin{align*}
[\wt F_{\omega _j}:k] =*+(j\cdot  N)\ m^{(j)}\;, \\
\text{(resp.} \quad
[\wt F_{\o\omega _j}:k] =*+(j\cdot  N)\ m^{(j)}\;);\end{align*}

but, as mentioned in sections 1.1.9 and 1.1.10, the  completions with ranks$[\wt F_{\omega _j}:k]=(j \cdot N)\  m^{(j)}$
(resp. $[\wt F_{\o \omega _j}:k]=(j\cdot N)\  m^{(j)} $~) and $[\wt F_{\omega _{j,m_j}}:k]=(j \cdot N)\  m^{(j)}$
(resp. $[\wt F_{\o\omega _{j,m_j}}:k]=(j\cdot N)\  m^{(j)}$~) will only be taken into account.

\item as the basic completion $F_{\omega _j}$ (resp. $F_{\o\omega _j}$~) is generated from an irreducible central completion $F_{\omega^1 _j}$ (resp.   $F_{\o\omega^1 _j}$~), the basic completions $F_{\omega _j}$ (resp. $F_{\o\omega _j}$~) and the equivalent completions $F_{\omega _{j,m_j}}$ (resp. $F_{\o\omega _{j,m_j}}$~) will be respectively decomposed into $j$ irreducible equivalent subcompletions $F_{\omega^{j'} _{j}}$~, $1\le j'\le j$~, 
(resp.  $F_{\o\omega^{j'} _{j}}$~) and $F_{\omega^{j'} _{j,m_j}}$ (resp. $F_{\o\omega^{j'} _{j,m_j}}$~) or rank $N\cdot m^{(j)} $~.	

\item a character $\chi _{\omega _j}$ (resp. $\chi _{\o\omega _j}$~) is associated to each basic completion $F_{\omega _j}$ (resp. $F_{\o\omega _j}$~); and an equivalent Dirichlet character $\chi _{\omega _{j,m_j}}$ (resp. $\chi _{\o\omega _{j,m_j}}$~) corresponds the equivalent completion $F_{\omega _{j,m_j}}$ (resp. $F_{\o\omega _{j,m_j}}$~).

\item the \lr $j$-th packet of complex pseudo-ramified extensions corresponds to the $j$-th \lr place $\omega _j$ (resp. $\o\omega _j$~) of $F_L$ (resp. $F_R$~).

\item a \lr infinite adele semiring was constructed from the tower $F_\omega $ (resp. $F_{\o\omega }$~) of these packets of equivalent completions characterized by increasing ranks according to:
\[ \Aa^\infty  _{F_\omega }=\prod_{j_p=1}F_{\omega _{j_p}}
\qquad \text{(resp.}\quad 
 \Aa^\infty  _{F_{\o\omega }}=\prod_{j_p=1}F_{\o\omega _{j_p}}\;),\]
where:
\Bi
\item $F_{\omega} =\{F_{\omega _1},\cdots,F_{\omega _{j,m_j}},\cdots,F_{\omega_r}\}$ (resp. $F_{\o\omega} =\{F_{\o\omega _1},\cdots,F_{\o\omega _{j,m_j}},\cdots,F_{\o\omega_r}\}$~);
\item $F_{\omega_\oplus} =
\bigoplus_j\bigoplus_{m_j}F_{\omega _{j,m_j}}$ (resp. $F_{\o\omega_\oplus} =
\bigoplus_j\bigoplus_{m_j}F_{\o\omega _{j,m_j}}$~) (see section. 1.1.6).
\Ei
\Ee
\vskip 11pt

\Bi
\item So, the set of packets of \lr complex pseudo-ramified extensions covering the associated completions is a \lr affine semigroup  $\SS^1_L$ (resp. $\SS^1_R$~) whose complex fibers are one-dimensional and to which a ``complex'' Picard semigroup $\Pic({\wt F_\omega })$ (resp. $\Pic({\wt F_{\o\omega }})$~) 
corresponds: $\Pic({\wt F_\omega })$ (resp. $\Pic({\wt F_{\o\omega }})$~) thus is the set of $r$ isomorphism classes of finitely generated semimodules of complex dimension 1 over ${\wt F_\omega }$ (resp. ${\wt F_{\o\omega }}$~).

\vskip 11pt 

\item These few considerations relative to complex pseudo-ramified extensions (real pseudo-ramified extensions can be handled similarly) deal with the abelian class field theory.  Now, the $n$-dimensional global Langlands program is based on nonabelian global class field theory \cite{Lan2}.  The challenge then consists in constructing the $n$-dimensional analog of the affine semigroup  $\SS^1_L$ (resp. $\SS^1_R$~) such that we get an injective  homomorphism:
\begin{align*}
\sigma _L : \quad W^{ab}_{F_L}&\To \GL_n({\wt F_{\omega_\oplus} })\\
\text{(resp.} \quad
\sigma _R : \quad W^{ab}_{F_R}&\To \GL_n({\wt F_{\o\omega_\oplus }})\ )\end{align*}
from the  Weil global group $W^{ab}_{F_L}$ (resp. $W^{ab}_{F_R}$~) to $ \GL_n({\wt F_{\omega _\oplus}})$ (resp. $ \GL_n({\wt F_{\o\omega_\oplus}})$~).

More concretely, as it was justified in the introduction and in \cite{Pie3}, bilinearity, instead of linearity, will be envisaged:
 so, enveloping (semi)algebras $A^e=A\otimes_R A^{\op}$ of a given $R$-(semi)algebra $A$~, where $A^{\op}$ denotes the opposite (semi)algebra of $A$~, as well as bisemialgebras will be considered.  Thus, we are interested in the products, right by left:
\Bi
\item  of infinite adele semirings $\Aa^\infty _{F_{\o\omega}}\times \Aa^\infty _{F_\omega }$ and of semigroups ${F_{\o\omega}}\times {F_\omega }$ and ${F_{\o\omega_\oplus}}\times {F_{\omega_\oplus} }$~.
\item of Picard semigroups $\Pic({\wt F_{\o\omega }})\times\Pic({\wt F_\omega })$~.
\item of affine semigroups $\SS^1_R\times \SS^1_L$~.
\item of  Weil global groups $W^{ab}_{F_R}\times W^{ab}_{F_L}$~.
\Ei
\vskip 11pt

\item The $n$-dimensional analog of the bilinear  affine semigroup  $\SS^1_R\times \SS^1_L$ will be a $2n$-dimensional bilinear affine semigroup  which is a reductive (and semisimple)
bilinear algebraic semigroup $G^{(2n)}({\wt F_{\o\omega }}\times {\wt F_\omega })$ considered as generated from the 
product $\GL_n({\wt F_{\o\omega }}\times {\wt F_\omega })\equiv T^t_n({\wt F_{\o\omega }})\times 
T_n({\wt F_\omega })$ of the group $T^t_n({\wt F_{\o\omega }})$ of lower triangular matrices with entries in 
${\wt F_{\o\omega }}$ by the group $T_n({\wt F_{\omega }})$ of upper triangular matrices with entries in 
${\wt F_\omega }$~.  So, the algebraic bilinear affine semigroup  $G^{(2n)}({\wt F_{\o\omega }}\times {\wt F_\omega })$ over ${\wt F_{\o\omega }}\times {\wt F_\omega }$~, resulting from the general bilinear algebraic
 semigroup of matrices $\GL_n({\wt F_{\o\omega }}\times{\wt F_\omega })$~, is a $\GL_n({\wt F_{\o\omega }}\times{\wt F_\omega })$-bisemimodule $\wt M_R\otimes \wt M_L$ (see the appendix).
 
 Recall that a reductive (resp. semisimple) group is a group having no unipotent (resp. solvable) infinite normal subgroup.
\vskip 11pt

\item The consideration of a  bilinear algebraic semigroup is justified by the fact that a bilinear algebraic (semi)group covers its linear equivalent as it is proved in proposition 2.1.7: this leads us to take into account  bialgebras and enveloping algebras.

As the bilinear algebraic semigroup $G^{(2n)}({\wt F_{\o\omega }}\times{\wt F_\omega })$ is built over 
$({\wt F_{\o\omega }}\times{\wt F_\omega })$~, it is composed of $r$ conjugacy classes, $1\le j\le r$~, having multiplicities $m^{(r)}=\sup(m_r+1)$ and corresponding to the $r$ biplaces of $({F_{\o\omega }}\times{F_\omega })$~.  Note that $m^{(r)}$ denotes the number of equivalent representatives in the $r$-th conjugacy class.

The algebraic representation of the bilinear algebraic semigroup of matrices $\GL_n({\wt F_{\o\omega }}\times{\wt F_\omega })$ into the 
$ \GL_n({\wt F_{\o\omega }}\times{\wt F_\omega })$ -bisemimodule 
$(\wt M_R\otimes \wt M_L)$ corresponds to an algebraic morphism from $\GL_n({\wt F_{\o\omega }}\times{\wt F_\omega })$ into $\GL(\wt M_R\otimes \wt M_L)$ where $\GL(\wt M_R\otimes \wt M_L)$ denotes the group of automorphisms of $\wt M_R\otimes \wt M_L$~.

Let $\wt M_{R_\oplus}\otimes \wt M_{L_\oplus}$ be the representation space of $\GL_n(\wt F_{\o\omega _\oplus}\times \wt F_{\omega_\oplus})$ decomposing into the direct sum of subbisemimodules representing its conjugacy classes as it will be seen in section 2.1.3.

Then, $\GL(\wt M_{R_\oplus}\otimes \wt M_{L_\oplus})$ constitutes the $2n$-dimensional equivalent of the product $W^{ab}_{F_R}\times W^{ab}_{F_L}$ of 
the global Weil groups and the bilinear algebraic semigroup $G^{(2n)}({\wt F_{\o\omega_\oplus }}\times{\wt F_{\omega_\oplus} })$ becomes 
naturally the $2n$-dimensional (irreducible) representation space\linebreak $(\Irr)\Rep^{(2n)}_{W_{F \RL}}(W^{ab}_{F_R}\times W^{ab}_{R_L})$ of 
$(W^{ab}_{F_R}\times W^{ab}_{R_L})$~.  The isomorphisms \cite{Bor2}:
\[\begin{matrix}
i_{\alg}&: \quad  &\GL_n({\wt F_{\o\omega }}\times{\wt F_\omega })
&\To  &G^{(2n)}({\wt F_{\o\omega }}\times{\wt F_\omega })\;,\\
i_{\aut}&: \quad  &\GL( \wt M_R \otimes \wt M_L)
&\To & \GL_n({\wt F_{\o\omega }}\times{\wt F_\omega })\;,\end{matrix}\]
leading to the (irreducible) representation 
\[(\Irr)\Rep^{(2n)}_{W_{\wt F \RL}}(W^{ab}_{F_R}\times W^{ab}_{F_L}):\quad \GL(\wt M_{R_\oplus}\otimes \wt M_{L_\oplus})\To G^{(2n)}(
{\wt F_{\o\omega_\oplus }}\times {\wt F_{\omega_\oplus} })\;,\]
implying the injective homomorphism:
\[\sigma _R\times \sigma _L :\quad W^{ab}_{F_R}\times W^{ab}_{F_L}\To \GL_n( {\wt F_{\o\omega_\oplus }}\times {\wt F_{\omega_\oplus} })\]
as announced precedingly.
\vskip 11pt 

\item So, the main purpose of chapter 2 will consist in introducing the bilinear algebraic semigroup $G^{(2n)}({\wt F_{\o\omega }}\times {\wt F_\omega })$ and in showing how it can be cut into conjugacy (bi)classes under the actions of Hecke operators.  On these basis, nonabelian global class field concepts will be taken up.
\vskip 11pt 
\Ei

\subsubsection{Notations}

The developments of this chapter will concern complex (algebraically closed) symmetric extension semifields $\widetilde F_L$ and $\widetilde F_R$ of a global number field $ k$~, taking into account that the ``real'' case can be handled similarly.

Let then
\begin{align*}
\wt F_{\omega_\oplus}
&= \bigoplus^r_{j=1}\wt F_{\omega _j}\ \bigoplus_{j,m_j}\wt F_{\omega _{j,m_j}}\\[11pt]
\text{(resp.} \quad
\wt F_{\o\omega_\oplus }
&= \bigoplus^r_{j=1}\wt F_{\o\omega _j}\ \bigoplus_{j,m_j}\wt F_{\o\omega _{j,m_j}}\ )\end{align*}
denote the sum of the \lr complex pseudo-ramified basic and equivalent decompositions of the semifield $\widetilde F_L$ (resp. $\widetilde F_R$~) at the set of complex places $\omega _j$ (resp. $\o\omega _j$~) as introduced in section 1.1.6 and let
\begin{align*}
\wt F^{nr}_{\omega _\oplus}
&= \bigoplus^r_{j=1}\wt F^{nr}_{\omega _j}\ \bigoplus_{j,m_j}\wt F^{nr}_{\omega _{j,m_j}}\\[11pt]
\text{(resp.} \quad
\wt F^{nr}_{\o\omega_\oplus }
&= \bigoplus^r_{j=1}\wt F_{\o\omega _j}\ \bigoplus_{j,m_j}\wt F^{nr}_{\o\omega _{j,m_j}}\ )\end{align*}
denote the sum of the respective pseudo-unramified extensions.
\vskip 11pt 

\subsubsection{The general reductive bilinear algebraic semigroups}

\Bean
\item
\Bi
\item If we refer to the introduction, the Gauss bilinear decomposition of the algebraic bilinear semigroup of matrices $\GL_n(\widetilde F_R\times \widetilde F_L)$ over the product of the extension semifields $\widetilde F_R$ and $\widetilde F_L$ can be developed according to:
\begin{align*}
\GL_n(\widetilde F_R\times \widetilde F_L)
& \equiv T^t_n(\widetilde F_R)\times T_n(\widetilde F_L)\\
&= [D_n(\widetilde F_R)\times UT^t_n(\widetilde F_R)] \times [UT_n(\widetilde F_L)\times D_n(\widetilde F_L)]
\end{align*}
in such a way that
\begin{align*}
T_n(\widetilde F_L) : \quad \widetilde F_L &\To T^{(2n)}(\widetilde F_L)\\
\text{(resp.} \quad 
T^t_n(\widetilde F_R) : \quad \widetilde F_R &\To T^{(2n)}(\widetilde F_R)\ )\end{align*}
can be viewed as an operator sending $\widetilde F_L $ (resp. $\widetilde F_R$~) into the affine semispace $T^{(2n)}(\widetilde F_L)$ (resp. $T^{(2n)}(\widetilde F_R)$~) of dimension $2n$~.

If $F_L$ (resp $F_R$~) denotes the set of equivalence classes of the completions associated with $\widetilde F_L$ (resp. $\widetilde F_R$~), then we have the commutative diagram:
\[ \begin{CD}
\widetilde F_L @>{T_n}>> T^{(2n)}(\widetilde F_L) \\ @VVV @VVV\\ F_L @>{T_n}>> T^{(2n)}(F_L)
\end{CD}
\qquad \mbox{\Huge(} \text{\ resp.}\quad 
\begin{CD}
\widetilde F_R @>{T^t_n}>> T^{(2n)}(\widetilde F_R) \\ 
@VVV @VVV\\ 
F_R @>{T^t_n}>> T^{(2n)}(F_R)
\end{CD}\quad \mbox{\Huge)}
\]
where  $T^{(2n)}(F_L)$ (resp. $T^{(2n)}(F_R)$~) is the locally compact linear complete semigroup isomorphic to $T^{(2n)}(\widetilde F_L)$ (resp. $T^{(2n)}(\widetilde F_R)$~).
\vskip 11pt

\item The question which arises now consists in knowing in what extent the generation of the affine semispace $T^{(2n)}(\widetilde F_L)$ (resp. $T^{(2n)}(\widetilde F_R)$~) from the ``action'' of the group of upper (resp. lower) triangular matrices $T_n(\centerdot)$ (resp. $T_n^t(\centerdot)$~) on the semifield $\widetilde F_L$ (resp. $\widetilde F_R$~) 
corresponds to the ``cross action'' of the group of upper (resp. lower) unitriangular matrices $UT_n(\centerdot)$ (resp. $UT_n^t(\centerdot)$~) by the group of diagonal matrices $D_n(\centerdot)$~.

The response lies on the {\bf existence of a lattice\/} decomposing the semifield $\widetilde F_L$ (resp. $\widetilde F_R$~) into conjugacy classes which are in one-to-one correspondence with the places of $F_L$ (resp. $F_R$~) in such a way that, at each place $\omega _j$ (resp. $\o\omega  _j$~), we have a basic completion $F_{\omega _j}$ (resp. $F_{\o\omega _j}$~) and a set of equivalent completions $\{F_{\omega _{j,m_j}}\}_{m_j\neq 0}$ (resp. $\{F_{\o\omega _{j,m_j}}\}_{m_j\neq 0}$~) as developed in section 1.1.5.
\vskip 11pt

\item Let then $\widetilde F_{b_L}$ (resp. $\widetilde F_{b_R}$~) denote the subfield of $\widetilde F_L$ (resp. $\widetilde F_R$~) composed of the extensions which are in one-to-one correspondence with the set $\{F_{\omega _j}\}_j$ (resp. $\{F_{\o\omega _j}\}_j$~) of basic completions of $F_L$ (resp. $F_R$~) and let 
$\widetilde F_{eq_L}=\widetilde F_L\setminus \widetilde F_{b_L}$ (resp. $\widetilde F_{eq_R}=\widetilde F_R\setminus \widetilde F_{b_R}$~) denote the subfield of $\widetilde F_L$ (resp. $\widetilde F_R$~) being in one-to-one correspondence with the set of equivalent completions $\{F_{\omega _{j,m_j}}\}_{m_j\neq 0}$ (resp. $\{F_{\o\omega _{j,m_j}}\}_{m_j\neq 0}$~).
\vskip 11pt

\item Then, the ``action'' of $T_n(\centerdot)$ (resp. $T^t_n(\centerdot)$~) will clearly correspond to the `` cross action'' of $UT_n(\centerdot)$ (resp. $UT^t_n(\centerdot)$~) by $D_n(\centerdot)$ if:
\Bi
\item $D_n(\centerdot)$ acts on the basic subfield  $\widetilde F_{b_L}$ (resp. $\widetilde F_{b_R}$~);
\item $UT_n(\centerdot)$ (resp. $UT^t_n(\centerdot)$~) acts on the complementary subfield 
$\widetilde F_{eq_L}$ (resp. $\widetilde F_{eq_R}$~).
\Ei

By this way, the diagonal group $D_n(\centerdot)$ generates an affine subsemispace\linebreak $T^{(2n)}(\widetilde F_{b_L})\subset T^{(2n)}(\widetilde F_L)$ (resp.
$T^{(2n)}(\widetilde F_{b_R})\subset T^{(2n)}(\widetilde F_R)$~) in one-to-one correspondence with the set of basic conjugacy class representatives of $T^{(2n)}(F_L)$ (resp. $T^{(2n)}(F_R)$~), generated from the set of basic completions $\{F_{\omega _j}\}_j$ (resp. $\{F_{\o\omega _j}\}_j$~), while the unitriagonal group $UT_n(\centerdot)$ (resp. $UT^t_n(\centerdot)$~) generates an affine subspace $T^{(2n)}(\widetilde F_{eq_L})\subset T^{(2n)}(\widetilde F_L)$ (resp.
$T^{(2n)}(\widetilde F_{eq_R})\subset T^{(2n)}(\widetilde F_R)$~), complementary of
$T^{(2n)}(\widetilde F_{b_L})\subset T^{(2n)}(\widetilde F_L)$ (resp.
$T^{(2n)}(\widetilde F_{b_R})\subset T^{(2n)}(\widetilde F_R)$~), in one-to-one correspondence with the set of equivalent conjugacy class representatives of $T^{(2n)}(F_L)$ (resp. $T^{(2n)}(F_R)$~), generated from the set of equivalent completions 
$\{F_{\omega _{j,m_j}}\}_{m_j\neq 0}$ (resp. $\{F_{\o\omega _{j,m_j}}\}_{m_j\neq 0}$~).
\vskip 11pt

\item So, we can state more precisely that:
\begin{align*}
T_n(\widetilde F_L) &= UT_n(\widetilde F_{b_L})\times D_n(\widetilde F_{b_L})\\
\text{(resp.} \quad
T^t_n(\widetilde F_R) &= D_n(\widetilde F_{b_R})\times UT^t_n(\widetilde F_{b_R})\ )\end{align*}
where:
\[ \widetilde F_L = \widetilde F_{b_L}\cup \widetilde F_{eq_L}
\qquad \text{(resp.} \quad 
\widetilde F_R = \widetilde F_{b_R}\cup \widetilde F_{eq_R}\ ),\]
but, the standard easier notation
\begin{align*}
T_n(\widetilde F_L) &= UT_n(\widetilde F_{L})\times D_n(\widetilde F_{L})\\
\text{(resp.} \quad
T_n(\widetilde F_R) &=D_n(\widetilde F_{R}) \times UT^t_n(\widetilde F_{R})\ )\end{align*}
will be commonly used in the following.

We can then summarize all that in the diagram:
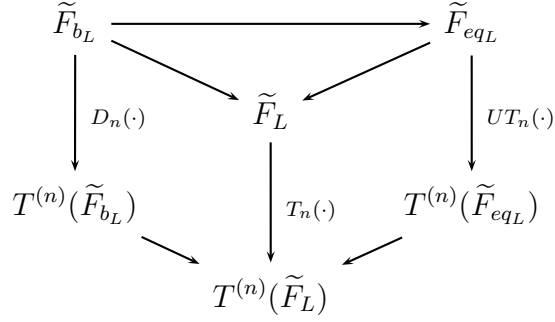
\[\begin{psmatrix}[colsep=1cm,rowsep=.6cm]
\widetilde F_{b_L} & & \widetilde F_{eq_L}\\
&\widetilde F_L \\
 T^{(n)}(\widetilde F_{b_L}) && T^{(n)}(\widetilde F_{eq_L})\\
&T^{(n)}(\widetilde F_L)
\psset{arrows=->,nodesep=5pt}
\everypsbox{\scriptstyle}
\ncline{1,1}{1,3}
\ncline{1,1}{2,2}
\ncline{1,3}{2,2}
\ncline{1,1}{3,1}>{D_n(\cdot)}
\ncline{1,3}{3,3}>{UT_n(\cdot)}
\ncline{2,2}{4,2}>{T_n(\cdot)}
\ncline{3,1}{4,2}
\ncline{3,3}{4,2}
\end{psmatrix}\]
(the right case being handled similarly).
\Ei

\item\Bi
\item Let $B_L$ and $B_R$ be two division semialgebras of dimension $2n$ respectively over the centers $\widetilde F_L$ and $\widetilde F_R$ such that 
$B_R$ be the opposite division semialgebra of $B_L$~.  Similarly, let $B_\omega $ and $B_{\o\omega }$ denote two division 
semialgebras of dimension $n$ respectively over $F_\omega $ and $F_{\o\omega }$~.  If we fix an isomorphism: $B_\omega \simeq T_n(F_\omega )$
(resp. $B_{\o\omega} \simeq T^t_n(F_\omega )$~), where $T_n(F_\omega )$ (resp. $T^t_n(F_\omega )$~) denotes the matrix algebra of upper (resp. lower) Borel triangular matrices, then we have that:
\[ B_{\o\omega }\otimes B_\omega \simeq T^t_n(F_{\o\omega })\times T_n(F_\omega )\equiv\GL_n(F_{\o\omega }\times F_\omega )\;.\]
$\GL_n(F_{\o\omega} \times F_\omega )$ thus is a reductive bilinear general complete semigroup of invertible $n\times n$ matrices over $F_{\o\omega }\times F_\omega $ having its representation space given by the tensor product $M_R\otimes M_L$ of a right $2n$-dimensional $B_{\o\omega }$-semimodule $M_R$ by a left $2n$-dimensional $B_\omega $-semimodule $M_L$ so that $M_L$ is provided with a left action of $T_n(F_\omega )$ while $M_R$ is provided with a right action of $T_n^t(F_{\o \omega })$~.
\vskip 11pt

\item We could also consider the bilinear complete semigroup $\GL_n(F_{\o\omega _\oplus}\times F_{\omega_\oplus})$ over product, right by left, of sums of complex pseudo-ramified completions $F_{\o\omega_\oplus}=\bigoplus_j\bigoplus_{m_j} F_{\o\omega_{j,m_j}}$ and
$F_{\omega_\oplus} =\bigoplus_j\bigoplus_{m_j} F_{\omega_{j,m_j}}$~.

Its representation space is given by the tensor product $M_{R_\oplus} \otimes M_{L_\oplus}$ of a right $n$-dimensional $T^t_n(F_{\o\omega _\oplus})$-semimodule $M_{R_\oplus}$ by its left equivalent 
$T_n(F_{\omega _\oplus})$-semimodule $M_{L_\oplus}$ in such a way that
 $M_{R_\oplus}\otimes M_{L_\oplus}$ decomposes according to:
\[M_{R_\oplus}\otimes M_{L_\oplus}=\bigoplus_j \bigoplus_{m_j}
(M_{\o\omega _{j,m_j}}\otimes M_{\omega _{j,m_j}})\]
where $M_{\omega _{j,m_j}} $ (resp. $M_{\o\omega _{j,m_j}}$~) are
$B_{\omega _{j,m_j}}$-subsemimodules (resp. $B_{\o\omega _{j,m_j}}$-sub\-semi\-modules).

This decomposition is a consequence of the decomposition of $F_{\o\omega_\oplus} \times F_{\omega_\oplus} $ according to:
\[F_{\o\omega_\oplus} \times F_{\omega_\oplus} 
=\bigoplus_j\bigoplus_{m_j}( F_{\o\omega_{j,m_j}} \times
F_{\omega_{j,m_j}} )\;.\]
\vskip 11pt

\item Similarly, the bilinear general complete semigroup of matrices over the product $\Aa^\infty _{F_{\o\omega }}\times \Aa^\infty _{F_\omega }$ of infinite adele semirings is given by:
\[ \GL_n(\Aa^\infty _{F_{\o\omega }}\times \Aa^\infty _{F_\omega })\equiv T^t_n( \Aa^\infty _{F_{\o\omega }} )\times T_n( \Aa^\infty _{F_\omega })\]
and verify
\[ B_{\Aa^\infty _{F_{\o\omega }}}\otimes B_{\Aa^\infty _{F_\omega }}\simeq \GL_n(\Aa^\infty _{F_{\o\omega }}\times \Aa^\infty _{F_\omega })\]
where   $B_{\Aa^\infty _{F_{\omega }}}$ (resp. $B_{\Aa^\infty _{F_{\o\omega }}}$~) is a division semialgebra of dimension $n$ over the adele 
semiring $\Aa^\infty _{F_{\omega }}$ (resp. $\Aa^\infty _{F_{\o\omega }}$~).

Its representation space is given by the tensor product 
$M_{R_{\tiny\textcircled{$\pi$}}}\otimes M_{L_{\tiny\textcircled{$\pi$}}}$ of the $T_n(\Aa^\infty _{F_\omega })$-semimodule $M_{L_{\tiny\textcircled{$\pi$}}}$ by the
$T^t_n(\Aa^\infty _{F_{\o\omega} })$-semimodule $M_{R_{\tiny\textcircled{$\pi$}}}$  which decompose respectively according to:
\[M_{L_{\tiny\textcircled{$\pi$}}}=\prod_{j_p} M_{\omega _{j_p}} \quad \and \quad 
M_{R_{\tiny\textcircled{$\pi$}}}=\prod_{j_p} M_{\o\omega _{j_p}} \]
where $M_{\omega _{j_p}} $ (resp. $M_{\o\omega _{j_p}} $~) is a $2n$-dimensional $B_{\omega_{j_p}}$-subsemimodule
(resp. $B_{\o\omega_{j_p}}$-subsemimodule) on the $j_p$-th Archimedean prime complex place $\omega _{j_p}$ (resp. $\o\omega _{j_p}$~) in the sense of section 1.1.6.
\vskip 11pt

\item The pseudo-unramified case can be introduced similarly. Indeed, let $B^{nr}_\omega $ and $B^{nr}_{\o\omega }$ denote two division 
semialgebras of dimension $2n$ respectively over the sums $F^{nr}_\omega $ and $F^{nr}_{\o\omega }$ of pseudo-unramified completions of the 
semifields $\wt  F_L$ and $\wt F_R$~.  Fixing the isomorphisms: $B^{nr}_\omega \simeq T_n(F^{nr}_\omega )$ and $B^{nr}_{\o\omega} 
\simeq T^t_n(F^{nr}_{\o\omega}) $~, we get:
\[ B^{nr}_{\o\omega }\otimes B_\omega ^{nr}\simeq T^t_n(F^{nr}_{\o\omega })\times T_n(F^{nr}_\omega )
\equiv \GL_n(F^{nr}_{\o\omega }\times F^{nr}_\omega )\]
and  $\GL_n(F^{nr}_{\o\omega} \times F^{nr}_\omega )$ has for representation space the tensor product $M^{nr}_R\otimes M^{nr}_L$ of a  pseudo-unramified right $B_{\o\omega }^{nr}$-semimodule $M^{nr}_R$ by its left equivalent $M^{nr}_L$~.

The bilinear general complete semigroup of invertible matrices over the product $ \Aa^{nr,\infty }_{F_{\o\omega} }\times \Aa^{nr,\infty }_{F_\omega}$ of pseudo-unramified infinite adele semirings is given by $\GL_n( \Aa^{nr,\infty }_{F_{\o\omega} }\times \Aa^{nr,\infty }_{F_\omega})$ and verifies:
\[ B_{\Aa^{nr,\infty }_{F_{\o\omega }}}\otimes B_{\Aa^{nr,\infty }_{F_\omega}} \simeq \GL_n( \Aa^{nr,\infty }_{F_{\o\omega} }
\times \Aa^{nr,\infty }_{F_\omega})\;.\]
\Ei\Ee
\vskip 11pt 

\subsubsection{Proposition}
{\em
\Bean
\item Any matrix $g_n(F_{\o\omega _j}\times F_{\omega_j }) \in \GL_n(F_{\o\omega }\times F_{\omega })$~, restricted over the product
$F_{\o\omega _j}\times F_{\omega _j}$ of the completions $F_{\o\omega _j}$ and $F_{\omega _j}$~, has the following {\bf Gauss bilinear decomposition\/}:
\[g_n(F_{\o\omega _j}\times F_{\omega _j})
= [(d_n(F_{\omega _j})\times d_n(F_{\o\omega _j})]\times [u^t_n(F_{\omega _j}]\times u_n(F_{\o\omega _j}]\]
where
\Bi
\item $d_n(F_{\omega _j})\in D_n(F_\omega )$ is an element of the group $D_n(F_\omega )$ of diagonal matrices of order $n$~.
\item $u_n(F_{\omega _j})\in UT_n(F_\omega )$ is an element of the group $UT_n(F_\omega )$ of upper unitriangular matrices.
\Ei

\item The representation space $(M_R\otimes M_L)$ of 
$\GL_n(F_{\o\omega }\times F_\omega )$ is (bi)-homomorphic to the 
representation space $M_{R_\oplus}\otimes M_{L_\oplus}$ of $\GL_n(F_{\o\omega _\oplus}\times F_{\omega _\oplus})$ and to the representation space 
$M_{R_{\tiny\textcircled{$\pi$}}}\otimes M_{L_{\tiny\textcircled{$\pi$}}}$ of 
$\GL_n(\Aa^\infty _{F_{\o\omega}}\times \Aa^\infty _{F_{\omega }})$ as it can be figured in the commutative diagram:
\vskip 11pt

\centerline{\begin{pspicture}(8,6)
\put(0,3){$\begin{CD}
\GL_n(F_{\o\omega }\times F_\omega ) @>>> 
\GL_n(F_{\o\omega _\oplus}\times F_{\omega _\oplus})\\
@VVV @VVV\\
M_R\otimes M_L @>{H_\oplus}>> M_{R_\oplus}\otimes M_{L_\oplus}\\
@VV{H_{\tiny\textcircled{$\pi$}}}V\\
M_{R_{\tiny\textcircled{$\pi$}}} \otimes M_{L_{\tiny\textcircled{$\pi$}}}\\
@AAA\\
\GL_n(\Aa^\infty _{F_{\o\omega }}\times \Aa^\infty _{F_\omega })
\end{CD}$}
\psline{->}(3.3,2.5)(4.8,3.2)
\psline{<-}(3.3,2.3)(4.8,3)
\put(3.3,2.7){\rotatebox{25}
{${}^{H_{\tiny\textcircled{$\pi$}\to\oplus}}$}}
\put(3.7,2){\rotatebox{25}
{${}^{H_{\oplus\to\tiny\textcircled{$\pi$}}}$}}
\end{pspicture}}
\Ee}
\vskip 11pt 

\bpr \Bean
\item The Gauss bilinear decomposition of $g_n(F_{\o\omega _j}\times F_{\omega _j})$ is an evident and natural generalization of the classical Gauss linear decomposition of the elements $g_n(F)$ of the general linear group $\GL_n(F)$ where $F$ denotes a number field.

The diagonal  $(B_{\o\omega _j}\otimes B_{\omega _j})$-subbisemimodules $M_{\o\omega _j}\otimes M_{\omega _j}$ correspond to the diagonal representation space of $D_n(F_{\o\omega }\times F_\omega )$~, while the off-diagonal
$(B_{\o\omega _{j,m_j}}\otimes B_{\omega _{j,m_j}})$-subbisemimodules $M_{\o\omega _{j,m_j}}\otimes M_{\omega _{j,m_j}}$ correspond to the off-diagonal representation space of $\GL_n(F_{\o\omega }\times F_\omega )$~.

Thus, the off-diagonal subsemimodules $M_{\o\omega _{j,m_j}}\otimes M_{\omega _{j,m_j}}$ are obtained from the diagonal subbisemimodules under the off-diagonal nilpotent action of $UT_n^t(F_{\o\omega })\times UT_n(F_\omega )$~.

\item As $\GL_n(F_{\o\omega _\oplus}\times F_{\omega _\oplus})$ and 
$\GL_n(\Aa^\infty _{F_{\o\omega }}\times \Aa^\infty _{F_\omega })$ are bilinear subgroups of $\GL_n(F_{\o\omega }\times F_\omega )$ by construction according to section 2.1.3, the homomorphisms:
\[ H_\oplus:\quad M_R\otimes M_L \To M_{R_\oplus}\otimes M_{L_\oplus} \qquad \and \qquad
H _{\tiny\textcircled{$\pi$}} :
\quad M_R\otimes M_L \To M_{R_{\tiny\textcircled{$\pi$}}}
\otimes M_{L_{\tiny\textcircled{$\pi$}}} \]
on their representation spaces follow directly.

On the other hand, the isomorphism
\[H_{\oplus\to {\tiny\textcircled{$\pi$} } }:\quad
M_{R_\oplus}\otimes M_{L_\oplus}\To
M_{R_{\tiny\textcircled{$\pi$} }}\otimes 
M_{L_{\tiny\textcircled{$\pi$} }}\;, \]
which is postutaled, results from the equality between the developments in series and in Eulerian products of the associated $L$-functions.\epr
\Ee
\vskip 11pt 

\subsubsection{Corollary}

{\em \Be
\item The $\GL_n(F_{\o\omega  _\oplus}\times F_{\omega _\oplus})$-bisemimodule $M_{R_\oplus}\otimes M_{L_\otimes}$ and the 
$\GL_n(\Aa^\infty _{F_{\o\omega 
}}\times \Aa^\infty _{F_{\omega 
}})$-bisemimodule $M_{R_{\tiny\textcircled{$\pi$} }}\otimes M_{L_{\tiny\textcircled{$\pi$} }}$ are representation spaces of the bilinear complete semigroup $\GL_n(F_{\o\omega }\times F_\omega )$~.

\item The $\GL_n(F_{\o\omega }\times F_\omega )$-bisemimodule $M_R\otimes M_L$ is a representation space of the bilinear algebraic semigroup $\GL_n(\Aa^\infty _{F_{\o\omega }}\times \Aa^\infty _{F_\omega })$~.\Ee
}
\vskip 11pt 

\bpr \Be
\item Part 1. of the corollary results from the homomorphism:
\begin{gather*}
\GL_n(F_{\o\omega }\times F_\omega ) \To M_R\otimes M_L \To M_{R_{\tiny\textcircled{$\pi$} }}\otimes M_{L_{\tiny\textcircled{$\pi$} }}\\
\begin{CD} @VVV \end{CD}\\
M_{R_\oplus}\otimes M_{L_\oplus}
\end{gather*}

\item Part 2. of the corollary results from the homomorphism:
\[ \GL_n(\Aa^\infty _{F_{\o\omega }}\times \Aa^\infty _{F_\omega }) \To M_{R_{\tiny\textcircled{$\pi$} }}\otimes M_{L_{\tiny\textcircled{$\pi$} }} \To
M_R\otimes M_L\]

\item and so on.\epr
\Ee
\vskip 11pt

\subsubsection{Linear and bilinear algebraic (semi)groups}

Let $\wt F_{\o\omega -\omega } =\wt F_{\o\omega }\cup \wt F_\omega $ denote the set of extensions decomposing   the algebraically (closed) symmetric splitting field $\widetilde F=\widetilde F_R\cup \widetilde F_L$~.

Let $\GL_n(\wt F_{\o\omega -\omega })$ denote the group of invertible $n\times n$ matrices with entries in $\wt F_{\o\omega -\omega }$~.  Then, the algebraic general linear group $\GL_n(\wt F_{\o\omega -\omega })$ has for representation space a vectorial space $W$ of dimension $4n^2$ isomorphic to $(\wt F_{\o\omega -\omega })^{4n^2}$ and has the Gauss decomposition:
\[ g_n(\wt F_{\o\omega -\omega_j }) = (u^t_n(\wt F_{\o\omega -\omega_j })\times u_n(\wt F_{\o\omega -\omega_j }))\times d_n
(\wt F_{\o\omega -\omega _j})]\]
for any matrix $g_n(\wt F_{\o\omega -\omega _j})\in \GL_n(\wt F_{\o\omega -\omega })$~.

On the other hand, let $\GL_n(\wt F_{\o\omega} \times \wt F_\omega )$ be the bilinear algebraic semigroup as introduced in section 2.1.3.  It has the Gauss bilinear decomposition introduced in proposition 2.1.4.

Then, the linear algebraic group $\GL_n(\wt F_{\o\omega -\omega })$ and the bilinear algebraic semigroup are in one-to-one correspondence under the conditions of proposition 2.1.7.
\vskip 11pt 

\subsubsection{Proposition}{\em

Let $ \wt F_{\o\omega -\omega } =\wt F_{\o\omega }\cup \wt F_\omega $ denote the set of extensions decomposing the algebraically closed symmetric splitting field $\widetilde F=\widetilde F_R\cup \widetilde F_L$~.

Let $g_n (\wt F_{\o\omega -\omega _j}) =(u^t_n(\wt F_{\o\omega -\omega_j })\times u_n(\wt F_{\o\omega -\omega_j }))
\times d_n(\wt F_{\o\omega -\omega_j })$ be the Gauss decomposition of the matrix $g_n(\wt F_{\o\omega -\omega_j })$ of the linear algebraic group $\GL_n(\wt F_{\o\omega -\omega_j })$~.

And, let
\[ g_n(\wt F_{\o\omega_j} \times \wt F_{\omega_j })
= [(d_n(\wt F_{\o\omega _j})\times d_n(\wt F_{\omega _j})]\times [u^t_n(\wt F_{\o\omega _j})\times u_n(\wt F_{\omega _j})]\]
be the Gauss decomposition of the matrix $g_n(\wt F_{\o\omega_j} \times \wt F_{\omega_j })$ of the bilinear algebraic semigroup $\GL_n(\wt F_{\o\omega} \times \wt F_{\omega})$~.

Then, if we take into account the maps:
\Bean
\item \quad $u_n(\wt F_{\o\omega -\omega_j } )\overset{\sim}{\To} u_n(\wt F_{\omega _j})$~,
\item \quad $u^t_n(\wt F_{\o\omega -\omega_j } )\overset{\sim}{\To} u^t_n(\wt F_{\o\omega _j})$~,
\item \quad $d_n(\wt F_{\o\omega -\omega_j } )\To d_n(\wt F_{\o\omega _j})\times d_n(\wt F_{\omega _j})$~,
\Ee
\Be
\item the linear algebraic group $\GL_n (\wt F_{\o\omega -\omega })$ is in one-to-one correspondence with the bilinear algebraic semigroup $\GL_n(\wt F_{\o\omega }\times \wt F_\omega )$~;
\item the $4n^2$-dimensional representation space $W$ of $\GL_n(\wt F_{\o\omega -\omega })$ coincides with the $4n^2$-dimensional representation space $\wt M_R\otimes \wt M_L$ of $\GL_n (\wt F_{\o\omega }\times \wt F_\omega )$~.
\Ee
}\vskip 11pt 

\bpr
\Be
\item the maps a) and b) send the unipotent matrices $u_n(\wt F_{\o\omega -\omega_j })$ and $u^t_n(\wt F_{\o\omega -\omega_j })$~, with entries in the completions of the symmetric splitting field $\wt F$~, respectively in the unipotent matrices
$u_n(\wt F_{\omega _j})$~, with entries in the extensions $\wt F_{\omega _j}$~, restricted to the upper half space, and in the transposed unipotent matrices 
$u^t_n(\wt F_{\o\omega _j})$~, with entries in the extensions $\wt F_{\o\omega _j}$~, restricted to the lower half space.

\item the map c) sends the centralizer $d_n (\wt F_{\o\omega -\omega_j }) $ of $g_n (\wt F_{\o\omega -\omega_j }) \in \GL_n
 (\wt F_{\o\omega -\omega }) $ into the ``bilinear'' centralizer $d_n (\wt F_{\o\omega _j}\times \wt F_{\omega_j }) =
d_n (\wt F_{\o\omega _j})\times d_n(\wt F_{\omega_j }) $ of the matrix $g_n (\wt F_{\o\omega_j} \times \wt F_{\omega_j }) \in \GL_n
 (\wt F_{\o\omega} \times \wt F_{\omega}) $~.

\item Consequently, the basis of the $4n^2$-dimensional representation space $W$ of the ``linear'' algebraic group $\GL_n
 (\wt F_{\o\omega -\omega }) $ must correspond to the basis of the $4n^2$-dimensional representation space $\wt M_R\otimes \wt M_L$ of the bilinear algebraic semigroup $\GL_n (\wt F_{\o\omega }\times \wt F_\omega ) $~, which implies that $W\simeq \wt M_R\otimes \wt M_L$~.  This means that the linear algebraic group $\GL_n (\wt F_{\o\omega -\omega }) $ is in fact in one-to-one correspondence with the bilinear algebraic semigroup $\GL_n (\wt F_{\o\omega }\times \wt F_\omega) $~.\epr
\Ee\vskip 11pt


\subsection{Lattices and Hecke algebras}

\subsubsection{Maximal orders and lattices}

Fix the maximal orders $\Os_{\wt F,\omega }$ of $\wt F_\omega $~, $\Os_{\wt F,\o\omega }$ of $\wt F_{\o\omega }$~, 
$\Os_{\wt F^{nr},\omega }$ of $\wt F^{nr}_\omega $ and  $\Os_{\wt F^{nr},\o\omega }$ of $\wt F^{nr}_{\o\omega }$~.  Then, the equivalent of the maximal order $\Lambda _\omega  $ (resp. $\Lambda _{\o\omega}  $~) in the division semialgebra $B_\omega $ (resp. $B_{\o\omega}$~) is a pseudo-ramified $\ZZ\big/N\ \ZZ$-lattice in the \lr $B_\omega $ (resp. $B_{\o\omega }$~)-semimodule $M_L$ (resp. $M_R$~).  Similarly, the equivalent of the maximal order 
$\Lambda ^{nr}_\omega  $ (resp. $\Lambda^{nr} _{\o\omega} $~) in $B_\omega ^{nr}$ (resp. $B_{\o\omega }^{nr}$~) is a  (pseudo-)unramified $\ZZ$-lattice in the (pseudo-)unramified \lr $B_\omega ^{nr}$ (resp. $B_{\o\omega }^{nr}$~)-semimodule $M_L^{nr}$ (resp. $M_R^{nr}$~).

Considering the decomposition of $F_\omega $ and $F_{\o\omega }$ into pseudo-ramified completions and of 
$F^{nr}_\omega $ and $F^{nr}_{\o\omega }$ into pseudo-unramified completions according to section 2.1.2, 
we have the expected decomposition of the pseudo-ramified $\ZZ\big/N\ \ZZ$-lattice $\Lambda _\omega $ (resp. $\Lambda _{\o\omega }$~) into:
\[ \Lambda _\omega  =\bigoplus_j\bigoplus_{m_j}\Lambda _{\omega _{j,m_j}}\qquad
\text{(resp.} \quad
\Lambda _{\o\omega } =\bigoplus_j\bigoplus_{m_j}\Lambda _{\o\omega _{j,m_j}}\ )\]
where $\Lambda _{\omega _{j,m_j}}\simeq \Os_{B_{\omega _{j,m_j}}}$ is a (pseudo-)ramified sublattice in the $B_{\omega _{j,m_j}}$-subsemimodule $M_{\omega _{j,m_j}}$~.

In the (pseudo-)unramified case, we have the similar decomposition of the pseudo-unramified $\ZZ$-lattice $\Lambda _\omega ^{nr}$ (resp. $\Lambda _{\o\omega }^{nr}$~) into:
\[ \Lambda^{nr} _\omega  =\bigoplus_j\bigoplus_{m_j}\Lambda^{nr} _{\omega _{j,m_j}}\qquad
\text{(resp.} \quad
\Lambda ^{nr}_{\o\omega } =\bigoplus_j\bigoplus_{m_j}\Lambda ^{nr}_{\o\omega _{j,m_j}}\ ).\]
On the other hand, we can fix the isomorphisms:
\[\begin{matrix}
& \Lambda _\omega  & \simeq & T_n(\Os_{F,\omega }) & \qquad & \text{(resp.}\quad &
  \Lambda _{\o\omega}  & \simeq & T^t_n(\Os_{F,\o\omega })\ )\\
\text{and} \quad 
& \Lambda^{nr} _\omega  & \simeq & T_n(\Os_{F^{nr},\omega }) & \qquad & \text{(resp.}\quad &
  \Lambda^{nr} _{\o\omega}  & \simeq & T^t_n(\Os_{F^{nr},\o\omega })\ )\end{matrix}\]
for the pseudo-ramified and pseudo-unramified lattices.
\vskip 11pt

\subsubsection{Proposition}
{\em Let $\Lambda _\omega $ and $\Lambda _{\o\omega }$ be pseudo-ramified $\ZZ\big/N\ \ZZ$-lattices and
$\Lambda ^{nr}_\omega $ and $\Lambda ^{nr}_{\o\omega }$ be the corresponding pseudo-unramified $\ZZ$-lattices as introduced in section 2.2.1.  Then, the pseudo-ramified bilattice  $\Lambda _{\o\omega }\otimes \Lambda _\omega $ in the $B_{\o\omega }\otimes B_\omega $-bisemimodule $M_R\otimes M_L$
\Bi
\item verifies:
\begin{align*}
\Lambda _{\o\omega }\otimes \Lambda _\omega 
&\simeq T^t_n(\Os_{F,\o\omega })\times T_n(\Os_{F,\omega} )\\
&\simeq \GL_n(\Os_{F,\o\omega }\times \Os_{F,\omega })\;.\end{align*}

\item has the decomposition:
\[\Lambda _{\o\omega }\otimes \Lambda _\omega
=\bigoplus_j \bigoplus_{m_j}( \Lambda _{\o\omega_{j ,m_j}}\otimes \Lambda _{\omega_{j,m_j}})\]
into subbilattices $\Lambda _{\o\omega_{j ,m_j}}\otimes \Lambda _{\omega_{j,m_j}}$~.
\Ei
And the pseudo-unramified bilattice
$\Lambda^{nr} _{\o\omega}\otimes \Lambda ^{nr}_{\omega} $ in the $B^{nr}_{\o\omega }\otimes B^{nr}_\omega $-bisemimodule $M^{nr}_R\otimes M^{nr}_L$
\Bi
\item satisfies:
\begin{align*}
\Lambda ^{nr}_{\o\omega }\otimes \Lambda^{nr} _\omega 
&\simeq T^t_n(\Os_{F^{nr},\o\omega })\times T_n(\Os_{F^{nr},\omega} )\\
&\simeq \GL_n(\Os_{F^{nr},\o\omega }\times \Os_{F^{nr},\omega })\;.\end{align*}

\item has the decomposition:
\[\Lambda ^{nr}_{\o\omega }\otimes \Lambda ^{nr}_\omega
=\bigoplus_j \bigoplus_{m_j}( \Lambda ^{nr}_{\o\omega_{j ,m_j}}\otimes \Lambda ^{nr}_{\omega_{j,m_j}})\;.\]
\Ei
}
\vskip 11pt 

\subsubsection{Proposition}

{\em
Consider that $\GL_n(\Os_{F,\o\omega }\times \Os_{F,\omega })$ has the following Gauss bilinear decomposition:
\[\GL_n(\Os_{F,\o\omega }\times \Os_{F,\omega })
= [D_n(\Os_{F,\omega })\times D_n(\Os_{F,\o\omega })]
\times [UT_n(\Os_{F,\omega })\times UT^t_n(\Os_{F,\o\omega })]\]
(see section 2.1.3).

Then, the nilpotent action of $[UT_n(\Os_{F,\omega })\times UT^t_n(\Os_{F,\o\omega })]$ on 
$ [D_n(\Os_{F,\omega })\times D_n(\Os_{F,\o\omega })]$ generates the off-diagonal subbilattices $\Lambda _{\o\omega _{j,m_j}}\otimes \Lambda _{\omega _{j,m_j}}$ equivalent to the basic subbilatice $\Lambda _{\o\omega _j,m_0}\otimes \Lambda _{\omega _j,m_0}\equiv 
\Lambda _{\o\omega _j}\otimes \Lambda _{\omega _j}$ for every $1\le j\le r$~.
}
\vskip 11pt 

\bpr This is a direct consequence of the generation of equivalent completions 
$F_{\o\omega _{j,m_j}}\times F_{\omega _{j,m_j}}$  from the basic completion
 $F_{\o\omega _j}\otimes F_{\omega _j}$ by the action of the nilpotent group element $u_{\omega ^2_j,m^2_j}$ as developed in section 1.1.5.\epr

\subsubsection{Definition: global decomposition groups and Frobenius automorphisms}  
\Be 
\item Let $g_n(\Os_{\o\omega _j}\times \Os_{\omega _j} )
=[d_n( \Os_{\omega _j} )\cdot d_n(\Os_{\o\omega _j})] [u_n(\Os_{\omega _j})\cdot u_n(\Os_{\omega _j})^t]$ denote the $j$-th component of $\GL_n( \Os_{F,\o\omega }\times \Os_{F,\omega } )$~.\par

The nilpotent part $[u_n(\Os_{\omega _j})\cdot u_n(\Os_{\o\omega _j})^t]$ of 
$g_n(\Os_{\o\omega _j}\times \Os_{\omega _j})$ is interpreted in the pseudo-ramified case as the element of the decomposition group of the $j$-th bisublattice in $(\wt M_R\otimes \wt M_L )$~.  The decomposition group element referring to the $(j,m_j)$-th bisublattice is denoted $D_{j^2;m_{j^2}}$ and has a representation given by
\[D_{j^2;m_{j^2}}=u_{n;m_j}(\Os_{\omega _j}) \cdot u_{n;m_j}(\Os_{\o\omega _j})^t
\in UT_n(\Os_{F,\omega })\times UT^t_n(\Os_{F,\o\omega })\;.\]
$D_{j^2;m_{j^2}}$ acts on the split Cartan subgroup element $\alpha _{n;j^2}= d_{n}(\Os_{\o\omega _j}) \cdot d_{n}(\Os_{\omega _j}) 
$ and is unimodular.\par

Let $\det(D_{j^2;m_{j^2}}\cdot \alpha _{n;j^2})_{ss}\simeq (\#\omega ^{2n}_j)\cdot N^{2n}\cdot (m^{(j)})^{2n}=j^{2n}\cdot N^{2n}\cdot (m^{(j)})^{2n}$ be the determinant of the semisimple form of $D_{j^2;m_{j^2}}\cdot \alpha _{n;j^2}$ where:
\Bi
\item $2n$ is the dimension;
\item the cardinality $\#\omega ^2_j=\#\o\omega _j\times \#\omega _j$ of the $j$-th quadratic place $\omega ^2_j$ is equal to $\#\omega ^2_j=f_{\o\omega _j}\cdot f_{\omega _j}=j^2$~;
\item $(m^{(j)}\cdot N)^{2n}$ is the rank of an irreducible quadratic complex completion according to section 1.1.4.
\Ei
Indeed, $\Lambda _\omega \simeq \Os^{2n}_{F,\omega }$ and $\Lambda _{\o\omega} \simeq \Os^{2n}_{F,\o\omega }$~.\par
Then, the unimodularity of the decomposition group $D_{j^2;m_{j^2}}$ leads to:
\[D_{j^2;m_{j^2}}:\det(\alpha _{2n;j^2}) \begin{CD} @>->-> \end{CD} \det(D_{j^2;m_{j^2}}\cdot \alpha _{2n;j^2})_{ss}\;.\]

\item in the pseudo-unramified case, i.e. when $N=1$~, we have that
\Bi
\item $D_{j^2;m_{j^2}}\equiv\Frob_{j^2;m_{j^2}}
=u_{n;m_j}(\Os^{nr}_{\omega _j})\cdot u_{n;m_j}(\Os^{nr}_{\o\omega _j})^t\in UT_n(\Os_{F^{nr},\omega })
\times UT^t_n(\Os_{F^{nr},\o\omega })
$~;\vskip 11pt

\item $\det(\Frob_{j^2;m_{j^2}}\cdot \alpha ^{nr}_{2n;j^2}) =j^{2n}$~, where $\alpha ^{nr}_{2n;j^2} =d_n(\Os^{nr}_{\omega _j})\cdot d_n(\Os^{nr}_{\o\omega _j})$~;\vskip 11pt

\item $\Frob_{j^2;m_{j^2}}:\det(\alpha _{n;j^2}^{nr}) \begin{CD} @>->-> \end{CD}\det(\alpha ^{nr}_{n;j^2}\cdot \Frob _{j^2;m_{j^2}})_{ss}$~.
\Ei\Ee
\vskip 11pt

\subsubsection{Proposition}  

{\em Let $\GL_n( \ZZ\big/N\ \ZZ)^2$ (resp. $\GL_n(\ZZ)^2)$~) be the ``pseudo-ramified'' (resp. ``pseudo-unramified'') general bilinear semigroup of matrices of order $n$ with entries in $\ZZ\big/N\ \ZZ$ (resp. $\ZZ$~).\newline
Then, $\Hs_{R\times L}(n)$ (resp. $\Hs^{nr}_{R\times L}(n)$~) 
will denote the pseudo-ramified (resp. pseudo-unramified) Hecke bialgebra   generated by all the pseudo-ramified (resp. pseudo-unramified) Hecke bioperators 
$(T_R(n;r)\otimes T_L(n;r))$ 
(resp. $(T_R^{nr}(n;r)\otimes T_L^{nr}(n;r))$~) having a representation 
in the subgroup of matrices of $\GL_n(\ZZ\big/N\ \ZZ)^2)$ (resp. 
$\GL_n( \ZZ )^2)$~), $1\le j\le r\le\infty$~, the index $j$ varying.}\vskip 11pt

\bpr
\Be
\item A \lr Hecke operator $T_L(n;r)$ (resp. $T_R(n;r)$~) is a \lr correspondence which associates to a \lr lattice 
$\Lambda _\omega $ (resp. $\Lambda _{\o\omega }$~) the sum of its \lr sublattices $\Lambda _{\omega_j}$ 
(resp. $\Lambda _{\o\omega _j} $~) of index $j$ and multiplicities $m^{(j)}=\sup (m_j)$~:
\begin{align*}
T_L(n;r)\ \Lambda _\omega 
&= \bigoplus^r_{j=1} \bigoplus_{m_j} \Lambda _{\omega _{j,m_j}}\\[11pt]
\text{(resp.}\quad 
T_R(n;r)\ \Lambda _{\o\omega }
&= \bigoplus^r_{j=1} \bigoplus_{m_j} \Lambda _{\o\omega _{j,m_j}}\ ).\end{align*}
A Hecke bioperator $T_R(n;r)\otimes T_L(n;r)$ is then defined by its action on the pseudo-ramified bilattice $\Lambda _\omega \otimes \Lambda _{\o\omega }$
\[(T_R(n;r)\otimes T_L(n;r))\ (\Lambda _\omega \otimes \Lambda _{\o\omega })
= \bigoplus^r_{j=1}\bigoplus_{m_j}(\Lambda _{\o\omega_{j,m_j}} \otimes \Lambda _{\omega _{j,m_j}})\;.\]
According to proposition 2.2.2, we have
\[\Lambda _{\o\omega }\otimes\Lambda _\omega \simeq \GL_n(\Os_{F,\o\omega }\times \Os_{F,\omega })\]
such that
\[\Lambda _{\o\omega_{j ,m_j}}\otimes\Lambda _{\omega_{j,m_j}}
 \simeq g_n( \Os_{F_{\o\omega_{j,m_j}} } \times \Os_{F_{\omega_{j,m_j}} })\in\GL_n( \Os_{F,\o\omega } \times\Os_{F,\omega })\]
where $\Os_{F_{\omega_{j,m_j}} }$ (resp. $\Os_{F_{\o\omega_{j,m_j}} }$~) is the maximal order $\Os_{F,\omega }$ (resp. $\Os_{F,\o\omega }$~) restricted to the extension $ \wt F_{\omega_{j,m_j}} $ (resp. $\wt F_{\o\omega_{j,m_j} }$~).

Remark that the entries of $g_n( \Os_{F_{\o\omega_{j,m_j}} } \times \Os_{F_{\omega_{j,m_j} }})$ are integers of
$( {\wt F_{\o\omega_{j,m_j}} } \times {\wt F_{\omega_{j,m_j} }})$ and the eigenvalues of 
$g_n( \Os_{F_{\o\omega_{j,m_j}} } \times \Os_{F_{\omega_{j,m_j} }})$ are algebraic integers.

The bisublattice $\Lambda  _{\o\omega_{j,m_j}}  \otimes \Lambda _{\omega_{j,m_j}} $ corresponds to the coset representative 
$g_n( \Os_{F_{\o\omega_{j,m_j} }} \times \Os_{F_{\omega_{j,m_j}} })$ of $\GL_n( \Os_{F,{\o\omega} } \times \Os_{F,{\omega_j} })$~, which involves that 
$g_n( \Os_{F_{\o\omega_{j,m_j}} } \times \Os_{F_{\omega_{j,m_j}} })$ can be chosen as a coset representative of the tensor product of Hecke operators \cite{Lang}.
\vskip 11pt 

\item Indeed, the ring of endomorphisms of the $(\ZZ\big/N\ \ZZ)^2$-bilattice $\Lambda _{\o\omega }\times \Lambda _\omega $ is generated over $(\ZZ\big/N\ \ZZ)^2$ by the products $(T_{j_R}\otimes T_{j_L})$ of Hecke operators $T_{j_R}$ and $T_{j_L}$ for the primes $j\nmid {N}$ and by the products $(U_{j_R}\otimes U_{j_L})$ of Hecke operators $U_{j_R}$ and $U_{j_L}$ for $j\mid {N}$ \cite{M-W}.

The coset representative of $U_{j_L}$~, referring to the upper half space, can be chosen to be upper-triangular and given by the integral matrix
\[t_{n;j,m_j}=\BM
j_{1_{N}} b^1_{2_{N}} & \cdots & \cdots& b^1_{n_{N}}\\
 & j_{2_{N}} &\cdots & b^2_{n_{N}}\\
 && \ddots & \vdots \\
\qquad  \raisebox{5mm}{\Huge{0}}&&&j_{n_{N}}\EM \qquad \text{of}\quad T_n(\Os_{F,\omega })\subset T_n(\ZZ\big/ N\ \ZZ)\]
such that:
\Bi
\item \quad $j_{n_{N}}=0\mod N$~,
\item \quad $b^1_{2_{N}}=0\mod N$~,
\item \quad $j_{1_{N}}\times \cdots j_{n_{N}} \simeq j^n\cdot N^n\cdot (m^{(j)})^n$
 \;\; such that \;\; $(f_{\o\omega _j}\cdot f_{\omega _j})^n\cdot N^{2n}\simeq j^{2n}\cdot N^{2n}\cdot (m^{(j)})^{2n}$ 

\quad (see definition 2.2.4),
\item \quad $t_{n;j,m_j}\equiv t_n(\Os_{F_{\omega _{j,m_j}}})\subset g_n(\Os_{F_{\o\omega _{j,m_j}}}\times 
\Os_{F_{\omega _{j,m_j}}})$~.
\Ei

Similarly, the coset representative of $U_{j_R}$~, referring to the lower half space can be chosen to be lower triangular and given by the integral matrix $t^t_{n;j,m_j}\equiv t^t_n(\Os_{F_{\o\omega _{j;m_j}}})$ of $T^t_n(\Os_{F,\o\omega })$~, i.e. by the transposed matrix of $t_{n;j,m_j}$~.

So, the coset representative of $(U_{j_R}\times U_{j_L})$ will be
\[ g_n(\Os_{F_{\o\omega _{j,m_j}}}\times \Os_{F_{\omega _{j,m_j}}})
= t^t_{n;j,m_j}\times t_{n;j,m_j} \in\GL_n(\Os_{F,\o\omega }\times \Os_{F,\omega })\subset \GL_n((\ZZ\big/ N\ \ZZ)^2)\;.\]
Note that $g_n(\Os_{F_{\o\omega _{j,m_j}}}\times \Os_{F_{\omega _{j,m_j}}})$ has the Gauss decomposition in diagonal and nilpotent part so that the nilpotent part is given by $u_n(\Os_{F_{\omega _{j,m_j}}})\times u^t_n( \Os_{F_{\o\omega _{j,m_j}}})$ and constitutes an element of the decomposition group $D_{j^2;m^2_j}$ of the corresponding bisublattice
$(\Lambda _{\o\omega _{j,m_j}}\otimes \Lambda _{\omega _{j,m_j}})$ as introduced in definition 2.2.4.
\vskip 11pt 

\item The pseudo-unramified Hecke bioperator $T_R^{nr}(n;r)\otimes T_L^{nr}(n;r)$ can also be envisaged throughout its action on the pseudo-unramified bilattice $\Lambda ^{nr}_{\o\omega }\otimes \Lambda _\omega ^{nr}$~:
\[ T_R^{nr}(n;r)\otimes T_L^{nr}(n;r))\ (\Lambda _{\o\omega }^{nr}\otimes \Lambda _\omega ^{nr})
= \bigoplus^r_{j=1}\bigoplus_{m_j} (\Lambda ^{nr}_{\o\omega _{j,m_j}}\otimes \Lambda ^{nr}_{\omega _{j,m_j}})\]
such that:
\[
(\Lambda ^{nr}_{\o\omega _{j,m_j}}\otimes \Lambda ^{nr}_{\omega _{j,m_j}})
\simeq g_n(\Os_{F^{nr}_{\o\omega _{j,m_j}}}\times \Os_{F^{nr}_{\omega _{j,m_j}}})
\in\GL_n(\Os_{F^{nr},\o\omega }\times \Os_{F^{nr},\omega })\subset \GL_n(\ZZ^2)\;.\]
The pseudo-unramified correspondent $U^{nr}_{j_L}$ (resp. $U^{nr}_{j_R}$~) of the coset representative of $U_{j_L}$ (resp. $U_{j_R}$~) is an upper (resp. lower) triangular integral matrix $t_n(\Os_{F^{nr}_{\omega _{j,m_j}}})$
(resp. $t^t_n(\Os_{F^{nr}_{\o\omega _{j,m_j}}})$~) of $T_n(\Os_{F^{nr},\omega })$ (resp. $T_n(\Os_{F^{nr},\o\omega })$~) $\subset T_n(\ZZ)$ (resp. $T^t_n(\ZZ)$~).\epr
\Ee
\vskip 11pt

\subsubsection{Proposition} 

{\em The Hecke bialgebra $\Hs_{R\times L}(n)$ (resp. $\Hs^{nr}_{R\times L}(n)$~) of the endomorphisms of the pseudo-ramified (resp. pseudo-unramified) $B_{\o\omega }\otimes B_\omega$-bisemimodule $M_R\otimes M_L$ (resp. $B^{nr}_{\o\omega }\otimes B^{nr}_\omega$-bisemimodule $M^{nr}_R\otimes M^{nr}_L$~) is generated by the pseudo-ramified (resp. pseudo-unramified) Hecke bioperators $(T_R(n;r)\otimes T_L(n;r))$ (resp. $(T^{nr}_R(n;r)\otimes T^{nr}_L(n;r))$~) having as representation $\GL_n(( \ZZ\big/ N\ \ZZ )^2)$ (resp. $\GL_n((\ZZ)^2)$~).}\vskip 11pt

\bpr as $\Lambda _{\o\omega }\otimes \Lambda _\omega $ is a pseudo-ramified bilattice in the $B_{\o\omega }\otimes B_\omega $-bisemimodule $M_R\otimes M_L$ and as the
decomposition of the $B_{\o\omega }\otimes B_{\omega }$-bisemimodule $M_R\otimes M_L$ into the set of $B_{\o\omega _j}\otimes B_{\omega _j}$-subbisemimodules $M_{\o\omega _j}\otimes M_{\omega _j}$ having multiplicities $m_j$ is obtained through the action of the Hecke bioperator $(T_R(n;r)\otimes T_L(n;r))$~, it is immediate to check that  $(T_R(n;r)\otimes T_L(n;r))$ has a representation in the  bilinear congruence subgroup $\GL_n(\ZZ\big/N\ \ZZ)^2)$~.\newline
The pseudo-unramified case can be handled similarly.\epr\vskip 11pt

\subsection{Lattice semispaces}

\subsubsection{Definition: pseudo-ramified lattice semispaces}  

The space $X=\GL_n(\RR)/\GL_n(\ZZ)$ corresponds to the set of lattices of  $\RR^n$~.  In this perspective, 
if $\widetilde F_L$ (resp. $\widetilde F_R$~) denotes a complex symmetric splitting semifield (see section 1.1.4),
we shall introduce a \lr (pseudo-)ramified lattice semispace:
\[X_{S_L}=T_n({\widetilde F_L})\Big/ T_n(\ZZ\big/N\ \ZZ) \qquad \text{(resp.}\quad 
X_{S_R}=T^t_n({\widetilde F_{L} })\Big/ T^t_n( \ZZ\big/N\ \ZZ)\ )\]
isomorphic to the \lr (pseudo-)ramified $B_\omega $-semimodule $M_L$  (resp. $B_{\o\omega} $-semimodule $M_R$~) such that the \lr cosets of $X_{S_L}$ (resp. $X_{S_R}$~) are isomorphic to the \lr pseudo-ramified $B_{\omega _j}$-subsemimodules $M_{\omega _j}$ (resp. $B_{\o\omega _j}$-subsemimodules $M_{\o\omega _j}$~).\par

Consequently, we can define a pseudo-ramified lattice bisemispace
\[X_{S_{R\times L}}=\GL_n({\widetilde F_{R}}\times {\widetilde F_L})\Big/ \GL_n( (\ZZ\big/N\ \ZZ)^2)\]
isomorphic to the pseudo-ramified $B_{\o\omega }\otimes B_\omega $-bisemimodule $M_R\otimes M_L$~.\par

As the \lr pseudo-ramified lattice  semispace $X_{S_L}$ (resp. $X_{S_R}$~) is defined with respect to the matrix algebra $T_n( \ZZ\big/N\ \ZZ)$ (resp. $T^t_n( \ZZ\big/N\ \ZZ)$~) in one-to-one correspondence with the pseudo-ramified \lr lattice $\Lambda _\omega $ (resp. $\Lambda _{\o\omega} $~), $X_{S_L}$ (resp. $X_{S_R}$~) is a \lr Hecke lattice semispace.
\par
Similarly, the lattice bisemispace $X_{S_{R\times L}}$ is a pseudo-ramified Hecke lattice bisemispace since $\GL_n( \ZZ\big/N\ \ZZ)$ is in one-to-one correspondence with the bilattice $\Lambda _{\o\omega \times \omega }=\Lambda _{\o\omega }\times \Lambda _\omega $~.
\vskip 11pt

\subsubsection{Definition: pseudo-unramified Hecke lattice semispaces}  

A \lr pseudo-unramified lattice semispace:
\[X^{nr}_{S_L}=T_n( \widetilde F^{nr}_L)\Big/ T_n(\ZZ)\qquad \text{(resp.}\quad
X^{nr}_{S_R}=T^t_n(\widetilde F^{nr}_{R} )\Big/ T^t_n( \ZZ )\]
can similarly be introduced; it is a \lr pseudo-unramified Hecke lattice semispace isomorphic to the \lr pseudo-unramified $B_\omega ^{nr}$-semimodule $M_L^{nr}$ (resp. $B_{\o\omega} ^{nr}$-semimodule $M_R^{nr}$~) verifying $M_L^{nr}\subset M_L$ (resp. $M_R^{nr}\subset M_R$~).
\par
A pseudo-unramified lattice bisemispace
\[X^{nr}_{S_{R\times L}}=\GL_n(\widetilde F^{nr}_{R}\times \widetilde F^{nr}_L)\Big/ 
\GL_n((\ZZ^2))\]
can also be defined: it is a pseudo-unramified Hecke lattice bisemispace isomorphic to the pseudo-unramified $B_{\o\omega} ^{nr}\otimes B_\omega ^{nr}$-bisemimodule $M_R^{nr}
\otimes M_L^{nr}$~.
\vskip 11pt

\subsubsection{Proposition}  

{\em \Bean
\item Any \lr coset 
\[X_{S_L}^{nr}(\wt  F^{nr}_{\omega _j })=t_n(\wt F_{\omega _j}^{nr})\Big/ 
t_n(\Os_{F_{\omega _j}^{nr}})\qquad
\text{(resp.} \quad X^{nr}_{S_R}(\wt F^{nr}_{\o\omega _j})=t_n(\wt F_{\o\omega _j}^{nr})^t\Big/ 
t_n(\Os_{F_{\o\omega _j}^{nr}})^t\ )\] of $X_{S_L}^{nr}$ (resp. $X_{S_R}^{nr}$~) having a multiplicity 
one is a \lr pseudo-unramified Hecke lattice semisubspace having a rank equal to $j^{2n}$ if $f_{\omega _j}=j$~.

\item the product $X_{S_{R\times L}}^{nr} (\wt F^{nr}_{\o\omega_j}\times \wt F^{nr}_{\omega _j}) =g_n(\wt F^{nr}_{\o\omega_j }\times \wt F_{\omega _j}^{nr})\Big/g_n( \Os_{F^{nr}_{\o\omega _j}}\times \Os_{F^{nr}_{\omega _j}})$ of a right coset 
$X_{S_R}^{nr}(\wt F^{nr}_{\o\omega _j})$ by a left coset 
$X_{S_L}^{nr}(\wt F^{nr}_{\omega _j})$ having multiplicities one is a  pseudo-unramified Hecke lattice bisemisubspace having a rank equal to $j^{2n} $~.
\Ee}
\vskip 11pt

\bpr  
\Bean
\item The \lr pseudo-unramified Hecke lattice semisubspace $X_{S_L}^{nr}(\wt F^{nr}_{\omega _j})$ (resp.\linebreak $X_{S_R}^{nr}(\wt F^{nr}_{\o\omega _j})$~) is isomorphic to the pseudo-unramified \lr $B^{nr}_{\omega _j}$-subsemimo\-dule $M^{nr}_{\omega _j}$ (resp. $B^{nr}_{\o\omega _j}$-subsemimodule $M^{nr}_{\o\omega _j}$~) which has a rank equal to the $n$-fold product  of its global class residue degree $f_{\omega _j}=[F_{\omega _j}^{nr}:k_L]=j$ (resp. $f_{\o\omega _j}=[F_{\o\omega _j}^{nr}:k_R]=j$~); so, we have that the rank of $X_{S_L}^{nr}(F^{nr}_{\omega _j})$ (resp. $X_{S_R}^{nr}(F^{nr}_{\o\omega _j})$~) is equal to $ j^{2n}$~.

\item the pseudo-unramified Hecke lattice bisemisubspace $X_{S_{R\times L}}^{nr}( \wt F^{nr}_{\o\omega _j}\times \wt F^{nr}_{\omega _j}  )$ is isomorphic to the pseudo-unramified $B_{\o\omega _j}^{nr}\otimes B_{\omega _j}^{nr}$-subbisemimodule $M_{\o\omega _j}^{nr}\otimes M_{\omega _j}^{nr}$ having a rank equal to $j^{2n}$~.
\epr\Ee\vskip 11pt

\subsubsection{Proposition}  

Let $f_{\omega _j}=j$ denote the global class residue degree of the $j$-th place having multiplicity $m^{(j)}$~.

{\em \Bean
\item The pseudo-unramified $B_{\o\omega }^{nr}\otimes B_{\omega }^{nr}$-bisemimodule $M_R^{nr}\otimes M_L^{nr}$ has a rank given by:
\[r_{\o\omega \times \omega }^{nr}=\txt\bigoplus\limits_{j=1}^r (f_{\o\omega _j}\cdot m^{(j)})^{2n}(f_{\omega _j}\cdot m^{(j)})^{2n}=\txt\bigoplus\limits_j (j\cdot m^{(j)})^{2n}\;.\]

\item the pseudo-ramified $B_{\o\omega }\otimes B_{\omega}$-bisemimodule $M_R\otimes M_L$ has a rank given by:
\[r_{\o\omega \times \omega }=\txt\bigoplus\limits_{j=1}^rN^{2n}(f_{\o\omega _j}\cdot m^{(j)})^{2n}(f_{\omega _j}\cdot m^{(j)})^{2n}=\txt\bigoplus\limits_jN^{2n} (j\cdot m^{(j)})^{2n} \;.\]
\Ee}\vskip 11pt

\subsubsection{Summary}  

The preceding developments can be summarized as follows:\par
Let
\Bi
\item $\Lambda _{\o\omega }\otimes \Lambda _\omega $ be a $( \ZZ\big/N\ \ZZ)^2$-bilattice (pseudo-ramified);\vskip 11pt
\item $\Lambda ^{nr}_{\o\omega }\otimes \Lambda^{nr} _\omega $ be a $\ZZ^2$-bilattice (pseudo-unramified);\vskip 11pt
\item $M ^{nr}_{R }\otimes M^{nr} _L $ be a $\GL_n(F^{nr}_{\o\omega }\times F_\omega ^{nr})$-bisemimodule (pseudo-unramified) isomorphic to the pseudo-unramified lattice bisemispace $X_{S_{R\times L}}^{nr}$~;\vskip 11pt
\item $M _{R }\otimes M _L $ be a $\GL_n(F_{\o\omega }\times F_\omega )$-bisemimodule (pseudo-ramified) isomorphic to the pseudo-ramified lattice bisemispace $X_{S_{R\times L}}$~.
\Ei
Then,  the following commutative diagram is evident:
\[
\begin{array}{cccccc}
\Lambda ^{nr}_{\o\omega}\otimes \Lambda _\omega^{nr}  & &\scalebox{1.5}{$\hookrightarrow$}   & \Lambda _{\o\omega }\otimes 
\Lambda _\omega  &:&  \text{bilattices}\\
 \Big\downarrow & \underset{\substack{pseudo-\\ramification}}{\mbox{{$\dashrightarrow$}}} & &
 \Big\downarrow &&
 \Big\downarrow \\
M^{nr}_{R}\otimes M _L^{nr}  && \scalebox{1.5}{$\hookrightarrow$}   & M_R\otimes M_L  &:&  \GL_n(F_{\o\omega}^{(nr)}\times F_\omega^{(nr)})\text{-bisemimodules}\end{array}
\]
where
\Bi
\item the vertical arrows stand for injective morphisms from bilattices to $\GL_n(F_{\o\omega}^{(nr)}\times F_\omega^{(nr)})$-bisemimodules;
\item the horizontal arrows stand for pseudo-ramification injective morphisms.
\Ei
\vskip 11pt

\subsection{Representations of global Weil groups}

\subsubsection{Nonabelian global class field concepts: a summary}

\Be
\item 
Let $G^{(2n)}({\wt F_{\o\omega }}\times {\wt F_\omega })$~, 
$\GL_n({\wt F_{\o\omega }}\times {\wt F_\omega })$ and $\GL(\wt M_R\otimes \wt M_L)$ denote respectively the bilinear algebraic semigroup over 
$({\wt F_{\o\omega }}\times {\wt F_\omega })$ which is a bilinear affine semigroup, its bilinear algebraic semigroup of matrices and the group of 
automorphisms of $(\wt M_R\otimes \wt M_L)$ which are isomorphic:
\[
G^{(2n)}({\wt F_{\o\omega }}\times {\wt F_\omega })
\simeq \GL_n({\wt F_{\o\omega }}\times {\wt F_\omega })\simeq
\GL(\wt M_R\otimes \wt M_L)\]
according to section 2.1.1: this implies the choice of a basis in $\wt M_R\otimes \wt M_L$~, which is evident from its construction from 
$({\wt F_{\o\omega }}\times {\wt F_\omega })$~.

The diagonal conjugacy classes of $G^{(2n)}({\wt F_{\o\omega }}\times {\wt F_\omega })$~, noted $g^{(2n)}\RL[j,m_j=0]$~, then correspond to the $\wt B_{\o\omega _j}\otimes \wt B_{\omega _j}$-subbisemimodules $\wt M_{\o\omega _j}\otimes \wt M_{\omega _j}$ and their off-diagonal equivalents $g^{(2n)}\RL[j,m_j]$~, $m_j>0$~, related to the multiplicity $m^{(j)}=\sup(m_j)$ of $g^{(2n)}\RL[j,m_j=0]$~, correspond to the $\wt B_{\o\omega _{j,m_j}}\otimes \wt B_{\omega _{j,m_j}}$-subbisemimodules $\wt M_{\o\omega _{j,m_j}}\otimes \wt M_{\omega _{j,m_j}}$~: they are obtained from $(\wt M_R\otimes \wt M_L)$ by the action of the Hecke bioperator $T_R(n;r)\otimes T_L(n;r)$~.

As $G^{(2n)} ({\wt F_{\o\omega }}\times {\wt F_\omega })$ is a smooth reductive bilinear affine semigroup, it
can be considered as the $n$-dimensional analog of the product
$(\SS^1_R\times \SS^1_L)$ of the one-dimensional affine complex semigroups $\SS^1_R$ and $\SS^1_L$~, introduced in section 2.1.1, we have that:
\[G^{(2n)}({\wt F_{\o\omega }}\times {\wt F_\omega })\simeq \prod_n(\SS^1_R\times \SS^1_L)\;.\]
Similarly, $\GL(\wt M_{R_\oplus}\otimes \wt M_{L_\oplus})$ can be viewed as the $n$-dimensional correspondent of the product
$W^{ab}_{F_R}\times W^{ab}_{F_L}$ of the Weil global groups included in the product
$ \Gal(F^{ac}_R\big/k)\times \Gal(F^{ac}_L\big/k)$ of the Galois groups as introduced in sections 1.1.7 and 1.1.9 which implies that:
\[ \GL(\wt M_{R_\oplus}\otimes \wt M_{L_\oplus})\simeq \prod _n(W^{ab}_{F_R}\times W^{ab}_{F_L})\subset
\prod _n(\Gal(F^{ac}_R\big/k_R)\times \Gal(F^{ac}_L\big/k_L))\;.\]
\vskip 11pt 

\item
Let us now introduce the normal bilinear algebraic semisubgroup $P^{(2n)}({\wt F_{\o\omega^1 }}\times {\wt F_{\omega^1} })$ 
of $G^{(2n)}({\wt F_{\o\omega }}\times {\wt F_\omega })$~: its $j$-th class noted $P ^{(2n)} (\wt F_{\o\omega ^1_j}\times \wt F_{\omega ^1_j})$ is defined over the product $(\wt F_{\o\omega ^1_j}\times \wt F_{\omega ^1_j})$ of irreducible central extensions 
$\wt F_{\o\omega ^1_j}$ and $\wt F_{\omega ^1_j}$ having a rank $N\cdot m^{(j)}$~, as given in section 1.1.5.  So, $P^{(2n)}
(\wt F_{\o\omega ^1_j}\times \wt F_{\omega ^1_j})$ has  a rank equal to $(N\cdot m^{(j)})^{2n}$~.  
$P^{(2n)}({\wt F_{\o\omega^1}}\times {\wt F_{\omega^1}})$ also acts on
$G^{(2n)}({\wt F_{\o\omega }}\times {\wt F_\omega })$ by conjugation.

Recall that the Galois subgroup of the irreducible extension $\widetilde F_{\omega ^1_j}$ associated to the completion 
$F_{\omega ^1_j}$ is the global inertia subgroup $I_{F_{\omega _j}}$~.

More generally, let
\begin{align*}
{\wt F_{\omega^1} }&
= \{\wt F_{\omega ^1_1},\cdots,\wt F_{\omega ^1_{j,m_j}},\cdots,
\wt F_{\omega ^1_r}\}\\
\text{(resp.} \quad
{\wt F_{\o\omega^1} }&
= \{\wt F_{\o\omega ^1_1},\cdots,\wt F_{\o\omega ^1_{j,m_j}},\cdots,
\wt F_{\o\omega ^1_r}\}\ )\end{align*}
denote the set of pseudo-ramified irreducible extensions 
$\wt F_{\omega ^1_j}$ and $\wt F_{\omega ^1_{j,m_j}}$ of ${\wt F_\omega }$ and ${\wt F_{\o\omega} }$ respectively.

Let $G^{(2n)}(\wt F^{nr}_{\o\omega }\times \wt F^{nr}_\omega )$ denote the pseudo-unramified bilinear algebraic semigroup over the product $\wt F^{nr}_{\o\omega }\times \wt F^{nr}_\omega $ of pseudo-unramified extensions.

Then, $P^{(2n)}({\wt F_{\o\omega^1 }}\times {\wt F_{\omega^1} })$~, which is in 
fact a bilinear (complex) minimal (not standard) parabolic subsemigroup, 
can be defined as being the kernel of the map:
\[ G^{(2n)}_{F\to F^{nr}}: \quad 
G^{(2n)}({\wt F_{\o\omega }}\times {\wt F_\omega })\To
G^{(2n)} (\wt F^{nr}_{\o\omega }\times \wt F^{nr}_\omega ) \;;\]
so:
\[ P^{(2n)}({\wt F_{\o\omega^1 }}\times {\wt F_{\omega^1} })=\Ker(G^{(n)}_{F\to F^{nr}})\]
is the smallest bilinear normal pseudo-ramified subgroup of $G^{(2n)}({\wt F_{\o\omega }}\times {\wt F_\omega })$~: it is the
 minimal parabolic subgroup viewed as the connected component of the identity in $G^{(2n)} ({\wt F_{\o\omega }}\times {\wt F_\omega }) $ on the Zarisky topology.
\vskip 11pt 

\item Referring to the decomposition of the product, right by left, $\wt F_{\o\omega _j}\times \wt F_{\omega _j}$ of the basic extensions $\wt F_{\o\omega _j}$ and $\wt F_{\omega _j}$ into $j'$ equivalent bisubextensions $\wt F_{\o\omega^{j'} _j}\times \wt F_{\omega ^{j'}_j}$~, $1\le j'\le j$~, as well as the decomposition of the equivalent biextensions 
$\wt F_{\o\omega _{j,m_j}}\times \wt F_{\omega _{j,m_j}}$ into $j'$ equivalent bisubextensions 
$\wt F_{\o\omega^{j'} _{j,m_j}}\times \wt F_{\omega ^{j'}_{j,m_j}}$~, we have that the $j$-th diagonal conjugacy class $g^{(2n)}\RL[j,m_j=0]$ of $G^{(2n)}({\wt F_{\o\omega }}\times {\wt F_\omega })$ will be cut into $j'$ equivalent diagonal conjugacy subclasses $g^{(2n)}\RL[j,m_j=0,j']$ and that the $j$-th off-diagonal equivalent conjugacy class
$g^{(2n)}\RL[j,m_j]$~, $m_j>0$~, will also be cut into $j'$ equivalent off-diagonal conjugacy subclasses $g^{(2n)}\RL[j,m_j,j']$~, $1\le j'\le j$~.

So, the $j$-th diagonal conjugacy class $g^{(2n)}\RL[j,m_j=0]$ and its off-diagonal equivalent $g^{(2n)}\RL[j,m_j]$~, $m_j>0$~, is decomposed into a set $\{g^{(2n)}\RL[j,m_j,j']\}^j_{j'=1}$ of $j'$ conjugacy subclasses $g^{(2n)}\RL[j,m_j,j']$~.
\Ee
\vskip 11pt 

Then, we can state the following propositions:
\vskip 11pt 

\subsubsection{Lemma}

{\em Let $G^{(2n)} ({\wt F_{\o\omega }}\times {\wt F_\omega }) $ be a smooth reductive bilinear affine  semigroup and let 
$P^{(2n)} ({\wt F_{\o\omega^1 }}\times {\wt F_{\omega^1} }) $ be its parabolic subgroup viewed as a locally  subgroup 
of $G^{(2n)} ({\wt F_{\o\omega }}\times {\wt F_\omega }) $~.

Then, it is clear that:
\Be
\item at every (bi)extension corresponding to the (bi)place $\o\omega _j\times \omega _j$ of $({F_{\o\omega }}\times {F_\omega }) $ corresponds a conjugacy class of $G^{(2n)} ({\wt F_{\o\omega }}\times {\wt F_\omega })$ whose number of representatives corresponds to the number of equivalent completions of
$\o\omega _j\times \omega _j$~.

\item $G^{(2n)} ({\wt F_{\o\omega }}\times {\wt F_{\omega} }) $ acts on 
$P^{(2n)} ({\wt F^1_{\o\omega }}\times {\wt F^1_\omega }) $ by conjugation in such a way that the number of conjugates of 
$P^{(2n)} ({\wt F_{\o\omega^1_j }}\times {\wt F_{\omega ^1_j}}) $ in 
$G^{(2n)} ({\wt F_{\o\omega_{j,m_j} }}\times {\wt F_{\omega _{j,m_j}}}) $ is the index of the normalizer of $P^{(2n)} ({\wt F_{\o\omega^1_j }}\times {\wt F_{\omega^1_j} }) $~,
\[ \L| G^{(2n)} ({\wt F_{\o\omega_{j,m_j} }}\times {\wt F_{\omega_{j,m_j}} }) 
: N_{G^{(2n)}} (P^{(2n)} ({\wt F_{\o\omega ^1_j}}\times {\wt F_{\omega^1_j} })) \R|=j\]
where
\Bi
\item $G^{(2n)} ({\wt F_{\o\omega_{j,m_j} }}\times {\wt F_{\omega_{j,m_j}} }) $ is
$G^{(2n)} ({\wt F_{\o\omega }}\times {\wt F_{\omega} }) $ restricted to the conjugacy class representative $g^{(2n)}\RL[j,m_j]$~;

\item $N_{G^{(n)}} (P^{(2n)} ({\wt F_{\o\omega^1_j }}\times {\wt F_{\omega^1_j} })) $ is the normalizer of
 $P^{(2n)} ({\wt F_{\o\omega ^1_j}}\times {\wt F_{\omega^1_j} }) $ restricted to the $j$-th irreducible subextension ${\wt F_{\o\omega^1_j }}\times {\wt F_{\omega^1_j} }$~;

\item $j$ is the global residue degree of the $j$-th extension class $\wt \omega _j$ (or $\wt {\o \omega} _j$~).
\Ei\Ee
}
\vskip 11pt

\bpr
\Be
\item As the bilinear affine semigroup  $G ^{(2n)} ({\wt F_{\o\omega }}\times {\wt F_\omega }) $ is assumed to be reductive, each conjugacy class representative of 
$G ^{(2n)} ({\wt F_{\o\omega }}\times {\wt F_\omega }) $ is isogeneous to a direct product of tori (or of simple groups).  But, taking into account the projective toroidal compactification of the conjugacy classes 
$g^{(2n)}\RL[j,m_j]\in G ^{(2n)} ({\wt F_{\o\omega }}\times {\wt F_\omega }) $~, introduced in section 3.3, we shall assume that $g^{(2n)}n\RL[j,m_j]$ is isomorphic to a direct product of tori, which is the product, right by left, of algebraic (semi)tori.

So, $G ^{(2n)} ({\wt F_{\o\omega }}\times {\wt F_\omega }) \simeq \prod_n (\SS^1_R\times \SS^1_L)$~, which is the product, right by left, of algebraic (semi)tori,is considered as the $n$-fold product of a right affine semigroup  $\SS^1_R$ localized in the lower half space by its left equivalent $\SS^1_L$ localized in the upper half space.

But, as $G ^{(2n)} ({\wt F_{\o\omega }}\times {\wt F_\omega }) $ is defined over the product 
$ {\wt F_{\o\omega }}\times {\wt F_\omega } $~, 
and as
$G ^{(2n)} ({\wt F_{\o\omega }}\times {\wt F_\omega }) $ 
is assumed to act on
$P ^{(2n)} ({\wt F_{\o\omega^1 }}\times {\wt F_{\omega^1} }) $ 
by conjugation, a conjugacy class representative $g^{2n}\RL[j,m_j]\in G ^{(2n)} ({\wt F_{\o\omega }}\times {\wt F_\omega }) $ corresponds to the biextension $ {\wt F_{\o\omega_{j,m_j }}}\times {\wt F_{\omega_{j,m_j}} } $~.  Then, we have the injective morphism:
\[ Fg: \quad \{ {\wt F_{\o\omega_{j,m_j} }}\times {\wt F_{\omega _{j,m_j}}}\}
\To \{g^{(2n)}\RL[j,m_j]\]
from the set of equivalent representatives $ {\wt F_{\o\omega _{j,m_j}}}\times \wt F_{\omega _{j,m_j}}$~, associated with 
 all places $\o\omega _j\times \omega _j$ 
of $F_R\times F_L$~, to the set of corresponding conjugacy class representatives 
$g^{(2n)}\RL[j,m_j]$ of $G ^{(2n)} ({\wt F_{\o\omega }}\times {\wt F_\omega }) $~, such that
$\prod_n ( {\wt F_{\o\omega_{j,m_j} }}\times {\wt F_{\omega_{j,m_j}} }) 
\simeq g^{(2n)}\RL[j,m_j]$~.

\item $G ^{(2n)} ({\wt F_{\o\omega_{j,m_j} }}\times {\wt F_{\omega_{j,m_j }}}) $  acts on
$P ^{(2n)} ({\wt F_{\o\omega^1_j }}\times {\wt F_{\omega^1_j} }) $ by conjugation, in the sense that, if 
$p_{j_R}\in P^{(2n)}(\wt F_{\o\omega ^1_j})$~, $p_{j_L}\in P^{(2n)}(\wt F_{\omega ^1_j})$ and 
$x_{g_{j_R}}\in G ^{(2n)} ({\wt F_{\o\omega_{j,m_j }}} )$~, 
$x_{g_{j_L}}\in G ^{(2n)} ({\wt F_{\omega_{j,m_j} }} ) $~, then the bielement 
$x_{g_{j_R}} \cdot x_{g_{j_L}}\in g ^{(2n)}\RL [j,m_j]$ is said to be conjugate of the bielement 
$x'_{g_{j_R}} \cdot x'_{g_{j_L}}\in g ^{(2n)}\RL [j,m_j]$ if:
\[ x'_{g_{j_R}} \cdot x'_{g_{j_L}}= p_{j_R}\cdot x_{g_{j_R}} \cdot  x_{g_{j_L}}
\cdot p_{j_L} \;.\]
On the other hand, as a extension $\wt F_{\omega _{j,m_j}}$ (resp. $\wt F_{\o\omega _{j,m_j}})$~, is cut into $j'$ irreducible equivalent subextensions 
$\wt F_{\omega^{j'} _{j,m_j}}$ (resp. $\wt F_{\o\omega ^{j'}_{j,m_j}})$ according to section 1.1.5, the conjugacy class representative $g^{(2n)}\RL[j,m_j]$ will also be cut into ``~$j'$~'' equivalent conjugacy class subrepresentatives $g^{(2n)}\RL[j,m_j,j']$~, $1\le j'\le j$~, with respect to 1.

So, we have the corresponding injective morphism:
\begin{align*}
F_{g_{j'}}:\quad
& \wt F_{\o\omega _{j,m_j}}  \times \wt F_{\omega _{j,m_j}}= \{ 
\wt F_{\o\omega ^{j'}_{j,m_j} } \times \wt F_{\omega^{j'} _{j,m_j}}\}^j_{j'=1}\\
& \qquad \To \quad g^{(2n)}\RL[j,m_j] = \{ g^{(2n)}\RL[j,m_j,j']^j_{j'=1}\}\;.\end{align*}

This implies that the conjugation of
$G ^{(2n)} ({\wt F_{\o\omega_{j,m_j} }}\times {\wt F_{\omega_{j,m_j }}}) $ 
by $P ^{(2n)} ({\wt F_{\o\omega^1_{j} }}\times {\wt F_{\omega^1_{j }}}) $ 
sends the bielement
$x_{g_{j_R}}\cdot x_{g_{j_L}}\in g^{(2n)}\RL[j,m_j,j^{''}]$~, $1\le j^{''}\le j$~, of the $j^{''}$-th conjugacy class 
subrepresentative of $g^{(2n)}\RL[j,m_j]$ into the bielement $x'_{g_{j_R}}\cdot
x'_{g_{j_L}}$ of the $j'$-th conjugacy class subrepresentative of $g^{(2n)}\RL[j,m_j]$~.

Finally, the number of conjugacy class subrepresentatives of $g^{(2n)}\RL[j,m_j]$ is exactly the number of conjugates of $P^{(2n)}(\wt F_{\o\omega ^1_j}\times \wt F_{\omega ^1_j})$ in
$G^{(2n)}(\wt F_{\o\omega _{j,m_j}}\times \wt F_{\omega _{j,m_j}})
\equiv g^{(2n)}\RL[j,m_j]$~.

So the index of the normalizer 
$N_{G^{(2n)}} (P^{(2n)}( \wt F_{\o\omega ^1_j}\times \wt F_{\omega ^1_j})  )$ of
$P^{(2n)}(\wt F_{\o\omega ^1_j}\times \wt F_{\omega ^1_j})$ in\linebreak
$G^{(2n)}(\wt F_{\o\omega _{j,m_j}}\times \wt F_{\omega _{j,m_j}})$ corresponds to the number of conjugacy class subrepresentatives of 
$G^{(2n)}(\wt F_{\o\omega _{j,m_j}}\times \wt F_{\omega _{j,m_j}})$~.

Consequently, we have
\[ \L| G^{(2n)} ({\wt F_{\o\omega_{j,m_j} }}\times {\wt F_{\omega_{j,m_j} }}) 
: N_{G^{(2n)}} (P^{(2n)} ({\wt F_{\o\omega^1_j }}\times{\wt F_{\omega ^1_j}})) \R|=j\]
where $j$ is also the global residue degree of the $j$-th biplace $\o\omega _j\times\omega _j$ taking into account the inverse morphism $F^{-1}_{g_{j'}}$~.\epr
\Ee
\vskip 11pt 


\subsubsection{Proposition}

{\em Let $g^{(2n)}\RL[j,m_j]$~, $m_j\in\nit$~, $m_j\ge 0$~, denote the $j$-th diagonal or off-diagonal conjugacy class representative  of the bilinear algebraic semigroup $G^{(2n)}({\wt F_{\o\omega }}\times {\wt F_\omega })$ and let $g^{(2n)}\RL[j,m_j,j']$~, $1\le j'\le j$~, be the $j'$-th conjugacy subclass representative of $g^{(2n)}\RL[j,m_j]$~.

Let $D^{(2n)} ({\wt F^D_{\o\omega }}\times {\wt F^D_\omega })$ denote the bilinear algebraic semigroup $G^{(2n)}({\wt F_{\o\omega }}\times {\wt F_\omega })$ restricted to the product of the sets:
\[ {\wt F^D_{\o\omega }}=
\{\wt F_{\omega _1},\cdots,\wt F_{\omega _j},\cdots,\wt F_{\omega _r}\}
\qquad \text{and} \qquad \wt F^D_\omega =
\{\wt F_{\o\omega _1},\cdots,\wt F_{\o\omega _j},\cdots,\wt F_{\o\omega _r}\}
\]
with $m_j=0$~.

Then:
\Be
\item the $j$-th equivalent off-diagonal conjugacy class representatives 
$g^{(2n)}\RL[j,m_j]\in G^{(2n)} ({\wt F_{\o\omega }}\times {\wt F_\omega })$ are generated from the diagonal 
conjugacy class representatives $g^{(2n)}\RL[j,m_j=0]\in D^{(2n)} ({\wt F^D_{\o\omega }}\times {\wt F^D_\omega })$ under the product of the actions $\WW_{\wt F_{\o\omega _j}}\times \WW_{\wt F_{\omega_j} }$ of the Weyl groups $\WW_{\wt F_{\o\omega }}$ and $\WW_{\wt F_{\omega} }$~.
\item the Weyl groups $\WW_{\wt F_{\o\omega }}$ and $\WW_{\wt F_\omega }$ acting in a nilpotent way respectively
 on the $j$-th right and left diagonal conjugacy class representatives $g^{(2n)}_R[j,m_j=0]$ and $g^{(2n)}_L[j,m_j=0]\in g^{(2n)}\RL[j,m_j=0]$ are the Coxeter groups $A_{j-1}$ (and possibly $B_{j-1}$ and $D_{j-1}$~).
\Ee}
\vskip 11pt 

\bpr
\Be
\item The equivalent off-diagonal conjugacy class representatives $g^{(2n)}\RL[j,m_j]$~, corresponding to the $\wt B_{\o\omega _{j,m_j}}\otimes \wt B_{\omega _{j,m_j}}$-subbisemimodules $\wt M_{\o\omega _{j,m_j}}\otimes \wt M_{\omega _{j,m_j}}$ according to section 2.4.1, proposition 2.1.4 and definition 2.2.4, are isomorphic to the products, right by left, of $2n$ dimensional maximal semitori as it will be seen below.

So, the generation of the off-diagonal conjugacy class representatives $g^{(2n)}\RL[j,m_j]$ from the diagonal conjugacy class representatives $g^{(2n)}\RL[j,m_j=0]$~, which are also isomorphic to the products, right by left, of $2n$-dimensional semitori, is obtained in a classical way under the product $\WW_{\wt F_{\o\omega _j}}\times \WW_{\wt F_{\omega _j}}$ of the Weyl groups.
\item From point 2), it results that the Coxeter group $A_{j-1}$ (and $B_{j-1}$ and $D_{j-1}$~) is directly related to the integer $j$ which corresponds to the number of equivalent conjugacy subclasses $g^{(n)}\RL[j,m_j=0,j']$ of the diagonal conjugacy class $g^{(n)}\RL[j,m_j=0]$~.\epr
\Ee
\vskip 11pt 


\subsubsection{Proposition}

{\em
Let 
\[ \out( G^{(2n)}({\wt F_{\o\omega }}\times {\wt F_\omega }) )
= \Aut (G^{(2n)}({\wt F_{\o\omega }}\times {\wt F_\omega }) 
\Big/ \Int (G^{(2n)}({\wt F_{\o\omega }}\times {\wt F_\omega }) )\]
be the group of Galois automorphisms of the reductive bilinear algebraic semigroup\linebreak
$G^{(2n)}({\wt F_{\o\omega }}\times {\wt F_\omega })  $ where $\Int (G^{(2n)}({\wt F_{\o\omega }}\times {\wt F_\omega }) )$ denotes the group of Galois inner automorphisms.

Let $\Aut( P^{(2n)}({\wt F_{\o\omega^1 }}\times {\wt F_{\omega^1} }) )$ denote the group of the Galois automorphisms of the\linebreak bilinear parabolic subsemigroup $ P^{(2n)}({\wt F_{\o\omega^1 }}\times {\wt F_{\omega ^1}})$~.

Then, it is proved that:
\[ \Int (G^{(2n)}({\wt F_{\o\omega }}\times {\wt F_\omega }) )=
\Aut (P^{(2n)}({\wt F_{\o\omega^1 }}\times {\wt F_{\omega^1} }) )\]
which implies that explicit irreducible $2n$-dimensional representations exist:
\begin{align*}
\Rep^{(2n)}_{I_{F\RL}}: \quad I_{F_R}\times I_{F_L}
&\To P^{(2n)}({\wt F_{\o\omega ^1}}\times {\wt F_{\omega^1} }) \;, \\
\Rep^{(2n)}_{I^D_{F\RL}}: \quad I^D_{F_R}\times I^D_{F_L}
&\To P^{(2n)}({\wt F^D_{\o\omega^1 }}\times {\wt F^D_{\omega^1} }) \;, \end{align*}
from the product, right by left, of inertia subgroups to the corresponding bilinear parabolic subsemigroups where:
\Bi
\item $I_{F_L}:\bigoplus_j I_{F^D_{\omega _j}}\bigoplus_{m_j}I_{F_{\omega _{j,m_j}}}$

\item $I^D_{F_L}=\bigoplus_j I_{F^D_{\omega _j}} $ according to section 1.1.7.

\Ei
}
\vskip 11pt 

\bpr
\Be
\item The bilinear parabolic subsemigroups $P^{(2n)}({\wt F_{\o\omega ^1}}\times {\wt F_{\omega^1} }) $
and $P^{(2n)}({\wt F^D_{\o\omega ^1}}\times {\wt F^D_{\omega^1}}) $ are defined respectively over the products
$\wt F_{\o\omega ^1_{j,m_j}}\times \wt F_{\omega ^1_{j,m_j}}$ and $\wt F_{\o\omega ^1_{j}}\times \wt F_{\omega ^1_{j}}$ of irreducible central extensions which constitute the set of fixed points respectively of the products 
$\wt F_{\o\omega _{j,m_j}}\times \wt F_{\omega _{j,m_j}}$ and 
$\wt F_{\o\omega _{j}}\times \wt F_{\omega _{j}}$ of the pseudo-ramified extensions of
$ G^{(2n)}({\wt F_{\o\omega }}\times {\wt F_\omega }) $~.

Then, the group of automorphisms $\Aut( P^{(2n)}({\wt F_{\o\omega ^1}}\times {\wt F_{\omega^1} }) )$ of the reductive bilinear semigroup $ P^{(2n)}({\wt F_{\o\omega^1 }}\times {\wt F_{\omega^1} }) $ is the group of automorphisms of the fixed points of $ G^{(2n)}({\wt F_{\o\omega }}\times {\wt F_\omega }) $~.

Indeed, the map 
\[ \sigma _{ P^{(2n)}} :
\quad  P^{(2n)}({\wt F_{\o\omega^1 }}\times {\wt F_{\omega ^1}}) \To
 P^{(2n)}({\wt F_{\o\omega ^1}}\times {\wt F_{\omega ^1}}) \]
given by
\[ \sigma _{P^{(2n)}}(x_{p_R}\cdot x_{p_L}) = p_R\ x_{p_R}\ x_{p_L}\ p_L=x'_{p_R}\ x'_{p_L}\]
for all  $x_{p_R}\in P^{(2n)}({\wt F_{\o\omega ^1}})$~, 
$x_{p_L}\in P^{(2n)}({\wt F_{\omega^1 }})$~, \\
is an automorphism of $P^{(2n)}({\wt F_{\o\omega^1 }}\times {\wt F_{\omega^1} }) $ induced by $p_R\in P^{(2n)}({\wt F_{\o\omega ^1}})$ and $p_L\in P^{(2n)}({\wt F_{\omega^1 }})$~, in such a way that $\sigma _{P^{(2n)}}\in 
\Aut (P^{(2n)}({\wt F_{\o\omega^1 }}\times {\wt F_{\omega^1} }) )$~.

Similarly, the map
\[ \sigma^{\rm int} _{ G^{(2n)}} :
\quad  G^{(2n)}({\wt F_{\o\omega }}\times {\wt F_\omega }) \To 
G^{2n)}({\wt F_{\o\omega }}\times {\wt F_\omega }) \]
given by
\[ \sigma ^{\rm int}_{G^{(2n)}}({p_R}\ x_{g_R}\ x_{g_L}\ p_L) = x_{g_R}\cdot x_{g_L}\;,\]
for $x_{g_R}\in G^{(2n)}({\wt F_{\o\omega }})$ and
$x_{g_L}\in G^{(2n)}({\wt F_{\omega }})$~, and verifying $x_{g_R}\cdot x_{g_L}=x'_{p_R}\cdot x'_{p_L}$~,\\
is an inner automorphism of $G^{(2n)}({\wt F_{\o\omega }}\times {\wt F_\omega }) $ induced by $p_R\in P^{(2n)}({\wt F_{\o\omega ^1}})$ and by $p_L\in P^{(2n)}({\wt F_{\omega^1 }})$~, in such a way that $\sigma ^{\rm int}_{G^{(2n)}}= \sigma _{P^{(2n)}}$ for $\sigma ^{\rm int}_{G^{(2n)}}\in
\Int (G^{(2n)}({\wt F_{\o\omega }}\times {\wt F_\omega }) ) $~.

This leads to:
\[ \Int (G^{(2n)}({\wt F_{\o\omega }}\times {\wt F_\omega }) )
= \Aut (P^{(2n)}({\wt F_{\o\omega^1 }}\times {\wt F_{\omega^1} }) )\;.\]

\item On the other hand, the bilinear parabolic subsemigroup
\[ P^{(2n)}({\wt F_{\o\omega ^1}}\times {\wt F_{\omega^1} }) =
T_n^t ({\wt F_{\o\omega ^1}}) \times T_n ({\wt F_{\omega^1 }})\]
is constructed from the product of Borel subgroups 
$T_n^t ({\wt F_{\o\omega^1 }}) $ and 
$T_n ({\wt F_{\omega ^1}})$ which are solvable.

And, as it is known \cite{Ser2} that every solvable finite group is Galois, it is clear that
$\Aut(P^{(2n)} ({\wt F_{\o\omega^1 }} \times  {\wt F_{\omega ^1}}))$ constitutes the $n$-dimensional analogue of the product, right by left, of inertia subgroups.  Finally, as
$\Aut(P^{(2n)} ({\wt F_{\o\omega^1 }} \times  {\wt F_{\omega^1 }}))$ is isomorphic to
$P^{(2n)} ({\wt F_{\o\omega ^1}} \times  {\wt F_{\omega^1 }})$~, the bilinear parabolic subsemigroup
$P^{(2n)} ({\wt F_{\o\omega^1 }} \times  {\wt F_{\omega ^1}})$ constitutes a $2n$-dimensional representation
$\Rep^{(2n)}_{I_{F\RL}}$ of the product $I_{F_R}\times I_{F_L}$ of inertia subgroups in $\Gal(\wt F^{ac}_R\big/k)\times \Gal(\wt F^{ac}_L\big/ k)$~.\epr
\Ee
\vskip 11pt 

\subsubsection{Proposition}

{\em Let $P^{(2n)} ({\wt F_{\o\omega^1 }} \times {\wt F_{\omega^1 }})$ be the bilinear parabolic affine subsemigroup  of the smooth reductive bilinear affine semigroup  $G^{(2n)} ({\wt F_{\o\omega }} \times  {\wt F_{\omega }})$~.  Then, $P^{(2n)} ({\wt F_{\o\omega ^1}} \times  {\wt F_{\omega^1 }})$ can be considered as the unitary irreducible $2n$-dimensional representation space of the algebraic bilinear semigroup of matrices
$\GL_n({\wt F_{\o\omega }}  \times  {\wt F_{\omega }})$~.
}
\vskip 11pt 

\bpr  According to proposition 2.4.4, 
\[ \Int(G^{(2n)} ({\wt F_{\o\omega }}  \times  {\wt F_{\omega }}))
= \Aut ( P^{(2n)} ({\wt F_{\o\omega ^1}}  \times {\wt F_{\omega^1 }}))\]
since we have:
\[ \sigma ^{\rm int}_{G^{(n)}}({p_R}\ x_{g_R}\ x_{g_L}\ p_L) = x_{g_R}\cdot x_{g_L}\;,\]
for ${p_R}\in P^{(2n)}({\wt F_{\o\omega^1 }})$~,
${p_L}\in P^{(2n)}({\wt F_{\omega^1 }})$~, $x_{g_R}\in 
G^{(2n)}({\wt F_{\o\omega }})$~,
$x_{g_L}\in G^{(2n)}({\wt F_{\omega }})$~.

So, $P^{(2n)} ({\wt F_{\o\omega ^1}}  \times {\wt F_{\omega ^1}})$ is the bilinear isotropy subgroup of
$G^{(2n)} ({\wt F_{\o\omega }}  \times {\wt F_{\omega }})$ since it is the subgroup fixing the bielements of
$G^{(2n)} ({\wt F_{\o\omega }}  \times  {\wt F_{\omega }})$~.

Consequently, $P^{(2n)} ({\wt F_{\o\omega^1 }}  \times  {\wt F_{\omega^1 }})$ is a unitary representation of
$\GL_n({\wt F_{\o\omega }}  \times  {\wt F_{\omega }})$~.\epr
\vskip 11pt


\subsubsection{Proposition}  

{\em Let $D^{(2n)}({\wt F^D_{\o\omega _\oplus}}\times {\wt F^D_{\omega_\oplus} })$ denote the diagonal bilinear algebraic subsemigroup of  $G^{(2n)}({\wt F_{\o\omega_\oplus }}\times {\wt F_{\omega_\oplus} })$~.
Then, there exist explicit irreducible $2n$-dimensional	representations:
\begin{alignat*}{3}
\Rep^{(2n)}_{\Gal_{\wt F\RL}}
&: \quad \Gal(\wt F^{ac}_R\big/k) \times \Gal(\wt F^{ac}_L\big/k) &\To &
G^{(2n)}({\wt F_{\o\omega_\oplus }}\times {\wt F_{\omega_\oplus} }) \;, \\
\Rep^{(2n)}_{\Gal^D\wt F\RL}
&: \quad \Gal^D(\wt F^{ac}_R\big/k) \times \Gal^D(\wt F^{ac}_L\big/k) &\To& 
D^{(2n)}({\wt F^D_{\o\omega_\oplus }}\times {\wt F^D_{\omega_\oplus} })\; 
%
%
 \end{alignat*}
where $\wt F^D_{\omega_\oplus }=\bigoplus\limits_j\wt F_{\omega _j}$
(resp. $\wt F^D_{\o\omega_\oplus }=\bigoplus\limits_j\wt F_{\o\omega _j}$~) are constructed from $\wt F_\omega $ (resp. $\wt F_{\o\omega }$~) introduced in section 1.1.6.
}
%
%
%
\vskip 11pt 

\bpr If a group $G$ is quasi split, in the sense that it contains a Borel subgroup, then, the homomorphism
\[ \Gal (\o\QQ\Big/ \QQ) \To \out (G)=\Aut(G)\Big/ \Int(G)\]
always exists \cite{Art3}.

Transposed into a bilinear context for the bilinear reductive algebraic semigroup $G^{(2n)}({\wt F_{\o\omega }}\times {\wt F_\omega })$~, this homomorphism becomes:
\[ \Gal (\wt F^{ac}_R\big/ k)\times \Gal (\wt F^{ac}_L\big/ k)
\To \out (G^{(2n)}({\wt F_{\o\omega_\oplus }}\times {\wt F_{\omega_\oplus} }))\]
where
\[ \out(G^{(2n)}({\wt F_{\o\omega_\oplus }}\times {\wt F_{\omega_\oplus} }))
= \Aut (G^{(2n)}({\wt F_{\o\omega_\oplus }}\times {\wt F_{\omega_\oplus} }))\Big/
\Int (G^{(2n)}({\wt F_{\o\omega_\oplus }}\times {\wt F_{\omega_\oplus} }))\;.\]
Indeed, as $G^{(2n)}({\wt F_{\o\omega_\oplus }}\times {\wt F_{\omega_\oplus } })$ is generated from the product of Borel subgroups 
which are solvable, it is of Galois type according to proposition 2.4.3.  And, 
$\Aut(G^{(2n)}({\wt F_{\o\omega_\oplus  }}\times {\wt F_{\omega_\oplus }}))$ 
constitutes the $2n$-dimensional analogue of the product 
$\Gal (\wt F^{ac}_R\big/k)\times \Gal (\wt F^{ac}_L\big/ k)$ of Galois groups.  But, as the bilinear algebraic semigroup 
$G^{(2n)}({\wt F_{\o\omega_\oplus }}\times {\wt F_{\omega_\oplus }})$ is pseudo-ramified, it is 
$\out (G^{(2n)}({\wt F_{\o\omega_\oplus }}\times {\wt F_{\omega_\oplus }}))$ which must be considered as the $2n$-dimensional analogue 
of $\Gal(\wt F^{ac}_R\big/ k)\times \Gal(\wt F^{ac}_L\big/k)$ since 
$\Int (G^{(2n)}({\wt F_{\o\omega_\oplus }}\times {\wt F_{\omega_\oplus }}))= 
\Aut(P^{(2n)}({\wt F_{\o\omega^1_\oplus}} \times {\wt F_{\omega ^1_\oplus }}))$ 
according to proposition 2.4.3.

So, \[\out (G^{(2n)}({\wt F_{\o\omega_\oplus  }}\times {\wt F_{\omega _\oplus }}))
\simeq \prod_n(\Gal(\wt F^{ac}_R\big/ k)\times \Gal (\wt F^{ac}_L\big/k))\]
which implies the commutative diagram:\\
\centerline{\unitlength=1cm
\begin{picture}(10.5,3)
\put(0,1.2){{$\begin{array}{ccc}
\Gal(\wt F^{ac}_R\big/k) \times \Gal(\wt F^{ac}_L\big/k)
&\To & \out(G^{(2n)}({\wt F_{\o\omega_\oplus  }}\times {\wt F_{\omega _\oplus }}))\\[-6pt]
&& \begin{CD} @VV{\wr}V\end{CD}\\[-6pt]
&& G^{(2n)}({\wt F_{\o\omega_\oplus  }}\times {\wt F_{\omega_\oplus } })\end{array}$}}
\put(5.6,1.4){\vector(2,-1){1.5}}
\end{picture}
}

Thus, as $\out ( G^{(2n)}({\wt F_{\o\omega_\oplus  }}\times {\wt F_{\omega_\oplus } }))$ is isomorphic to 
$G^{(2n)}({\wt F_{\o\omega_\oplus  }}\times {\wt F_{\omega_\oplus } })$~, the irreducible $2n$-dimensional representation 
$\Rep^{(2n)}_{\Gal_{\wt F\RL}}$ of
$\Gal(\wt F^{ac}_R\big/k) \times \Gal(\wt F^{ac}_L\big/k)$ is given by the reductive bilinear algebraic semigroup
$G^{(2n)}({\wt F_{\o\omega_\oplus  }}\times {\wt F_{\omega_\oplus } })$~.
\epr\vskip 11pt

\subsubsection{Proposition} 

{\em
Let 
\[ W^{ab}_{F_R}\times W^{ab}_{F_L}=
\Gal(\wt {\dot F}^{ac}_R\big/k) \times \Gal(\wt {\dot F}^{ac}_L\big/k)
\subset \Gal(\wt F^{ac}_R\big/k) \times \Gal(\wt F^{ac}_L\big/k)\]
be the product of global Weil groups as introduced in definition 1.1.9.

Then, the $2n$-dimensional irreducible representation of the product of global Weil groups is given by:
\[ \Irr\ W^{(2n)}_{F\RL}:\quad W^{ab}_{F_R}\times W^{ab}_{F_L}
\To G^{(2n)}({\wt F_{\o\omega_\oplus  }}\times {\wt F_{\omega_\oplus } })\;.\]
}\vskip 11pt 

\bpr This directly follows from $\Rep^{(2n)}_{\Gal_{\wt F\RL}}$ introduced in proposition 2.4.6, implying that $\wt F_{\o\omega _\oplus}$ and $\wt F_{\omega _\oplus}$
bear here on extensions characterized by extension degrees\linebreak $d=0\mod N$~.\epr
\vskip 11pt 


\subsection{Trace formulas}

\subsubsection{Lemma}

{\em Let $ G^{(2n)}({F_{\o\omega }}\times {F_\omega })$ be a reductive bilinear complete semigroup.  As it is solvable, the set of its conjugacy classes forms an increasing sequence:
\[ g\RL^{(2n)}[1,m_1] \subset g\RL^{(2n)}[2,m_1]\subset \cdots \subset 
g\RL^{(2n)}[j,m_j]\subset \cdots \subset g\RL^{(2n)}[r,m_r]\]
characterized by the increasing ranks of their conjugacy class representatives $g\RL^{(2n)}[j,m_j]$~, $1\le j\le r$~.
}
\vskip 11pt 

\bpr According to definition 2.3.1 and proposition 2.1.4, the pseudo-ramified lattice bisemispace 
\[X_{S\RL}= \GL_n({\widetilde F_{R }}\times {\widetilde F_L})\Big/ \GL_n((\ZZ\big/N\  \ZZ)^2)\]
is homomorphic to the pseudo-ramified $B_{\o\omega _\oplus }\times B_{\omega_\oplus } $-bisemimodule $M_{R_\oplus }\otimes M_{L_\oplus }$ developing (see proposition 2.1.4) as follows:
\[ M_{R_\oplus }\otimes M_{L_\oplus } = \bigoplus_j \bigoplus_{m_j} \L(M_{\o\omega _{j,m_j}}\otimes M_{\omega _{j,m_j}}\R)\;.\]
As the cosets of $X\SRL$ correspond to the conjugacy class representatives of 
$G^{(2n)}({F_{\o\omega }}\times {F_\omega })$  and as $(M_R\otimes M_L)$ constitutes the representation space
 of $\GL_n({F_{\o\omega }}\times {F_\omega })$~,
each conjugacy class representative $g^{(2n)}\RL[j,m_j]$ is in fact a 
$\GL_n({F_{\o\omega_{j,m_j} }}\times {F_{\omega_{j,m_j}} })$-subbisemimodule $M_{\o\omega _{j,m_j}}\otimes M_{\omega _{j,m_j}}$ characterized by a rank $r_{\o\omega _j\times \omega _j}=(j\cdot m^{(j)}\cdot N)^{2n}$ according to proposition 2.3.4.

If we consider the $m_{j+1}$-th representative $g^{(2n)}\RL[j+1,m_{j+1}]$ of the $(j+1)$-th conjugacy class of 
$G^{(2n)}({F_{\o\omega }}\times {F_\omega })$~, we see that it has a rank given by
\[ r_{\o\omega _{j+1}\times \omega _{j+1}}=(j+1)^{2n}\cdot N^{2n}\cdot (m^{(j+1)})^{2n}\]
which is superior to the rank $r_{\o\omega _j\times\omega _j}$ of $g^{(2n)}\RL[j,m_j]$ ~.\epr
\vskip 11pt

\subsubsection{Tower of sums of conjugacy class representatives}

Remark first that a conjugacy class representative $g^{(2n)}\RL[j,m_j]$ of $G^{(2n)}({F_{\o\omega }}\times {F_\omega })$ is in fact a bilinear complete subsemigroup  noted 
$G^{(2n)}({F_{\o\omega_{j,m_j} }}\times {F_{\omega_{j,m_j} }})$~.  Then, we can introduce a tower of embedded bilinear 
complete subsemigroups  on $G^{(2n)}(F_{\o\omega _\oplus}\times F_{\omega _\oplus})$~:
\[G^{(2n)}(F_{\o\omega _1}\times F_{\omega _1}) \subset \cdots \subset 
\bigoplus_{j=1}^j\bigoplus_{m_j}G^{(2n)}(F_{\o\omega _{j,m_j}}\times F_{\omega _{j,m_j}} )\subset \cdots
\subset \bigoplus_{j=1}^r\bigoplus_{m_j}G^{(2n)}(F_{\o\omega _{j,m_j}}\times F_{\omega _{j,m_j}} )\]
corresponding to increasing sums of conjugacy class representatives of
$G^{(2n)}({F_{\o\omega }}\times {F_{\omega }} )$~.
\vskip 11pt 

\subsubsection{{\boldmath The bialgebra of bifunctions on $G^{(2n)}({F_{\o\omega }}\times {F_\omega })$}}

Now, we consider the set of all smooth continuous bifunctions on $G^{(2n)}({F_{\o\omega }}\times {F_\omega })$ which is a bilinear complete semigroup.  Consequently, we have to envisage the set $\widehat G^{(2n)}_R({F_{\o\omega }})$ of smooth continuous functions $\phi ^{(2n)}_{G_R}(x_{g_R})$~, $x_{g_R}\in G^{(2n)}({F_{\o\omega }}) $~, on
$G^{(2n)}({F_{\o\omega }})\simeq T^t_n ({F_{\o\omega }})$ and localized in the lower half space as well as the corresponding set $\widehat G^{(2n)}_L({F_{\omega }})$ of smooth continuous functions 
$\phi ^{(2n)}_{G_L}(x_{g_L})$~, $x_{g_L}\in G^{(2n)}({F_{\omega }}) $~, on
$G^{(2n)}({F_{\omega }})\simeq T_n ({F_{\omega }})$ and localized in the upper half space in such a way that we have on $G^{(2n)}({F_{\o\omega }}\times {F_{\omega }})$ 
the product
$\widehat G^{(2n)}({F_{\o\omega }})\times \widehat G^{(2n)}({F_{\omega }})$ 
of  the tensor products
$\phi ^{(2n)}_{G_R}(x_{g_R})\otimes \phi ^{(2n)}_{G_L}(x_{g_L})$ of the smooth continuous functions 
$\phi ^{(2n)}_{G_R}(x_{g_R})$ and $\phi ^{(2n)}_{G_L}(x_{g_L})$~.  

$\widehat G^{(2n)}_R(F_{\o\omega })$ (resp. $\widehat G^{(2n)}_L(F_{\omega })$~) is the coordinate (semi)ring $k[G^{(2n)}(F_{\o\omega })]$ (resp. $k[G^{(2n)}(F_{\omega })]$~) of $G^{(2n)}(F_{\o\omega })$ (resp. $G^{(2n)}(F_{\o\omega })$~).

But, as $G^{(2n)}({F_{\o\omega }}\times  {F_{\omega }})$ is partitioned into conjugacy classes, we can consider bifunctions on the conjugacy class representatives of 
$G^{(2n)}({F_{\o\omega }} \times {F_{\omega }})$~.
\\
So, let $\phi ^{(2n)}_{G_{j_R}}(x_{g_{j_R}})\otimes \phi ^{(2n)}_{G_{j_L}}(x_{g_{j_L}})$ denote a bifunction on
$G^{(2n)}({F_{\o\omega_{j,m_j}}}\times {F_{\omega_{j,m_j} }})$~.  Then, the set of smooth continuous bifunctions
$\{\phi ^{(2n)}_{G_{j_R}}(x_{g_{j_R}})\otimes \phi ^{(2n)}_{G_{j_L}}(x_{g_{j_L}})\}^r_{j=1}$ on
$G^{(2n)}({F_{\o\omega}}\times {F_{\omega }})$ forms a bialgebra
\[ \widehat G^{(2n)}({F_{\o\omega}}\times {F_{\omega }})=
\widehat G^{(2n)}({F_{\o\omega }})\otimes \widehat G^{(2n)}( {F_{\omega }})\]
where $\widehat G^{(2n)}({F_{\o\omega}})$ is the coalgebra of $\widehat G^{(2n)}({F_{\omega}})$~.
\vskip 11pt 

\subsubsection{{\boldmath The bialgebras $L^{1-1}\RL(G^{(2n)}({F_{\o\omega }}\times {F_\omega }))$ and
$L^{2}_{L\times L}(G^{(2n)}({F_{\omega }}\times {F_\omega }))$}}

The algebra (resp. coalgebra) $\widehat G^{(2n)}({F_{\omega }})$ (resp. $\widehat G^{(2n)}({F_{\o\omega }})$~) is the set of all continuous complex-valued measurable functions on $G^{(2n)}({F_{\omega }})$  (resp.
$G^{(2n)}({F_{\o \omega }})$~) satisfying
\begin{align*}
\int_{G^{(2n)}({F_{\omega }})}\L| \phi ^{(2n)}_{G_L}(x_{g_L})\R|\ dx_{g_L}<\infty \;, \\[11pt]
\text{(resp.}\quad
\int_{G^{(2n)}({F_{\o\omega }})}\L| \phi ^{(2n)}_{G_R}(x_{g_R}) \R|\ dx_{g_R}<\infty \;), \end{align*}
with respect to a unique Haar measure on
${G^{(2n)}({F_{\omega }})}$ (resp. ${G^{(2n)}({F_{\o\omega }})}$~): it is noted 
$L^1_L(G^{(2n)}({F_{\omega }}))$ (resp. $L^1_R(G^{(2n)}({F_{\o\omega }}))$~).

The bialgebra $\widehat G^{(2n)}({F_{\o\omega}}\times {F_{\omega }})$ of all continuous complex-valued measurable bifunctions $ \phi ^{(2n)}_{G_R}(x_{g_R}) \otimes  \phi ^{(2n)}_{G_L}(x_{g_L}) $ on 
$G^{(2n)}({F_{\o\omega}}\times {F_{\omega }})$ satisfying
\[ \int_{G^{(2n)}({F_{\o\omega }}\otimes {F_\omega })}
\L| \phi ^{(2n)}_{G_R}(x_{g_R}) \otimes
 \phi ^{(2n)}_{G_L}(x_{g_L}) \R|\ dx_{g_R}\ 
 dx_{g_L}<\infty \]
is also noted $L^{1-1}\RL(G^{(2n)}({F_{\o\omega}}\times {F_{\omega }}))$~.

On the other hand, let
\begin{align*}
i_{G_{R\to L}}\quad
L^1_R( G^{(2n)}({F_{\o\omega}})) & \To L^1_L( G^{(2n)}( {F_{\omega }}))\\[11pt]
\phi ^{(2n)}_{G_R}(x_{g_R}) &\To \phi ^{(2n)}_{G_L}(x_{g_L})\end{align*}
be the involution mapping each continuous function
$\phi ^{(2n)}_{G_R}(x_{g_R}) \in L^1_R( G^{(2n)}({F_{\o\omega}}))$~, restricted to the lower half space, to the corresponding continuous function $\phi ^{(2n)}_{G_L}(x_{g_L}) \in L^1_L( G^{(2n)}({F_{\omega}}))$~, restricted to the upper half space.

Then, the involution:
\begin{align*}
i_{G_{R\to L}} : \qquad
&L^1_R( G^{(2n)}( {F_{\o\omega }})) \times L^1_L( G^{(2n)}( {F_{\omega }}))=L^{1-1}\RL( G^{(2n)}( {F_{\o\omega }}\times {F_{\omega }}))\\
& \qquad \To L^2_{L\times L}( G^{(2n)}( {F_{\omega }}\times  {F_{\omega }}))\\[11pt]
& \phi ^{(2n)}_{G_R}(x_{g_R})\otimes \phi ^{(2n)}_{G_L}(x_{g_L}) \To 
\phi ^{(2n)}_{G_L}(x_{g_L})\otimes \phi ^{(2n)}_{G_L}(x_{g_L}) \end{align*}
has the property that the continuous bifunctions
$\phi ^{(2n)}_{G_R}(x_{g_R})\otimes \phi ^{(2n)}_{G_L}(x_{g_L}) \in L^{1-1}\RL( G^{(2n)}( {F_{\o\omega }}
\times  {F_{\omega }}))$ are transformed into the bifunctions
$ \phi ^{(2n)}_{G_L}(x_{g_L})\otimes \phi ^{(2n)}_{G_L}(x_{g_L}) \in  L^{2}_{L\times L}( G^{(2n)}( {F_{\omega }}
\times  {F_{\omega }}))$ satisfying
\[\int_{G^{(2n)}({F_\omega }\times {F_\omega })}
\L|  \phi ^{(2n)}_{G_L}(x_{g_L})\R|^2\ dx_{g_L}<\infty \;.\]
Thus, $L^2_{L\times L}( G^{(2n)}( {F_{\omega }}\times  {F_{\omega }})) $ is the standard space of square integrable continuous functions on $G^{(2n)}({F_{\omega}}\times {F_{\omega }})$
restricted to the upper half space and defined with respect to a fixed Haar measure.
\vskip 11pt 

\subsubsection{Proposition}

{\em The bialgebra $\widehat G^{(2n)}({F_{\o\omega}}\times {F_{\omega }})$ of smooth continuous complex-valued measurable bifunctions on the locally compact bilinear complete semigroup $G^{(2n)}({F_{\o\omega}}\times {F_{\omega }})$ is composed of the increasing sequence:
\[ \phi ^{(2n)}_{G_{1_R}}( x_{g_{1_R}} ) \otimes \phi ^{(2n)}_{G_{1_L}}( x_{g_{1_L}} )
\subset \cdots \subset
\phi ^{(2n)}_{G_{j_R}}( x_{g_{j_R}} )\otimes \phi ^{(2n)}_{G_{j_L}}( x_{g_{j_L}} )
\subset \cdots \subset
\phi ^{(2n)}_{G_{r_R}}( x_{g_{r_R}} )\otimes \phi ^{(2n)}_{G_{r_L}}( x_{g_{r_L}} )\]
of these bifunctions in such a way that the $j$-th bifunction 
$\phi ^{(2n)}_{G_{j_R}}( x_{g_{j_R}} )\otimes \phi ^{(2n)}_{G_{j_L}}( x_{g_{j_L}} )$ 
is defined on the closed bilinear  subsemigroup  $G^{(2n)}( F_{\o\omega _{j,m_j}}\times
F_{\omega _{j,m_j}})\equiv g^{(2n)}\RL[j,m_j]$ corresponding to the $j$-th conjugacy  class of 
$G^{(2n)}( {F_{\o\omega }}\times  {F_{\omega }})$~.}
\vskip 11pt 

\paragraph{Sketch of proof}:  This results from lemma 2.5.1 and from section 2.5.3.\epr
\vskip 11pt 

\subsubsection{Proposition}

{\em The bialgebra $\widehat G^{(2n)}({F_{\o\omega_\oplus }}\times {F_{\omega _\oplus }})
\equiv L^{1-1}\RL (
G^{(2n)}({F_{\o\omega_\oplus }}\times {F_{\omega_\oplus  }}))$
of bifunctions on $G^{(2n)}({F_{\o\omega_\oplus }}\times {F_{\omega_\oplus  }})$ is  composed of a tower:
\begin{align*}
\phi ^{(2n)}_{G_{1_R}}( x_{g_{1_R}} )\otimes \phi ^{(2n)}_{G_{1_L}}( x_{g_{1_L}} )
\subset \cdots &\subset 
\bigoplus^j_{j=1} \L( \phi ^{(2n)}_{G_{j_R}}( x_{g_{j_R}} )\otimes \phi ^{(2n)}_{G_{j_L}}( x_{g_{j_L}} )\R)\\[11pt]
& \subset \cdots \subset
\bigoplus^r_{j=1} \L( \phi ^{(2n)}_{G_{j_R}}( x_{g_{j_R}} )\otimes \phi ^{(2n)}_{G_{j_L}}( x_{g_{j_L}} )\R)\end{align*}
of the partial sums of these bifunctions
$\bigoplus^j_{j=1} \L( \phi ^{(2n)}_{G_{j_R}}( x_{g_{j_R}} )\otimes \phi ^{(2n)}_{G_{j_L}}( x_{g_{j_L}} )\R)$
on the corresponding sums of conjugacy class representatives of
$G^{(2n)}({F_{\o\omega}}\times {F_{\omega }})$~.
}
\vskip 11pt 

\paragraph{Sketch of proof}: This is a consequence of the introduction of a tower of sums of conjugacy class representatives of
$G^{(2n)}({F_{\o\omega_\oplus }}\times {F_{\omega_\oplus  }})$ given in section 2.5.2.
\vskip 11pt 

\subsubsection{{\boldmath The bialgebra $  L^{1-1}\RL (
G^{(2n)}(F^{nr}_{\o\omega}\times F^{nr}_{\omega }))$
of bifunctions on the pseudo-unramified bilinear complete semigroup $G^{(2n)}(F^{nr}_{\o\omega}\times F^{nr}_{\omega })$}}

Similarly, on the set of the conjugacy class representatives of the pseudo-unramified bilinear complete semigroup 
$G^{(2n)}(F^{nr}_{\o\omega}\times F^{nr}_{\omega })$~, an increasing sequence
\begin{align*}
\phi ^{(2n)}_{G^{nr}_{1_R}}( x^{nr}_{g_{1_R}} )\otimes \phi ^{(2n)}_{G^{nr}_{1_L}}( x^{nr}_{g_{1_L}} )
\subset \cdots &\subset 
 \phi ^{(2n)}_{G^{nr}_{j_R}}( x^{nr}_{g_{j_R}} )\otimes \phi ^{(2n)}_{G^{nr}_{j_L}}( x^{nr}_{g_{j_L}} )\\[11pt]
& \subset \cdots \subset
 \phi ^{(2n)}_{G^{nr}_{r_R}}( x^{nr}_{g_{r_R}} )\otimes \phi ^{(2n)}_{G^{nr}_{r_L}}( x^{nr}_{g_{r_L}} )\end{align*}
of smooth continuous complex-valued measurable  bifunctions can be introduced in such a way that they satisfy:
\[\int_{G^{(2n)}(F^{nr}_{\o\omega}\times F^{nr}_{\omega })}\
\L|
 \phi ^{(2n)}_{G^{nr}_{j_R}}( x^{nr}_{g_{j_R}} )\otimes \phi ^{(2n)}_{G^{nr}_{j_L}}( x^{nr}_{g_{j_L}} )\R|
dx^{nr}_{g_{j_R}}\ dx^{nr}_{g_{j_L}}<\infty \;.\]
Then, the bialgebra of these bifunctions on
$G^{(2n)}(F^{nr}_{\o\omega}\times F^{nr}_{\omega })$ is written
$\widehat G^{(2n)}(F^{nr}_{\o\omega}\times F^{nr}_{\omega })$ or
$L^{1-1}\RL( G^{(2n)}(F^{nr}_{\o\omega}\times F^{nr}_{\omega }))$~.
\vskip 11pt 

\subsubsection{{\boldmath The bialgebra $  L^{1-1}\RL (
P^{(2n)}({F_{\o\omega^1}}\times {F_{\omega ^1}}))$
of bifunctions on the complete parabolic semigroup $P^{(2n)}({F_{\o\omega^1}}\times {F_{\omega ^1}})$}}

As
$P^{(2n)}({F_{\o\omega^1}}\times {F_{\omega ^1}})$ is defined over the products
$F_{\o\omega^1 _{j,m_j}}\times F_{\omega ^1_{j,m_j}}$ of the irreducible completions of the completions
$F_{\o\omega _{j,m_j}}\times F_{\omega _{j,m_j}}$ at all infinite places
${\o\omega_j} \times {\omega _j}$ of $F_R\times F_L$~, a bialgebra 
$L^{1-1}\RL( P^{(2n)}({F_{\o\omega^1}}\times {F_{\omega^1 }}))$ of smooth continuous complex-valued measurable bifunctions
$
\phi ^{(2n)}_{P_{j_R}}( x_{p_{j_R }} )\otimes \phi ^{(2n)}_{P_{j_L}}( x_{p_{j_L}} )$~, $1\le j\le r$~, on
$P^{(2n)}({F_{\o\omega^1}}\times {F_{\omega^1 }})$ 
can be introduced in such a way that:
\[\int_{P^{(2n)}({F_{\o\omega^1}}\times {F_{\omega ^1}})}\
\L|
 \phi ^{(2n)}_{P_{j_R}}( x_{p_{j_R}} )\otimes \phi ^{(2n)}_{P_{j_L}}( x_{p_{j_L}} ) \R|
dx_{p_{j_R}}\ dx_{p_{j_L}}<\infty \;.\]
Each bifunction $ \phi ^{(2n)}_{P_{j_R}}( x_{p_{j_R}} )\otimes \phi ^{(2n)}_{P_{j_L}}( x_{p_{j_L}} ) $
is defined on the irreducible equivalence class representative 
$P^{(2n)}({F_{\o\omega^1_{j,m_j}}}\times {F_{\omega ^1_{j,m_j}}})$ characterized by a rank
$r_{P_{\o\omega _j\times \omega _j}}=(N\cdot m^{(j)})^{2n}$~.

The normalized rank $r_{/N_{P_{\o\omega _j\times \omega _j}}}$ of
$P^{(2n)}({F_{\o\omega^1_{j,m_j}}}\times {F_{\omega^1_{j,m_j} }})$ is defined by:
\[r_{/N_{P_{\o\omega _j\times \omega _j}}}
= r_{P_{\o\omega _j\times \omega _j}}\Big/ (N\cdot m^{(j)})^{2n}=1\;.\]
A character on $G^{(2n)}({F_{\o\omega}}\times {F_{\omega }})$ is then defined to be a bifunction
$\chi _{j\RL}= \phi ^{(2n)}_{P_{j_R}}( x_{p_{j_R}} )\otimes \phi ^{(2n)}_{P_{j_L}}( x_{p_{j_L}} ) $ from the irreducible equivalence class representative
$P^{(2n)}({F_{\o\omega^1_{j,m_j}}}\times {F_{\omega^1_{j,m_j} }})$ of
$G^{(2n)}({F_{\o\omega}}\times {F_{\omega }})$~, having as normalized rank 
$r_{/N_{P_{\o\omega _j\times \omega _j}}}$~, into $\cit$~.

And, the group of characters $\chi _{j\RL}$ of 
$P^{(2n)}({F_{\o\omega^1}}\times {F_{\omega^1 }})$ is written
$X( P^{(2n)}({F_{\o\omega^1}}\times {F_{\omega^1 }}),\cit)$~.
\vskip 11pt 

\subsubsection{Two types of trace formulas: an introduction}

Two kinds of trace formulas will now be developed \cite{Kot1,Art1}.  The first one is an adaptation of the Selberg-Arthur trace formula \cite{Sel,Art2} of an operator associated with the quotient of a semisimple Lie group and a discrete subgroup in such a way that this operator corresponds to a unitary representation of the semisimple Lie group $G(\Aa)$~.  The adaptation considered here concerns the trace of the operator associated with the bilinear parabolic semigroup
$P^{(2n)}({F_{\o\omega^1}}\times {F_{\omega^1 }})$ envisaged as the unitary representation of the complete bilinear semigroup $\GL_n({F_{\o\omega}}\times {F_{\omega }})$~.  This operator acts by convolution on the bialgebra 
$L^{1-1}\RL( G^{(2n)}(F^{nr}_{\o\omega}\times F^{nr}_{\omega }))$ of bifunctions on the pseudo-unramified bilinear complete semigroup 
$G^{(2n)}(F^{nr}_{\o\omega}\times F^{nr}_{\omega })$~.  As a result, this operator decomposes according to the unitary conjugacy classes of the pseudo-ramified bilinear algebraic semigroup
$G^{(2n)}({F_{\o\omega}}\times {F_{\omega }})$ and the corresponding trace formula occurs in the bialgebra
$L^{1-1}\RL( G^{(2n)}({F_{\o\omega_\oplus }}\times {F_{\omega _\oplus }}))$ of bifunctions on 
$G^{(2n)}({F_{\o\omega_\oplus }}\times {F_{\omega _\oplus }})$~.  Notice that the complete trace formula referring to the set of irreducible representations of
$G^{(2n)}({F_{\o\omega}}\times {F_{\omega }})$ will be envisaged in chapter 4.

Remark that a trace formula, being an operation of bilinear type since it refers to the matrix representation of an operator, is well adapted to the bilinear case considered here and, especially, to the second kind of trace formulas occurring directly in the bialgebra 
$L^{1-1}\RL( G^{(2n)}({F_{\o\omega_\oplus }}\times {F_{\omega_\oplus  }}))$ of bifunctions on the pseudo-ramified complete bilinear semigroup 
$ G^{(2n)}({F_{\o\omega}}\times {F_{\omega }})$~.  This latter case then leads to the Plancherel formula.

Consequently, the two types of trace formulas correspond to each other.
\vskip 11pt

Let us now start with the first trace formula referring to the Arthur-Selbeg trace formula.
\vskip 11pt 

\subsubsection{{\boldmath Unitary representation of  $ \GL_n({F_{\o\omega}}\times {F_{\omega }})$}}

Let $\GL_n(\ZZ^2)$ be the discrete subgroup of the pseudo-unramified bilinear complete semigroup 
$ \GL_n(F^{nr}_{\o\omega}\times F^{nr}_{\omega })$~.  We then consider the bialgebra
$L^{1-1}\RL( \GL_n(F^{nr}_{\o\omega}\times F^{nr}_{\omega })\big/ \GL_n(\ZZ^2))$ of bifunctions
$\phi ^{(2n)}_{G^{nr}_{j_R}}( x^{nr}_{g_{j_R }} )\otimes \phi ^{(2n)}_{G^{nr}_{j_L}}( x^{nr}_{g_{j_L}} )$ on the conjugacy class representatives of $ \GL_n(F^{nr}_{\o\omega}\times F^{nr}_{\omega })$ according to section 2.5.7.

The conjugacy class representatives 
$ \GL_n(F^{nr}_{\o\omega_{j,m_j}}\times 
F^{nr}_{\omega_{j,m_j }})$ correspond to the cosets of\linebreak 
$ \GL_n(F^{nr}_{R}\times F^{nr}_{L})\big/ \GL_n(\ZZ^2)$ in such a way that the space of these cosets is locally compact.

On the other hand, let $ \GL_n({F_{\o\omega}}\times {F_{\omega }})$ be the reductive pseudo-ramified complete bilinear semigroup of matrices and let $ G^{(2n)}({F_{\o\omega}}\times {F_{\omega }})= M_R\otimes M_L$ be its representation space, i.e. a bilinear  semigroup.  Then, a unitary representation of
$ \GL_n({F_{\o\omega}}\times {F_{\omega }})$ in this bilinear  semigroup 
$ G^{(2n)}({F_{\o\omega}}\times {F_{\omega }})$ is any unitary mapping:                                                                           %
\begin{align*}
U_{P^{(2n)}}: \quad 
 P^{(2n)}({F_{\o\omega^1}}\times {F_{\omega^1 }})
& \To R(P^{(2n)}({F_{\o\omega^1}}\times {F_{\omega ^1}}))\\[11pt]
 P^{(2n)}(F_{\o\omega^1_{j,m_j}}\times F_{\omega^1_{j,m_j} })
& \To R(P^{(2n)}({F_{\o\omega^1_{j,m_j}}}\times{F^1_{\omega_{j,m_j} }}))\end{align*}
from the bilinear parabolic semigroup 
 $P^{(2n)}({F_{\o\omega^1}}\times {F_{\omega ^1}})$ to the unitary operator\linebreak
$ R(P^{(2n)}({F_{\o\omega^1}}\times {F_{\omega^1 }}))$ according to proposition 2.4.5.
\vskip 11pt 

\subsubsection{{\boldmath Kernel of the integral operator restricted to the $(j,m_j)$-th conjugacy class representative of
$P^{(2n)}({F_{\o\omega^1}}\times {F_{\omega ^1}})$}}

Let $ \phi ^{(2n)}_{P_{j_R}}( x_{p_{j_R}} )\otimes \phi ^{(2n)}_{P_{j_L}}( x_{p_{j_L}} ) \in
L^{1-1}\RL( P^{(2n)}({F_{\o\omega^1}}\times {F_{\omega ^1}})) $ be the bifunction on the $j$-th irreducible equivalence class representative of 
$P^{(2n)} ({F_{\o\omega^1}}\times {F_{\omega^1 }})$ and let
 $ \phi ^{(2n)}_{G^{nr}_{j_R}}( x^{nr}_{g_{j_R}} )\otimes \phi ^{(2n)}_{G^{nr}_{j_L}}( x^{nr}_{g_{j_L}} ) \in 
L^{1-1}\RL( G^{(2n)}(F^{nr}_{\o\omega}\times F^{nr}_{\omega })) $ be the bifunction on the corresponding $j$-th equivalence class representative of the pseudo-unramified complete bilinear semigroup
 $G^{(2n)}(F^{nr}_{\o\omega}\times F^{nr}_{\omega })$~.

Then, the unitary operator $ R(P^{(2n)}({F^1_{\o\omega}}\times {F^1_{\omega }}))$ acts by convolution on every bifunction $ \phi ^{(2n)}_{G ^{nr}_{j_R}}( x^{nr}_{g_{j_R}} )\otimes \phi ^{(2n)}_{G^{nr}_{j_L}}( x^{nr}_{g_{j_L}} )$ on the pseudo-unramified bilinear complete semigroup  $G^{(2n)}(F^{nr}_{\o\omega}\times F^{nr}_{\omega })$ according to:
\begin{align*}
&R(P^{(2n)}(F_{\o\omega^1_{j,m_j}}\times F_{\omega^1_{j,m_j} }))\cdot
(\phi ^{(2n)}_{G^{nr}_{j_R}}( x^{nr}_{g_{j_R}} )\otimes \phi ^{(2n)}_{G^{nr}_{j_L}}( x^{nr}_{g_{j_L}} ))\\[11pt]
& \quad 
= \int_{G^{(2n)} ({F_{\o\omega}}\times {F_{\omega }})}
(\phi ^{(2n)}_{P_{j_R}}( x_{p_{j_R}} )\otimes \phi ^{(2n)}_{P_{j_L}}( x_{P_{j_L}} ))\\[11pt]
& \qquad \times
(\phi ^{(2n)}_{G^{nr}_{j_R}}( x^{nr}_{g_{j_R}} \cdot x_{p_{j_R}})\otimes \phi ^{(2n)}_{G^{nr}_{j_L}}
( x^{nr}_{g_{j_L}}\cdot x_{p_{j_L}} ))\ dx_{p_{j_R}}\ dx_{p_{j_L}}\\[11pt]
& \quad 
= \int_{G^{(2n)} ({F_{\o\omega}}\times {F_{\omega }})}
(\phi ^{(2n)}_{P_{j_R}} ( x_{p_{j_R}}\cdot ( x^{nr}_{g_{j_R}} )^{-1}))
\otimes \phi ^{(2n)}_{P_{j_L}}( x_{P_{j_L}}\cdot( x^{nr}_{g_{j_L}})^{-1}
 )))\\[11pt]
& \qquad \times
(\phi ^{(2n)}_{G^{nr}_{j_R}}( x_{p_{j_R}} ) \otimes \phi ^{(2n)}_{G^{nr}_{j_L}}
( x_{p_{j_L}}  ))\ dx_{p_{j_R}}\ dx_{p_{j_L}}\\[11pt]
& \quad 
= \int_{G^{(2n)} ({F_{\o\omega}}\times {F_{\omega }})}
( N_{P_{j_R}}(x_{p_{j_R}},x^{nr}_{g_{j_R}})
\otimes N_{P_{j_L}}(x_{p_{j_L}},x^{nr}_{g_{j_L}}))\\[11pt]
& \qquad \times
(\phi ^{(2n)}_{G^{nr}_{j_R}}( x_{p_{j_R}} ) \otimes \phi ^{(2n)}_{G^{nr}_{j_L}}
( x_{p_{j_L}}  )) \ dx_{p_{j_R}}\ dx_{p_{j_L}}\end{align*}
where
\begin{align*}
N_{P_{j_R}}(x_{p_{j_R}},x^{nr}_{g_{j_R}}) &= 
\phi ^{(2n)}_{P_{j_R}} ( x_{p_{j_R}}\cdot ( x^{nr}_{g_{j_R}} )^{-1})\\[11pt]
\text{(resp.} \quad 
N_{P_{j_L}}(x_{p_{j_L}},x^{nr}_{g_{j_L}}) &= 
\phi ^{(2n)}_{P_{j_L}}\ ( x_{p_{j_L}}\cdot ( x^{nr}_{g_{j_L}} )^{-1})\ )\end{align*}
is the \rl kernel of the integral operator
\begin{align*}
R_{P^{(2n)}({F_{\o\omega^1 }})}
(\phi ^{(2n)}_{P_{j_R}} )
&= \int_{P^{(2n)}} R_{P^{(2n)}({F_{\o\omega ^1}})} 
(x_{p_{j_R}})
\ \phi ^{(2n)}_{P_{j_R}}\  (x_{p_{j_R}})\ dx_{p_{j_R}}\\[11pt]
\text{(resp.} \quad 
R_{P^{(2n)}({F_{\omega^1 }})}
(\phi ^{(2n)}_{P_{j_L}} )
&= \int_{P^{(2n)}} R_{P^{(2n)}({F_{\omega ^1}})}
(x_{p_{j_L}})
\ \phi ^{(2n)}_{P_{j_L}} \ (x_{p_{j_L}})\ dx_{p_{j_L}}\ )\end{align*}
restricted to the $(j,m_j)$-th conjugacy class representative of 
$P^{(2n)}({F_{\o\omega^1}}\times {F_{\omega^1} })$~.
\vskip 11pt 

\subsubsection{Proposition}

{\em The unitary operator $ R(P^{(2n)}({F_{\o\omega^1}}\times {F_{\omega ^1}}))$ corresponding to the bilinear parabolic semigroup $ P^{(2n)}({F_{\o\omega^1}}\times {F_{\omega^1 }})$ acts by convolution on the bifunction(s)
$(\phi ^{(2n)}_{G^{nr}_{R}}( x^{nr}_{g_{R}} ) \otimes \phi ^{(2n)}_{G^{nr}_{L}}
( x^{nr}_{p_{L}}  )$ on the pseudo-unramified bilinear complete semigroup $ G^{(2n)}(F^{nr}_{\o\omega}\times F^{nr}_{\omega })$ according to
\begin{align*}
& R( P^{(2n)}({F_{\o\omega^1}}\times {F_{\omega^1 }}) )
(\phi ^{(2n)}_{G^{nr}_{R}}( x^{nr}_{g_{R}} ) \otimes \phi ^{(2n)}_{G^{nr}_{L}}
( x^{nr}_{p_{L}}  ))\\[11pt]
&\quad =
\int_{G^{(2n)}({F_{\o\omega}}\times {F_{\omega }}) }
( N_{P_R}(x_{p_R},x^{nr}_{g_R})\otimes N_{P_L}(x_{p_L},x^{nr}_{g_L}))
(\phi ^{(2n)}_{G^{nr}_{R}}( x_{p_{R}} ) \otimes \phi ^{(2n)}_{G^{nr}_{L}}
( x_{p_{L}}  ))\ dx_{p_R}\ dx_{p_L}\\[11pt]
& \quad 
=\bigoplus_{j,m_j} \int_{G^{(2n)}({F_{\o\omega}}\times {F_{\omega }})\big/ \GL_n(\ZZ\big/N\ZZ)^2)}
( N_{P_{j_R}}(x_{p_{j_R}},x^{nr}_{g_{j_R}})
\otimes N_{P_{j_L}}(x_{p_{j_L}},x^{nr}_{g_{j_L}}))\\[11pt]
& \qquad \qquad  \times 
( \phi ^{(2n)}_{G^{nr}_{j_R}}( x_{p_{j_R}} ) \otimes \phi ^{(2n)}_{G^{nr}_{j_L}}
( x_{p_{j_L}}  ))\ dx_{p_{j_R}}\ dx_{p_{j_L}}\end{align*}
where
\Bi
\item $\phi ^{(2n)}_{G^{nr}_{R}}( x_{p_{R}} ) \otimes \phi ^{(2n)}_{G^{nr}_{L}}
( x_{p_{L}}  ))=\bigoplus_j \bigoplus_{m_j}
(\phi ^{(2n)}_{G^{nr}_{j_R}}( x_{p_{j_R}} ) \otimes \phi ^{(2n)}_{G^{nr}_{j_L}}
( x_{p_{j_L}}  ))$~;

\item the (bi)kernel of the integral (bi)operator $R_{P^{(2n)}( {F_{\o\omega^1 }}\times {F_{\omega ^1}}  )}
( \phi ^{(2n)}_{P_{R}}\otimes  \phi ^{(2n)}_{P_{L}}  )$ decomposes:
\[ N_{P_R}(x _{p_R},x^{nr}_{g_R}) \otimes N_{P_L}(x _{p_L},x^{nr}_{g_L}) 
= \bigoplus_{j,m_j} ( N_{P_{j_R}}(x_{p_{j_R}},x^{nr}_{g_{j_R}}) \otimes
N_{P_{j_L}}(x_{p_{j_L}},x^{nr}_{g_{j_L}})\]
according to  the conjugacy class representatives of $P^{(2n)}( {F_{\o\omega^1 }}\times {F_{\omega^1 }}  )$~.
\Ei
}
\newpage

\bpr \Be
\item By this way, the (bi)kernel $(N_{P_R}(x _{p_R},x^{nr}_{g_R}) \otimes N_{P_L}(x _{p_L},x^{nr}_{g_L})) $
is separable since it decomposes according to the conjugacy class representatives of 
$P^{(2n)}( {F_{\o\omega^1 }}\times {F_{\omega^1 }}  )$  which corresponds to the cosets of
$G^{(2n)}({F_{\o\omega}}\times {F_{\omega }})\big/\GL_n((\ZZ/N\ \ZZ)^2)$~.

\item The integral 
$\int_{G^{(2n)}({F_{\o\omega}}\times {F_{\omega }}) }
( N_{P_R}(x_{p_R},x^{nr}_{g_R})\otimes N_{P_L}(x_{p_L},x^{nr}_{g_L}))
(\phi ^{(2n)}_{G^{nr}_{R}}( x_{p_{R}} ) \otimes \phi ^{(2n)}_{G^{nr}_{L}}
( x_{p_{L}}  ))\linebreak dx_{p_R}\ dx_{p_L}$~, decomposing according to the conjugacy class representatives of 
$G^{(2n)}({F_{\o\omega}}\times {F_{\omega }}) $~, converges since the integral
$ \int_{G^{(2n)}({F_{\o\omega}}\times {F_{\omega }})\big/ \GL_n((\ZZ\big/N\ZZ)^2)}
( N_{P_{j_R}}(x_{p_{j_R}},x^{nr}_{g_{j_R}})
\otimes\linebreak N_{P_{j_L}}(x_{p_{j_L}},x^{nr}_{g_{j_L}}))
( \phi ^{(2n)}_{G^{nr}_{j_R}}( x_{p_{j_R}} ) \otimes \phi ^{(2n)}_{G^{nr}_{j_L}}
( x_{p_{j_L}}  ))\ dx_{p_{j_R}}\ dx_{p_{j_L}}$  restricted to the $(j,m_j)$-th conjugacy class representative is bounded.
\epr
\Ee
\vskip 11pt 

\subsubsection{Proposition}

{\em The trace of the unitary operator $R(P^{(2n)}( {F_{\o\omega^1 }}\times {F_{\omega^1 }}  ))$ 
can be expressed as a sum of integrals of $ \phi ^{(2n)}_{P_{j_R}}( x_{p_{j_R}} \ (x^{nr}_{g_{j_R}})^{-1}) 
\otimes_D \phi ^{(2n)}_{P_{j_L}}
( x_{p_{j_L}} \ (x^{nr}_{g_{j_L}})^{-1}) $ according to:
\begin{align*}
& {\rm tr}( R_{P^{(2n)}({F_{\o\omega^1}}\times {F_{\omega^1 }})}
(\phi ^{(2n)}_{P_{R}}\otimes \phi ^{(2n)}_{P_{L}} ) )\\[11pt]
&\quad =\bigoplus_{j,m_j}\ {\rm vol}(P^{(2n)}({F_{\o\omega^1}}\times {F_{\omega^1 }})
 \int_{P^{(2n)}}  \phi ^{(2n)}_{P_{j_R}} \ (x_{p_{j_R}}\ (x^{nr}_{g_{j_R}})^{-1})  \otimes_D 
 \phi ^{(2n)}_{P_{j_L}} \ (x_{p_{j_L}}\ (x^{nr}_{g_{j_L}})^{-1})  \ dx_{p_{j_R}}\ dx_{p_{j_L}}\\[11pt]
&\quad =\bigoplus_{j,m_j} R_{P^{(2n)}(F_{\o\omega^1_{j,mj}}\times F_{\omega^1_{j,m_j }})}
(\phi ^{(2n)}_{P_{j_R}}, \phi ^{(2n)}_{P_{j_L}} ) 
\end{align*}
where $(\cdot,\cdot)$ is a bilinear form.
}
\vskip 11pt 

\bpr According to proposition 2.5.12, the kernel $N_{P_R}(x_{p_R},x^{nr}_{g_R})\otimes
N_{P_L}(x_{p_L},x^{nr}_{g_L})$ of the integral operator 
$R(P^{(2n)}( {F_{\o\omega^1 }}\times {F_{\omega ^1}}  ))$ corresponds precisely to the action of the unitary parabolic bilinear semigroup $P^{(2n)}( {F_{\o\omega ^1}}\times {F_{\omega ^1}}  )$ on the bifunction
$\phi ^{(2n)}_{G^{nr}_R}(x_{p_R}) \otimes \phi ^{(2n)}_{G^{nr}_L}(x_{p_L}) $ on the pseudo-unramified complete bilinear semigroup
$G^{(2n)}( F^{nr}_{\o\omega }\times F^{nr}_{\omega }  )$~.

As a result, we have the thesis which is an adaptation of the Arthur's trace formula for reductive linear groups \cite{Art2}.\epr
\vskip 11pt 

\subsubsection{{\boldmath Trace formula for the operator $R(G^{(2n)}( {F_{\o\omega }}\times {F_{\omega }}  ))$}}

A representation of the pseudo-ramified bilinear complete semigroup 
$G^{(2n)}( {F_{\o\omega }}\times {F_{\omega }}  )$ is any mapping
\begin{align*}
R_{G^{(2n)}} : \quad
G^{(2n)}( {F_{\o\omega }}\times {F_{\omega }}  )
&\To R(G^{(2n)}( {F_{\o\omega }}\times {F_{\omega }}  ))\\[11pt]
G^{(2n)}( F_{\o\omega_{j,m_j} }\times F_{\omega_{j,m_j }}  )
&\To R(G^{(2n)}( F_{\o\omega_{j,m_j}}\times F_{\omega _{j,m_j}}  ))\end{align*}
from $G^{(2n)}( {F_{\o\omega }}\times {F_{\omega }}  )$ to the operator
$ R(G^{(2n)}( {F_{\o\omega }}\times {F_{\omega }}  ))\simeq \GL_n( {F_{\o\omega }}\times {F_{\omega }}  )$ in such a way that the matrix 
$\GL_n( F_{\o\omega_{j,m_j} }\times F_{\omega_{j,m_j }}  )$ corresponds to the conjugacy class representative
$g^{(2n)}\RL[j,m_j]\equiv
G^{(2n)}( F_{\o\omega_{j,m_j} }\times F_{\omega_{j,m_j }}  )$ of the pseudo-ramified bilinear  semigroup 
$G^{(2n)}({ F_{\o\omega }}\times {F_{\omega }}  )$~.

Then, a trace formula of the operator $R(G^{(2n)}( {F_{\o\omega }}\times {F_{\omega }}  ))$ can be reached from the trace of the unitary parabolic operator 
$R(P^{(2n)}( {F_{\o\omega ^1}}\times {F_{\omega^1 }}  ))$ by considering the analytic continuation action of the (bi)kernel  $N_{P_R}(x_{p_R},x^{nr}_{g_R})\otimes
N_{P_L}(x_{p_L},x^{nr}_{g_L})$ according to the following proposition.
\vskip 11pt 

\subsubsection{Proposition}

{\em The trace of the unitary parabolic operator $R(P^{(2n)}( {F_{\o\omega^1 }}\times {F_{\omega^1 }}  ))$ can be extended to the trace of the pseudo-ramified bilinear complete semigroup 
$R(G^{(2n)}( {F_{\o\omega }}\times {F_{\omega }}  ))$ if the action of the (bi)kernel
$N_{P_R}(x_{p_R},x^{nr}_{g_R})\otimes
N_{P_L}(x_{p_L},x^{nr}_{g_L})$  of $R(P^{(2n)}( {F_{\o\omega ^1}}\times {F_{\omega^1 }}  ))$ on the (pseudo-)unramified bifunction
$\phi ^{(2n)}_{G^{nr}_R}(x_{p_R}) \otimes \phi ^{(2n)}_{G^{nr}_L}(x_{p_L}) $ consists in the analytic continuation of this bifunction in such a way that the (pseudo-)ramified bifunction
$ \phi ^{(2n)}_{G_R}(x_{g_R}) \otimes \phi ^{(2n)}_{G_L}(x_{g_L}) \in L^{1-1}\RL
(G^{(2n)}( {F_{\o\omega }}\times {F_{\omega }}  ))$ on the pseudo-ramified bilinear complete semigroup
$G^{(2n)}( {F_{\o\omega }}\times {F_{\omega }}  )$ verifies:
\[
\phi ^{(2n)}_{G_R}(x_{g_R}) \otimes \phi ^{(2n)}_{G_L}(x_{p_L}) 
= (N_{P_R}(x_{p_R},x^{nr}_{g_R}) \otimes N_{P_L}(x_{p_L},x^{nr}_{g_L}) )
(\phi ^{(2n)}_{G^{nr}_R}(x^{nr}_{p_R}) \otimes \phi ^{(2n)}_{G^{nr}_L}(x^{nr}_{p_L}) )\;.\]
So, we get the following trace formula:
\begin{align*}
& {\rm tr}( R_{G^{(2n)}({F_{\o\omega}}\times {F_{\omega }})}
(\phi ^{(2n)}_{G_{R}}\otimes \phi ^{(2n)}_{G_{L}} ) )\\[11pt]
&\quad =\bigoplus_{j,m_j}\ {\rm vol}(G^{(2n)}({F_{\o\omega_\oplus }}\times {F_{\omega_\oplus  }})
 \int_{G^{(2n)}}  \phi ^{(2n)}_{G_{j_R}} \ (x_{g_{j_R}})  \otimes_D 
 \phi ^{(2n)}_{G_{j_L}} \ (x_{g_{j_L}})\ dx_{g_{j_R}}\ dx_{g_{j_L}} \end{align*}
where the sum ranges over the diagonal part of the conjugacy class representatives of 
$G^{(2n)}( {F_{\o\omega }}\times {F_{\omega }}  )$~.
}
\vskip 11pt 

\bpr The action of the bikernel 
\[
N_{P_R}(x_{p_R},x^{nr}_{g_R}) \otimes N_{P_L}(x_{p_L},x^{nr}_{g_L}) 
= \bigoplus_{j,m_j} ( \phi ^{(2n)}_{P_{j_R}}(x_{p_{j_R}}\ (x^{nr}_{g_{j_R}})^{-1})\otimes
\phi ^{(2n)}_{P_{j_L}}(x_{p_{j_L}}\ (x^{nr}_{g_{j_L}})^{-1}))\]
on the (pseudo-)unramified bifunction
\[
\phi ^{(2n)}_{G^{nr}_{R}}(x_{p_{R}}) \otimes \phi ^{(2n)}_{G^{nr}_{L}}(x_{p_{L}})
= \bigoplus_{j,m_j} 
(\phi ^{(2n)}_{G^{nr}_{j_R}}(x_{p_{j_R}}) \otimes \phi ^{(2n)}_{G^{nr}_{j_L}}(x_{p_{j_L}}))\;,\]
defined on the pseudo-unramified complete bilinear semigroup
$G^{(2n)}( F^{nr}_{\o\omega }\times F^{nr}_{\omega }  )$~, consists in the pseudo-ramification of this bifunction according to:
\begin{align*}
I^{nr\to r}\RL: \quad 
\phi ^{(2n)}_{G^{nr}_{R}}(x_{p_{R}}) \otimes \phi ^{(2n)}_{G^{nr}_{L}}(x_{p_{L}})
&\To \phi ^{(2n)}_{G_{R}}(x_{g_{R}}) \otimes \phi ^{(2n)}_{G_{L}}(x_{g_{L}})\end{align*}
where $I^{nr\to r}\RL$ is an injective morphism to the (pseudo-)ramified bifunction
$\phi ^{(2n)}_{G_{R}}(x_{g_{R}}) \otimes \phi ^{(2n)}_{G_{L}}(x_{g_{L}})$ on  the pseudo-ramified bilinear complete semigroup 
$G^{(2n)}( {F_{\o\omega }}\times {F_{\omega }}  )$~.  Consequently, the injective morphism of pseudo-ramification
$I^{nr\to r}\RL$ of the pseudo-unramified bifunction 
$\phi ^{(2n)}_{G^{nr}_{R}}(x_{p_{R}}) \otimes \phi ^{(2n)}_{G^{nr}_{L}}(x_{p_{L}})$ corresponds to an analytic continuation of this bifunction since the injective morphism $I^{nr\to r}\RL$ corresponds to the inverse of the map:
\[ G^{(2n)}_{F\to F^{nr}}:\quad G^{(2n)}( {F_{\o\omega }}\times {F_{\omega }})
\To G^{(2n)}( F^{nr}_{\o\omega }\times F^{nr}_{\omega })\]
introduced in section 2.4.1.

The volume ${\rm vol}(G^{(2n)}( F^{nr}_{\o\omega }\times F^{nr}_{\omega }))$ of the pseudo-unramified complete 
bilinear semigroup $G^{(2n)}( F^{nr}_{\o\omega }\times F^{nr}_{\omega })$ with respect to a fixed Haar 
measure is then inflated by the volume ${\rm vol}(P^{(2n)}( {F_{\o\omega^1 }}\times {F_{\omega^1 }}))$ 
of the kernel $\Ker (G^{(2n)}_{F\to F^{nr}})=
 P^{(2n)}( {F_{\o\omega ^1}}\times {F_{\omega^1 }})$ of $G^{(2n)}_{F\to F^{nr}}$ into the volume
${\rm vol}(G^{(2n)}( {F_{\o\omega }}\times {F_{\omega }}))$~.

And, the trace formula ${\rm tr}( R_{G^{(2n)}({F_{\o\omega}}\times {F_{\omega }})}
(\phi ^{(2n)}_{G_{R}}\otimes \phi ^{(2n)}_{G_{L}} ) )$ then follows from the developments of propositions 2.5.12 and 2.5.13.
\epr
\vskip 11pt 

\subsubsection{Corollary}

{\em Let $L^{1-1}\RL
(G^{(2n)}( F^{nr}_{\o\omega }\times F^{nr}_{\omega }  ))$ denote the bialgebra of (pseudo-)unramified bifunctions on the (pseudo-)unramified bilinear complete semigroup $G^{(2n)}( F^{nr}_{\o\omega }\times F^{nr}_{\omega }  )$  and let
$L^{1-1}\RL
(G^{(2n)}( {F_{\o\omega }}\times {F_{\omega }}  ))$ be the bialgebra of (pseudo-)ramified bifunctions on the pseudo-ramified bilinear complete semigroup $
G^{(2n)}( {F_{\o\omega }}\times {F_{\omega }}  )$~.

Then, the injective morphism:
\begin{align*}
I^{nr\to r}\RL : \quad
L^{1-1}\RL
(G^{(2n)}( F^{nr}_{\o\omega }\times F^{nr}_{\omega }  ))
& \To L^{1-1}\RL
(G^{(2n)}( {F_{\o\omega }}\times {F_{\omega }}  ))\\[11pt]
\phi ^{(2n)}_{G^{nr}_{R}}(x^{nr}_{g_R}) \otimes \phi ^{(2n)}_{G^{nr}_{L}}(x^{nr}_{g_L})
& \To \phi ^{(2n)}_{G_{R}}(x_{g_R})\otimes \phi ^{(2n)}_{G_{L}}(x_{g_L})\end{align*}
corresponds to a pseudo-ramification morphism of the bialgebra
$ L^{1-1}\RL
(G^{(2n)}( F^{nr}_{F_{\o\omega }}\times F^{nr}_{F_{\omega }}  ))$~.
}
\vskip 11pt


\subsubsection{Proposition}

{\em Let the bifunction
\[
\phi ^{(2n)}_{G_{j_R}}(x_{g_{j_R}})\otimes \phi ^{(2n)}_{G_{j_L}}(x_{g_{j_L}})
\in L^{1-1}\RL( G^{(2n)}( {F_{\o\omega }}\times {F_{\omega }}  ))
\]
 be developed according to:
\[
\phi ^{(2n)}_{G_{j_R}}(x_{g_{j_R}})\otimes \phi ^{(2n)}_{G_{j_L}}(x_{g_{j_L}})
= (\phi ^{(2n)}_{P_{j_R}}(x_{p_{j_R}})\otimes \phi ^{(2n)}_{P_{j_L}}(x_{p_{j_L}}))
\times (\phi ^{(2n)}_{G^{nr}_{j_R}}(x^{nr}_{g_{j_R}})
\otimes \phi ^{(2n)}_{G^{nr}_{j_L}}(x^{nr}_{g_{j_L}}))\;.
\]
Then, the following trace formula
\begin{align*}
& {\rm Tr}(R_{G^{(2n)}( {F_{\o\omega }}\times {F_{\omega }})}
(\phi ^{(2n)}_{G_{R}}  \otimes \phi ^{(2n)}_{G_{L}} ))\\[11pt]
& \quad =\bigoplus_{j,m_j} \L[
{\rm vol}(P^{(2n)}( {F_{\o\omega^1 }}\times {F_{\omega ^1}}  ))
\int_{P^{(2n)}} \phi ^{(2n)}_{P_{j_R}}(x_{p_{j_R}})\otimes_D \phi ^{(2n)}_{P_{j_L}}(x_{p_{j_L}})
\ dx_{p_{j_R}}\ dx_{p_{j_L}}\R]\\[11pt]
& \qquad \qquad \times \L[
{\rm vol}(G^{(2n)}( F^{nr}_{\o\omega}\times F^{nr}_{\omega }  ))
\int_{G^{(2n)}} \phi ^{(2n)}_{G^{nr}_{j_R}}(x^{nr}_{g_{j_R}})\otimes_D \phi ^{(2n)}_{G^{nr}_{j_L}}(x^{nr}_{g_{j_L}})
\ dx^{nr}_{g_{j_R}}\ dx^{nr}_{g_{j_L}}\R]\end{align*}
follows and relies on the Lefschetz trace formula.
}\vskip 11pt

\bpr The considered trace formula results from the development of the bifunction
$\phi ^{(2n)}_{G_{j_R}}(x_{{j_R}})\otimes \phi ^{(2n)}_{G_{j_L}}(x_{{j_L}})$
given in proposition 2.5.15 and relies on the Lefschetz trace formula because $\phi ^{(2n)}_{P_{j_R}}(x_{p_{j_R}})\otimes \phi ^{(2n)}_{P_{j_L}}(x_{p_{j_L}})$
is defined on the bilinear parabolic semigroup 
$P^{(2n)}( {F_{\o\omega^1}}\times {F_{\omega^1}})$ which is the isotropy subgroup of
$G^{(2n)}( {F_{\o\omega }}\times {F_{\omega }})$ fixing the bielements of
$G^{(2n)}( {F_{\o\omega }}\times {F_{\omega }})$ (see proposition 2.4.5 and \cite{Gro3}).
\epr
\vskip 11pt

\subsubsection{Proposition}

{\em Let \quad ${\Tr}(R_{P^{(2n)}( {F^1_{\o\omega }}\times {F^1_{\omega }})}
(\phi ^{(2n)}_{P_{R}}  \otimes \phi ^{(2n)}_{P_{L}} ))$ \quad be given by:
\begin{align*}
&{\Tr}(R_{P^{(2n)}( {F_{\o\omega^1}}\times {F_{\omega^1}})}
(\phi ^{(2n)}_{P_{R}}  \otimes \phi ^{(2n)}_{P_{L}} ))\\[11pt]
& \quad = \bigoplus_j\bigoplus_{m_j}
{\rm vol}(P^{(2n)}( {F_{\o\omega^1}}\times {F_{\omega ^1 }} ))
\int_{P^{(2n)}} \phi ^{(2n)}_{P_{j_R}}(x_{p_{j_R}})\otimes _D
\phi ^{(2n)}_{P_{j_L}}(x_{p_{j_L}})
\ dx_{p_{j_R}}\ dx_{p_{j_L}}\\[11pt]
& \quad = \bigoplus_j\bigoplus_{m_j}
R_{P^{(2n)}( 
{F_{\o\omega^1_{j,m_j} }}\times {F_{\omega^1_{j,m_j} }})}
(\phi ^{(2n)}_{P_{j_R}} , \phi ^{(2n)}_{P_{j_L}} )\end{align*}
according to proposition 2.5.13.

Then, the trace formula ${\rm Tr}(R_{G^{(2n)}( {F_{\o\omega }}\times {F_{\omega }})}
(\phi ^{(2n)}_{G_{R}}  \otimes \phi ^{(2n)}_{G_{L}} ))$ of proposition 2.5.17, relied on the Lefschetz trace formula, can be developed as follows:
\begin{align*}
& {\rm Tr}
(R_{G^{(2n)}( {F_{\o\omega }}\times {F_{\omega }})}
(\phi ^{(2n)}_{G_{R}}  \otimes \phi ^{(2n)}_{G_{L}} ))\\[11pt]
& \quad =
{\rm Tr}
(R_{P^{(2n)}( {F_{\o\omega^1}}\times {F_{\omega ^1)}})}
(\phi ^{(2n)}_{P_{R}}  \otimes \phi ^{(2n)}_{P_{L}} )) \times
{\rm Tr}(R_{G^{(2n)}( F^{nr}_{\o\omega }\times F^{nr}_{\omega })}
(\phi ^{(2n)}_{G^{nr}_{R}}  \otimes \phi ^{(2n)}_{G^{nr}_{L}} )) \end{align*}
where 
\begin{align*}
& {\rm Tr}(R_{G^{(2n)}( F^{nr}_{\o\omega }\times F^{nr}_{\omega })}
(\phi ^{(2n)}_{G^{nr}_{R}}  \otimes \phi ^{(2n)}_{G^{nr}_{L}} ))\\[11pt]
& \quad = \bigoplus_{j,m_j}
{\rm vol}(G^{(2n)}( F^{nr}_{\o\omega }\times F^{nr}_{\omega }  ))
\int_{G^{(2n)}} \phi ^{(2n)}_{G^{nr}_{j_R}}(x^{nr}_{g_{j_R}})\otimes_D 
\phi ^{(2n)}_{G^{nr}_{j_L}}(x^{nr}_{g_{j_L}})
\ dx^{nr}_{g_{j_R}}\ dx^{nr}_{g_{j_L}}\end{align*}
in such a way that 
$R_{P^{(2n)}( {F_{\o\omega^1_{j,m_j} }}\times {F_{\omega^1_{j,m_j} }})}
(\phi ^{(2n)}_{P_{j_R}}  \otimes \phi ^{(2n)}_{P_{j_L}} )$ is an invariant of\linebreak
${\rm Tr}(R_{G^{(2n)}( {F_{\o\omega }}\times {F_{\omega }})}
(\phi ^{(2n)}_{G_{R}}  \otimes \phi ^{(2n)}_{G_{L}} ))$~.
}\vskip 11pt

\bpr The parabolic bifunctions $
\phi ^{(2n)}_{P_{j_R}}(x_{p_{j_R}})  \otimes \phi ^{(2n)}_{P_{j_L}} (x_{p_{j_L}})$~, $1\le j\le r$~,
are in one-to-one correspondence and are defined on the bilinear parabolic semigroup
$P^{(2n)}({F_{\o\omega ^1}}\times {F_{\omega^1} })$~, being the isotropy subgroup of
$G^{(2n)}( {F_{\o\omega }}\times {F_{\omega} })$~: so, the
$R_{P^{(2n)}( {F_{\o\omega^1_{j,m_j} }}\times {F_{\omega^1_{j,m_j} }})}
(\phi ^{(2n)}_{P_{j_R}}  \otimes \phi ^{(2n)}_{P_{j_L}} )$~, $1\le j\le r$~, are the invariants of
${\rm Tr}(R_{G^{(2n)}( {F_{\o\omega }}\times {F_{\omega }})}
(\phi ^{(2n)}_{G_{R}}  \otimes \phi ^{(2n)}_{G_{L}} ))$~.\epr
\vskip 11pt

\subsubsection{Remark}

If $P^{(2n)}( {F_{\o\omega^1 }}\times {F_{\omega ^1}}  )$ is considered as the unitary representation of the bilinear  semigroup  
$G^{(2n)}( {F_{\o\omega }}\times {F_{\omega }}  )$~, then the trace formula
${\rm Tr}(R_{G^{(2n)}( {F_{\o\omega }}\times {F_{\omega }})}
(\phi ^{(2n)}_{G_{R}}  \otimes \phi ^{(2n)}_{G_{L}} ))$ of proposition 2.5.17 refers to the Arthur's trace formula.

But, if 
$P^{(2n)}( {F_{\o\omega^1 }}\times {F_{\omega ^1}}  )$ is viewed as the isotropy subgroup of
$G^{(2n)}( {F_{\o\omega }}\times {F_{\omega }}  )$~, then
${\Tr}(R_{G^{(2n)}( {F_{\o\omega }}\times {F_{\omega }})}
(\phi ^{(2n)}_{G_{R}}  \otimes \phi ^{(2n)}_{G_{L}} ))$ relies on the Lefschetz trace formula.
\vskip 11pt


\subsubsection{{\boldmath Second trace formula for the operator $R(G^{(2n)}( {F_{\o\omega }}\times {F_{\omega }}  ))$}}

A trace formula of the operator $R(G^{(2n)}( {F_{\o\omega }}\times {F_{\omega }}  ))$ can be developed directly from the set of bifunctions 
$\phi ^{(2n)}_{G_{j_R}}(x_{g_{j_R}})\otimes \phi ^{(2n)}_{G_{j_L}}(x_{g_{j_L}})\in 
L^{1-1}\RL
(G^{(2n)}( {F_{\o\omega }}\times {F_{\omega }}  ))$ on the conjugacy class representatives of
$G^{(2n)}( {F_{\o\omega }}\times {F_{\omega }}  )$ and not from an extension of the trace formula of the unitary parabolic operator $R(P^{(2n)}( {F_{\o\omega^1 }}\times {F_{\omega^1 }}  ))$~.

As the bifunction $\phi ^{(2n)}_{G_{R}}(x_{g_{R}})\otimes \phi ^{(2n)}_{G_{L}}(x_{g_{L}})$~, defined on the pseudo-ramified bilinear complete semigroup $G^{(2n)}( {F_{\o\omega_\oplus }}\times {F_{\omega_\oplus }}  )$~, develops
\[\phi ^{(2n)}_{G_{R}}(x_{g_{R}})\otimes \phi ^{(2n)}_{G_{L}}(x_{g_{L}})
= \bigoplus_{j,m_j} \L(\phi ^{(2n)}_{G_{j_R}}(x_{g_{j_R}})\otimes \phi ^{(2n)}_{G_{j_L}}(x_{g_{j_L}})\R)\]
according to the bifunctions 
$\phi ^{(2n)}_{G_{j_R}}(x_{g_{j_R}})\otimes \phi ^{(2n)}_{G_{j_L}}(x_{g_{j_L}})$ on the conjugacy class representatives
$g^{(2n)}\RL[j,m_j]$ of $G^{(2n)}( {F_{\o\omega }}\times {F_{\omega }}  )$~, which correspond to the cosets of
$\GL_n( {F_{R}}\times {F_{L }}  )\big/\linebreak \GL_n((\ZZ\big/N\ \ZZ)^2)$~, we get clearly the searched trace formula:
\begin{align*}
&{\rm tr}(R_{\GL_n( {F_{\o\omega }}\times {F_{\omega }}  )\big/\GL_n((\ZZ\big/N\ \ZZ)^2)}
(\phi ^{(2n)}_{G_{R}}(x_{g_{R}})\otimes \phi ^{(2n)}_{G_{L}}(x_{g_{L}}))\\[11pt]
&\quad =\bigoplus_{j,m_j}
{\rm vol}(G^{(2n)}( {F_{\o\omega }}\times {F_{\omega }}  ))\int_{G^{(2n)}}
\phi ^{(2n)}_{G_{j_R}}(x_{g_{j_R}})\otimes_D \phi ^{(2n)}_{G_{j_L}}(x_{g_{j_L}})\ dx_{g_{j_R}}\ dx_{g_{j_L}}\;.\end{align*}
This trace formula of
$R(G^{(2n)}( {F_{\o\omega }}\times {F_{\omega }}  ))$ corresponds to the trace formula obtained in proposition 2.5.15 from the trace of $R(P^{(2n)}( {F_{\o\omega^1 }}\times {F_{\omega^1 }}  ))$.  It results from the action of the Hecke bialgebra $\Hs\RL(n)$ on the pseudo-ramified bilinear  semigroup  
$G^{(2n)}( {F_{\o\omega }}\times {F_{\omega }}  )$ since it is generated by the Hecke bioperators
$(T_R(n;r)\otimes T_L(n;r))$ having as representation
$\GL_n((\ZZ\big/N\ \ZZ)^2)$ according to proposition 2.2.6.
\vskip 11pt

This trace formula will be proved in chapter 3 to correspond to the Plancherel formula if the considered pseudo-ramified bilinear algebraic semigroup has undergone a projective toroidal isomorphism of compactification.
\vskip 11pt

\section[{\boldmath Langlands global correspondences for irreducible representations of\protect\linebreak $\GL(n)$}]{{\boldmath Langlands global correspondences for irreducible representations of $\GL(n)$}}

\setcounter{subsection}{0}
\subsubsection*{Introductive keys}
\addcontentsline{toc}{subsection}{Introductive keys}

In chapter 2, nonabelian global class field concepts were introduced, i.e. essentially the construction of an algebraic bilinear semigroup  $G^{(2n)}({\wt F_{\o\omega }}\times {\wt F_\omega })$ constituting the $2n$-dimensional irreducible representation $\Irr W^{(2n)}\FRL(W^{ab}_{F_R}\times W^{ab}_{F_L})$ of the product, right by left, of global Weil groups.

It was shown that $G^{(2n)}({\wt F_{\o\omega }}\times {\wt F_\omega })$ decomposes into diagonal and off-diagonal
$r$ conjugacy (bi)classes whose equivalent representatives $m_r$ are noted $g^{(2n)}\RL[j,m_j]$~, $1\le j\le r$~.  
The equivalent representative $m_j$ of the $j$-th \lr class of $T^{(2n)}({\wt F_{\omega }})$ (resp. $T^{(2n)}({\wt F_{\o\omega }})^t\subset G^{(2n)}({\wt F_{\o\omega }}\times {\wt F_\omega })$ is noted $g^{(2n)}_L[j,m_j]$ (resp. $g_R^{(2n)}[j,m_j]$~).

The algebraic representation space of the algebraic bilinear semigroup of matrices\linebreak $\GL_n({\wt F_{\o\omega }}\times {\wt F_\omega })$ is noted $\Reps (\GL_n({\wt F_{\o\omega }}\times {\wt F_\omega }))$ and is precisely the algebraic bilinear semigroup $G^{(2n)}({\wt F_{\o\omega }}\times {\wt F_\omega })$ isomorphic to the $\GL_{n}({\wt F_{\o\omega }}\times {\wt F_\omega })$-bisemimodule $(\wt M_R\otimes \wt M_L)$~.

Then, the content of this chapter will concern:
\Bean

\item {\bf the introduction of a general bilinear cohomology\/} defined as a contravariant bifunctor $H^*$ from smooth abstract (resp. algebraic) bisemivarieties together with a natural transformation
$n_{H^*\to H^{[*,*]}}$ from $H^*$ to the associated the de Rham bilinear cohomology $H^{[*,*]}$~;
\item the Fulton-McPherson compactification of bilinear algebraic semigroups;
\item a toroidal compactification of the lattice bisemispace $X\SRL$~;
\item a double coset decomposition of $\GL_n({F_{\o\omega }}\times {F_\omega })$ leading to the equivalent 
of a bilinear Shimura variety $\o S^{P_n}_{K_n}$~;
\item the Langlands global correspondence on $\Reps (\GL_n({\wt F_{\o\omega_\oplus }}\times {\wt F_{\omega_\oplus} }))$ based 
upon the study of the ring of regular (bi)functions on $\Reps (\GL_n({\wt F_{\o\omega }}\times {\wt F_\omega }))$ and the corresponding Langlands ``real'' correspondence.
\Ee
\vskip 11pt

\subsection{Compactifications of bilinear algebraic semigroups and global holomorphic correspondences}

\subsubsection{The compactification of bilinear algebraic semigroups}

The compactification of the bilinear algebraic semigroup $G^{(n)}(\widetilde F^+_{\o v}\times\widetilde F^+_v)$ over the product, right by left, of real pseudo-ramified extension semirings leads to the compactified bilinear algebraic semigroup
$G^{(n)}({F^+_{\o v}}\times {F^+_v})$~, taking into account the inclusions
\[\begin{CD}
 G^{(n)}({F^+_{\o v}}\times {F^+_v})@. \quad \scalebox{1.5}{$\hookrightarrow$}\quad  @. G^{(2n)}({F_{\o\omega }}\times {F_\omega })\\
@VV\wr V @. @VV\wr V \\
G^{(n)}(\widetilde F^+_{\o v}\times \widetilde F^+_v) @. \scalebox{1.5}{$\hookrightarrow$}@. G^{(2n)}( \widetilde F_{\o\omega }\times \widetilde F_\omega )\end{CD}\]
of the ``real'' bilinear algebraic semigroups into the corresponding ``complex'' bilinear algebraic semigroups as justified in section 3.5.1.
\vskip 11pt

We show that this problem is closely associated to the generation of a smooth bilinear general  bisemivariety 
$\tau _\sigma ^{(n)}({F^+_{\o v}}\times {F^+_v})$ characterized by polyhedral convex bicones
$\sigma \RL^{(n)}(F^+_{\o v_{j_\delta ,m_{j_\delta }}} \times F^+_{v_{j_\delta ,m_{j_\delta }}})$~.
\vskip 11pt

And, this bilinear general  bisemivariety $\tau _\sigma ^{(n)}({F^+_{\o v}}\times {F^+_v})$ is in one-to-one correspondence with the bilinear general ``toroidal'' bisemivariety $\tau _T^{(n)}({F^{+,T}_{\o v}}\times {F^{+,T}_v})$ in such a way that the polyhedral convex bicones 
$\sigma \RL^{(n)}(F^+_{\o v_{j_\delta ,m_{j_\delta }}} \times F^+_{v_{j_\delta ,m_{j_\delta }}})$ are transformed into products, right by left, of $n$-dimensional real semitori $T^n_R[j_\delta ,m_{j_\delta }]\times T^n_L[j_\delta ,m_{j\delta }]
\in G ^{(n)}({F^{+,T}_{\o v}}\times {F^{+,T}_v})$~.
\vskip 11pt 

\subsubsection{An adaptation of the Fulton-McPherson compactification}

The set of \lr real pseudo-ramified completions
\[ {F^+_{v}}= \{ F^+_{v_{1_\delta}},\cdots,F^+_{v_{j_\delta ,m_{j_\delta }}},\cdots,F^+_{v_{r_\delta }}\}
\qquad \text{(resp.} \quad 
{F^+_{\o v}}= \{ F^+_{\o v_{1_\delta}},\cdots,F^+_{\o v_{j_\delta ,m_{j_\delta }}},\cdots,F^+_{\o v_{r_\delta }}\}\ )\]
introduced in section 1.1.5, proceeds from the set of \lr real pseudo-ramified Galois extensions of $k$~:
\[ \widetilde F^+_{v}= \{ \widetilde F^+_{v_{1_\delta}},\cdots,\widetilde F^+_{v_{j_\delta ,m_{j_\delta }}},\cdots,\widetilde F^+_{v_{r_\delta }}\}
\qquad \text{(resp.} \quad 
\widetilde F^+_{\o v}= \{ \widetilde F^+_{\o v_{1_\delta}},\cdots,\widetilde F^+_{\o v_{j_\delta ,m_{j_\delta }}},\cdots,\widetilde F^+_{\o v_{r_\delta }}\} \ )\]
(see sections 1.1.4 and 1.1.7).
\vskip 11pt 

Let $G ^{(n)}(\widetilde F^{+}_{\o v}\times \widetilde F^{+}_v)$ denote the non-compact real bilinear algebraic semigroup with entries in $ (\widetilde F^{+}_{\o v}\times \widetilde F^{+}_v)$~.
\vskip 11pt

Several types of compactification can be envisaged for it: we shall introduce here an adaptation of the compactification of W. Fulton and R. McPherson \cite{F-M}, which allows to consider a compactification by unit blocks consisting in irreducible real extensions $\widetilde F_{v^{j'_\delta }_{j_\delta }}$~, $1\le j'_\delta \le j_\delta $~, of rank $N$ in one-to-one correspondence with the irreducible real completions $F_{v^{j'_\delta }_{j_\delta }}$~, as defined in section 1.1.5.
\vskip 11pt

Let $ \{\tilde g\RL^{(n)}[j_\delta ,m_{j_\delta }]\}^t_{j_\delta =1}$ denote the set of conjugacy class representatives of the real non compact bilinear algebraic semigroup 
$G ^{(n)}(\widetilde F^{+}_{\o v}\times \widetilde F^{+}_v)$ and let 
\begin{align*}
\{\tilde g_L^{(n)}[j_\delta ,m_{j_\delta }]\}^t_{j_\delta =1}
&\subset G ^{(n)}(\widetilde F^{+}_{v})\simeq  T_n(\widetilde F^{+}_{v})\\
\text{(resp.} \quad
\{\tilde g_R^{(n)}[j_\delta ,m_{j_\delta }]\}^t_{j_\delta =1}
&\subset G ^{(n)}(\widetilde F^{+}_{\o v})\simeq  T^t_n(\widetilde F^{+}_{\o v})\ )\end{align*}
be the corresponding set of \lr conjugacy class representatives of the real non compact linear algebraic semigroup
$ G ^{(n)}(\widetilde {F^{+}_{v}}) $ (resp.
$ G ^{(n)}(\widetilde {F^{+}_{\o v}}) $~).
\vskip 11pt

Referring to proposition 2.3.4, it is easy to see that the rank $r_{v_{j_\delta }}$ (resp. $r_{\o v_{j_\delta }}$~)of the conjugacy class representative
$\tilde g_L^{(n)}[j_\delta ,m_{j_\delta }]$ (resp. $\tilde g_R^{(n)}[j_\delta ,m_{j_\delta }]$~) is given by:
\[ r_{v_{j_\delta }}=N^n\cdot j^n \qquad \text{(resp.} \quad  r_{\o v_{j_\delta }}=N^n\cdot j^n\ )\]
where:
\Bi
\item $N$ is the rank of the real irreducible extensions $\widetilde F_{v^{j'_\delta }_{j_\delta }}$~,
\item $j$ is the global class residue degree of $\widetilde F_{v_{j_\delta }}$~.
\Ei
\vskip 11pt

If $n_F$ is the number of non-units of $\widetilde F_{v^{j'_\delta }_{j_\delta }}$~,
then:
\Bi
\item the number of points of 
$\tilde g_L^{(n)}[j_\delta ,m_{j_\delta }]$ (resp. $\tilde g_R^{(n)}[j_\delta ,m_{j_\delta }]$~) is:
\begin{align*}
n_{ \tilde g_L^{(n)}[j_\delta ,m_{j_\delta }] } &=j^n\cdot N^n\cdot n^n_F\\
\text{(resp.} \quad 
n_{\tilde g_R^{(n)}[j_\delta ,m_{j_\delta }]} &=j^n \cdot N^n\cdot n^n_F\ ),\end{align*}

\item the number of irreducible subsets $\widetilde F_{v^{j'_\delta}_{j_\delta }}(n)$ of
$ \tilde g_L^{(n)}[j_\delta ,m_{j_\delta }] $ (resp. $ \tilde g_R^{(n)}[j_\delta ,m_{j_\delta }] $~) is:
\[ n_{\widetilde F_{v^{j'_\delta }_{j_\delta }}}=j^n \qquad \text{(resp.} \quad
n_{\widetilde F_{\o v_{j_\delta }}}=j^n\ ).\]
\Ei
\vskip 11pt 

The compactification of the conjugacy class representatives 
$\tilde g_L^{(n)}[j_\delta ,m_{j_\delta }] $ and 
$ \tilde g_R^{(n)}[j_\delta ,m_{j_\delta }] $ is realized symmetrically respectively in the upper and in the lower half spaces according to the method developed by W. Fulton and R. McPherson in \cite{F-M}: this leads to the following propositions.
\vskip 11pt

\subsubsection{Proposition}
{\em At each irreducible subset $\widetilde F_{v^{j'_\delta}_{j_\delta }}(n) \in 
\tilde g_L^{(n)}[j_\delta ,m_{j_\delta }] $ (resp. 
$\widetilde F_{\o v^{j'_\delta}_{j_\delta }}(n) \in 
\tilde g_R^{(n)}[j_\delta ,m_{j_\delta }] $~) at $N\cdot n_F$ elements, corresponds a non-singular irreducible compact completion 
$F_{v^{j'_\delta}_{j_\delta }}(n) $ (resp. $F_{\o v^{j'_\delta}_{j_\delta }}(n) $~), i.e. a closed irreducible compact one-dimensional subset of 
$g_L^{(n)}[j_\delta ,m_{j_\delta }] \in G^{(n)}({F^+_v})$ (resp.
$g_R^{(n)}[j_\delta ,m_{j_\delta }] \in G^{(n)}({F^+_{\o v}})$~) such that:
\Bean
\item the union of the irreducible completions 
$F_{v^{j'_\delta}_{j_\delta }}(n) $ (resp. $F_{\o v^{j'_\delta}_{j_\delta }}(n) $~) is
$g_L^{(n)}[j_\delta ,m_{j_\delta }] $ (resp. $g_R^{(n)}[j_\delta ,m_{j_\delta }] $~);
\item these irreducible completions meet transversally.
\Ee}
\vskip 11pt

\bpr \Be
\item As the conjugacy class representatives $ \widetilde g_L^{(n)}[j_\delta ,m_{j_\delta }]$ (resp. $\widetilde g_R^{(n)}[j_\delta ,m_{j_\delta }] $~) are built from irreducible subsets 
 $\widetilde F_{v^{j'_\delta}_{j_\delta }}(n) $ (resp. 
$\widetilde F_{\o v^{j'_\delta}_{j_\delta }}(n) $~), it is normal that the compactification acts on these unit blocks.

\item A	set of blowups is then envisaged on the $N\cdot n_F$ elements of the irreducible subsets
$\widetilde F_{v^{j'_\delta}_{j_\delta }}(n) $ and $\widetilde F_{\o v^{j'_\delta}_{j_\delta }}(n) $in such a way that the $N\cdot n_F$ points of
$\widetilde F_{v^{j'_\delta}_{j_\delta }}(n) $ form an irreducible divisor (or completion) $F_{v^{j'_\delta}_{j_\delta }}(n) $ cent red on a point in the upper half space and that the 
$N\cdot n_F$ points of
$\widetilde F_{\o v^{j'_\delta}_{j_\delta }}(n) $ form a symmetric divisor $D(F_{\o v^{j'_\delta}_{j_\delta }}(n) )$ in the lower half space centered on a symmetric point.

A sequence of maximum $N\cdot n_F$ blowups \cite{F-M} is envisaged such that the $N\cdot n_F$ points become together a point and generate a completion (i.e. a closed irreducible compact one-dimensional subset) which constitutes a limiting configuration in the compactification of 
$\widetilde F_{v^{j'_\delta}_{j_\delta }}(n) $ (resp.
$\widetilde F_{\o v^{j'_\delta}_{j_\delta }}(n) $~), in the sense that equivalent compactifications of
$\widetilde F_{v^{j'_\delta}_{j_\delta }}(n) $ (resp.
$\widetilde F_{\o v^{j'_\delta}_{j_\delta }}(n) $~) exist (they are called degenerated in \cite{F-M}).

\item More concretely:
\Bi
\item a first irreducible subset 
$\widetilde F_{v^{1'_\delta}_{j_\delta }}(n) $ is compactified by blowups of its $N\cdot n_F$ points into the completion
$F_{v^{1'_\delta}_{j_\delta }}(n) $~;

\item a second irreducible subset $\widetilde F_{v^{2'_\delta}_{j_\delta }}(n) $ is compactified in the completion
$F_{v^{2'_\delta}_{j_\delta }}(n)$ by blowing up its $(N\cdot n_F)$ points in a compact neighbourhood of
$F_{v^{1'_\delta}_{j_\delta }}(n) $ in such a way that
$F_{v^{2'_\delta}_{j_\delta }}(n) $ be transversal to $F_{v^{1'_\delta}_{j_\delta }}(n) $~;

\item it is proceeded in that manner for all irreducible subsets of 
$\tilde g_L^{(n)}[j_\delta ,m_{j_\delta }]$ in such a way that the compactification of these gives a polyhedral convex cone
$g_L^{(n)}[j_\delta ,m_{j_\delta }]$ with  $j^n_\delta $ completions or closed one-dimensional compact subsets meeting transversally (or having normal crossings);

\item to the polyhedral convex cone
$g_L^{(n)}[j_\delta ,m_{j_\delta }]$ localized in the upper half space corresponds a symmetric polyhedral convex cone
$g_R^{(n)}[j_\delta ,m_{j_\delta }]$ localized in the lower half space.\epr
\Ei\Ee
\vskip 11pt 

\subsubsection{Definitions}

\Bi
\item A {\bf \lr polyhedral convex cone\/} 
$\sigma _L^{(n)}(F^+_{v_{j_\delta ,m_{j_\delta }}} )$ (resp. 
$\sigma _R^{(n)}(F^+_{\o v_{j_\delta ,m_{j_\delta }}} )$~) is a conjugacy class representative 
$g_L^{(n)}[j_\delta ,m_{j_\delta }]$ (resp. $g_R^{(n)}[j_\delta ,m_{j_\delta }]$~) of the linear complete algebraic semigroup
$G^{(n)}({F^+_v})$ (resp. $G^{(n)}({F^+_{\o v}})$~) on the $(\ZZ \big/N\ \ZZ )$-lattice 
$\Lambda _v$ (resp. $\Lambda _{\o v}$~) (corresponding to the 
$(\ZZ \big/N\ \ZZ )$-lattice 
$\Lambda _\omega $ (resp. $\Lambda _{\o \omega }$~) introduced in section 2.2 in such a way that 
$g_L^{(n)}[j_\delta ,m_{j_\delta }]$ (resp. $g_R^{(n)}[j_\delta ,m_{j_\delta }]$~) be included into the  subsemigroup 
 $g^{(2n)}_L[j,m_j]$
(resp. $g^{(2n)}_R[j,m_j]$)
 localized in the upper (resp. lower) half space) \cite{Dan}.
\vskip 11pt

\item A {\bf polyhedral convex bicone\/} $\sigma \RL^{(n)}(F^+_{\o v_{j_\delta ,m_{j_\delta }}} \times 
F^+_{v_{j_\delta ,m_{j_\delta }}} )$ is a conjugacy class representative 
$g\RL^{(n)}[j_\delta ,m_{j_\delta }]$ of the bilinear complete semigroup 
$G^{(n)}({F^+_{\o v}} \times {F^+_{v}})$ on the 
$(\ZZ \big/N\ \ZZ )^2$-bilattice $\Lambda _{\o v} \otimes\Lambda _{v }$~.
\vskip 11pt

\item A {\bf smooth linear general  semivariety\/} \cite{K-K-M-S}
$\tau _\sigma ^{(n)}({F^+_v})$ (resp.\linebreak $\tau _\sigma ^{(n)}({F^+_{\o v}})$~) is a linear complete algebraic semigroup
$G ^{(n)}({F^+_v})$ (resp. $G^{(n)}({F^+_{\o v}})$~) composed of the family
$\{g_L^{(n)}[j_\delta ,m_{j_\delta }]\}$ (resp. $\{g_R^{(n)}[j_\delta ,m_{j_\delta }]\}$~) of disjoint conjugacy class representatives together with a collection of charts from real conjugacy class representatives 
$g^{(n)}_L[j_\delta ,m_{j_\delta }]$ (resp. $g^{(n)}_R[j_\delta ,m_{j_\delta }]$~) 
to their complex equivalents
$g^{(2n)}_L[j,m_{j}]$ (resp. $g^{(2n)}_R[j,m_{j}]$~):
\begin{alignat*}{5}
c z^{j } &:& \quad \{g_L^{(n)}[j_\delta ,m_{j_\delta }] \}_{m_{j_\delta }}
&\To & g^{(2n)}_L[j,m_{j}] \\
\text{(resp.} \quad 
c z^{*j } &:& \quad \{g_R^{(n)}[j_\delta ,m_{j_\delta }] \}_{m_{j_\delta }}
&\To & g^{(2n)}_R[j,m_{j}] 
\ )\end{alignat*}
where $c z^{j}$ (resp.  $c z^{*j}$~) are coordinate functions on the corresponding conjugacy class representatives.

The $cz^{j}$ (resp.  $cz^{*j}$~) can be represented as Laurent monomials in functions of the complex variables $z^{j}_1,\cdots,z^{j}_n$ (resp. $z^{*j}_1,\cdots,z^{*j}_n$~).
\vskip 11pt

\item A {\bf compactified smooth linear general complex semivariety\/} 
$\tau _{c_\sigma} ^{(2n)}({F_\omega })$ (resp.\linebreak $\tau _{c_\sigma }^{(2n)}({F_{\o \omega }})$~) is obtained from
$\tau _\sigma ^{(n)}({F^+_v})$ (resp. $\tau _\sigma ^{(n)}({F^+_{\o v}})$~) by gluing together the conjugacy class representatives of 
$G ^{(2n)}({F_\omega })$ (resp. $G ^{(2n)}({F_{\o \omega }})$~) in such a way that, if
\begin{alignat*}{3}
f_{\omega _{j,m_j}}(z^{j}): & \quad g_L^{(2n)}[j ,m_{j }] & \To & \CC\\
\text{(resp.} \quad 
f_{\o \omega _{j,m_j}}(z^{*j}): & \quad g_R^{(2n)}[j ,m_{j }]& \To & \CC^*\ )\end{alignat*}
is a  function on $g_L^{(2n)}[j ,m_{j }]$ (resp. $g_R^{(2n)}[j ,m_{j }]$~),\\
then a Laurent polynomial
\begin{align*}
 f_\omega (z): \quad &G^{(2n)}(F_\omega )\To {\CC} \\[11pt]
\mbox{(resp.} \quad
&f_{\o \omega }(z^*): \quad G^{(2n)}(F_{\o \omega })\To {\CC^*} \ )\end{align*}
given by:
\begin{align*}
f_\omega (z) &= \sum^{r }_{j,m_j =1}f_{\omega_{j,m_j}} (z^{j})\;, \quad 1\le j \le r \le \infty\;, \\
\text{(resp.} \quad 
f_{\o \omega }(z^*) &= \sum^{r }_{j,m_j=1}f_{\o \omega _{j,m_j}}(z^{*j})\ )\end{align*}
can be introduced on 
$\tau _{c_\sigma} ^{(2n)}({F_{\omega _\oplus}})$ (resp. $\tau _{c_\sigma }^{(2n)}({F_{\o \omega _\oplus}})$~) as a(n) (in)finite linear combination of Laurent monomials corresponding to the \lr polyhedral convex cones 
$\sigma _L^{(n)}(F^+_{v_{j_\delta ,m_{j_\delta }}}) $ (resp. 
$\sigma _R^{(n)}(F^+_{\o v_{j_\delta ,m_{j_\delta }}}) $~).
\Ei
\vskip 11pt 

\subsubsection{Proposition}

{\em Let $\tau _{c_\sigma} ^{(2n)}({F_\omega })$ (resp. $\tau _{c_\sigma }^{(2n)}({F_{\o \omega }})$~) 
be a compactified smooth linear general  semivariety whose conjugacy class representatives 
$g_L^{(2n)}[j ,m_{j }] $ (resp. $g_R^{(2n)}[j ,m_{j }] $~) are glued together and on which the 
 functions $f_{\omega }(z^{j })$ (resp. $f_{\o \omega }(z^{*j })$~) are considered.
\vskip 11pt

Then, on the compactified  semivariety $\tau _{c_\sigma} ^{(2n)}({F_\omega })$ (resp. $\tau _{c_\sigma }^{(2n)}({F_{\o \omega }})$~), the (differentiable) function $f_\omega (z)$ (resp. $f_{\o \omega }(z^*)$~), defined in a neighbourhood of a point $z_0$ (resp. $z_0^*$~) of $\CC^n$~, is holomorphic at $z_0$ (resp. $z_0^*$~) if we have the following multiple power series development:
\begin{align*}
 f_\omega (z)&=\sum_{j =1}^\infty \sum_{m_{j }}c_{j ,m_{j }}\ ( z'_1-z_{0_1})^{j} 
\cdots (z'_n-z_{0_n})^{j}\;, \\[11pt]
 \text{(resp.} \quad 
f_{\o \omega }(z^*)&=\sum_{j =1}^\infty \sum_{m_{j }}c^*_{j ,m_{j }}\ (z^{'*}_1-z^{'*}_{0_1})^{j}
\cdots (z^{'*}_n-z^{'*}_{0_n})^{j}\ ) \end{align*}
where 
$z'_1,\cdots,z'_n$ are (functions of) complex variables, on unitary closed supports
\par
taking into account that:
\Bean
\item to each conjugacy class representative $g_L^{(2n)}[j ,m_{j }] \in 
\tau _{c_\sigma} ^{(2n)}({F_\omega })$  (resp. $g_R^{(2n)}[j ,m_{j }] \in 
\tau _{c_\sigma} ^{(2n)}({F_{\o \omega }})$~) corresponds a term of the multiple power series development of $f_\omega (z)$ (resp. $f_{\o \omega }(z^*)$~).
These terms have the same structure by construction of $g_L^{(2n)}[j,m_j]$
(resp. $g_R^{(2n)}[j,m_j]$~), $1\le j\le r\le \infty $~.
\item the coefficients $c_{j ,m_{j }}$ (resp. $c^*_{j ,m_{j }}$~) have an inflation action from functions on the conjugacy class representatives $P^{(2n)}(F_{\omega ^1_{j }})$ (resp. $P^{(2n)}(F_{\o \omega ^1_{j }})$~) of the ``complex'' parabolic subsemigroup to the corresponding functions on the conjugacy class representatives $g^{(2n)}_L[j ,m_{j }]$ (resp. $g^{(2n)}_R[j ,m_{j }]$~) of the algebraic semigroup $G^{(2n)}({F_{\omega }})$ (resp. $G^{(2n)}({F_{\o \omega }})$~).
\item $j \to\infty$~.
\Ee}
\vskip 11pt

\bpr
\Bean
\item According to definitions 3.1.4, the coordinate functions $cz^{j }$ (resp. $cz^{*j }$~) on the conjugacy classes of $\tau _{c_\sigma} ^{(2n)}({F_\omega })$  (resp. $\tau _{c_\sigma} ^{(2n)}({F_{\o \omega }})$~) can be represented as Laurent monomials.  So, each term of the multiple power series development of $f_\omega (z)$ (resp. $f_{\o \omega }(z^*)$~) corresponds to a conjugacy class representative $g_L^{(2n)}[j ,m_{j }] $ (resp. $g_R^{(2n)}[j ,m_{j }] $~).

\item If $j \to \infty$~, the number of conjugacy classes of $\tau _{c_\sigma} ^{(2n)}({F_\omega })$  (resp. $\tau _{c_\sigma} ^{(2n)}({F_{\o \omega }})$~)  must tend to infinity in order that the multiple power series, converging to $z$ in some neighbourhood of $z_0$~, be equal there to $f_\omega (z)$ (resp. $f_{\o \omega }(z^*)$~).\epr
\Ee
\pagebreak
 
\subsubsection{Proposition}

{\em The coefficients $c_{j ,m_{j }}$ of the holomorphic function $f_\omega (z)$ (resp. $f_{\o \omega }(z^*)$~) on the compactified algebraic semivariety $\tau _{c_\sigma} ^{(2n)}({F_\omega })$  (resp. $\tau _{c_\sigma} ^{(2n)}({F_{\o \omega }})$~) are given by the products:
\[\lambda ^{\half}(2n,j ,m_{j })=\prod^{2n}_{d=1}\lambda ^{\half}_d(2n,j ,m_{j })\]
of the square roots of the eigenvalues $\lambda _d(2n,j ,m_{j })$ of the 
$(j ,m_{j })$-th coset representatives 
$U_{j ,m_{j_{R }}} \times U_{j ,m_{j_{L }}} $ of the product $T_R(n;t)\otimes T_L(n;t)$ of the Hecke operators.
}
\vskip 11pt

\bpr As the holomorphic function $f_\omega (z)$ (resp. $f_{\o \omega }(z^*)$~) decomposes into a multiple power series whose terms are (sub)functions on the conjugacy class representatives $g_L^{(2n)}[j ,m_{j }] $ (resp. $g_R^{(2n)}[j ,m_{j }] $~) of the  semivariety $\tau _{c_\sigma} ^{(2n)}({F_\omega })$  (resp. $\tau _{c_\sigma} ^{(2n)}({F_{\o \omega }})$~), $f_\omega (z)$ (resp. $f_{\o \omega }(z^*)$~) corresponds to an endomorphism of 
 $\tau _{c_\sigma} ^{(2n)}({F_\omega })$  (resp. $\tau _{c_\sigma} ^{(2n)}({F^+_{\o \omega }})$~) into itself.

On the other hand, as the  semivariety  $\tau _{c_\sigma} ^{(2n)}({F_\omega })$  (resp. $\tau _{c_\sigma} ^{(2n)}({F_{\o \omega }})$~) is defined on the $(\ZZ\big/N\ \ZZ)$-lattice $\Lambda _\omega $ (resp. $\Lambda _{\o \omega }$~) (see definitions 3.1.4), on which acts the Hecke operator $T_L(n;t)$ (resp. $T_R(n;t)$~) whose product $T_R(n;t)\otimes T_L(n;t)$ has a representation in the subgroup of matrices $\GL_n((\ZZ\big/N\ \ZZ)^2)$~, the product of the square roots of the eigenvalues $\lambda (2n,j ,m_{j })$ of the $(j ,m_{j })$-th coset representative $U_{j ,m_{j_{R} }}\times U_{j ,m_{j_{L }}}$ of $T_R(n;t)\otimes T_L(n;t)$ will naturally constitute the searched coefficient $c_{j_ ,m_{j }}$ of the holomorphic function $f_\omega (z)$~: this results from proposition 2.2.5.\epr
\vskip 11pt 

\subsubsection{Laurent polynomial on the general  bisemivariety}

Let $G ^{(n)}(\widetilde F^+_{\o v} \times \widetilde F^+_{v})$ denote the non-compact real bilinear algebraic semigroup and let $G ^{(n)}({F^+_{\o v}} \times {F^+_{\o v}})$ be its locally compact equivalent obtained by the Fulton-McPherson compactification as introduced in proposition 3.1.4.  To the set of polyhedral convex bicones
$\sigma \RL^{(n)}(F^+_{\o v_{j_\delta ,m_{j_\delta }}} \times F^+_{v_{j_\delta ,m_{j_\delta }}})
\equiv g^{(n)}\RL[j_\delta ,m_{j_\delta }]$ of $G ^{(n)}({F^+_{\o v}} \times {F^+_{v}})$  corresponds the smooth general  bisemivariety $\tau _\sigma  ^{(n)}({F^+_{\o v}} \times {F^+_{v}})$ which leads to its compactified complex   equivalent $\tau _{c_\sigma } ^{(2n)}({F_{\o \omega }} \times {F_{\omega }})$ obtained by gluing together the conjucacy class representatives of\linebreak $G ^{(n)}({F^+_{\o v }} \times {F^+_{v }})$~.  And,  a Laurent (bi)polynomial
\[f_{\o \omega }(z^*)\otimes_D f_\omega (z)=\sum^{r }_{j =1}\sum_{m_{j }} c_{j ,m_{j }}
\cdot c^*_{j ,m_{j }}\ z^{*j}\cdot z^{j}\]
can be introduced on $\tau _{c_\sigma } ^{(2n)}({F_{\o \omega _\oplus}} \times {F_{\omega _\oplus}})$~.
\vskip 11pt 
 
\subsubsection{Proposition}

{\em On the compactified complex bisemivariety $\tau _{c_\sigma } ^{(2n)}({F_{\o \omega _\oplus}} \times {F_{\omega _\oplus}})$~, the bifunction $f_{\o \omega }(z^*)\otimes f_\omega (z)$~, defined in the neighbourhood of a bipoint $(z^*_0\times z_0)$ of $\CC^n\times \CC^n$~, is holomorphic at $(z^*_0\times z_0)$ if we have the multiple power series development:
\[ f_{\o \omega }(z^*)\otimes_{(D)} f_\omega (z)=\sum_{j =1}^\infty\sum_{m_{j }}c^*_{j ,m_{j }}\ 
c_{j ,m_{j }}\ (z^{*'}_1z'_1-z^*_{0_1}z_{0_1})^{j}\cdots(z^{*'}_nz'_n-z^*_{0_n}z_{0_n})^{j}\]
 verifying:
\Bean
\item to each conjugacy class representative 
$g^{(2n)}\RL[j ,m_{j }]\in \tau  ^{(2n)}_{c_\sigma }({F_{\o \omega }} \times {F_{\omega }})$  
corresponds a term of $f_{\o \omega }(z^*)\otimes f_\omega (z)$~;
\item $j \to\infty$~;
\item each term of $f_{\o \omega }(z^*)\otimes_{(D)} f_\omega (z)$ generates a function subspace of dimension $j^{2n}$~.
\Ee
}
\vskip 11pt

\bpr This is an adaptation of proposition 3.1.6 to the bilinear case.\epr
\vskip 11pt 

\subsubsection{Corollary}

{\em The holomorphic bifunction $f_{\o \omega }(z^*)\otimes_{(D)} f_\omega (z)$ on the compactified algebraic bisemivariety
$\tau  ^{(2n)}_{c_\sigma }({F_{\o \omega _\oplus}} \times {F_{\omega _\oplus}})$  constitutes an irreducible holomorphic representation $\Irr\hol(\GL_n({F_{\o \omega }} \times {F_{\omega }}))$  of the bilinear complete semigroup $G ^{(2n)} ({F_{\o \omega }} \times {F_{\omega }})\approx \GL_n ({F_{\o \omega }} \times {F_{\omega }})$~.
}
\vskip 11pt  

\subsubsection{Power series development on the real compactified semigroups}

If the inclusion $G^{(n)}(F^+_{\o v}\times F^+_{v}\hookrightarrow G^{(2n)}(F_{\o\omega} \times F_\omega )$ of the real compactified bilinear complete semigroup 
$G^{(n)}(F^+_{\o v}\times F^+_v)$ into its complex equivalent $ G^{(2n)}(F_{\o\omega} \times F_\omega )$ is not taken into account, a power series development can be envisaged on
$G^{(n)}(F^+_{v})=\{g_L^{(n)}[j_\delta ,m_{j_\delta }]\}$
(resp. $G^{(n)}(F^+_{\o v})=\{g_R^{(n)}[j_\delta ,m_{j_\delta }]\}$~).

Consider the compactified smooth linear general real semivariety
$\tau ^{(n)}_{c_\sigma }(F^+_{v})$
(resp. $\tau ^{(n)}_{c_\sigma }(F^+_{\o v})$~) obtained from
$\tau ^{(n)}_{\sigma }(F^+_{v})$
(resp. $\tau ^{(n)}_{\sigma }(F^+_{\o v})$~) by gluing together the conjugacy class representatives of
$G^{(n)}(F^+_{v})$
(resp. $G^{(n)}(F^+_{\o v})$~) in such a way that:
\begin{align*}
f_{v_{j_\delta ,m_{j_\delta }}}(x^{j_\delta }) : & \quad g^{(n)}_L[j_\delta ,m_{j_\delta }] \To \rit\;, \quad &x&=(x_1,\dots,x_n)\in\rit^n\;, \\
\text{(resp.} \quad
f_{\o v_{j_\delta ,m_{j_\delta }}}(-x^{j_\delta }) : & \quad g^{(n)}_R[j_\delta ,m_{j_\delta }] \To \rit\;) \end{align*}
be a function on
$g_L^{(n)}[j_\delta ,m_{j_\delta }]$
(resp. $g_R^{(n)}[j_\delta ,m_{j_\delta }]$~).
\vskip 11pt

\subsubsection{Proposition}
{\em On the compactified real semivariety
$\tau ^{(n)}_{c_\sigma }(F^+_{v})$
(resp. $\tau ^{(n)}_{c_\sigma }(F^+_{\o v})$~), the (differentiable) function
$f_v(x)$ (resp. $f_{\o v}(x)$) has a power series development if it can be expressed according to
\begin{align*}
f_v(x) &=\sum_{j_\delta }\sum_{m_{j_\delta }}c_{j_\delta ,m_{j_\delta }}(x')^{j_\delta }\\
\text{(resp.} \quad 
f_{\o v}(-x) &=\sum_{j_\delta }\sum_{m_{j_\delta }}\o c_{j_\delta ,m_{j_\delta }}(-x')^{j_\delta }\;),
\end{align*}
where $x'$ is a function of the real variables $x=(x_1,\dots,x_n )$.

This implies that:
\begin{align*}
f_{v_{j_\delta ,m_{j_\delta }}}(x^{j_\delta }) &= c_{j_\delta ,m_{j_\delta }}(x')^{j_\delta }\\
\text{(resp.}\quad
f_{\o v_{j_\delta ,m_{j_\delta }}}(-x^{j_\delta }) &= \o c_{j_\delta ,m_{j_\delta }}(-x')^{j_\delta }\;).
\end{align*}
}

\bpr The power series development of $f_v(x)$ (resp. $f_{\o v}(-x)$) on the complete algebraic semigroup $G^{(n)}(F^+_v)$ (resp. $G^{(n)}(F^+_{\o v})$~) is rather ``natural'' since its conjugacy class representatives have the same structure by construction.  So,  the terms of $f_v(x)$ (resp. $f_{\o v}(-x)$~) on the conjugacy class representatives have the same  functional structure.\epr
\vskip 11pt

\subsubsection{Corollary}

{\em
The power series bifunction $f_{\o v}(-x)\otimes_{(D)}f_v(x)$ on the compactified algebraic real bisemivariety
$\tau ^{(n)}_{c_\sigma }(F^+_{\o v_\oplus}\times F^+_{v_\oplus})$ constitutes an irreducible functional representation space\linebreak
$\Irr\Repsp (\GL_n(F^+_{\o v}\times F^+_v))$ of the bilinear complete semigroup
$G^{(n)}(F^+_{\o v}\times F^+_v)$.
}
\vskip 11pt


\subsubsection{Proposition}

{\em Let 
$\widetilde F^+_{v_\oplus}=\bigoplus_{j_\delta }\widetilde F^+_{v_{j_\delta }}\bigoplus_{m_{j_\delta }}
\widetilde F^+_{v_{j_\delta ,m_{j_\delta }}}$ (resp.
$\widetilde F^+_{\o v_\oplus}=\bigoplus_{j_\delta }\widetilde F^+_{\o v_{j_\delta }}\bigoplus_{m_{j_\delta }}
\widetilde F^+_{\o v_{j_\delta ,m_{j_\delta }}}$~) denote the sum of the \lr real pseudo-ramified extensions of $k$~.  Then, the non compact real bilinear algebraic semigroup 
$G^{(n)}(\widetilde F^+_{\o v_\oplus}\times \widetilde F^+_{v_\oplus})$~, associated with 
$G ^{(n)} ({F^+_{\o v}} \times {F^+_{v}}) $~, constitutes a $n$-dimensional non compact irreducible representation of the product $W^{ab}_{F^+_R}\times W^{ab}_{F^+_L}$ of global Weil groups according to:
\[ \Irr\Rep W^{(n)}_{F^+_R\times F^+_L}: \quad 
W^{ab}_{F^+_R}\times W^{ab}_{F^+_L}\To  G^{(n)} ({ \widetilde F^+_{\o v_\oplus}} \times {\widetilde F^+_{v_\oplus}}) \;.\]
}
\vskip 11pt

\bpr This is an adaptation of proposition 2.4.7 to the real case.\epr
\vskip 11pt 

\subsubsection{Proposition}

{\em On the non compact real bilinear algebraic semigroup $G^{(n)} ({\widetilde F^+_{\o v_\oplus }} \times {\widetilde F^+_{v_\oplus }}) $ there exists the global functional correspondence:
\[\begin{array}[t]{c}
 \Irr\Rep W^{(n)}_{F^+_R\times F^+_L}(
W^{ab}_{F^+_R}\times W^{ab}_{F^+_L})\\
\| \\
G^{(n)} (\widetilde F^+_{\o v_\oplus } \times \widetilde F^+_{v_\oplus }) \end{array}\To
\begin{array}[t]{c}
 \Irr\Repsp (G^{(n)}(F^+_{\o v_\oplus }\times F^+_{v_\oplus }))\\
\| \\
f_{\o v}(-x)\otimes f_v(x) \end{array}\]
\[\begin{array}{rcl}
\searrow  && \nearrow  \\
& G^{(n)}(F^+_{\o v_\oplus }\times F^+_{v_\oplus })& \end{array}\]
\Bi
\item from the sum of the products, right by left, of the equivalence classes of the irreducible $n$-dimensional non-compact 
representation $\Irr\Rep W^{(n)}_{F^+_R\times F^+_L}(
W^{ab}_{F^+_R}\times W^{ab}_{F^+_L})$ of the product
$W^{ab}_{F^+_R}\times W^{ab}_{F^+_L}$ of global Weil groups, given by the non-compact real bilinear algebraic 
semigroup 
$G^{(n)}(\widetilde F^+_{\o v_\oplus }\times \widetilde F^+_{v_\oplus })$

\item to the sum of the products, right by left, of the equivalence classes of the irreducible  representation space
$ \Irr\Repsp (G^{(n)}(F^+_{\o v}\times F^+_v))$ of $G^{(n)}(F^+_{\o v}\times F^+_v)$ given by the power series bifunction 
$f_{\o v}(-x)\otimes f_v(x) $ on the compactified algebraic bisemivariety 
$\tau _{c_\sigma }^{(n)} ({F^+_{\o v_\oplus }} \times {F^+_{v_\oplus }}) $ or on the compact bilinear algebraic semigroup 
$G^{(n)} ({F^+_{\o v}} \times {F^+_{v}}) $~.
\Ei
}
\vskip 11pt

\bpr \Be
\item It has been seen in proposition 3.1.13 that 
$\Irr\Rep W^{(n)}_{F^+_R\times F^+_L}(
W^{ab}_{F^+_R}\times W^{ab}_{F^+_L}) = G^{(n)}(\widetilde F^+_{\o v_\oplus }\times \widetilde F^+_{v_\oplus })$ 
where $G^{(n)}(\widetilde F^+_{\o v_\oplus }\times \widetilde F^+_{v_\oplus })$ is the non-compact real bilinear algebraic semigroup with entries in $\widetilde F^+_{\o v_\oplus }\times \widetilde F^+_{v_\oplus }$~.
\item The Fulton-McPherson compactification by blowups transforms $G^{(n)}(\widetilde F^+_{\o v}\times \widetilde F^+_v)$  into the compact real bilinear algebraic semigroup $G^{(n)}(F^+_{\o v}\times F^+_v)$~.
\item The power series development of $f_{\o v}(-x)\otimes f_v(x) $ on $G^{(n)}(F^+_{\o v_\oplus }\times F^+_{v_\oplus })$ constitutes an irreducible functional representation space of $G^{(n)}(F^+_{\o v}\times F^+_v)$ according to corollary\linebreak 3.1.12.\epr
\Ee
\vskip 11pt 

\subsection{General bilinear cohomology theory}

\subsubsection{Introduction of bilinear cohomology}

As the envisaged mathematical structures are bilinear, we have to introduce a bilinear cohomology in one-to-one correspondence with its linear equivalent.

The existence of a bilinear cohomology is justified by:
\Bi
\item the cycle map;
\item the Tannakian category of representations of affine groups schemes.
\Ei
\vskip 11pt

The considered bilinear cohomology will be ``general'' in the sense that it will be a cohomology theory on abstract  bisemivarieties in one-to-one correspondence with the associated affine algebraic bisemivarieties covering those \cite{Hart}.

{\bbf The affine algebraic bisemivarieties will be affine bilinear algebraic semigroups 
$G^{(n)}(\wt F^+_{\o v}\times \wt F^+_v)$ embedded isomorphically in their complex equivalents
$G^{(n)}(\wt F_{\o \omega }\times \wt F_\omega )$ and covering the complete algebraic bilinear semigroups
$G^{(n)}( F^+_{\o v}\times  F^+_v )$\/}~.  Note that the conjugacy classes of
$G^{(n)}( F^+_{\o v}\times  F^+_v )$  are isomorphic to these of
$G^{(n)}(\wt F^+_{\o v}\times \wt F^+_v )$ by a suitable compactification described in section 3.1.
\vskip 11pt

The following proposition will define and justify {\bbf the general bilinear cohomology theory on smooth abstract (resp. algebraic) bisemivarieties\/}
$G^{(n)}( F^+_{\o v}\times  F^+_v )$
(resp. $G^{(n)}(\wt F^+_{\o v}\times \wt F^+_v )$)
{\bbf as the (graded) right derived bifunctor
$H^{2*}(G^{(n)}(F^+_{\o v}\times  F^+_v ),\FREPSP(\GL_{2*}( F^+_{\o v}\times  F^+_v )))$
(resp. $H^{2*}(G^{(n)}(\wt F^+_{\o v}\times  \wt F^+_v ),\FREPSP(\GL_{2*}( \wt F^+_{\o v}\times \wt F^+_v )))$\/}
of the functional representation spaces of the complete (resp. algebraic) bilinear parabolic subsemigroups
$P^{2*}(F^+_{\o v^1}\times F^+_{v^1})$
(resp. $P^{2*}(\wt F^+_{\o v^1}\times \wt F^+_{v^1})$) having trivial actions where
$\FREPSP(\GL_{2*}(\wt F^+_{\o v}\times \wt F^+_{v})$ denotes the functional representation spaces of the graded bilinear algebraic semigroups
\[\GL_{2*}(\wt F^+_{\o v}\times \wt F^+_v )=\bigoplus_{i\le n}\GL_{2i}(\wt F^+_{\o v}\times \wt F^+_v )\;.\]
\vskip 11pt

In connection with the bilinear cohomology semigroup
$H^{2*}(G^{(n)}(F^+_{\o v}\times  F^+_v ),\FREPSP(\GL_{2*}( F^+_{\o v}\times  F^+_v )))$
on the smooth abstract bisemivariety
$G^{(n)}(F^+_{\o v}\times  F^+_v )$~, there exists a morphism
\begin{multline*}
\HH^{2*\to[*,*]}: \quad
H^{2*}(G^{(n)}(F^+_{\o v}\times  F^+_v ),\FREPSP(\GL_{2*}( F^+_{\o v}\times  F^+_v )))\\
\To
H^{[*,*]}(G^{(n)}_d(F^+_{\o v}\times  F^+_v ),\Omega ^{*+*}_{G^{(n)}_d( F^+_{\o v}\times  F^+_v )})\end{multline*}
into the bilinear cohomology of the de Rham type
$H^{[*,*]}(G^{(n)}_d(F^+_{\o v}\times  F^+_v ),\Omega ^{*+*}_{G^{(n)}_d( F^+_{\o v}\times  F^+_v )})$ where:
\Bi
\item the dimension $*$ on the left in $[*,*]$ is a covariant cohomological dimension and the dimension $*$ on the right is a contravariant cohomological dimension;

\item $G^{(n)}_d(F^+_{\o v}\times  F^+_v )$ is a differentiable smooth abstract bisemivariety;

\item $\Omega ^{*+*}_{G^{(n)}_d( F^+_{\o v}\times  F^+_v )}=\bigoplus_i\Omega ^{i+i}_{G^{(n)}_d(F^+_{\o v}\times F^+_v)}$ are the graded bisemisheaves of differential $(i+i)$-forms on
$G^{(n)}_d( F^+_{\o v}\times  F^+_v )$~, i.e. the wedge product of right semisheaves of differential $(i,0)$-forms by left semisheaves of differential $(0,i)$-forms.
\Ei
\vskip 11pt

\subsubsection{Proposition (General bilinear cohomology)}

{\em
{\bbf A general bilinear cohomology theory is a contravariant bifunctor}
\begin{multline*}
H^{2*}: \quad \{\text{smooth abstract (resp. algebraic) bisemivarieties}
\begin{aligned}[t]
G^{(n)}(F^+_{\o v}\times F^+_v) \qquad &\\
\text{(resp.\; } G^{(n)}(\wt F^+_{\o v}\times \wt F^+_v)\; )\}\end{aligned}\\
\To \begin{aligned}[t]
\{ \text{graded (functional) representation spaces of the complete}&\\
\text{(resp. algebraic) bilinear semigroups\; } \GL_{2*}(F^+_{\o v}\times F^+_v)&\\
\text{(resp.\; } \GL_{2*}(\wt F^+_{\o v}\times \wt F^+_v)\ )\}&\end{aligned}
\end{multline*}
{\bbf together with a natural transformation:}
\begin{multline*}
n_{H^{2*}\to H^{[*,*]}}: \quad H^{2*}(G^{(n)}(F^+_{\o v}\times F^+_v),\FREPSP (\GL_{2*}(F^+_{\o v}\times F^+_v)))\\
\qquad \qquad \text{(resp.\; } H^{2*}(G^{(n)}(\wt F^+_{\o v}\times \wt F^+_v),\FREPSP (\GL_{2*}(\wt F^+_{\o v}\times \wt F^+_v)))\ )\\
\To H^{[*,*]}(G^{(n)}_d(F^+_{\o v}\times F^+_v),\Omega ^{*+*}_{G^{(n)}_d(F^+_{\o v}\times F^+_v)})
\end{multline*}
into the bilinear cohomology of the de Rham type of dimension $*+*$~.
\vskip 11pt

Let $H^{2i}(G^{(n)}( F^+_{\o v}\times  F^+_v),\FREPSP (\GL_{2i}( F^+_{\o v}\times  F^+_v)))$ be the general bilinear cohomology semigroup of dimension $2i$~.

Let $H^{2i}(G^{(n)}(\wt F^+_{\o v-v}))$ be the associated linear cohomology group where
$\wt F^+_{\o v-v}=\wt F^+_{\o v}\cup \wt F^+_v$ is the set of algebraic extensions corresponding to
$ F^+_{\o v}\times  F^+_v$~.
\vskip 11pt

The linear cohomology $H^{2i}(G^{(n)}(\wt F^+_{\o v-v}))$ is characterized by the cycle map:
\[\gamma ^i_{G^{(n)}_{\o v-v}}: \quad \hZ^i(G^{(n)}(\wt F^+_{\o v-v}))\To H^{2i}(G^{(n)}(\wt F^+_{\o v-v}))\]
where $\hZ^i(G^{(n)}(\wt F^+_{\o v-v}))$ is the group of algebraic cycles of codimension $i$ over the linear algebraic group $G^{(n)}(\wt F^+_{\o v-v})$~.

The bilinear cohomology $H^{2i}(G^{(n)}(F^+_{\o v}\times F^+_v),\FREPSP(\GL_{2i}(F^+_{\o v}\times F^+_v))$ is characterized by the bisemicycle map
\[ \gamma ^i_{G^{(n)}_{\o v\times v}}: \quad
\hZ^i(G^{(n)}(F^+_{\o v}\times F^+_v))\To
H^{2i}(G^{(n)}(F^+_{\o v}\times F^+_v)),\FREPSP(\GL_{2i}(F^+_{\o v}\times F^+_v))\]
where $\hZ^i(G^{(n)}(F^+_{\o v}\times F^+_v))$ is the bilinear semigroup of compactified bisemicycles of codimension $i$ on the bilinear complete semigroup $G^{(n)}(F^+_{\o v}\times F^+_v)$ in such a way that\linebreak 
$FREPSP(\GL_{2i}(F^+_{\o v}\times F^+_v))$ be embedded isomorphically in its complex equivalent\linebreak
$FREPSP(\GL_{2i}(F^+_{\o \omega }\times F^+_\omega ))$~: this is at the origin of this cycle map.

The associated the de Rham bilinear cohomology is characterized by the bisemicycle map:
\[ \gamma ^{[i,i]}_{G^{(n)}_{\o v-v}}:
\quad \hZ^{[i,i]}(G^{(n)}_d(F^+_{\o v}\times F^+_v))\To H^{[i,i]}(G^{n)}_d(F^+_{\o v}\times F^+_v),\Omega ^{*+*}_{G^{(n)}_d(F^+_{\o v}\times F^+_v)}\]
where $\hZ^{[i,i]}(G^{(n)}_d(F^+_{\o v}\times F^+_v))$ is the bilinear semigroup of compactified bisemicycles of codimension $2i=i+i$ on the bilinear complete semigroup $G^{(n)}_d(F^+_{\o v}\times F^+_v)$~.
\vskip 11pt

{\bbf The correspondence between the linear cohomology and the associated bilinear cohomologies is given by the commutative diagram\/}:
\[
\begin{CD}
\text{{\bf Linear}} @. \text{{\bf Bilinear}}  @. \text{{\bf Bilinear}} 
\\
\hZ^i(G^{(n)}(\wt F^+_{\o v-v})) 
@>>> 
\hZ^{[i,i]}(G^{(n)}_d( F^+_{\o v}\times F^+_v))
@>>> 
\hZ^i(G^{(n)}(F^+_{\o v}\times F^+_v))
\\
@VV{\gamma ^i_{G^{(n)}_{\o v-v}}}V
@VV{\wr \; \gamma ^{[i,i]}_{G^{(n)}_{\o v-v}}}V
@VV{\wr \; \gamma ^{i}_{G^{(n)}_{\o v\times v}}}V
\\
H^{2i}(G^{(n)}(\wt F^+_{\o v-v})) 
@>\sim>>
\begin{aligned}[t]
H^{[i,i]}(G^{(n)}_d( F^+_{\o v}\times F^+_v),&\\ \Omega ^{i+i}_{G^{(n)}_d(F^+_{\o v}\times F^+_v))} &\end{aligned}
@>>>
\begin{aligned}[t]
&H^{2i}(G^{(n)}( F^+_{\o v}\times F^+_v),\\ &\FREPSP (\GL_{2i}(F^+_{\o v}\times F^+_v)))\end{aligned}
\end{CD}
\]
}

\bpr Remark first that the general bilinear cohomology theory is induced from:
\Bean
\item 
The cycle map $\gamma ^i_{G^{(n)}_{\o v\times v}}$ which is related to the linear cycle map $\gamma ^i_{G^{(n)}_{\o v-v}}$ and to the bilinear the de Rham cycle map $\gamma ^{[i,i]}_{G^{(n)}_{\o v-v}}$~.

\item the Tannakian category of representations of affine group schemes \cite{Del4}.
\Ee
\vskip 11pt

The commutative diagram of this proposition clearly shows that the linear cycle map $\gamma ^i_{G^{(n)}_{\o v-v}}$ proceeds from the bilinear cycle map of the de Rham type $\gamma ^{[i,i]}_{G^{(n)}_{\o v-v}}$ and from the isomorphism:
\[H^{2i}(G^{(n)}(\wt F^+_{\o v-v}))\overset{\sim}{\To} H^{[i,i]}(G^{(n)}_d(F^+_{\o v}\times F^+_v),
\Omega ^{i+i}_{G^{(n)}_d(F^+_{\o v}\times F^+_v)})\;.\]
\vskip 11pt

Note that {\bbf the bilinear cohomology of the de Rham type\/}
$ H^{[i,i]}(G^{(n)}_d(F^+_{\o v}\times F^+_v)$ {\bbf has dimension $2i=i+i$} (as the corresponding bilinear semigroup of compactified bisemicycles $\hZ^{[i,i]}(G^{(n)}_d(F^+_{\o v}\times F^+_v))$ from the consideration of bisemisheaves of differential $(i+i)$-forms {\bbf while the general bilinear cohomology} $H^{2i}(G^{(n)}F^+_{\o v}\times F^+_v),\FREPSP(\GL_{2i}(F^+_{\o v}\times F^+_v))$ {\bbf has dimension ``$2i$'' due to its ``diagonal'' isomorphic embedding into\/} $H^{2i}(G^{(n)}(F_{\o \omega }\times F_\omega ),\FREPSP(\GL_{2i}(F_{\o\omega} \times F_\omega ))$~.
\epr
\vskip 11pt

\subsubsection{Proposition}

{\em 
The general bilinear cohomology is a general cohomology theory in the sense that:
\Bi
\item it is a motivic (bilinear) cohomology theory or a Weil (bilinear) cohomology theory;
\item  it is directly related to the standard cohomology theories like the singular, de Rham and Betti cohomologies.
\Ei
}

\bpr \Be
\item The general bilinear cohomology $H^{2*}(G^{(n)}(F^+_{\o v}\times F^+_v)),\FREPSP(\GL_{2*}(F^+_{\o v}\times F^+_v))$ is directly related to:
\Bi
\item the Betti cohomology because the complete (algebraic) bilinear semigroup\linebreak $G^{(n)}(F^+_{\o v}\times F^+_v)$ is covered by the affine bilinear algebraic semigroup $G^{(n)}(\wt F^+_{\o v}\times \wt F^+_v)$ which is embedded isomorphically in its complex equivalent $G^{(n)}(\wt F_{\o \omega }\times \wt F_\omega )$~.

Thus, we have an embedding $\sigma \RL:F^+_{\o v}\times F^+_v\hookrightarrow \CC\times \CC$~.

This allows to define a bilinear Betti cohomology which must be in one-to-one correspondence with the classical linear Betti cohomology according to proposition 3.2.2.

\item the de Rham cohomology if the considered complete (algebraic) bilinear semigroup is differentiable, i.e. is $G^{(n)}_d(F^+_{\o v}\times F^+_v)$~.

The obtained de Rham cohomology $H^{[*,*]}(G^{(n)}_d(F^+_{\o v}\times F^+_v), \Omega ^{*+*}_{G^{(n)}_d(F^+_{\o v}\times F^+_v})$ is bilinear, but, by the commutative diagram of proposition 3.2.2, it is in one-to-one correspondence with the classical linear de Rham cohomology.

\item the singular bilinear cohomology which is in one-to-one correspondence with the de Rham bilinear cohomology.
\Ei

\item The general bilinear cohomology is a motivic or Weil (bilinear) cohomology theory if it verifies the standard conjectures on algebraic (bi)cycles as developed in the following proposition.\epr
\Ee
\vskip 11pt

\subsubsection{Proposition}

{\em The general bilinear cohomology theory, defined by the contravariant bifunctor $\HH^{2*}$ and the natural transformation $n_{H^{2*}\to H^{[*,*]}}$~, is characterized by:
\Be
\item Hodge bisemicycles
$H^{i+i}(G^{(n)}(F_{\o \omega }\times F_\omega ),\FREPSP(\GL_{i+i}(F^+_{\o v}\times F^+_v))$
(resp.\linebreak $H^{i+i}(G^{(n)}(\wt F^+_{\o \omega }\times \wt F^+_\omega ),\FREPSP(\GL_{i+i}(\wt F^+_{\o v}\times \wt F^+_v))$~), $2i=i+i$~, on abstract (resp. algebraic) complex bisemivarieties 
into ``real'' functional representation spaces of $\GL_{i+i}(F^+_{\o v}\times F^+_v)$
in such a way that there is a bifiltration $F^p\RL$ on the right and left cohomology semigroups of $H^{2i}(G^{(n)}(\pt\times\pt),-)$ given by:
\[ F^p\RL H^{2i}( G^{(n)}(\pt\times \pt ),- )=
\bigoplus_{i=p+q} H^{2(p+q)}(G^{(n)}(\pt\times\pt),-)\;.\]

\item a Künneth standard conjecture implying that the projectors 
\[H^{2*}(G^{(n)}(\pt\times\pt),-)\To
H^{2i}(G^{(n)}(\pt\times\pt),-)\]
are induced by algebraic bisemicycles $\CY^i(G^{(n)}(\wt F^+_{\o v}\times\wt F^+_v))\subset
\hZ^i(G^{(n)}(\wt F^+_{\o v}\times \wt F^+_v))$ decomposing into rational subbisemicycles according to the conjugacy class representatives of $\GL_{2i}(\wt F^+_{\o v}\times\wt F^+_v)$~.

\item a Künneth biisomorphism:
\begin{multline*}
H^{2*}(G^{(n)}(F^+_{\o v}),\FREPSP(\GL_{2*}(F^+_{\o v})))\otimes_{F^+_{\o v}\times F^+_v}
H^{2*}(G^{(n)}(F^+_{v}),\FREPSP(\GL_{2*}(F^+_{v})))\\
\To
H^{2*}(G^{(n)}(F^+_{\o v}\times F^+_v),\FREPSP(\GL_{2*}(F^+_{\o v}\times F^+_v)))
\end{multline*}
associated with the exact sequence:
\begin{multline*}
0\To
\sum_{i=p+q} H^{2i}(G^{(n)}(F^+_{\o v}\times F^+_v),\FREPSP ( \GL_{2p+2q}(F^+_{\o v}\times F^+_v)))\\
\overset{\Hs_{i-0}}{\To}
\sum_{i=p+q} H^{2i}(G^{(n)}(F^+_{\o v}\times F^+_v),\FREPSP(\GL_{2p+2q}(F^+_{\o v}\times F^+_v)))\\
\qquad \qquad \qquad + H^0(G^{(n)}(F^+_{\o v}\times F^+_v),\FREPSP(\GL_{1}(F^+_{\o v}\times F^+_v)))\\
\To
H^0(G^{(n)}(F^+_{\o v}\times F^+_v),\FREPSP(\GL_{1}(F^+_{\o v}\times F^+_v)))\To 0
\end{multline*}
in such a way that
\begin{multline*}
H^{2p}(G^{(n)}(F^+_{\o v}\times F^+_v),\FREPSP(\GL_{2p}(F^+_{\o v}\times F^+_v)))\\ \otimes
H^{2i-2p}(G^{(n)}(F^+_{\o v}\times F^+_v),\FREPSP(\GL_{2i-2p}(F^+_{\o v}\times F^+_v)))\\
\To
H^{0}(G^{(n)}(F^+_{\o v}\times F^+_v),\FREPSP(\GL_{1}(F^+_{\o v}\times F^+_v)))\end{multline*}
is the bilinear version of the intersection cohomology according to the Goresky-Mac Pherson approach {\em \cite{G-MP}} and referring to the Poincare duality.
\Ee
}

\bpr \Be
\item The Hodge bisemicycles $H^{i+i}(G^{(n)}(F_{\o \omega }\times F_\omega ),\FREPSP(\GL_{i+i}(F^+_{\o v}\times F^+_v))$ are justified because 
$G^{(n)}(F_{\o \omega }\times F_\omega )$ is an abstract complex bisemivariety embedding its real equivalent
$G^{(n)}(F^+_{\o v }\times F^+_v )$~.  So,
$G^{(n)}(F_{\o \omega }\times F_\omega )$ is a Kähler bisemivariety on which a Hodge decomposition \cite{D-M-O-S}an be performed.

Since abstract bisemivarieties are concerned, the Hodge decomposition must apply on the right and left semivarieties
$G^{(n)}(F_{\o \omega })$ and 
$G^{(n)}(F_{\omega })$ of
$G^{(n)}(F_{\o \omega }\times F_\omega )$ and Hodge bisemicycles follow themselves.

Similarly, a Hodge bifiltration on right and left (abstract) bisemivarieties of
$G^{(n)}(F^+_{\o v}\times F^+_v)$ is justified.

\item The Künneth standard conjecture dealing with the projectors
$H^{2*}(G^{(n)}(\pt\times\pt),-)\to\linebreak
H^{2i}(G^{(n)}(\pt\times\pt),-)$ are induced from algebraic bisemicycles
$\CY^i(G^{(n)}(\wt F^+_{\o v}\times \wt F^+_v))$ because the concerned cohomologies are bilinear cohomologies of abstract bisemivarieties.

On the other hand, these projectors refer to mappings interpreted as inverse quantum deformations of Galois representations as introduced in \cite{Pie6} and applied in \cite{Pie2}.

\item The considered Künneth biisomorphism takes into account the intersection cohomology which is a pairing between the $2p$-th and $2(i-p)$-th bilinear cohomologies in  the $0$-th bilinear cohomology $H^0(G^{(n)}(\pt\times\pt),-)$ yielding a $0$-dimensional bisemicycle because the intersection is composed of a finite set of bipoints.\epr
\Ee

\subsection{Borel-Serre toroidal compactification}
\subsubsection{General bilinear toroidal bisemivariety}

A following step is the transformation of the bilinear general toric bisemivariety 
$\tau _{\sigma }^{(n)} ({\wt F^+_{\o v}} \times {\wt F^+_{v}}) $ 
(or $\tau _{c_\sigma }^{(n)} ({\wt F^+_{\o v}} \times {\wt F^+_{v}}) $~) into the bilinear general toroidal bisemivariety 
$\tau _{T }^{(n)} ({F^{+,T}_{\o v}} \times {F^{+,T}_{v}}) $~.  This can be realized by considering the projective emergent isomorphism 
$\gamma ^c_R\times \gamma ^c_L$ from the bilinear algebraic semigroup $G^{(n)}({\wt F^+_{\o v}}\times {\wt F^+_v})$~, corresponding to 
$\tau _{\sigma }^{(n)} ({\wt F^+_{\o v}} \times {\wt F^+_{v}}) $ or to $\tau _{c_\sigma }^{(n)} ({\wt F^+_{\o v}} \times {\wt F^+_{v}}) $~, to its toroidal equivalent 
$G ^{(n)} ({F^{+,T}_{\o v}} \times {F^{+,T}_{v}}) $~, corresponding to
 $\tau _{T }^{(n)} ({F^{+,T}_{\o v}} \times {F^{+,T}_{v}}) $ as considered in the following section: in fact, the projective emergent isomorphism $\gamma ^c_R\times \gamma ^c_L$ will be more precisely developed for the complex bilinear algebraic semigroup $G ^{(2n)} ({\wt F_{\o \omega }} \times {\wt F_{\omega }}) $~, taking into account the inclusion
\[ G ^{(n)} ({\wt F^{+}_{\o v}} \times {\wt F^{+}_{v}}) \hookrightarrow 
G ^{(2n)} ({\wt F_{\o \omega }} \times {\wt F_{\omega }}) \;.\]
\vskip 11pt



 A toroidal compactification will been envisaged for the lattice bisemispace $X_{S\RL}$~: it will correspond to the Borel-Serre toroidal compactification which will be decomposed into a two step sequence which will lead to the equivalent of a Shimura (bisemi)variety.
This compactification is in fact a projective mapping of $X_{S_{R\times L}}$~, generating products of pairs of embedded complex semitori (which are the only complex compact Lie (semi)groups).\vskip 11pt

\subsubsection{Definition}  

The \lr toroidal projective emergent (iso)morphism $\gamma ^C_L:X_{S_L}\to\o X_{S_L}$ (resp. $\gamma ^C_R:X_{S_R}\to\o X_{S_R}$~) can be decomposed into the two steps sequence \cite{Pie1}:
\Bean
\item the points $P_{a_{L[j,m_j]}}\in \tilde g^{(2n)}_L[j,m_j]$ (resp. $P_{a_{R[j,m_j]}}\in \tilde g^{(2n)}_R[j,m_j]$~) are mapped onto the origin of $\wt F^{2n}_\omega $ (resp. $\wt F^{2n}_{\o\omega }$~).
\item these points $P_{a_{L[j,m_j]}} $ (resp. $P_{a_{R[j,mj]}}
$~) are then projected symmetrically from the origin of 
$\wt F^n_\omega $ (resp. $\wt F^n_{\o\omega}$~) into a connected compact complete semivariety which is a complex semitorus $T^{2n}_L[j,m_j]$ (resp. $T^{2n}_R[j,m_j]$~) in 
$({F^T_\omega })^n$ (resp. $({F^T_{\o\omega }})^n$~) where 
$F^T_\omega =\{F^T_{\omega _1},\cdots,F^T_{\omega _{j,m_j}},\cdots,F^T_{\omega _r}\}$ (resp.
$F^T_{\o\omega} =\{F^T_{\o\omega _1},\cdots,F^T_{\o\omega _{j,m_j}},\cdots,F^T_{\o\omega _r}\}$~)
 is the set of toroidal (compactified) completions
$F^T_{\omega _{j,m_j}}$ (resp. $F^T_{\o\omega _{j,m_j}}$~)
 corresponding to $F_{\omega _{j,m_j}}$(resp. $F_{\o\omega _{j,m_j}}$~).
\Ee
This \lr projective emergent morphism $\gamma ^C_L$ (resp. $\gamma ^C_R$~) is an isomorphism because it is characterized by its representatives $h_{a_{L[j,m_j]}}$ (resp. $h_{a_{R[j,m_j]}}$~) given by the triple
\begin{align*}
& h_{a_{L[j,m_j]}} = \{ P_{a_{L[j,m_j]}},r(P_{a_{L[j,m_j]}}), \gamma _{a_{L[j,m_j]}} \}\\
\text{(resp.} \qquad & h_{a_{R[j,m_j]}} = \{ P_{a_{R[j,m_j]}},r(P_{a_{R[j,m_j]}}),\gamma _{a_{R[j,m_j]}}\}\ )\end{align*}
where:
\Bi
\item $r(P_{a_{L[j,m_j]}})$ (resp. $r(P_{a_{R[j,m_j]}})$~) is the euclidian distance, from the origin, of a point $P_{a_{L[j,m_j]}}$ (resp. $P_{a_{R[j,m_j]}}$~) projected into the $2n$-dimensional semitorus $T^{2n}_{{L[j,m_j]}}$ (resp. $T^{2n}_{{R[j,m_j]}}$~);
\item $\gamma _{a_{L[j,m_j]}} $ (resp. $\gamma _{a_{R[j,m_j]}}$~) is a one-to-one correspondence between the point $P _{a_{L[j,m_j]}}$ (resp. $P _{a_{R[j,m_j]}}$~) and its projective localization given by $r(P _{a_{L[j,m_j]}})$ (resp. $r(P _{a_{R[j,m_j]}})$~).  
\Ei
\vskip 11pt

\subsubsection{Proposition} 

{\em The 
projective emergent isomorphism $\gamma ^C_{R\times L}:X_{S_{R\times L}} 
\to \o X_{S_{R\times L}}$~, mapping the pseudo-ramified lattice bisemispace $X_{S_{R\times L}}=
\GL_n({\widetilde F_{R}}\times {\widetilde F_L})\Big/ \GL_n( ( \ZZ\big/N\ \ZZ)^2 )$ into the corresponding 
toroidal compactified lattice bisemispace: $\o X_{S_{R\times L}}=\GL_n({F^T_{R}}\times 
{F^T_L}\Big/ \GL_n( ( \ZZ\big/N\ \ZZ)^2)$~, is such  that $X_{S_{R\times L}}$ may be viewed as the 
interior of its compactified corresponding $\o X_{S_{R\times L}}$ in the sense of the Borel-Serre compactification.}
\vskip 11pt

\bpr according to definition 3.3.2, the projective emergent isomorphism sends all equivalent representatives of conjugacy classes into their toroidal compactified equivalents\linebreak $g^{(2n)}_{T\RL}[j,m_j]$ consisting in products of $n$-dimensional complex semitori $T^{2n}_R[j,m_j]\times T^{2n}_L[j,m_j]$~.
The inclusion $X_{S_{R\times L}}\hookrightarrow \o X_{S_{R\times L}}$ is a homotopy equivalence so that 
$X_{S_{R\times L}}$ may be considered as the interior of $\o X_{S_{R\times L}}$~.\epr
\vskip 11pt

\subsubsection{Double coset decomposition}

A double coset decomposition of the bilinear toroidal semigroup $\GL_n({F^T_{R }}\times {F^T_L })$  gives rise to the compactified bisemispace:
\[ \o S^{P_n}_{K_n}
=P_n( {F^{T}_{\o\omega^1 }}\times {F^{T}_{\omega^1 }})\setminus
\GL_n({F^{T}_{R}}\times {F^{T}_{L}}
\Big/ \GL_n((\ZZ\big/N\ \ZZ)^2)\]
where
$P_n({F^{T}_{\o\omega^1 }}\times {F^{T}_{\omega ^1}})$ is the bilinear parabolic subsemigroup of matrices associated with the complete (algebraic) bilinear parabolic subsemigroup $P^{(2n)}({F^{T}_{\o\omega^1 }}\times {F^{T}_{\omega^1 }})$ introduced in section 2.4.1.

$P_n({F^{T}_{\o\omega^1 }}\times {F^{T}_{\omega^1 }})$ has the Gauss decomposition as developed in the following section.


Remark that
\[  \GL_n({F^{T}_{R}}\times {F^{T}_{L}})
\Big/ \GL_n((\ZZ\big/N\ \ZZ)^2) \simeq
P_n({F^{T}_{\o\omega ^1}}\times {F^{T}_{\omega^1 }})\setminus \GL_n({F^{T}_{R}}\times {F^{T}_{L}})\]
since $\GL_n({F^{T}_{R}}\times {F^{T}_{L}})
\Big/ \GL_n((\ZZ\big/N\ \ZZ)^2) $
corresponds to the set of lattices of\linebreak $({F^{T}_{\o\omega }})^n\times ({F^{T}_{\omega }})^n$ \cite{Bor1}, \cite{Del3}.
\vskip 11pt 

\subsubsection{Definition:  Levi and Gauss decompositions for (complex) parabolic subgroups}

 $P_n( {F^{T}_{\omega^1 }} )$ (resp. $P_n( {F^{T}_{\o\omega^1 }} )$~) has the standard Levi decomposition: 
\[P_n( {F^{T}_{\omega^1 }} )=\pZ_n( {F^{T}_{\omega^1 }} )\cdot U_n( {F^{T}_{\omega^1 }} )\qquad \text{(resp.}\quad P_n( {F^{T}_{\o\omega^1 }} )=\pZ_n( {F^{T}_{\o\omega^1 }} )\cdot U_n( {F^{T}_{\o\omega^1 }} )\ )\] where:
\Bi
\item $\pZ_n( {F^{T}_{\omega^1 }} )$ (resp. $\pZ_n( {F^{T}_{\o\omega^1 }} )$~) is the \lr centralizer of $ {F^{T}_{\omega^1 }}  $
(resp. ${F^{T}_{\o\omega ^1}} $~) in $T_n( {F^{T}_{\omega^1 }} )$ (resp. $T_n( {F^{T}_{\o\omega^1 }} )^t\subset 
\GL_n({F^T_{\o\omega }}\times {F^T_\omega })$ and is represented by diagonal matrices $d_n( {F^{T}_{\omega^1 }} )$ (resp. $d_n( {F^{T}_{\o\omega^1 }} )$~);
\item $U_n( {F^{T}_{\omega^1 }} )$ (resp. $U_n( {F^{T}_{\o\omega^1 }} )$~) is the \lr unipotent radical of $P_n( {F^{T}_{\omega^1 }} )$ (resp. $P_n( {F^{T}_{\o\omega^1 }} )$~) and is represented by upper (resp. lower) unitriangular matrices 
$u_n( {F^{T}_{\omega^1 }} )$ (resp. $u_n( {F^{T}_{\o\omega^1 }} )^t$~).
\Ei
Similarly, the complex bilinear parabolic subsemigroup $P_n( {F^{T}_{\o\omega^1 }} \times {F^{T}_{\omega^1 }} )\equiv P_n( {F^{T}_{\o\omega^1 }} )\times P_n( {F^{T}_{\omega^1 }} )$~) has the Gauss decomposition:
\[P_n( {F^{T}_{\o\omega^1 }} \times {F^{T}_{\omega^1 }} )=(D_n( {F^{T}_{\o\omega^1 }} )\cdot D_n( {F^{T}_{\omega^1 }} ))(UT_n( {F^{T}_{\o\omega^1 }} )^t\cdot UT_n( {F^{T}_{\omega^1 }} ))\]
according to proposition 2.1.4.
$P_n( {F^{T}_{\o\omega^1 }} \times {F^{T}_{\omega^1 }} )$ is a bilinear  normal subgroup of $\GL_n({F^T_{\o\omega }}\times {F^T_\omega })$ because it is represented by the product of the subgroups of diagonal matrices, isomorphic to maximal (semi)tori, by the subgroups of unipotent matrices \cite{Bor2}.\vskip 11pt

\subsubsection{Representation spaces of bilinear algebraic semigroups}

\Bean
\item The double coset decomposition of $ \GL_n({F^{T}_{R}}\times {F^{T}_{L}})
$ is a compactified bisemispace $\o S^{P_n}_{K_n}$ which corresponds to the  representation space $\Reps (\GL_n({F^{T}_{\o\omega }}\times {F^{T}_{\omega }}))$ of $\GL_n({F^{T}_{\o\omega }}\times {F^{T}_{\omega }})$~.  So we have
\[\o S^{P_n}_{K_n} \simeq
\Reps (\GL_n({F^{T}_{\o\omega }}\times {F^{T}_{\omega }})
)= G^{(2n)}({F^{T}_{\o\omega }}\times {F^{T}_{\omega }})\]
since $G^{(2n)}({F^{T}_{\o\omega }}\times {F^{T}_{\omega }})$ decomposes into conjugacy classes having equivalent representatives $g^{(2n)}_{T\RL}[j,m_j]$~.
\vskip 11pt 

\item Considering the Gauss decomposition of $\GL_n({F^{T}_{\o\omega }}\times {F^{T}_{\omega }})$ and of $P_n({F^{T}_{\o\omega^1 }}\times {F^{T}_{\omega^1 }})$~, we can introduce the double coset decomposition of the diagonal part of $\GL_n({F^{T}_{R}}\times {F^{T}_{L}})$ which is a compactified diagonal bisemispace noted $\o S^{Z[P_n]}_{Z[K_n]}$ and given by:
\[\o S^{Z[P_n]}_{Z[K_n]}= D_n({F^{T}_{\o\omega^1 }}\times {F^{T}_{\omega^1 }})
\setminus  D_n({F^{T}_{R}}\times {F^{T}_{L}})\Big/D_n( (\ZZ\big/N\ \ZZ)^2)\]
where:
\Bi
\item $D_n({F^{T}_{\o\omega^1 }}\times {F^{T}_{\omega^1 }})
\subset P_n({F^{T}_{\o\omega^1 }}\times {F^{T}_{\omega^1 }})$ following definition 3.3.5;
\item $D_n({F^{T}_{R}}\times {F^{T}_{L}})
\subset \GL_n({F^{T}_{R}}\times {F^{T}_{L}})$~.
\Ei
So, as in a), we can state:
\[\o S^{Z[P_n]}_{Z[K_n]}
\simeq \Reps (D_n({F^{T}_{\o\omega }}\times {F^{T}_{\omega }}))
= D^{(n)}({F^{T}_{\o\omega }}\times {F^{T}_{\omega }})\]
with evident notations.\vskip 11pt

\item In the following sections, bilinear cohomology semigroups of $\o S^{P_n}_{K_n}$ will be considered with coefficients in complete algebraic representation spaces $\Repsp(\GL_n({F^T_{\o\omega }\times F^T_\omega }))$ of 
$\GL_n({F^T_{\o\omega }\times F^T_\omega })$~, $1\le n\le \infty $~.

These representation spaces $\Repsp(\GL_n({F^T_{\o\omega }\times F^T_\omega }))$ are $\GL_n({F^T_{\o\omega_\oplus  }\times F^T_{\omega_\oplus } })$-bisemimo\-dules $M_{T_{R_\oplus }}\otimes M_{T_{L_\oplus }}$ (see proposition 2.1.4).  Now, generally, the {\bf{coefficients of the}} {\bbf{(bilinear) cohomology are envisaged in (bisemi)sheaf of rings\/}.}  In this purpose, a (semi)sheaf of rings  $\widehat M^{2i}_{T_R}$ (resp. $\widehat M^{2i}_{T_L}$~) of complex-valued continuous functions $\phi ^{(2i)}_{G^T_{j_R}}(x_{g^T_{j_R}})$ (resp.
  $\phi ^{(2i)}_{G^T_{j_L}}(x_{g^T_{j_L}})$~) on the basic conjugacy class representatives
$g^{(2i)}_{T_R}[j,m_j=0]$ (resp. $g^{(2i)}_{T_L}[j,m_j=0]$~) of 
$\Repsp(T^t_i ({F^T_{\o\omega }}))$ 
(resp. $\Repsp(T_i ({F^T_{\omega }}))$~) can thus be envisaged.

On $M^{2i}_{T_R}\otimes M^{2i}_{T_L}$~, it is then the bisemisheaf of rings $\widehat M^{2i}_{T_R}\otimes \widehat M^{2i}_{T_L}$ which has to be taken into account (in this context, see also section 3.4.10).

So, the developments dealing {\bf{with bilinear cohomology with coefficients\/}} in\linebreak 
$\GL_i({F^T_{\o\omega }\times F^T_\omega })$-bisemimodules $M^{2i}_{T_R}\otimes M^{2i}_{T_L}$~, $1\le i\le \infty $~, can  be easily transposed to those with coefficients {\bbf{ in bisemisheaves of rings
$\widehat M^{2i}_{T_R}\otimes \widehat M^{2i}_{T_L}$ over $M^{2i}_{T_R}\otimes M^{2i}_{T_L}$\/}}~.
\Ee
\vskip 11pt 

\subsubsection{Definition: Weil algebra of the bilinear nilpotent Lie algebra} 

{\parindent=0pt
Let
\[p_{G\to \pZ[G]}:\quad 
\Reps(\GL_i({F^T_{\o\omega }}\times {F^T_\omega }))\To
\Reps(D_i({F^T_{\o\omega }}\times {F^T_\omega }))\]
denote the projective mapping of the smooth principal bundle 
$\HH^{2i}$
given by the triple 
$(\Reps(\GL_i({F^T_{\o\omega }}\times {F^T_\omega })),\Reps(D_i({F^T_{\o\omega }}\times {F^T_\omega })),p_{G\to \pZ[G]})$
and having as fiber\linebreak 
$\Fs^{(2i)}_{R\times L}= UT_i( {F^T_{\o\omega}} )^t \times UT_i({F^T_\omega})$
the product of the unipotent complete subgroups\linebreak 
$UT_i({F^T_{\o\omega}})^t $
and
$UT_i({F^T_{\omega}}) $~.
\par
Let 
$\Lie(\Fs^{(2i)}_{R\times L})= \Us_i({F^T_{\o\omega}}) \times \Us_i({F^T_\omega})$
be the Lie algebra of the fiber
$\Fs^{(2i)}_{R\times L}$~.
\par
The exterior algebra 
$A(\Lie(\Fs^{(2i)}_{R\times L}))$ on $\Lie(\Fs^{(2i)}_{R\times L})$  
is the algebra of products of differential forms of all degrees. 
$\wedge(\Reps(\GL_i({F^T_{\o\omega }}\times {F^T_\omega })))$
will denote the graded differential algebra of differential forms of 
$\Reps(\GL_i({F^T_{\o\omega }}\times {F^T_\omega }))$
and 
$\wedge(\Reps(D_i({F^T_{\o\omega }}\times {F^T_\omega })))$
will denote the graded differential algebra of differential forms of 
$\Reps(D_i({F^T_{\o\omega }}\times {F^T_\omega }))$~.
\par
It is evident that:
\[ \wedge(\Reps(D_i({F^T_{\o\omega }}\times {F^T_\omega })))
\subset \Reps(\GL_i({F^T_{\o\omega }}\times {F^T_\omega })))\;.\]
Let 
$S(\Lie(\Fs^{(2i)}_{R\times L}))$ 
denote the symmetric algebra on 
$\Lie(\Fs^{(2i)}_{R\times L})$
which can be identified to the algebra of symmetric multilinear forms on 
$\Lie(\Fs^{(2i)}_{R\times L})$~.
\par
Then, the Weil algebra of the Lie algebra 
$\Lie(\Fs^{(2i)}_{R\times L})$
is by definition the graded algebra \cite{Car}:
\[W(\Lie(\Fs^{(2i)}_{R\times L}))= A(\Lie(\Fs^{(2i)}_{R\times L}))\times S(\Lie(\Fs^{(2i)}_{R\times L})) \;.\]
}\vskip 11pt

\subsubsection{Definition} 

{\parindent=0pt
{\em A connection on the fibered space\/}
$\Reps(\GL_i({F^T_{\o\omega }}\times {F^T_\omega }))$
consists in a bilinear mapping 
$f_{R\times L}$ of $A^1(\Lie(\Fs^{(2i)}_{R\times L}))$ in the subspace of the elements of degree one of the algebra $\wedge(\Reps(\GL_i({F^T_{\o\omega }}\times {F^T_\omega })))$~.
\par
$f_{R\times L}$ can be extended as a homomorphism:
\[f'_{R\times L}:\quad A(\Lie(\Fs^{(2i)}_{R\times L}))\To  
\wedge(\Reps(\GL_i({F^T_{\o\omega }}\times {F^T_\omega })))\;.\]
On the other hand, according to \cite{G-H-V} and \cite{Car}, the connection $f_{R\times L}$
 generates a homomorphism from $I_s(\Lie(\Fs^{(2i)}_{R\times L}))$~,
which is the subalgebra of invariant elements of $S(\Lie(\Fs^{(2i)}_{R\times L}))$
identified to the algebra of symmetric multilinear forms $\vee(\Lie(\Fs^{(2i)}_{R\times L}))$
on $\Lie(\Fs^{(2i)}_{R\times L})$~, into 
$H^{2n}(\wedge(\Reps(D_i({F^T_{\o\omega }}\times {F^T_\omega }))))$~: this is the Chern-Weil homomorphism
\[h_{R\times L} :\quad \vee(\Lie(\Fs^{(2i)}_{R\times L}))\To H^{2i}(\wedge(\Reps(D_i({F^T_{\o\omega }}\times {F^T_\omega }))))\]  
were $H^{2i}(\cdot)$ defines a Hodge bilinear $(i, i)$ class also written $i+i$ in section 3.2.
}\vskip 11pt

\subsubsection{Proposition} 

{\em\parindent=0pt

Let $(M^{2i}_{T_R}\otimes M^{2i}_{T_L})$ be the representation space of $\GL_i(F^T_{\o\omega }\times F^T_\omega )$~.\par
 Then, the existence of a biconnection associated with the homomorphism $\Is_s(\Lie(\Fs^{(2i)}\RL))\to \wedge(\Reps(D_i({F^T_{\o\omega }}\times {F^T_\omega })))$
is equivalent to the existence of a nilpotent fiber $\Fs^{(2i)}_{R\times L}$
in the smooth principal bundle $\HH^{2i}(\Reps(\GL_i({F^T_{\o\omega }}\times {F^T_\omega })),\Reps(D_i({F^T_{\o\omega }}\times {F^T_\omega })),p_{G\to \pZ[G]})$
which implies that:
\[H^{2i}(\o S^{P_n}_{K_n},M^{2i}_{T_R }\otimes M^{2i}_{T_L})\simeq
\Reps(\GL_i({F^T_{\o\omega }}\times {F^T_\omega })) \]
only if $\Reps(\GL_i({F^T_{\o\omega }}\times {F^T_\omega }))$~, which is the fibered space of the principal bundle $\HH^{2i}(\cdot,\cdot,\cdot)$~,  is an irreducible representation space of 
$\GL_i({F^T_{\o\omega }}\times {F^T_\omega })$~.
}\vskip 11pt

\bpr  as the connection $f_{R\times L}$
can be extended to the homomorphism $f'_{R\times L}$
given in definition 3.3.8 and generates the Chern-Weil homomorphism 
$h_{R\times L}$~, it appears that the knowledge of  $\Reps(D_i({F^T_{\o\omega }}\times {F^T_\omega }))$~,
together with the connection $f_{R\times L}$
 is sufficient to know the cohomology $H^{2i}(\o S^{P_n}_{K_n},M^{2i}_{T_R }\otimes M^{2i}_{T_L})$ 
\cite{Car}. Thus, the action of a connection on 
$\Reps(\GL_i({F^T_{\o\omega }}\times {F^T_\omega })))$
is equivalent to the action of the nilpotent fiber $\Fs^{(2i)}_{R\times L}$
on $\Reps(D_i({F^T_{\o\omega }}\times {F^T_\omega }))$~.
And, we have that the cohomology
$H^{2i}(\o S^{P_n}_{K_n},M^{2i}_{T_R }\otimes M^{2i}_{T_L})$~, which is the cohomology of the fibered 
space, is in one-to-one correspondence with $\Reps(\GL_i({F^T_{\o\omega }}\times 
{F^T_\omega }))$~.\epr\vskip 11pt

\subsubsection{Proposition}

{\em Let $M^{2i}_{T_R}\otimes M^{2i}_{T_L}$ be an irreducible 
$\GL_{i}({F^{T}_{\o\omega }}\times {F^{T}_{\omega }})$-subbisemimodule of 
$(M_{T_R}\otimes M_{T_L})$~, $i\le n$~.  Then, the $2i$-dimensional irreducible representation
$\Irr W^{(2i)}\FRL(W^{ab}_{F_R}\times W^{ab}_{F_L})$ of the product of the global Weil groups is given by:
\begin{align*}
\Irr \Rep W^{(2i)}\FRL(W^{ab}_{F_R}\times W^{ab}_{F_L})
&= G^{(2i)}({\wt F_{\o\omega_\oplus  }}\times {\wt F_{\omega_\oplus  }})\\
\simeq\Reps (\GL_{i}({F^{T}_{\o\omega_\oplus  }}\times {F^{T}_{\omega_\oplus  }})
&= H^{2i}(\o S^{P_n}_{K_n},M^{2i}_{T_{R_\oplus }}\otimes M^{2i}_{T_{L_\oplus }})\;.\end{align*}
}\vskip 11pt 

\bpr indeed, according to proposition 3.1.13, we have that 
\[ H^{2i}(\o S^{P_n}_{K_n},M^{2i}_{T_{R_\oplus }}\otimes M^{2i}_{T_{L_\oplus }})
= \Reps( \GL_{i}({F^{T}_{\o\omega_\oplus  }}\times {F^{T}_{\omega _\oplus }})\]
and, according to proposition 2.4.7, we have that:
\[\Irr W^{(2i)}\FRL(W^{ab}_{F_R}\times W^{ab}_{F_L})
=G^{(2i)}({\wt F_{\o\omega_\oplus  }}\times {\wt F_{\omega_\oplus  }})\;.\]

So, we are led to the isomorphism:
\[ \Irr W^{(2i )}\FRL(W^{ab}_{F_R}\times W^{ab}_{R_L})
\simeq H^{2i}(\o S^{P_n}_{K_n},M^{2i}_{T_{R_\oplus }}\otimes M^{2i}_{T_{L_\oplus }})\]
because $ \Irr W^{(2i )}\FRL(W^{ab}_{F_R}\times W^{ab}_{F_L})$ refers to not necessarily compact algebraic extensions while $H^{2i}(\o S^{P_n}_{K_n},M^{2i}_{T_R}\otimes M^{2i}_{T_L})$ is the cohomology of the toroidal compactified bisemispace $\o S^{P_n}_{K_n}$~.\epr
\vskip 11pt 

\subsection[Langlands global correspondence on the irreducible algebraic bilinear semigroup 
$\GL_n({\wt F_{\o\omega }}\times {\wt F_{\omega }})$]{{\boldmath Langlands global correspondence on the irreducible algebraic bilinear semigroup 
$\GL_n({\wt F_{\o\omega }}\times {\wt F_{\omega }})$}}

\subsubsection{Proposition}

{\em The cohomology $H^{2i}(\o S^{P_n}_{K_n},M^{2i}_{T_{R_\oplus }}\otimes M^{2i}_{T_{L_\oplus }})$ has a decomposition according to the equivalent representatives $g^{(2i)}_{T\RL}[j,m_j]$ of the conjugacy classes of the complete bilinear semigroup\linebreak $G^{(2i)}({F^{T}_{\o\omega }}\times {F^{T}_{\omega }})$ according to:
\[H^{2i}(\o S^{P_n}_{K_n},M^{2i}_{T_{R_\oplus }}\otimes M^{2i}_{T_{L_\oplus }})
= \bigoplus^r_{j=1} \bigoplus_{m_j}g^{(2i)}_{T\RL}[j,m_j]\;.\]
}
\vskip 11pt 

\bpr 
\[ H^{2i}(\o S^{P_n}_{K_n},M^{2i}_{T_{R_\oplus }}\otimes M^{2i}_{T_{L_\oplus }})=  G^{(2i)}({F^{T}_{\o\omega_\oplus  }}\times {F^{T}_{\omega _\oplus }})\]
where $G^{(2i)}({F^{T}_{\o\omega _\oplus }}\times {F^{T}_{\omega_\oplus  }})$ decomposes according to the equivalent representatives $g^{(2i)}_{T\RL}[j,m_j]$ of $G^{(2i)}(F^T_{\o\omega }\times F^T_{\omega })$~.  So, we get the thesis.\epr
\vskip 11pt 

\subsubsection{Corollary}

{\em Let 
\[ \o S^{Z[P_n]}_{Z[K_n]}
= D_n({F^{T}_{\o\omega^1 }}\times {F^{T}_{\omega^1 }})\setminus
D_n({F^{T}_{R}}\times {F^{T}_{L}})\Big/ D_n((\ZZ\big/N\ \ZZ)^2)\]
be the diagonal compactified bisemispace and let $D^{(2i)}({F^{T}_{\o\omega_\oplus  }}\times {F^{T}_{\omega_\oplus  }})$ be the diagonal bilinear  semigroup decomposing into conjugacy classes having representatives $g^{(2i)}_{T\RL}[j]$~.  Then, we have the following decomposition of the cohomology:
\[ H^{2i}(\o S^{Z[P_n]}_{Z[K_n]},M^{2i}_{T_{R_\oplus }}(D)\otimes M^{2i}_{T_{L_\oplus }}(D))
=\simeq \bigoplus_j g^{(2i)}_{T\RL}[j]\]
where $(M^{2i}_{T_{R_\oplus }}(D)\otimes M^{2i}_{T_{L_\oplus }}(D))=\bigoplus_j(M^{2i}_{T_{\o\omega _j}}\otimes M^{2i}_{T_{\omega _j}})$ is a $D_i({F^{T}_{\o\omega_\oplus  }}\times {F^{T}_{\omega_\oplus  }})$-bisemimodule.
}
\vskip 11pt 

\bpr referring to chapter 1, we note that the conjugacy classes of $D^{(2i)}({F^{T}_{\o\omega }}\times {F^{T}_{\omega }})$ have unique representatives $g^{(2i)}_{T\RL}[j]$~, since there is no nilpotent action in the diagonal case, the nilpotent action being generated by $UT^t_i({F^{T}_{\o\omega }})\times UT_i({F^{T}_{\omega }})\subset 
\GL_i({F^{T}_{\o\omega }}\times {F^{T}_{\omega }})$~.\epr
\vskip 11pt 

\subsubsection[Analytic development of the bilinear cohomology $H^{2i}(\o S^{P_n}_{K_n},M^{2i}_{T_R }\otimes M^{2i}_{T_L})$]{{\boldmath Analytic development of the bilinear cohomology $H^{2i}(\o S^{P_n}_{K_n},M^{2i}_{T_R }\otimes M^{2i}_{T_L})$}}

The objective consists now in proving that the decomposition of the bilinear cohomology
$H^{2i}(\o S^{P_n}_{K_n},M^{2i}_{T_R }\otimes M^{2i}_{T_L})$
corresponding to the representation space of the irreducible bilinear general semigroup 
$\GL_i({F^T_{\o\omega }}\times {F^T_\omega })$
has an analytic development given by the development
in truncated Fourier series  of the product of a ``right" 
$2i$-dimensional cusp form by its left equivalent.\vskip 11pt

As the decomposition of the cohomology $H^{2i}(\o S^{P_n}_{K_n},M^{2i}_{T_R }\otimes M^{2i}_{T_L})$
is given by the sum over $j$ and $m_j$  of equivalent representatives $ g^{(2i)}_{T\RL}[j,m_j]=
g^{(2i)}_{T_R}[j,m_j]\times g^{(2i)}_{T_L}[j,m_j]$
consisting in products, right by left, of
$i$-dimensional complex semitori $T^{2i}_R[j,m_j]\times T^{2i}_L[j,m_j]$ according to proposition 3.3.3, an analytic development must be given to these semitori as follows:\vskip 11pt

\subsubsection{Proposition}  

{\em Let \quad $\vec z=\sum\limits^{2i}_{d =1}z_d \vec e_d $\quad be a vector of $(F_{\omega _j})^i$ and, more precisely, a point of $G^{(2i)}(F_{\omega_j})$~. 
\par
Then, every left (resp. right)
$i$-dimensional complex semitorus has the analytic development:
\begin{align*}
T^{2i}_L[j,m_j] &\simeq \lambda ^{\half}(2i,j,m_j) e^{2\pi ijz}\\
(\text{resp.} \quad 
T^{2i}_R[j,m_j] &\simeq \lambda ^{\half}(2i,j,m_j)e^{-2\pi ijz}\ )\end{align*}
where \quad \bt[t]{p{14cm}}
\textbullet\  $\lambda ^{\half}(2i,j,m_j)\simeq j^{2i}\cdot N^{2i}\cdot (m^{(j)})^{2i}$\quad 
can be considered as a ``global'' Hecke  character;\\[-6pt]
\textbullet\  $z=\sum\limits_d z_d \ |\vec e_d |$~.\te
}\vskip 11pt

\bpr
\Be
\item The cohomology $H^{2i}(\o S^{P_n}_{K_n},M^{2i}_{T_R }\otimes M^{2i}_{T_L})$
corresponds to an endomorphism of $\o S^{P_n}_{K_n}$
 into itself. Its decomposition into classes ``~$j$~'' with representatives ``~$m_j$~'' refers to the cosets of $\GL_i({F^T_{R}}\times {F^T_L})\big/\GL_i( (\ZZ\big/N\ \ZZ)^2)$~.
So, the scalar $\lambda (2i, j,mj)$ 
will correspond to the eigenvalues of $g_i(\Os_{F^T_{\o\omega _{j,m_j}}}\times \Os_{F^+_{\omega _{j,m_j}}}) \in \GL_i( (\ZZ\big/N\ \ZZ)^2 )$
which is the $j$-th coset representative $U{j_R}\times U_{j_L}$
of the product $T_R(i;r)\otimes T_L(i;r)$
of the Hecke operators according to proposition 2.2.5.
\par
More precisely, let $\{\lambda_d(2i,j,m_j)\}^{2i}_{d=1}$
be the set of eigenvalues of $ U_{j_R}\times U_{j_L} $ and let $\lambda(2i,j,m_j)=\prod\limits_{d=1}^{2i}\lambda_d(2i,j,m_j)$ 
be the product of these eigenvalues. 
\par
According to definition 2.2.4, we have that
\[ \lambda(2i,j,m_j)=\txt\prod\limits_{d=1}^{2i}\lambda_d(2i,j,m_j)
= \det(\alpha _{i^2;j^2}\times D_{j^2;m_{j^2}})_{ss}\simeq j^{2i}\cdot N^{2i} \cdot (m^{(j)})^{2i}
\]
where $D_{j^2;m_{j^2}}$
is the decomposition group of the $j$-th bisublattice with representative $m_j$ and where 
$\alpha _{i^2;j^2}$
is the $j$-th split Cartan subgroup.
\par
So, it appears that the square root $\lambda ^{\half}(2i,j,m_j)$ 
of $\lambda (2i,j,m_j)$ can be considered as  a global Hecke character having an inflation action on $e^{2\pi ijz}$~.

\item  On the other hand, as we are concerned with the $j$-th global or infinite complex place 
$\omega _j$ or $\o\omega _j$~, we have to consider the global Frobenius substitution for the left (resp. right) place $\omega _j$ (resp. $\o\omega _j$~) 
as given by the mapping:
\begin{align*}
e^{2\pi iz} & \To e^{2\pi ijz}\; ,\\
(\text{resp.} \quad 
e^{-2\pi iz} & \To e^{-2\pi ijz}\; )\tag*{$z\in G^{(2i)}(F_{\omega _j})\equiv g^{(2i)}_{L}[j,m_j]$~,}\end{align*}
$e^{2\pi iz}$ being a ``circular'' function on the $j^{2i}_\delta $ completions of $g^{(2i)}_L[j,m_j]$ according to proposition 3.1.3.

So, the left (resp. right) semitorus corresponding to the equivalent representative $g^{(2i)}_{T_L}[j,m_j]$ (resp. $g^{(2i)}_{T_R}[j,m_j]$~)
has the following analytic development:
\begin{align*}
T^{2i}_L[j,m_j] &\simeq \lambda^{\half}(2i,j,m_j)e^{2\pi ijz} \\
(\text{resp.} \quad 
T^{2i}_R[j,m_j] &\simeq \lambda^{\half}(2i,j,m_j)e^{-2\pi ijz} \ )
\end{align*}
taking into account that the $n$-dimensional complex \lr semitorus $T^{2i}_L[j,m_j]$ (resp. $T^{2i}_R[j,m_j]$~) is diffeomorphic to the $i$-fold product
\begin{align*}
T^{2i}_L[j,m_j]
&= T^2_L[j,m_j]\times \cdots \times T^2_L[j,m_j]\\
\text{(resp.\quad}
T^{2i}_R[j,m_j]
&= T^2_R[j,m_j]\times \cdots \times T^2_R[j,m_j]\ )\end{align*}
of $1$-dimensional complex \lr semitori $T^2_L[j,m_j]$ (resp. $T^2_R[j,m_j]$~).  
These semitori have the analytic development:
\begin{align*}
T^2_L[j,m_j]
&= \lambda _{T^2_L}\ e^{2\pi ijz_d}\\
&= S^1_{d_1}[j,m_j] \times S^1_{d_2}[j,m_j]  \tag*{(~$ z_d\in\cit\;, $~)} \\
&= r_{S^1_{d_1}}\ e^{2\pi ijx_{d_1}} \times  r_{S^1_{d_2}}\ e^{2\pi ij(iy_{d_2})} 
\tag*{(~$x_{d_1}\in\rit\ , \; y_{d_2}\in\rit$~)}\\
&\simeq \lambda _{d_1}(1,j,m_j)\cdot\lambda _{d_2}(1,j,m_j)\ e^{2\pi ijx_{d_1}}\times e^{2\pi ij(iy_{d_2})}
\end{align*}
where
\Bi
\item $\lambda _{d_1}(1,j,m_j)\cdot \lambda _{d_2}(1,j,m_j)\simeq r_{S^1_{d_1}}\cdot r_{S^1_{d_2}}$ is a product of two eigenvalues of $U_{j_R}\times U_{j_L}$~,
\item $r_{S^1_{d_1}}$ and $r_{S^1_{d_2}}$ are radii.
\Ei
Indeed, the left $1D$-complex semitorus $T^2_L[j,m_j]$ is diffeomorphic to the product of two circles \[S^1_{d_1}[j,m_j]=r_{S^1_{d_1}}\ e^{2\pi ijx_{d_1}}\] and 
\[S^1_{d_2}[j,m_j]=r_{S^1_{d_2}}\ e^{2\pi ij(iy_{d_2})}\] localized in perpendicular planes with $\cos (2\pi ijy_{d_2})$ and $\sin (2\pi (ijy_{d_2})$ of $e^{2\pi i(ijy_{d_2})}$ defined over $i \rit$~.  

This is justified by the fact that a rotation of $90^{0}$ of the circle
$S^1_{d_2}[j,m_j]=\linebreak r_{S^1_{d_2}}\ e^{2\pi ij(iy_{d_2})}$ over $i\rit$ transforms it into the circle 
$S^1_{d_{2_\perp}}[j,m_j]=r_{S^1_{d_2}}\ e^{2\pi ij(y_{d_2})}$ over $\rit$~, 
localized in the same plane as the circle $S^1_{d_1}[j,m_j]$~, according to:
\begin{align*}
{\rm rot}(90^0): \quad S^1_{d_2}[j,m_j] &\To S^1_{d_{2_\perp}}[j,m_j]\\
r_{S^1_{d_2}}\ e^{2\pi ij(iy_{d_2})}\big/i\rit 
&\To 
r_{S^1_{d_2}}\ e^{2\pi ij(y_{d_2})}\big/\rit \end{align*}
where:
\Bi
\item $S^1_{d_2}[j,m_j]=r_{S^1_{d_2}}[\cos(2\pi ijy_{d_2}) + i\sin (2\pi ijy_{d_2})]$~,
\item $S^1_{d_{2_\perp}}[j,m_j]=r_{S^1_{d_2}}[\cos(2\pi jy_{d_2}) + i\sin (2\pi jy_{d_2})]$~,
\Ei

Finally, recall that 
\[ e^{2\pi ijz}=e^{2\pi ijz_1}\times \cdots\times e^{2\pi ijz_d}\times \cdots\times e^{2\pi ijz_{2n}}\]
is a Laurent monomial in the variables $z_1,\cdots,z_{2n}$~.\epr
\Ee
\vskip 11pt

\subsubsection{Proposition} 

{\em  The cohomology $H^{2i}(\o S^{P_n}_{K_n},\widehat M^{2i}_{T_{R_\oplus } }\otimes \widehat M^{2i}_{T_{L_\oplus }})$
has the following analytic
development:
\[H^{2i}(\o S^{P_n}_{K_n},\widehat M^{2i}_{T_{R_\oplus } }\otimes \widehat M^{2i}_{T_{L_\oplus }})= 
\txt\sum\limits_{j=1}^r\sum\limits_{m_j}(\lambda^{\half}(2i,j,m_j)e^{-2\pi ijz})
\otimes_{(D)} (\lambda^{\half}(2i,j,m_j)e^{2\pi ijz})\] 
where
\begin{align*}
\Eis_L(2i,j,m_j) &= \txt\sum\limits_{j=1}^r \sum\limits_{m_j} \lambda^{\half}(2i,j,m_j)e^{2\pi ijz} \\
(\text{resp.} \quad 
\Eis_R(2i,j,m_j) &= \txt\sum\limits_{j=1}^r \sum\limits_{m_j} \lambda^{\half}(2i,j,m_j)e^{-2\pi ijz} \ )\end{align*}
is the (truncated) Fourier development of the equivalent of a normalized $2i$-dimensional left (resp. right) Eisenstein series of weight $k = 2$ restricted to the upper (resp. lower) half space.}\vskip 11pt
\bpr
\Be\item  First we recall that the Fourier development of a normalized Eisenstein series is:
\[G_k(y)\simeq \txt\sum\limits_{n=1}^\infty \sigma _{k-1}(n) \ q^n\]
 with \quad $y\in\CC$~, \quad $q=e^{2\pi iy}$~, \quad $\sigma _{k-1}(n) =\sum\limits_{d\mid n}d^{k-1}$ \quad
and \quad $k $ the weight.
\item Referring to the development of the cohomology $H^{2i}(\o S^{P_n}_{K_n},M^{2i}_{T_R }\otimes M^{2i}_{T_L})$
according to the
sum over the products of equivalent representatives 
$g^{(2i)}_{T\RL}[j,m_j]=g^{(2i)}_{T_R}[j,m_j] \times 
g^{(2i)}_{T_L}[j,m_j] = T^{2i}_R[j,m_j]\times T^{2i}_L[j,m_j]$ of the complete bilinear semigroup $G^{(2i)}({F^T_{\o\omega }}\times {F^+_\omega })$ such
that an $i$-dimensional left
(resp. right) complex semitorus has the analytic development:
\begin{align*}
T^{2i}_L[j,m_j] &\simeq \lambda^{\half}(2i,j,m_j)e^{2\pi ijz} \\
(\text{resp.} \quad 
T^{2i}_R[j,m_j] &\simeq \lambda^{\half}(2i,j,m_j)e^{-2\pi ijz} \ ),\end{align*}
we get the searched analytic decomposition of the cohomology in terms of the product $\Eis_R(2i,j,m_j)\times \Eis_L(2i,j,m_j)$ of the (truncated) Fourier development of the equivalent of a normalized $2i$-dimensional right Eisenstein series of weight two
by its left equivalent. 
\par
Remark that
 the sum \quad $\sigma _{k-1}(n) =\sum\limits_{d\mid n}d^{k-1}$ \quad
in $G_k(y)$
corresponds in the $2i$-dimen\-sio\-nal envisaged
case to the sum $\sum\limits_{m_j}$ over the equivalent representatives ${m_j}$ of the classes $j$~.
\epr
\Ee\vskip 11pt

\subsubsection{Proposition} 

{\em
 The analytic development $\Eis_R(2i,j,m_j)\otimes \Eis_L(2i,j,m_j)$
of the bilinear cohomology\linebreak $H^{2i}(\o S^{P_n}_{K_n},\widehat M^{2i}_{T_R }\otimes \widehat M^{2i}_{T_L})$
is an eigenbifunction of the product of Hecke operators $T_R(i;r)\otimes T_L(i;r)$
whose $j$-th bicoset representative is
\[ U_{j_R}\times U_{j_L}=t^t_i(\Os_{F_{\o\omega_{j,m_j}}})\times t_i(\Os_{F_{\omega_{j,m_j}}})\;.\]
}\vskip 11pt

\bpr {\parindent=0pt
 the cohomology $H^{2i}(\o S^{P_n}_{K_n},M^{2i}_{T_R }\otimes M^{2i}_{T_L})$
having coefficients in the 
$\GL_i({F^T_{\o\omega}}\times {F^T_\omega})$-bisemimodule $M^{2i}_{T_R }\otimes M^{2i}_{T_L}$
is in one-to-one correspondence with the functional representation space $\FReps(\GL_i({F^T_{\o\omega }}\times {F^T_\omega }))$
of the general bilinear semigroup $\GL_i({F^T_{\o\omega }}\times {F^T_\omega })$~.  Now, 
$\GL_i({F^T_{\o\omega }}\times {F^T_\omega })$ is a solvable bilinear semigroup such that we have a chain of embedded subsemigroups:
\begin{align*}
\GL_i(F^T_{\o\omega_1}\times F^T_{\omega _1}) &\subset
\GL_i(F^T_{\o\omega_{1+2}}\times F^T_{\omega _{1+2}})\subset \cdots\\
&\quad \subset \GL_i(F^T_{\o\omega_{1+\cdots+j}}\times F^T_{\omega _{1+\cdots+j}})\\
&\qquad \subset \GL_i(F^T_{\o\omega_{1+\cdots+j+\cdots+r}}\times F^T_{\omega _{1+\cdots+j+\cdots+r}})\end{align*}
corresponding to the sums of the bicosets of the quotient semigroup
\[\o X_{S_{R\times L}}=\GL_i({F^T_{R}}\times {F^T_L})\Big/\GL_i((\ZZ\big/N\ \ZZ)^2)=G^{(2i)}(F^T_{\o\omega }\times F^T_\omega )\;.\]
The conjugacy class representative $g^{(2i)}_{T\RL}[j,m_j]$ of $G^{(2i)}({F^T_{\o\omega}}\times {F^T_{\omega }})$
is obtained throughout the action of the $j$-th bicoset representative $U_{j_R}\times U_{j_L}$
of the product of Hecke operators $T_R(i;r)\otimes T_L(i;r)$
according to proposition 2.2.5 since the Hecke bialgebra $\Hs_{R\times L}(i)$~,
generated by all the pseudo-ramified Hecke bioperators, has a representation in 
$\GL_i((\ZZ\big/N\ \ZZ)^2)$~.
\par
According to propositions 3.4.1 and 3.4.5, 3.3.10 and 2.1.4, we have that:
\[\FReps(\GL_i({F^T_{\o\omega }}\times {F^T_\omega })) 
= G^{(2i)}(F^T_{\o\omega _\oplus }\times F^T_{\omega _\oplus })
=
\Eis_R(2i,j,m_j)\otimes \Eis_L(2i,j,m_j)\]
for the upper value of $j=r\le \infty$~.
\par
So, considering the quotient semigroup $\o X_{S_{R\times L}}=\GL_n({F^T_{R}}\times {F^T_L})\big/\GL_n( (\ZZ\big/N\ \ZZ)^2)$
is equivalent to take into account the action of the Hecke bioperator $T_R(n;r)\otimes T_L(n;r)$
on $\Reps(\GL_n({F^T_{R}}\times {F^T_L})) $~. It then results that 
$\Eis_R(2n,j,m_j)\otimes \Eis_L(2n,j,m_j)$
 is an eigenbifunction of $T_R(n;r)\otimes T_L(n;r)$~.
\par
On the other hand, it is well-known that the Eisenstein series are eigenfunctions of the Hecke operators. 
\par
But, for more details on the algebraic spectral representations, see \cite{Pie2} and \cite{Clo}.\epr}\vskip 11pt

\subsubsection{Corollary} 

{\em
 The truncated  (bi)series $\Eis_{R\times L}(2i,j,m_j)=\Eis_R(2i,j,m_j)\otimes \Eis_L(2i,j,m_j)$
is:
\Be
\item a truncated supercuspidal (bi)form over $\CC^{i}$~;

\item a solvable power (bi)series.\Ee
}\vskip 11pt

\bpr  {\parindent=0pt
\Be
\item As $\Eis_L(2i,j,m_j)$ (resp. $\Eis_R(2i,j,m_j)$~) constitutes the $i$-dimensional equivalent representative of the one-dimensional complex truncated Fourier development of a normalized Eisenstein series of weight two, which is a (quasi-) modular form, 
$\Eis_L(2i,j,m_j)$ (resp. $\Eis_R(2i,j,m_j)$~) is a truncated (super)cuspidal form over $\CC^i$~.
\item
As the general bilinear semigroup $\GL_i({F^T_{\o\omega }}\times {F^T_\omega })$
 is solvable and as we have that:
\[\FReps(\GL_i({F^T_{\o\omega_\oplus  }}\times 
{F^T_{\omega_\oplus } })) 
=\Eis_{R\times L}(2i,j,m_j)\;, \quad 
F^T_{\omega _\oplus }=\bigoplus_j \bigoplus_{m_j}F^T_{\omega _{j,m_j}}\;,\]
according to propositions 3.4.6 and 2.1.4, the ``Eisenstein'' biseries $\Eis_{R\times L}(2i,j,m_j)$
will be said solvable.
\par
We thus have a chain of embedded (truncated)  biseries:

\begin{multline*} \Eis_{R\times L}(2i,j^{\text{up}}=1,m_j)\subset\cdots \subset
\Eis_{R\times L}(2i,j^{\text{up}}=j,m_j)\\ \subset\cdots\subset
\Eis_{R\times L}(2i,j^{\text{up}}=r,m_j)
\end{multline*}
where $j^{\text{up}}$
is the upper value of the integer ``~$j$~'' labeling the conjugacy classes. 
\par
This implies the injective mapping:
\[E_{j^{\text{up}}}:
\quad \Eis_{R\times L}(2i,j^{\text{up}}=j,m_j)\To
\Eis_{R\times L}(2i,j^{\text{up}}=j+1,m_j)\]
between two (truncated)  biseries such that $\Eis_{R\times L}(2i,j^{\text{up}}=j,m_j)
$ is an eigenbifunction of the Hecke bioperator $T_R(i;j)\otimes T_L(i;j)$ while 
$\Eis_{R\times L}(2i,j^{\text{up}}=j+1,m_j)$ is an eigenbifunction of $T_R(i;j+1)\otimes T_L(i;j+1)$~. \epr
\Ee}\vskip 11pt


\subsubsection{Proposition} 

{\em The cohomology $H^{2i}(\o S^{\pZ[P_n]}_{\pZ[K_n]}, \widehat M^{2i}_{T_{R_\oplus } }(D)\otimes \widehat M^{2i}_{T_{L_\oplus }} (D))$
of the compactified diagonal bisemispace $\o S^{\pZ[P_n]}_{\pZ[K_n]}$ has the analytic development: 
\[ H^{2i}(\o S^{\pZ[P_n]}_{\pZ[K_n]},\widehat M^{2i}_{T_{R_\oplus } }(D)\otimes \widehat M^{2i}_{T_{L_\oplus }}(D))
= \sum^r_{j=1} (\lambda ^{\half}(2i,j,1)\ e^{-2\pi ijz})
\otimes (\lambda ^{\half}(2i,j,1)\ e^{2\pi ijz})
\]
where
\begin{align*}\Eis_L(2i,j,0)
&=\txt\sum\limits_{j=1}^r \lambda^{\half}(2i,j,0)\ e^{2\pi ijz}\\
(\text{resp.}\quad 
\Eis_R(2i,j,0)
&=\txt\sum\limits_{j=1}^r \lambda^{\half}(2i,j,0)\ e^{-2\pi ijz}\ )\end{align*}
is  a solvable $2i$-dimensional left (resp. right) cusp form of weight $k = 2 $~.
}\vskip 11pt

\bpr this is a particular case of proposition 3.4.5 where the number of equivalent representatives of the left and right conjugacy classes is restricted to one.\epr
\vskip 11pt

%

\subsubsection[Bialgebra of cusp forms $L^{1-1}_{\rm cusp}(
G^{(2n)}({F^T_{\o\omega}} \times{F^T_{\omega}}))$]{{\boldmath Bialgebra of cusp forms $L^{1-1}_{\rm cusp}(
G^{(2n)}({F^T_{\o\omega}} \times{F^T_{\omega}}))$}}

Let $L^{1-1}\RL(
G^{(2n)}({F_{\o\omega}} \times{F_{\omega}}))$ be the bialgebra of complex-valued smooth continuous bifunctions
$ \phi ^{(2n)}_{G_{j_R}}(x_{g_{j_R}})\otimes \phi ^{(2n)}_{G_{j_L}}(x_{g_{j_L}})$ on the pseudo-ramified bilinear   semigroup   $G^{(2n)}({F_{\o\omega}} \times{F_{\omega}})$ as developed in section 2.5.

Let 
\[ \gamma ^C\RL : \quad G^{(2n)}({ \widetilde F_{\o\omega}} \times{\widetilde F_{\omega}})
\To G^{(2n)}({F^T_{\o\omega}} \times{F^T_{\omega}})\]
be the projective toroidal isomorphism of compactification of 
$G^{(2n)}({\widetilde F_{\o\omega}} \times{\widetilde F_{\omega}})$ mapping each conjugacy class representative
$\tilde g^{(2n)}\RL[j,m_j]$ of $G^{(2n)}({\widetilde F_{\o\omega}} \times{\widetilde F_{\omega}})$ into its toroidal equivalent
$g^{(2n)}_{T\RL}[j,m_j]$ according to proposition 3.3.3.

It is evident that $\gamma ^C\RL$ generates
$L^{1-1}_{\rm cusp}(
G^{(2n)}({F^T_{\o\omega}} \times{F^T_{\omega}}))$
in such a way that the bialgebra of smooth continuous bifunctions
$ \phi ^{(2n)}_{G^T_{j_R}}(x_{g^T_{j_R}})\otimes \phi ^{(2n)}_{G^T_{j_L}}(x_{g^T_{j_L}})$ on the pseudo-ramified toroidal bilinear  semigroup  
$G^{(2n)}({F^T_{\o\omega}} \times{F^T_{\omega}})$ is the bialgebra of cusp forms
$L^{1-1}_{\rm cusp}(
G^{(2n)}({F^T_{\o\omega}} \times{F^T_{\omega}}))$~.

On the other hand, the projective toroidal isomorphism 
maps each bifunction
$\phi ^{(2n)}_{G_{j_R}}(x_{g_{j_R}})\otimes \phi ^{(2n)}_{G_{j_L}}(x_{g_{j_L}})
\in L^{1-1}\RL(
G^{(2n)}({F_{\o\omega}} \times{F_{\omega}}))
$ on the conjugacy class representative $g^{(2n)}\RL[j,m_j]\in
G^{(2n)}({F_{\o\omega}} \times{F_{\omega}})$ into its toroidal equivalent 
$\phi ^{(2n)}_{G^T_{j_R}}(x_{g^T_{j_R}})\otimes \phi ^{(2n)}_{G^T_{j_L}}(x_{g^T_{j_L}})
\in L^{1-1}_{\rm cusp}(
G^{(2n)}({F^T_{\o\omega}} \times{F^T_{\omega}}))$ on the conjugacy class representative
$g^{(2n)}_{T\RL}[j,m_j]\in G^{(2n)}({F^T_{\o\omega}} \times{F^T_{\omega}})$ in such a way that:
\Be
\item \quad $\begin{array}[t]{ll}
\phi ^{(2n)}_{G^T_{j_R}}(x_{g^T_{j_R}})\otimes_D \phi ^{(2n)}_{G^T_{j_L}}(x_{g^T_{j_L}})
&= T^{2n}_R[j,m_j]\otimes_D T^{2n}_L[j,m_j]\\
&= \lambda ^{\half}(2n,j,m_j)\ e^{-2\pi ijz}\otimes_D \lambda ^{\half}(2n,j,m_j)\ e^{2\pi ijz}\end{array}$

according to proposition 3.4.4.

\item \quad $\begin{array}[t]{ll}
\phi ^{(2n)}_{G^T_{R}}(x_{g^T_{R}})\otimes_D \phi ^{(2n)}_{G^T_{L}}(x_{g^T_{L}})
&= \bigoplus_{j=1}^r\bigoplus_{m_j}
(\phi ^{(2n)}_{G^T_{j_R}}(x_{g^T_{j_R}})\otimes_D \phi ^{(2n)}_{G^T_{j_L}}(x_{g^T_{j_L}}))\\
&= \bigoplus_{j=1}^r\bigoplus_{m_j}(\lambda ^{\half}(2n,j,m_j)\ e^{-2\pi ijz}\otimes_D \lambda ^{\half}(2n,j,m_j)
\ e^{2\pi ijz}\\
&= \Eis_R(2n,j,m_j)\otimes_{(D)} \Eis_L(2n,j,m_j)\end{array}$

is the (truncated) supercuspidal biform on the pseudo-ramified toroidal bilinear  semigroup 
$G^{(2n)}({F^T_{\o\omega}} \times {F^T_{\omega}})$ given by the product, right by left, of the equivalent of $2n$-dimensional Eisenstein series.
\Ee
\vskip 11pt 

\subsubsection{Proposition}

{\em Let $R_{G^{(2n)}({F^T_{\o\omega}} \times{F^T_{\omega}})}
(\phi ^{(2n)}_{G^T_{R}}  \otimes \phi ^{(2n)}_{G^T_{L}} )$ denote  the integral operator on cusp biforms
$(\phi ^{(2n)}_{G^T_{R}}   \otimes \phi ^{(2n)}_{G^T_{L}}) $ over the pseudo-ramified toroidal bilinear  semigroup  $G^{(2n)}({F^T_{\o\omega}} \times{F^T_{\omega}})$~.

Then, the trace formula of this integral operator is given by:
\[
\Tr(R_{G^{(2n)}({F^T_{\o\omega}} \times{F^T_{\omega}})}
(\phi ^{(2n)}_{G^T_{R}}  \otimes \phi ^{(2n)}_{G^T_{L}} )
= {\rm vol}(G^{(2n)}({F^T_{\o\omega}} \times{F^T_{\omega}}))
\sum_j\sum_{m_j}\lambda (2n,j,m_j)\]
where $\lambda (2n,j,m_j)=\prod\limits_{d=1}^{2n}\lambda _d(2n,j,m_j)$ is the product of the eigenvalues
$\lambda _d(n,j,m_j)$ of the $j$-th coset representative of the product, right by left, of Hecke operators.
}\vskip 11pt 

\bpr Indeed, according to section 2.5.18, the trace of this integral operator can be developed as follows:
\begin{align*}
&\Tr(R_{G^{(2n)}({F^T_{\o\omega}} \times{F^T_{\omega}})}
(\phi ^{(2n)}_{G^T_{R}}  \otimes \phi ^{(2n)}_{G^T_{L}} )\\[11pt]
& \quad = \bigoplus_{j,m_j}
{\rm vol}(G^{(2n)}({F^T_{\o\omega}} \times{F^T_{\omega}}))
\int_{G^{(2n)}} (\phi ^{(2n)}_{G^T_{j_R}}(x_{g^T_{j_R}})\otimes_D \phi ^{(2n)}_{G^T_{j_L}}(x_{g^T_{j_L}})
\ dx_{g^T_{j_R}}\ dx_{g^T_{j_L}}\\[11pt]
& \quad = \bigoplus_{j,m_j}
{\rm vol}(G^{(2n)}({F^T_{\o\omega}} \times{F^T_{\omega}}))
\int_{G^{(2n)}} \lambda ^{\half}(2n,j,m_j)\ e^{-2\pi ijz}\otimes_D  \lambda ^{\half}(2n,j,m_j)\ e^{2\pi ijz}
\ dx_{g^T_{j_R}}\ dx_{g^T_{j_L}}\\[11pt]
& \quad \simeq {\rm vol}(G^{(2n)}({F^T_{\o\omega}} \times{F^T_{\omega}}))
(\Eis_R(2n,j,m_j),\Eis_L(2n,j,m_j))\\[11pt]
& \quad = {\rm vol}(G^{(2n)}({F^T_{\o\omega}} \times{F^T_{\omega}}))
\sum_{j,m_j}\lambda (2n,j,m_j)\end{align*}
where $(\Eis_R(2n,j,m_j),\Eis_L(2n,j,m_j))$ is a bilinear form from
$\Irr\cusp (\GL_n({F^T_{\o\omega}} \times{F^T_{\omega}}))
= \Eis\RL(2n,j,m_j)$ to $\cit$~.\epr
\vskip 11pt 

\subsubsection{Proposition}

{\em Let $L^{1-1}\RL(
P^{(2n)}({F_{\o\omega^1}} \times{F_{\omega^1}}))$ be the bialgebra of complex-valued bifunctions on the bilinear parabolic semigroup 
  $P^{(2n)}({F_{\o\omega^1}} \times{F_{\omega^1}})$ and let
$L^{1-1}\RL(
P^{(2n)}({F^{T}_{\o\omega^1}} \times{F^{T}_{\omega^1}}))$ be its toroidal equivalent obtained by the projective 
toroidal isomorphism
$ \gamma ^C\RL$~.  

Let $R_{P^{(2n)}({F^{T}_{\o\omega^1}} \times{F^{T}_{\omega^1}})} 
(\phi ^{(2n)}_{P^T_{R}}\otimes \phi ^{(2n)}_{P^T_{L}})$ denote  the integral operator on 
bifunctions $\in\linebreak L^{1-1}\RL(
P^{(2n)}({F^{T}_{\o\omega^1}} \times{F^{T}_{\omega^1}}))$~.

Then, the trace formula of this integral operator is given by:
\begin{align*}
& \Tr( R_{P^{(2n)}({F^{T}_{\o\omega^1}} \times{F^{T}_{\omega^1}})} 
(\phi ^{(2n)}_{P^T_{R}}\otimes \phi ^{(2n)}_{P^T_{L}}))\\[11pt]
&\quad = {\rm vol}(P^{(2n)}({F^{T}_{\o\omega^1}} \times{F^{T}_{\omega^1}})
\sum_j\sum_{m_j} \int \chi ^T_{j,m_{j\RL}}\ dx_{p_{j_R}}\ dx_{{j_L}}\\[11pt]
&\quad = {\rm vol}(P^{(2n)}({F^{T}_{\o\omega^1}} \times{F^{T}_{\omega^1}}))
\sum_j\sum_{m_j}\lambda _I(2n,j,m_j)\end{align*}
where
\Bi
\item $\chi ^T_{j,m_{j\RL}}:x_{p_{j_R}}\times x_{p_{j_L}}\to
\phi ^{(2n)}_{P^T_{j_R}} (x_{p_{j_R}}) \otimes_D \phi ^{(2n)}_{P^T_{j_L}} (x_{p_{j_L}})$ is a 
(bi)character on the irreducible equivalence class representative 
$P^{(2n)}(F^{T}_{\o\omega^1_{j,m_j}} \times F^{T}_{\omega^1_{j,m_j}})$~;

\item $\lambda _I(2n,j,m_j)$ is the ``unitary'' Hecke  character restricted to 
$P^{(2n)}(F^{T}_{\o\omega^1_{j,m_j}} \times F^{T}_{\omega^1_{j,m_j}})$~.
\Ei
}
\vskip 11pt 

\bpr This formula corresponds to the trace formula of proposition 2.5.13 submitted to the toroidal projective isomorphism $\gamma ^C\RL$~.\epr\vskip 11pt

\subsubsection{Proposition}

{\em Let 
\begin{align*}
& \Tr( R_{G^{(2n)}({F^T_{\o\omega}} \times{F^T_{\omega}})}))
(\phi ^{(2n)}_{G^T_{R}}  \otimes \phi ^{(2n)}_{G^T_{L}} )\\[11pt]
&\quad = {\rm vol}(G^{(2n)}({F^T_{\o\omega}} \times{F^T_{\omega}})\sum_j\sum_{m_j}
\lambda (2n,j,m_j)\end{align*}
be the trace of the integral operator on cusp biforms 
$(\phi ^{(2n)}_{G^T_{R}}  \otimes \phi ^{(2n)}_{G^T_{L}} )$ over the pseudo-ramified toroidal bilinear   semigroup 
$ G^{(2n)}({F^T_{\o\omega}} \times{F^T_{\omega}})$ and let 
\begin{align*}
& \Tr( R_{P^{(2n)}({F^{T}_{\o\omega^1}} \times{F^{T}_{\omega^1}})}
(\phi ^{(2n)}_{P^T_{R}}  \otimes \phi ^{(2n)}_{P^T_{L}} )\\[11pt]
&\quad = {\rm vol}(P^{(2n)}({F^{T}_{\o\omega^1}} \times {F^{T}_{\omega^1}})\sum_j\sum_{m_j}
\int \chi _{j,m_{j\RL}}\ dx_{p_{j_R}}\ dx_{p_{j_L}}
\end{align*}
be the trace of the integral operator on bifunctions over the bilinear parabolic semigroup
$P^{(2n)}({F^{T}_{\o\omega^1}} \times{F^{T}_{\omega^1}})$~.

Then, the trace of the integral operator on pseudo-unramified cusp biforms $(\phi ^{(2n)}_{G^{T,nr}_{R}}  \otimes \phi ^{(2n)}_{G^{T,nr}_{L}} )$ over the pseudo-unramified toroidal bilinear  semigroup  
$ G^{(2n)}(F^{T,nr}_{\o\omega_\oplus } \times F^{T,nr}_{\omega_\oplus })$ is given by:
\begin{align*}
& \Tr( R_{G^{(2n)}(F^{T,nr}_{\o\omega} \times F^{T,nr}_{\omega})}))
(\phi ^{(2n)}_{G^{T,nr}_{R}}  \otimes \phi ^{(2n)}_{G^{T,nr}_{L}} )\\[11pt]
&\quad = \L({\rm vol}(G^{(2n)}({F^T_{\o\omega}} \times{F^T_{\omega}})\Big/
{\rm vol}(P^{(2n)}({F^{T}_{\o\omega^1}} \times{F^{T}_{\omega^1}})\R)
\sum_j\sum_{m_j} \lambda (2n,j,m_j)\Big/ \lambda _I(2n,j,m_j)\;.\end{align*}
}
\vskip 11pt 

\paragraph{Sketch of proof:}  This results from the action by convolution of the unitary operator
$R({P^{(2n)}({F^{T}_{\o\omega^1}} \times{F^{T}_{\omega^1}}))}$ on the bifunctions of the pseudo-unramified 
bilinear   semigroup $G^{(2n)}(F^{T,nr}_{\o\omega} \times F^{T,nr}_{\omega})
)$ as developed in proposition 2.5.12 and, more precisely, from the map:
\[ G^{(2n)}_{F^T\to F^{T,nr}}: \quad
G^{(2n)}({F^T_{\o\omega}} \times{F^T_{\omega}}))
\To G^{(2n)}(F^{T,nr}_{\o\omega} \times F^{T,nr}_{\omega})\]
whose kernel is $\Ker (G^{(2n)}_{F^T\to F^{T,nr}})=
P^{(2n)}(F^{T,nr}_{\o\omega^1} \times F^{T,nr}_{\omega^1})$ as developed in section 2.4.1.\epr
\vskip 11pt 

\subsubsection{Proposition}

{\em Let
\begin{align*}
i_{G^T_{R\to L}} : \quad
 L^{1-1}_{\rm cusp}(G^{(2n)}({F^{T}_{\o\omega}} \times{F^{T}_{\omega}}))
&\To
 L^{2}_{\rm cusp}(G^{(2n)}({F^{T}_{\omega}} \times{F^{T}_{\omega}}))\;, \\[11pt]
\phi ^{(2n)}_{G^{T}_{j_R}} (x_{g^T_{j_R}}) \otimes
\phi ^{(2n)}_{G^{T}_{j_L}} (x_{g^T_{j_L}})
&\To
\phi ^{(2n)}_{G^{T}_{j_L}} (x_{g^T_{j_L}}) \otimes
\phi ^{(2n)}_{G^{T}_{j_L}} (x_{g^T_{j_L}})\;,\end{align*}
be the involution map as introduced in section 2.5.4 and transforming the cusp biforms of
$ L^{1-1}_{\rm cusp}(G^{(2n)}({F^{T}_{\o\omega}} \times{F^{T}_{\omega}}))$ 
into the cusp biform $\phi ^{(2n)}_{G^{T}_{j_L}} (x_{g^T_{j_L}}) \otimes
\phi ^{(2n)}_{G^{T}_{j_L}} (x_{g^T_{j_L}})$ of 
$L^{2}_{\rm cusp}(G^{(2n)}({F^{T}_{\omega}} \times{F^{T}_{\omega}}))$ localized in the upper half space.

Then, the trace formula
\[  \Tr( R_{G^{(2n)}({F^T_{\omega}} \times{F^T_{\omega}})}
(\phi ^{(2n)}_{G^{T}_{L}} \otimes
\phi ^{(2n)}_{G^{T}_{L}} ))\]
of the integral operator on cusp biforms 
$(\phi ^{(2n)}_{G^{T}_{L}} \otimes
\phi ^{(2n)}_{G^{T}_{L}} )$ leads to the Plancherel formula:
\[ \int_{G^{(2n)}}\L| 
\phi ^{(2n)}_{G^{T}_{L}} (x_{g^T_{L}}) \R|^2\ dx_{g^T_L}
= \sum_j\sum_{m_j}\L| \widehat \phi ^{(2n)}_{G^{T}_{j_L,m_{j_L}}}\R|^2\]
where
\Bi
\item $ \widehat \phi ^{(2n)}_{G^{T}_{j_L,m_{j_L}}}=\lambda ^{\half}(2n,j,m_j)$ can be considered as a global Hecke character;
\item the sum over $j$ and $m_j$ runs over the conjugacy class representatives of the pseudo-ramified algebraic   semigroup
$G^{(2n)}({F^T_{\o\omega}} \times{F^T_{\omega}})$ localized in the upper half space.
\Ei
}
\vskip 11pt 

\bpr
According to proposition 3.4.10, we have from the following trace formula:
\begin{align*}
&  \Tr( R_{G^{(2n)}({F^T_{\o\omega}} \times{F^T_{\omega}})}
(\phi ^{(2n)}_{G^{T}_{L}} \otimes
\phi ^{(2n)}_{G^{T}_{L}} ))\\[11pt]
& \quad = \bigoplus_{j,m_j} {\rm vol}(G^{(2n)}({F^T_{\o\omega}} \times{F^T_{\omega}})
\int_{G^{(2n)}} \L| (\phi ^{(2n)}_{G^{T}_{j_L}} (x_{g^T_{j_L}})\R|^2\ dx_{g^T_{j_L}}\\[11pt]
 & \quad = \bigoplus_{j,m_j} {\rm vol}(G^{(2n)}({F^T_{\o\omega}} \times{F^T_{\omega}})
 \L| \lambda ^{\half}(2n,j,m_j)\R|^2\;, \end{align*}
the equality
\[
 \bigoplus_{j,m_j} 
\int_{G^{(2n)}} \L| (\phi ^{(2n)}_{G^{T}_{j_L}} (x_{g^T_{j_L}})\R|^2\ dx_{g^T_{j_L}}
=  \bigoplus_{j,m_j} 
 \L| \lambda ^{\half}(2n,j,m_j)\R|^2\;.\]
If we take into account that
\Be
\item $\int_{G^{(2n)}} \L| (\phi ^{(2n)}_{G^{T}_{L}} (x_{g^T_{L}})\R|^2\ dx_{g^T_{L}}
=  \bigoplus_{j,m_j} 
\int_{G^{(2n)}} \L| (\phi ^{(2n)}_{G^{T}_{j_L}} (x_{g^T_{j_L}})\R|^2\ dx_{g^T_{j_L}}
$ expressing the development of the integral of the square of the cusp biform
$ \phi ^{(2n)}_{G^{T}_{L}} (x_{g^T_{L}}) $ according to the squares of the cusp subbiforms
$ \phi ^{(2n)}_{G^{T}_{j_L}} (x_{g^T_{j_L}}) $ on the conjugacy class representatives
$g^{(2n)}_{T_L}[,m_j]$ of the pseudo-ramified   bilinear semigroup
$G^{(2n)}({F^T_{\o\omega}} \times{F^T_{\omega}})$~;

\item $\begin{array}[t]{ll}
 \widehat \phi ^{(2n)}_{G^{T}_{j_L,m_{j_L}}}
&= \int_{G^{(2n)}}  \phi ^{(2n)}_{G^{T}_{L}} (x_{g^T_{L}})\ e^{-2\pi ijz}\ dx_{g^T_L}\\[11pt]
&= \lambda ^{\half}(2n,j,m_j)\;, \end{array}$

we get the searched Plancherel formula
\begin{equation}
\int_{G^{(2n)}} \L| (\phi ^{(2n)}_{G^{T}_{L}} (x_{g^T_{L}})\R|^2\ dx_{g^T_{L}}
= \sum_j\sum_{m_j}
\L| \widehat \phi ^{(2n)}_{G^T_{{j_L,m_{j_L}}}}\R|^2
\;.\tag*{\eop}\end{equation}
\Ee
\vskip 11pt 

We can then consider the following Langlands global correspondence.
\vskip 11pt 

\subsubsection{Proposition} 

{\em  On the compactified lattice bisemispace
\[\o X_{S_{R\times L}}=\GL_n({F^T_{R}}\times {F^T_L})) \Big/ \GL_n((\ZZ\big/N\ \ZZ)^2)\;,\]
there exists the {\bf global Langlands correspondence\/}:

\[\begin{array}{ccc}
\Irr\Rep^{(2i)}_{W_{F_{R\times L}}}(W^{ab}_{F_{R}}\times W^{ab}_{F_{L}})  &\To &
\Irr\cusp (\GL_i({F^T_{\o\omega }}\times {F^T_\omega })) \\
\| && \|\\
G^{(2i)}({\wt F_{\o\omega_\oplus  }}\times {\wt F_{\omega_\oplus } }) &\To& \Eis_{R\times L}(2i,j,m_j)\\
\begin{CD} @VV\wr V \end{CD} &&\begin{CD} @AAA \end{CD}\\
 G^{(2i)}({F^T_{\o\omega_\oplus  }}\times {F^T_{\omega_\oplus } }) & \To & \widehat G^{(2i)}({F^T_{\o\omega  }}\times {F^T_{\omega } }) \\
 \| & &\\
 H^{2i}(\o S^{P_n}_{K_n},M^{2i}_{T_{R_\oplus }}\otimes M^{2i}_{T_{L_\oplus }}) & &\end{array}\]
\Bi
\item from the sum of the products, right by left, of the equivalence classes of the irreducible $2i$-dimensional Frobenius semisimple Weil-Deligne representation $\Irr\Rep^{(2i)}_{W_{F_{R\times L}}}(W^{ab}_{F_{R}}\times W^{ab}_{F_{L}})  $ of the bilinear global Weil group 
$(W^{ab}_{F_{R}}\times W^{ab}_{F_{L}})   $ given by the algebraic bilinear semigroup 
$ G^{(2i)}({\wt F_{\o\omega_\oplus  }}\times {\wt F_{\omega_\oplus } }))$

\item to the sum of the products, right by left, of the equivalence classes of the irreducible cuspidal representation $\Irr\cusp (\GL_i({F^T_{\o\omega }}\times {F^T_\omega }))$ of 
$\GL_i({F^T_{\o\omega }}\times {F^T_\omega })$
given by  $2i$-dimensional cusp biforms 
$\Eis_{R\times L}(2i,j,m_j)$~, in such a way that 
$\Irr\cusp (\GL_i({F^T_{\o\omega }}\times {F^T_\omega }))=
\Irr\cusp (\GL_i({F^T_{\o\omega_\oplus  }}\times {F^T_{\omega_\oplus } }))$~.
\Ei}
\vskip 11pt

\bpr we start with $\Irr\Rep^{(2i)}_{W\FRL}(W^{ab}_{F_R}\times W^{ab}_{F_L})$ which was proved to be equal to the algebraic bilinear semigroup $G^{(2i)}({\wt F_{\o\omega_\oplus  }}\times {\wt F_{\omega_\oplus  }})$ in proposition 3.3.10. Then, a toroidal compactification maps $G^{(2i)}({\wt F_{\o\omega_\oplus  }}\times {\wt F_{\omega_\oplus  }})$ into its compactified toroidal equivalent 
$G^{(2i)}({F^{T}_{\o\omega_\oplus  }}\times {F^{T}_{\omega_\oplus  }})$ according to sections 3.3.2 to 3.3.3.  The analytic development of the compactified   bilinear semigroup $G^{(2i)}({F^{T}_{\o\omega_\oplus  }}\times {F^{T}_{\omega_\oplus  }})$ is given by the product of the (truncated) Fourier development of  a normalized $2i$-dimensional cusp form by its left equivalent: $\Eis\RL(2i,j,m_j)$~.  So, we get the searched bijection:
\begin{equation}
 \Irr\Rep^{(2i)}_{W\FRL}(W^{ab}_{F_R}\times W^{ab}_{F_L})\simeq
\Irr \cusp (\GL_i({F^{T}_{\o\omega }}\times {F^{T}_{\omega }}))\;.\tag*{\eop}\end{equation}
\vskip 11pt 

\subsection[Langlands global correspondence on the boundary $\partial\o X^{(\rit)}_{S_{R\times L}}$ of
 the irreducible compactified lattice bisemispace $\o X_{S_{R\times L}}$]{{\boldmath Langlands global correspondence on the boundary $\partial\o X^{(\rit)}_{S_{R\times L}}$ of
 the irreducible compactified lattice bisemispace $\o X_{S_{R\times L}}$}}

\subsubsection{Definition:
the boundary of the Borel-Serre compactification} 

{\parindent=0pt
 Let $\gamma^c_{R\times L} : X_{S_{R\times L}}\to \o X_{S_{R\times L}}$\quad be the projective emergent isomorphism, as introduced in proposition 3.3.2, which maps the lattice bisemispace $ X_{S_{R\times L}}$
 into its compactified toroidal correspondent $\o X_{S_{R\times L}}$ such that 
$X_{S_{R\times L}}$
can be considered as the interior of $\o X_{S_{R\times L}}$ in the context of the Borel-Serre compactification.
\par
Then, the boundary $\partial \o X_{S_{R\times L}}$ of $\o X_{S_{R\times L}}$ will be defined as resulting from the inclusion morphism:
\[\gamma^\delta _{R\times L} : \quad \o X_{S_{R\times L}}\scalebox{1.5}{$\hookrightarrow$}  \partial \o X_{S_{R\times L}}\]
sending $\o X_{S_{R\times L}} =\GL_n({F^T_{R }}\times {F^T_L })\big/\GL_n( (\ZZ\big/N\ \ZZ)^2 )= GL_n(F^T_{\o\omega }\times F^T_\omega )$
to its real boundary $\partial \o X^{(\rit)}_{S_{R\times L}}
=\GL_n({F^{+,T}_{R }}\times {F^{+,T}_L })\big/\GL_n( (\ZZ\big/N\ \ZZ)^2 )= GL_n(F^{+,T}_{\o v }\times F^{+,T}_v )$~,
where $F^{+,T}_v=\{F^{+,T}_{v_{1_\delta }},\cdots,F^{+,T}_{v_{j_\delta ,m_{j_\delta }}},\cdots,F^{+,T}_{v_{r_\delta }}\}$ covers its complex equivalent $F^T_\omega $~.

By this way, each complex diagonal conjugacy class representative 
$g^{(2n)}_{T_{L}}[j,m_j=0]$ (resp.\linebreak $g^{(2n)}_{T_{R}}[j,m_j=0]$~) of $ \o X\SRL$ is covered by the set of $m_{j_\delta }$ real conjugacy class representatives 
$\{g^{(n)}_{T_{L}} [j_\delta ,m_{j_\delta }]\}_{m_{j_\delta }}$ (resp.
$\{g^{(n)}_{T_{R}} [j_\delta ,m_{j_\delta }]\}_{m_{j_\delta }}$~) of $\partial \o X^{(\rit)}\SRL$~.

The morphism $\gamma^\delta_{R\times L}:\quad 
\o X_{S_{R\times L}} \to \partial \o X_{S_{R\times L}} $ is such that there is a one-to-one correspondence between the left (resp. right) points 
$P^T_{a_{L[j,m_j]}} $ (resp. $P^T_{a_{R[j,m_j]}} $~) of $\partial \o X_{S_{L}}$  (resp. $\partial\o X_{S_{R}} $~) 
and their  ``real" equivalents on $\partial \o X^{(\rit)}_{S_{L}}$  (resp. $\partial \o X^{(\rit)}_{S_{R}} $~).
}

\subsubsection[Definition: the double coset decomposition of the equivalent of the Shimura bisemivariety $\partial\o S^{P_n}_{K_n}$]{{\boldmath Definition: the double coset decomposition of the equivalent of the Shimura bisemivariety $\partial\o S^{P_n}_{K_n}$ }}

The inclusion morphism $\gamma^\delta_{R\times L}$ has for result that the double coset decomposition
\[ \o S^{P_n}_{K_n}=P_n( {F^{T}_{\o\omega ^1}}\times {F^{T}_{\omega^1 }} )\setminus 
\GL_n ({F^{T}_{R   }} \times {F^{T}_L  })\big/\GL_n( (\ZZ\big/N\ \ZZ)^2 )\]
of $\GL_n( {F^{T}_{R  }} \times {F^{T}_L  })$
corresponding to $\o X_{S_{R\times L}} $
is transformed into the following double coset decomposition
\[  \partial\o S^{P_n}_{K_n}=P_n( {F^{+,T}_{\o v^1}} \times {F^{+,T}_{v^1}} )\setminus 
\GL_n( {F^{+,T}_{R}} \times {F^{+,T}_L  })\Big/\GL_n( (\ZZ\big/N\ \ZZ)^2)\]
where 
$P_n( {F^{+,T}_{v^1}} )$ (resp. $P_n( {F^{+,T}_{\o v^1}} )$~) is the standard left (resp. right) parabolic subgroup over 
the set of real irreducible completions 
 ${F^{+,T}_{v^1_{j_\delta,m_{j_\delta}}}} $ (resp. ${F^{+,T}_{\o v^1_{j_\delta,m_{j_\delta}}}} $~) 
(see section 2.4.1).
\par
$\partial\o S^{P_n}_{K_n}$  is the equivalent of a Shimura (bisemi)variety because \cite{Del2}, \cite{Har1}:
\Bi
\item  it is dependent on the morphism $\gamma^\delta_{R\times L}$ sending a ``complex" 
bisemivariety $ \o X^{(\rit)}_{S_{R\times L}} $
into a ``real" bisemivariety $\partial \o X^{(\rit)}_{S_{R\times L}} $~;
\item  it is defined with respect to an  open compact subgroup $P_n(F^{+,T}_{\o v^1}\times F^{+,T}_{v^1})$~.
\Ei
\vskip 11pt 

\subsubsection{Proposition} 

{\em\parindent=0pt
 Let $M^{2i}_{T_{\o v_R}} \otimes M^{2i}_{T_{ v_L}}$ be a $(B^{2i}_{\o  v^T}\otimes B^{2i}_{v^T})$-bisemimodule where
$B^{2i}_{\o  v^T}\otimes B^{2i}_{v^T}\simeq T_i(F^{+,T}_{\o v})^t \times T_i(F^{+,T}_{ v})$ is a tensor 
product of division semialgebras and 
let $M^{2i}_{T_{\o v_{R_\oplus }}}\otimes M^{2i}_{T_{v_{L_\oplus }}}$ be a $\GL_i(F^{+,T}_{\o v_\oplus }\times F^{+,T}_{v_\oplus })$-bisemimodule where $F^{+,T}_{v_\oplus }$ is given by
\[F^{+,T}_{v_\oplus }= \bigoplus_{j_\delta}F^{+,T}_{v_{j_\delta}}\bigoplus_{j_\delta ,m_{j_\delta }}F^{+,T}_{v_{j_\delta,m_{j_\delta }}}\;.\]
Then, the bilinear cohomology of the equivalent of the Shimura bisemivariety $\partial\o S^{P_n}_{K_n}$   is the Eisenstein cohomology:
 \[ H^{2i}(\partial \o S^{P_n}_{K_n},M^{2i}_{T_{\o v_{R_\oplus }}} \otimes M^{2i}_{T_{ v_{L_\oplus }}})
=G^{(2i)}( {F^{+,T}_{\o v_\oplus }} \times {F^{+,T}_{v_\oplus } })\]
in bijection with the $2i$-dimensional real irreducible representation $\Irr\Rep^{(2i)}_{W_{F^+\RL}}(W^{ab}_{F^+_{\o v}}\times W^{ab}_{F^+_v})$ of the product, right by left, of global Weil groups.}\vskip 11pt

\bpr
\Bi
\item This proposition is a transposition of proposition 3.3.10 to $H^{2i}(\partial \o S^{P_n}_{K_n},\cdot)$~.
\item The fact that $H^{2i}(\partial \o S^{P_n}_{K_n},M^{2i}_{T_{\o v_{R_\oplus }}} \otimes M^{2i}_{T_{ v_{L_\oplus }}})$ is the Eisenstein (bilinear) cohomology results from the works of G. Harder and J. Shimura \cite{Har2}, \cite{Sch}.
\item the global Weil groups $W^{ab}_{F^+_v}$ and $W^{ab}_{F^+_{\o v}}$ 
referring to finite real algebraic extensions can be defined as they were introduced for complex algebraic extensions in 
definition 1.1.9.
\Ei
Then, referring to proposition 3.3.10, we can state that:
\be\Irr \Rep^{(2i)}_{W_{F^+\RL}}(W^{ab}_{F^+_{\o v}}\times W^{ab}_{F^+_v})
= G^{(2i)}({\wt F^{+}_{\o v_\oplus }}\times {\wt F^{+}_{v_\oplus }})\tag*{\eop}\ee

\vskip 11pt 

\subsubsection{Proposition}

{\em The cohomology $H^{2i}(\partial \o S^{P_n}_{K_n},M^{2i}_{T_{\o v_{R_\oplus }}}\otimes M^{2i}_{T_{v_{L_\oplus }}})$ has a decomposition according to the equivalent representatives $g^{(2i)}_{T\RL}[j_\delta,m_{j_\delta}]$ of the conjugacy classes of the real   bilinear semigroup $G^{(2i)}({F^{+,T}_{\o v}}\times {F^{+,T}_{v}})$ according to:
\[ H^{2i}(\partial \o S^{P_n}_{K_n},M^{2i}_{T_{\o v_{R_\oplus }}}\otimes M^{2i}_{T_{v_{L_\oplus }}})
= \bigoplus_{j_\delta}\bigoplus_{m_{j_\delta}}g^{(i)}_{T\RL}[j_\delta,m_{j_\delta}]\;.\]
}
\vskip 11pt 

\bpr this is an adaptation of the real case of proposition 3.4.1 since
\begin{equation} H^{2i}(\partial \o S^{P_n}_{K_n},M^{2i}_{T_{\o v_{R_\oplus }}}\otimes M^{2i}_{T_{v_{L_\oplus }}})
= G^{(2i)}({F^{+,T}_{\o  v}}\times {F^{+,T}_{v}})\;.\tag*{\eop}
\end{equation}
\vskip 11pt 

\subsubsection{Proposition} 

{\em  The Eisenstein cohomology $H^{2i}(\partial \o S^{P_n}_{K_n}, \widehat M^{2i}_{T_{\o v_{R_\oplus }}} \otimes \widehat M^{2i}_{T_{ v_{L_\oplus }}})$
has a power series representation given by the product of a right $2i$-dimensional global elliptic semimodule $\ELLIP_R(2i, j_\delta, m_{j_\delta})$ by its left equivalent $\ELLIP_L(2i, j_\delta, m_{j_\delta})$ according to:
\[ H^{2i}(\partial \o S^{P_n}_{K_n},\widehat M^{2i}_{T_{\o v_{R_\oplus }}} \otimes \widehat M^{2i}_{T_{ v_{L_\oplus }}})=
\ELLIP_R(2i, j_\delta, m_{j_\delta})\otimes_{(D)} \ELLIP_L(2i, j_\delta, m_{j_\delta})\] 
where
\begin{align*}
\ELLIP_L(2i, j_\delta, m_{j_\delta}) &=\txt\bigoplus\limits_{j_\delta=1}^t \bigoplus\limits_{m_ {j_\delta}}\lambda^{\half}(i,j_\delta,m_{j_\delta})e^{2\pi ij_\delta x}\;,\\
\ELLIP_R(2i, j_\delta, m_{j_\delta}) &=\txt\bigoplus\limits_{j_\delta=1}^t \bigoplus\limits_{m_ {j_\delta}} \lambda^{\half}(i,j_\delta,m_{j_\delta}) e^{-2\pi ij_\delta x}\;,\end{align*}
with
\Bi
\item $\vec x=\sum\limits_{c=1}^ix_c\vec e_c$ a vector of $(F^{+}_{v_{j_\delta }})^i$ and more precisely a point of $G^{(i)}(F^{+}_{v_{j_\delta }})\equiv g^{(i)}_{L}[j_\delta ,m{j_\delta }]$ and 
$x=\sum\limits^i_{c=1}x_c\ |\vec e_c|$~;\vskip 11pt
\item  $\lambda(i,j_\delta,m_{j_\delta}) =\prod\limits_{c=1}^i \lambda_c(i,j_\delta,m_{j_\delta}) $ a product of eigenvalues $\lambda_c(i,j_\delta,m_{j_\delta}) $ of 
$g_i(\Os_{F^{+,T}_{\o v_{j_\delta,m_{j_\delta}}}}\times \Os_{F^{+,T}_{v_{j_\delta,m_{j_\delta}}}})
\in  \GL_i((\ZZ\big/N\ \ZZ)^2)$~.
\Ei}\vskip 11pt

\bpr {\parindent=0pt
since $g_{T_L}^{(i)} [j_\delta,m_{j_\delta}]$ (resp. $g^{(i)}_{T_R}[j_\delta,m_{j_\delta}]$~) is a $i$-dimensional left (resp. right) real semitorus, it has the following analytical development:
\begin{align*}
T^i_L[j_\delta,m_{j_\delta}] &= \lambda^{\half}(i,j_\delta,m_{j_\delta})e^{2\pi ij_\delta x}\\ 
(\text{resp.} \quad
T^i_R[j_\delta,m_{j_\delta}] &= \lambda^{\half}(i,j_\delta,m_{j_\delta})e^{-2\pi ij_\delta x}\ )\end{align*} 
where
\[ \lambda(i,j_\delta,m_{j_\delta})=\txt\prod\limits_{c=1}^i \lambda_c(i,j_\delta,m_{j_\delta})=\det(\alpha _{i^2;(j_\delta)^2}\times D_{(j_\delta)^2;m^2_{j_\delta}})_{ss} \simeq j^{2i}\cdot (N)^{2i}\]
is a global Hecke character since the $\lambda_c(i,j_\delta,m_{j_\delta})$
are eigenvalues of $g_i(\Os_{F^{+,T}_{\o v_{j_\delta,m_{j_\delta}}}}\times \Os_{F^{+,T}_{v_{j_\delta,m_{j_\delta}}}})$ which is the $j_\delta$-th coset representative of the product of Hecke operators (see propositions 3.4.4 and 3.4.5). 

$\lambda ^\half(i,j_\delta ,m_{j_\delta })$ has an inflation action on $e^{2\pi ij_\delta x }$ leading to the existence of a global elliptic {\bf semimodule\/} $\ELLIP(\centerdot,\centerdot, \centerdot)$~.

On the other hand, 
$g_{T_L}^{(i)} [j_\delta,m_{j_\delta}]$ (resp. $g^{(i)}_{T_R}[j_\delta,m_{j_\delta}]$~) is defined on the real left (resp. right) place $v_{j_\delta}$ (resp. $\o v_{j_\delta}$~) implying the global Frobenius substitution:
\begin{align*}
e^{2\pi ix} &\To e^{2\pi ij_\delta x}\\
(\text{resp.} \quad
e^{-2\pi ix} &\To e^{-2\pi ij_\delta x}\ ).\end{align*}
Note that these $2i$-dimensional left and right global elliptic semimodules $\ELLIP_L(2i, j_\delta, m_{j_\delta}) $ and $\ELLIP_R(2i, j_\delta, m_{j_\delta}) $ are the $2i$-dimensional equivalents of the $2$-dimensional global elliptic semimodules introduced in \cite{Pie3} (see also \cite{Ande}).\epr
}\vskip 11pt

\subsubsection{Characters and extended bicharacters}

The functions $x_c\to e^{2\pi ij_\delta x_c}$ from the real line, i.e. the irreducible central real completion $F^+_{v^1_{1_\delta }}$~, to the complex numbers of modulus one,   are characters if they 
satisfy the identity $\chi (x_cy_c)=\chi (x_c)\cdot \chi (y_c)$~.  In the bilinear context considered here, the characters are 
replaced by extended bicharacters defined as follows:
\vskip 11pt 

If $x^+_c\cdot x^-_c\equiv -(x_c)^2$ denote a bipoint of the real irreducible completion $F^+_{\o v^1_{1_\delta }}\times F^+_{v^1_{1_\delta }}$~, i.e. the product of a point $x^-_c\equiv x_c$ by its 
symmetrical $x^+_c\equiv -x_c$~, $x_c\in F^+_{v^1_{1_\delta }}$~, then an extended bicharacter $\chi _+(x^+_c)\cdot \chi _-(x^-_c)$ is a continuous bifunction
\begin{align*}
\chi _+\cdot \chi _- : \quad &x^+_c\cdot x^-_c\\
& \quad \To \lambda ^{\half}(1,j_\delta ,m_{j_\delta })\cdot e^{-2\pi ij_\delta x_c}\times 
 \lambda ^{\half}(1,j_\delta ,m_{j_\delta })\cdot e^{2\pi ij_\delta x_c}=
 \lambda (1,j_\delta ,m_{j_\delta }) \end{align*}
from $F^+_{\o v^1_{1_\delta }}\times F^+_{v^1_{1_\delta }}$ to the product of the complex numbers 
$ \lambda ^{\half}(1,j_\delta ,m_{j_\delta })\cdot e^{2\pi ij_\delta x_c}$ of modulus $\neq 1$ by their complex conjugates where $ \lambda  (1,j_\delta ,m_{j_\delta })$ is generally an eigenvalue of a product $U_{j_{\delta _R}}\times U_{j_{\delta _L}}$ of Hecke operators.
\vskip 11pt

\subsubsection{Proposition} 

{\em  The $2i$-dimensional global elliptic bisemimodule $\ELLIP_R(2i, j_\delta, m_{j_\delta}) \otimes_{(D)} \ELLIP_L(2i, j_\delta, m_{j_\delta})$  is:
\Be
\item  a solvable power biseries, constituting an irreducible (super)cuspidal representation of the (truncated) (super)cuspidal biform $\Eis\RL(n,j,m_j=0)$~;
\item  an eigenbifunction of the Hecke bioperator $T_R(2i; t) \otimes T_L (2i; t)$~.\Ee
}\vskip 11pt

\bpr
\Be
\item The fact that $\ELLIP_R(2i, j_\delta, m_{j_\delta})\otimes \ELLIP_L(2i, j_\delta, m_{j_\delta}) $ is a  solvable power biseries results  from: 
\[\Reps(\GL_i( {F^{+,T}_{\o v_\oplus }} \times {F^{+,T}_{v_\oplus }  }))=\ELLIP_{R\times L}(2i, j_\delta, m_{j_\delta})\]
where
$\ELLIP_{R\times L}\equiv \ELLIP_R\otimes \ELLIP_L$ (see proposition 3.4.7).

$ \ELLIP_R(2i,j_\delta,m_{j_\delta}) \otimes \ELLIP_L(2i,j_\delta,m_{j_\delta}) $ is an irreducible supercuspidal representation of $\Eis\RL(2i,j,m_j=0)$ due to the morphism $\gamma ^\delta \RL:\o X\SRL\to\partial\o X\SRL$ related to the generation of the boundary of the Borel-Serre compactification according to definition 3.5.1.

\item  The proof that $\ELLIP_{R\times L}(2i, j_\delta, m_{j_\delta})$ is an eigenbifunction of the Hecke bioperator  $T_R(2i; t) \otimes T_L (2i; t)$ can be handled as in proposition 3.4.6.  
\par
Let us yet note that the $j_\delta $-th bicoset representative of $T_R(2i; j_\delta ) \otimes T_L(2i; j_\delta )$ is given, with evident notations, by~:
\begin{equation}
U_{j_{\delta _R}} \times U_{j_{\delta _L}} = t^t_i( \Os_{F^{+,T}_{\o v_{j_\delta,m_{j_\delta}}}}
)\times t_i( \Os_{F^{+,T}_{v_{j_\delta,m_{j_\delta}}}})\;.\tag*{\eop}\end{equation}
\Ee \vskip 11pt

\subsubsection[Bialgebra of elliptic semimodules $L^{1-1}_{\ELLIP}(G^{(2i)}( {F^{+,T}_{\o v}}
\times {F^{+,T}_{v}}))$]{{\boldmath Bialgebra of elliptic semimodules $L^{1-1}_{\ELLIP}(G^{(2i)}( {F^{+,T}_{\o v}}
\times {F^{+,T}_{v}}))$}}

As in section 3.4.9, we can introduce a bialgebra $L^{1-1}_{\ELLIP}(G^{(2i)}( {F^{+,T}_{\o v}}
\times {F^{+,T}_{v}}))$ of complex-valued smooth continuous bifunctions
$ \phi ^{(i)}_{G^T_{j^\delta _R}} (x_{g^T_{j^\delta _R}}) \otimes
 \phi ^{(i)}_{G^T_{j^\delta _L}} (x_{g^T_{j^\delta _L}})$ on the pseudo-ramified toroidal bilinear   semigroup 
$G^{(2i)}( {F^{+,T}_{\o v}}
\times {F^{+,T}_{v}})$ from the corresponding bialgebra\linebreak
$L^{1-1}\RL(G^{(2i)}( {F^{+}_{\o v}}
\times {F^{+}_{v}}))$ of smooth continuous bifunctions on the pseudo-ramified bilinear   semigroup 
$ G^{(2i)}({F^{+}_{\o v}}
\times {F^{+}_{v}})$ under the action of a projective toroidal isomorphism.

Then, each bifunction on the conjugacy class representative
$g^{(i)}_{{T\RL}}[j_\delta ,m_{j_\delta }]\in G^{(i)} ( {F^{+,T}_{\o v}}
\times {F^{+,T}_{v}})$ is given, according to proposition 3.4.5, by
\begin{align*}
 \phi ^{(i)}_{G^T_{j^\delta _R}} (x_{g^T_{j^\delta _R}}) \otimes_D
 \phi ^{(i)}_{G^T_{j^\delta _L}} (x_{g^T_{j^\delta _L}}) 
&= T^i_R[j_\delta ,m_{j_\delta }]\otimes_D T^i_L[j^\delta ,m_{j^\delta }]\\[11pt]
&= \lambda ^{\half}(i,j_\delta ,m_{j_\delta })\ e^{-2\pi ij_\delta x} \otimes _{(D)}
 \lambda ^{\half}(i,j_\delta ,m_{j_\delta })\ e^{2\pi ij_\delta x} \end{align*}
and the sum of all bifunctions on the conjugacy class representatives of
$G^{(2i)}( {F^{+,T}_{\o v}}
\times {F^{+,T}_{v}})$~:
\begin{align*}
 \phi ^{(2i)}_{G^T_{\delta _R}} (x_{g^T_{R}}) \otimes_D
 \phi ^{(2i)}_{G^T_{\delta _L}} (x_{g^T_{L}}) 
&= \bigoplus^t_{j_\delta =1}\bigoplus_{m_{j_\delta }}
 \L(\phi ^{(i)}_{G^T_{j^\delta _R}} (x_{g^T_{j^\delta _R}}) \otimes_D
 \phi ^{(i)}_{G^T_{j^\delta _L}} (x_{g^T_{j^\delta _L}}) \R)\\[11pt]
&= \bigoplus_{j_\delta}\bigoplus_{m_{j_\delta }}
( \lambda ^{\half}(i,j_\delta ,m_{j_\delta })\ e^{-2\pi ij_\delta x} \otimes_{(D)}
 \lambda ^{\half}(i,j_\delta ,m_{j_\delta })\ e^{2\pi ij_\delta x} )\\[11pt]
&= \ELLIP_R(2i,j_\delta ,m_{j_\delta }) \otimes_{(D)} \ELLIP_L(2i,j_\delta ,m_{j_\delta })
\end{align*}
is the (diagonal) product of a right $2i$-dimensional global elliptic semimodule\linebreak
$\ELLIP_R(2i,j_\delta ,m_{j_\delta })$ by its left equivalent $\ELLIP_L(2i,j_\delta ,m_{j_\delta })$~.
\vskip 11pt 

\subsubsection{Proposition}

{\em
Let $R_{G^{(2i)}( {F^{+,T}_{\o v}}
\times {F^{+,T}_{v}})} 
 (\phi ^{(2i)}_{G^T_{\delta _R}} \otimes
 \phi ^{(2i)}_{G^T_{\delta _L}})$  be the integral operator on elliptic bisemimodules
$( \phi ^{(2i)}_{G^T_{\delta _R}}\otimes
 \phi ^{(2i)}_{G^T_{\delta _L}} )$ over the pseudo-ramified toroidal bilinear real   semigroup 
$G^{(2i)}( {F^{+,T}_{\o v_\oplus }}
\times {F^{+,T}_{v_\oplus }})$~.

Then,  the trace formula of this integral operator is given by:
\begin{align*}
&\Tr( R_{G^{(2i)}( {F^{+,T}_{\o v}}
\times {F^{+,T}_{v}})} 
 (\phi ^{(2i)}_{G^T_{\delta _R}} \otimes
 \phi ^{(2i)}_{G^T_{\delta _L}})\\[11pt]
&\quad = {\rm vol} (G^{(2i)}( {F^{+,T}_{\o v}}
\times {F^{+,T}_{v}})) \sum_{j^\delta }\sum_{m_{j^\delta }}\lambda (i,j_\delta ,m_{j_\delta })\;.\end{align*}
}
\vskip 11pt 

\paragraph{Sketch of proof:}  This is an adaptation of the trace formula
$\Tr( R_{G^{(2i)}( {F^{T}_{\o \omega }}
\times {F^{T}_{\omega }})} 
 (\phi ^{(2i)}_{G^T_{R}} \otimes
 \phi ^{(2i)}_{G^T_{L}})$ developed in proposition 3.4.10 to the real case.\epr
\vskip 11pt


\subsubsection{Proposition} 

{\em On the boundary $\partial \o X^{(\rit)}_{S_{R\times L}}= \GL_i( {F^{+,T}_{R}} \times {F^{+,T}_L  })\big/\GL_i( (\ZZ\big/N\ \ZZ)^2)$ of the compactified lattice bisemispace $\o X_{S_{R\times L}}$~, there exists the Langlands global correspondence:
\[\begin{array}{ccc}
\Irr\Rep^{(2i)}_{W^+_{F_{R\times L}}}(W^{ab}_{F^{+}_{\o v}}\times W^{ab}_{F^{+}_{v }}) &\To &
\Irr\ELLIP (\GL_i({F^{+,T}_{\o v }}\times {F^{+,T}_v })) \\
\| && \|\\
G^{(2i)}({\wt F^+_{\o v_\oplus }}\times {\wt F^+_{v_\oplus }}) &\To& \ELLIP_{R\times L}(2i,j_\delta,m_{j_\delta})\\
\begin{CD}@VVV \end{CD} &&\begin{CD}@AAA \end{CD}
\\ 
G^{(2i)}({F^{+,T}_{\o v_\oplus }}\times {F^{+,T}_{v_\oplus }})&
\To &
\widehat G^{(2i)}({F^{+,T}_{\o v }}\times {F^{+,T}_{v }})\\
 \| &&\\
 H^{2i}(\partial \o S^{P_n}_{K_n},M^{2i}_{T_{\o v_{R_\oplus }}}\otimes M^{2i}_{T_{v_{L_\oplus }}})&&\end{array}\]
\Bi
\item from the sum of products, right by left, of the equivalence classes of the irreducible $2i$-dimensional Frobenius semisimple 
Weil-Deligne representation $\Irr\Rep^{(2i)}_{W^+_{F_{R\times L}}}(W^{ab}_{F^{+}_{\o v}}
\times W^{ab}_{F^{+}_{v}})  $ of the bilinear global Weil group 
$(W^{ab}_{F^{+}_{\o v}}\times W^{ab}_{F^{+}_{v}})  $
given by the algebraic bilinear real semigroup
 $G^{(2i)}({\wt F^{+}_{\o v_\oplus  }}\times {\wt F^{+}_{v_\oplus } })) $
 \item to the sum of the products, right by left, of
 the equivalence classes of the irreducible elliptic representation
$\Irr\ELLIP
(\GL_i({F^{+,T}_{\o v }}\times {F^{+,T}_v })) $
of $\GL_i({F^{+,T}_{\o v }}\times {F^{+,T}_v }) $ given by the $2i$-dimensional solvable
 global elliptic bisemimodule $\ELLIP_{R\times L}(2i,j_\delta,m_{j_\delta})$
in such a way that $\Irr\ELLIP
(\GL_i({F^{+,T}_{\o v }}\times {F^{+,T}_v })) =
\ELLIP_{R\times L}(2i,j_\delta,m_{j_\delta}) $~.
\Ei
}\vskip 11pt

\stepcounter{subsubsection}\noindent{\bf{\thesubsubsection. The $n$-dimensional irreducible global Langlands correspondences\/}} can be summarized in the following diagram:
\[\begin{array}{ccc}
\Irr\Rep^{(2n)}_{W_{F_{R\times L}}}(W^{ab}_{F_{R}}\times W^{ab}_{F_{ L}})) &\To &
\Irr\cusp (\GL_n({F^{T}_{\o \omega  }}\times {F^{T}_\omega  })) \\
\Big\downarrow && \Big\downarrow\\
\Irr\Rep^{(2n)}_{W_{F^+_{R\times L}}}(W^{ab}_{F^{+}_{\o v }}\times W^{ab}_{F^{+}_{ v }})) &\To &
\Irr\ELLIP (\GL_n({F^{+,T}_{\o v  }}\times {F^{+,T}_v  })) \end{array}\]
\section[{{\boldmath Langlands global correspondences for reducible representations of $\GL(n)$}}]{{\boldmath Langlands global correspondences for reducible\protect\newline representations of $\GL((2)n)$}}
\subsection[Reducibility of $\GL_{(2)n}({F_{\o \omega   }}\times {F_\omega   }) $]{{\boldmath Reducibility of $\GL_{(2)n}({F_{\o \omega   }}\times {F_\omega   }) $}}

\subsubsection[Definition: Partial reducibility of $\GL_n({F_{\o \omega   }}\times {F_\omega   }) $]{{\boldmath Definition: Partial reducibility of $\GL_n({F_{\o \omega   }}\times {F_\omega   }) $}}

Another way of constructing admissible representations of the general bilinear semigroup $\GL_n({F_{\o \omega   }}\times {F_\omega   }) $
 is to consider a partition $n=n_1+n_2+\cdots +n_\ell+\cdots+n_s$
of the integer $n$~, $n_\ell\in\NN$~, $n_\ell<n\in\NN$~, leading to the representations $\Rep(\GL_{n_\ell}({F_{\o \omega   }}\times {F_\omega   })) $
of $\GL_{n_\ell}({F_{\o \omega   }}\times {F_\omega   }) $~,
for $1 \le \ell \le s$~, in such a way that 
$\Rep(\GL_{n_\ell}({F_{\o \omega   }}\times {F_\omega   })) $
 constitutes an irreducible representation of 
$\GL_{n_\ell}({F_{\o \omega   }}\times {F_\omega   }) $~.

The tensor product 
\[\Rep(\GL_{n_1}({F_{\o \omega   }}\times {F_\omega   })) 
\otimes\cdots\otimes
\Rep(\GL_{n_\ell}({F_{\o \omega   }}\times {F_\omega   })) 
\otimes\cdots\otimes
\Rep(\GL_{n_s}({F_{\o \omega   }}\times {F_\omega   })) \]
may not be irreducible but it has an irreducible quotient given by the formal sum:
\begin{align*}
&\Rep(\GL_{n=n_1+\cdots+n_s}({F_{\o \omega   }}\times {F_\omega   })) \\
&\quad = \Rep(\GL_{n_1}({F_{\o \omega   }}\times {F_\omega   })) \boxplus \cdots\\
&\qquad \quad \boxplus \Rep(\GL_{n_\ell}({F_{\o \omega   }}\times {F_\omega   })) 
\boxplus\cdots\boxplus
\Rep(\GL_{n_s}({F_{\o \omega   }}\times {F_\omega   })) \end{align*}
 which constitutes a partially reducible representation of 
$\GL_{n}({F_{\o \omega   }}\times {F_\omega   }) $~.\vskip 11pt

\subsubsection[Definition: complete reducibility of $\Rep(\GL_{2n}({F_{\o \omega   }}\times {F_\omega   })) $]{{\boldmath Definition: complete reducibility of $\Rep(\GL_{2n}({F_{\o \omega   }}\times {F_\omega   })) $}}

If the partition $2n = 2_1 + 2_2 + \cdots+2_\ell+\cdots+2_n$ of $2n$~, $2_\ell=|2|$~, is considered, then the tensor product
\[\Rep(\GL_{2_1}({F_{\o \omega   }}\times {F_\omega   })) \otimes \cdots\otimes
\Rep(\GL_{2_n}({F_{\o \omega   }}\times {F_\omega   })) \]
of representations of $\GL_{2_\ell}({F_{\o \omega   }}\times {F_\omega   }) $~,
$1 \le \ell \le n$~, must be taken into account. This tensor product has an irreducible quotient given by the formal sum:
\begin{align*}
&\Rep(\GL_{2n=2_1+\cdots+2_\ell+\cdots+2_n}({F_{\o \omega   }}\times {F_\omega   })) \\
&\quad = \Rep(\GL_{2_1}({F_{\o \omega   }}\times {F_\omega   })) \boxplus \cdots\\
&\qquad \boxplus \Rep(\GL_{2_\ell}({F_{\o \omega   }}\times {F_\omega   })) 
\boxplus\cdots\boxplus
\Rep(\GL_{2_n}({F_{\o \omega   }}\times {F_\omega   })) \end{align*}
which constitutes a completely reducible representation of $\GL_{2n}({F_{\o \omega   }}\times {F_\omega   })) $
decomposing into the direct sum of irreducible representations 
$\Rep(\GL_{2_\ell}({F_{\o \omega   }}\times {F_\omega   }))$~, $1\le \ell\le n$~.\vskip 11pt 

\subsubsection{Proposition} 

{\em Let $2n_L = 2_{1_L} + 2_{2_L} + \cdots + 2_{k_L} + \cdots + 2_{\ell_L} + \cdots + 2_{n_L}$ (resp. $2n_R = 2_{1_R} + 2_{2_R} + \cdots + 2_{k_R} + \cdots + 2_{\ell_R} + \cdots + 2_{n_R}$~) be a   partition of $2n_L$ (resp. $2n_R$~) labeling the irreducible representations of $T_{2n_L}({F_\omega})$ (resp. $T^t_{2n_R}({F_{\o \omega}})$~).
Then:
\Be
\item $\Rep(\GL_{2n=2_1+\cdots+2_\ell+\cdots+2_n}({F_{\o \omega   }}\times {F_\omega   })) = \mathop{\boxplus}\limits^{n}_{\ell=1} \Rep(\GL_{2_\ell}({F_{\o \omega   }}\times {F_\omega   })) $
 constitutes a
completely reducible orthogonal bilinear representation of $\GL_{2n}({F_{\o \omega   }}\times {F_\omega   }) $~.

\item $\Rep(\GL_{2n_{R\times L}}({F_{\o \omega   }}\times {F_\omega   })) = 
\mathop{\boxplus}\limits^{2n}_{2_{\ell_R}\equiv2_{\ell_L}=2} \Rep(\GL_{2_{\ell_{R\times L}}}
({F_{\o \omega   }}\times {F_\omega   })) \mathop{\boxplus}\limits_{2_{k_R}\neq 
2_{\ell_L}}$\linebreak
$ \Rep(T^t_{2_{k_R}}({F_{\o \omega   }})\times T_{2_{\ell_L}}({F_\omega   }))$~,
where $\GL_{2_{\ell_{R\times L}}} $ is another notation for $\GL_{2_\ell}$~, constitutes a completely reducible nonorthogonal bilinear representation of $\GL_{2n}({F_{\o \omega   }}\times {F_\omega   }))$~.\Ee}\vskip 11pt

\bpr if we consider the decomposition of $\GL_{2n}({F_{\o \omega   }}\times {F_\omega   })$ into a product of trigonal semigroups according to:
\[ \GL_{2n}({F_{\o \omega   }}\times {F_\omega   })
= T^t_{2{n_R}}({F_{\o \omega   }})\times T_{2{n_L}}({F_\omega   })\]
and the partitions of $2n_L$ and $2n_R$ as envisaged in this proposition, we have that:
\begin{align*}
& \Rep(\GL_{2n_{R\times L}}({F_{\o \omega   }}\times {F_\omega   }))\\
&\quad = \left(\mathop{\boxplus}\limits^{2n}_{2_{\ell_R}=
2} \Rep(T^t_{2_{\ell_R}}({F_{\o \omega   }}))\right)\txt\bigotimes \left(\mathop{\boxplus}\limits^{2n}_{2_{\ell_L}=
2} \Rep(T_{2_{\ell_L}}({F_{ \omega   }}))\right)\\
&\quad =
\mathop{\boxplus}\limits^{2n}_{2_{\ell_R}\equiv2_{\ell_L}=2} \Rep(\GL_{2_{\ell}}
({F_{\o \omega   }}\times {F_\omega   })) \mathop{\boxplus}\limits_{2_{k_R}\neq 
2_{\ell_L}}  \Rep(T^t_{2_{k_R}}({F_{\o \omega   }})\times T_{2_{\ell_L}}({F_\omega   }))
\end{align*}
where \quad $\GL_{2_{\ell}}
({F_{\o \omega   }}\times {F_\omega   }) =T^t_{2_{\ell_R}}({F_{\o \omega   }})\times T_{2_{\ell_L}}({F_\omega   })$~.  

$\Rep(\GL_{2n_{R\times L}}({F_{\o \omega   }}\times {F_\omega   })) $
then decomposes diagonally according to the direct sum of irreducible representations\linebreak 
$\Rep(\GL_{2_\ell}({F_{\o \omega   }}\times {F_\omega   })) $ leading to a completely reducible orthogonal bilinear representation of $\Rep(\GL_{2n}({F_{\o \omega   }}\times {F_\omega   })) $
if $\mathop{\boxplus}\limits_{2_{k_R}\neq 
2_{\ell_L}}  \Rep(T^t_{2_{k_R}}({F_{\o \omega   }})\times T_{2_{\ell_L}}({F_\omega   }))=0$~.\epr\vskip 11pt

\subsubsection{Corollary} 

{\em \Be
\item The representation of the general bilinear semigroup $\GL_{n}({F_{\o \omega   }}\times {F_\omega   })$ 
is {\bf partially reducible\/} if it decomposes according to the direct sum of irreducible bilinear representations $\Rep(\GL_{n_\ell}({F_{\o \omega   }}\times {F_\omega   }))$~. 

\item The representation of the general bilinear semigroup $\GL_{2n}({F_{\o \omega   }}\times {F_\omega   })$
is {\bf orthogonally completely reducible\/} if it decomposes diagonally according to the direct sum of irreducible bilinear representations $\Rep(\GL_{2_\ell}({F_{\o \omega   }}\times {F_\omega   })) $~.

\item The representation of the general bilinear semigroup $\GL_{2n}({F_{\o \omega   }}\times {F_\omega   })) $
is {\bf nonorthogonally completely reducible\/} if it decomposes diagonally according to the direct sum of irreducible bilinear representations $\Rep(\GL_{2_\ell}({F_{\o \omega   }}\times {F_\omega   })) $
and off diagonally according to the direct sum of irreducible bilinear representations 
$ \Rep(T^t_{2_{k_R}}({F_{\o \omega   }})\times T_{2_{\ell_L}}({F_\omega   }))$~.
\Ee}\vskip 11pt

\subsubsection{Lemma} 

{\em  Let $X^{2n}_{\SRL}= \GL_{n}({ \widetilde F_{R  }}\times {\widetilde F_L   })) \big/\GL_n( (\ZZ\big/N\ \ZZ)^2 )$
denote the lattice bisemispace
as introduced in definition 2.3.1, and let
\begin{align*}
X^{2n=2n_1+\cdots+2n_s}_{\SRL}&= \GL_{n=n_1+\cdots+n_s}({\widetilde F_{R   }}\times {\widetilde F_L   }) \Big/\GL_{n=n_1+\cdots+n_s}( (\ZZ\big/N\ \ZZ)^2 )\\
&=
\Reps(\GL_{n=n_1+\cdots+n_s}({\wt F_{\o \omega   }}\times {\wt F_\omega   })\end{align*}
be the partially reducible lattice bisemispace. 

Then, we have that the bilinear cohomology of 
$X^{2n=2n_1+\cdots+2n_s}_{\SRL}$ decomposes according to:
\[H^*(X^{2n=2n_1+\cdots+2n_s}_{\SRL})
=\txt\bigoplus\limits^{2n_s}_{2n_\ell=2n_1}H^{2n_\ell}(
X^{2n}_{\SRL},\widetilde M^{2n_\ell}_{R}\otimes \widetilde M^{2n_\ell}_{L})\] 
where
$\widetilde M^{2n_\ell}_{L}$ (resp. $\widetilde M^{2n_\ell}_{R})$~)  is the left (resp. right) $T_{n_\ell} ({\wt F_\omega}$)-subsemimodule $M^{2n_\ell}_L$ (resp. $T^t_{n_\ell} ({\wt F_{\o\omega}}$)-subsemimodule $M^{2n_\ell}_R$~)  of complex dimension $n_\ell$~.
}\vskip 11pt

\bpr this results from the partial reducibility of $X^{2n=2n_1+\cdots+2n_s}_{\SRL}$ and from the equality:
\[H^{2n_\ell}(
X^{2n}_{\SRL},\wt M^{2n_\ell}_{R}\otimes \wt M^{2n_\ell}_{L})=
\Reps(\GL_{n_\ell}({\wt F_{\o \omega   }}\times {\wt F_\omega   }))  \]
if we refer to section 3.2.\epr\vskip 11pt

\subsubsection{Proposition}  

{\em Let
\[
X^{2n=2_1+\cdots+2_n}_{\SRL}= 
\GL_{2n=2_1+\cdots+2_n}({\wt F_{\o \omega   }}\times {\wt F_\omega   }) \Big/\GL_{2n=2_1+\cdots+2_n}( (\ZZ\big/N\ \ZZ)^2 )\]
be the orthogonal completely reducible lattice bisemispace and let $X_{\SRL}^{2n_R\times 2n_L}$
 be its
nonorthogonal correspondent. Then, we have the following decompositions of the bilinear 
cohomologies:
\Bean
\item $H^{2n}(
X^{2n}_{\SRL},\widetilde M^{2n}_{R}\otimes \widetilde M^{2n}_{L})= \bigoplus\limits_{2_\ell}
H^{2}(
X^{2n}_{\SRL},\widetilde M^{2_{\ell_R}}_{R}\otimes \widetilde M^{2_{\ell_L}}_{L})$~;\vskip 11pt

\item $H^{2n}(
X^{2n_{R\times L}}_{\SRL},\widetilde M^{2n_R}_{R}\otimes \widetilde M^{2n_L}_{L})= \bigoplus\limits_{2_\ell}
H^{2}(
X^{2_{{n_{R\times L}}}}_{\SRL},\widetilde M^{2_{\ell_R}}_{R}\otimes \widetilde M^{2_{\ell_L}}_{L})$

\hspace{2cm} $\bigoplus\limits_{2_{k_R}, 2_{\ell_L}}
H^{2}(
X^{2_{k_R-\ell_L}}_{\SRL},\widetilde M^{2_{k_R}}_{R}\otimes \widetilde M^{2_{\ell_L}}_{L})$

where
\Bi
\item[\textbullet]
$ \widetilde M^{2n}_{R}\otimes \widetilde M^{2n}_{L}= \bigoplus\limits_{2_{\ell_R}= 2_{\ell_L}}
(\widetilde M^{2_{\ell_R}}_{R}\otimes \widetilde M^{2_{\ell_L}}_{L})$ \quad  is   the completely reducible orthogonal
$\GL_{2n=2_1+\cdots+2_n}({\wt F_{\o \omega   }}\times {\wt F_\omega   }) $-bisemimodule $\wt M^{2n}_R\otimes \wt M^{2n}_L$~;

\item[\textbullet]
$ \widetilde M^{2n_R}_{R}\otimes \widetilde 2M^{n_L}_{L}= \bigoplus\limits_{2_{\ell_R}= 2_{\ell_L}}
(\widetilde M^{2_{\ell_R}}_{R}\otimes \widetilde M^{2_{\ell_L}}_{L})
\bigoplus\limits_{2_{k_R}\neq 2_{\ell_L}}
(\widetilde M^{2_{k_R}}_{R}\otimes \widetilde M^{2_{\ell_L}}_{L}) $ \quad  is   the completely reducible nonorthogonal
$\GL_{2n_{R\times L}}({\wt F_{\o \omega   }}\times {\wt F_\omega   }) $-bisemimodule $\wt M^{2n_R}_R\otimes \wt M^{2n_L}_L$~.
\epr
\Ei\Ee}\vskip 11pt

\subsubsection{Proposition}

 {\em \parindent=0pt
Let $\CY^{2n_\ell}(X_{R})$ (resp. $\CY^{2n_\ell}(X_{L})$~) be an algebraic semicycle of complex codimension $n_\ell$ over 
a \rl semischeme $ X_{R}$ (resp. $ X_{L}$~) isomorphic to
the pseudo-ramified lattice semispace $X^{2n}_{S_R}$ (resp. $X^{2n}_{S_L}$~).\par
Let
\begin{align*}
\CY^{2n_\ell}(X_{R}) \times \CY^{2n_\ell}(X_{L})
&= \txt\bigoplus\limits_{j} \txt\bigoplus\limits_{m_j}
(\CY^{2n_\ell}(X_{R}[j,m_j]) \times \CY^{2n_\ell}(X_{L}[j,m_j]))\\
&\simeq \Reps(\GL_{n_\ell}({\wt F_{\o \omega _\oplus  }}\times {\wt F_{\omega_\oplus}   }))
\end{align*}
 be the product, right by left, of these algebraic semicycles decomposed according to their 
conjugacy classes (see section 1.3) in such a way that
$\CY^{2n_\ell}(X_{R})\times \CY^{2n_\ell}(X_{L})$ is isomorphic to an algebraic bilinear subsemigroup $ G^{(2n_\ell)}( {\wt F_{\o\omega_\oplus }}\times {\wt F_{\omega _\oplus}})$~.\par
Then, we have the following decompositions of the bilinear cohomologies:
\Bi
\item $H^{*}(
X^{2n=2n_1+\cdots+2n_s}_{\SRL},\widetilde M^{2n}_{R_\oplus}\otimes \widetilde M^{2n}_{L_\oplus})= 
\bigoplus\limits_{2n_\ell=n_1}^{2n_s}
H^{2n_\ell}(
X^{2n}_{\SRL},\widetilde M^{2n_\ell}_{R_\oplus}\otimes \widetilde M^{2n_\ell}_{L_\oplus})$

\hspace*{2cm} $\simeq  \txt\bigoplus\limits_{n_\ell} \txt\bigoplus\limits_{j} \txt\bigoplus\limits_{m_j}
(\CY^{2n_\ell}(X_{R}[j,m_j]) \times \CY^{2n_\ell}(X_{L}[j,m_j]))$\vskip 11pt

\item $H^{2n}(
X^{2n}_{\SRL},\widetilde M^{2n}_{R_\oplus}\otimes \widetilde M^{2n}_{L_\oplus})= 
\bigoplus\limits_{2_\ell}
H^{2}(
X^{2n}_{\SRL},\widetilde M^{2_{\ell_R}}_{R_\oplus}\otimes M^{2_{\ell_L}}_{L_\oplus})$

\hspace*{2cm} $\simeq \txt\bigoplus\limits_{2_\ell} \txt\bigoplus\limits_{j} \txt\bigoplus\limits_{m_j}
(\CY^{2_\ell}(X_{R}[j,m_j]) \times \CY^{2_\ell}(X_{L}[j,m_j]))$\vskip 11pt

\item $H^{2n} (
X^{2n_{R\times L}}_{\SRL},\widetilde M^{2n_R}_{R_\oplus}\otimes \widetilde M^{2n_L}_{L_\oplus})= 
\bigoplus\limits_{2_\ell}
H^{2}(
X^{2{n_{R\times L}}}_{\SRL},\widetilde M^{2_{\ell_R}}_{R_\oplus}\otimes \widetilde M^{2_{\ell_L}}_{L_\oplus})$

\hspace*{4cm} $\bigoplus\limits_{2_{k_R}\neq 2_{\ell_L}}
H^{2}(
X^{2{n\RL}}_{\SRL},\widetilde M^{2_{k_R}}_{R_\oplus}\otimes \widetilde M^{2_{\ell_L}}_{L_\oplus})$

\hspace*{2cm} $\simeq  \txt\bigoplus\limits_{2_{\ell_{R\times L}}} \txt\bigoplus\limits_{j} \txt\bigoplus\limits_{m_j}
(\CY^{2_\ell}(X_{R}[j,m_j]) \times \CY^{2_\ell}(X_{L}[j,m_j]))$

\hspace*{2cm} \quad $ \txt\bigoplus\limits_{2_{k_R}\neq 2_{\ell_L}} 
\txt\bigoplus\limits_{j_R}  \txt\bigoplus\limits_{m_{j_R}}
\txt\bigoplus\limits_{j_L}  \txt\bigoplus\limits_{m_{j_L}}
(\CY^{2_{k_R}}(X_{R}[j_R,m_{j_R}]) \times \CY^{2_{\ell_L}}(X_{L}[j_L,m_{j_L}]))$
\Ei
being equal to the sums of products, right by left, of the conjugacy classes, counted with\linebreak  their multiplicities, of the irreducible semicycles.
}\vskip 11pt

{\noindent \bf Sketch of proof:} indeed, 
$H^{2n_\ell}(X^{2n}_{\SRL},\widetilde M^{2n_\ell}_{R_\oplus}\otimes \widetilde M^{2n_\ell}_{L_\oplus})=
\Reps(\GL_{n_\ell}({\wt F_{\o \omega _\oplus  }}\times {\wt F_{\omega _\oplus}  }))
=\bigoplus_j\bigoplus_{m_j}(\widetilde M^{2n_\ell}_{\o\omega_{j,m_j}}
\otimes \widetilde M^{2n_\ell}_{\omega_{j,m_j}})$ where 
$\widetilde M^{2n_\ell}_{\o\omega_{j,m_j}}$ (resp. 
$\widetilde M^{2n_\ell}_{\omega_{j,m_j}}$~) is  
the conjugacy class representative $g^{2n_\ell}_{R}[j,m_j]$ 
(resp. $g^{2n_\ell}_{L}[j,m_j]$~) $\in G^{2n_\ell}(\wt F_{\o\omega}\times \wt F_{\omega})$~.
\epr\vskip 11pt

\subsection[Langlands global correspondences on the reducible compactified bisemispaces $\o X^{(2)n}_{\SRL}$ and $\partial\o X^{(2)n}_{\SRL}$]{{\boldmath Langlands global correspondences on the reducible\protect\newline  compactified bisemispaces $\o X^{(2)n}_{\SRL}$ and $\partial \o X^{(2)n}_{\SRL}$}}

\subsubsection{Definition: reducible compactified-lattice bisemispaces} 

Let
\[ X^{2n=2n_1+\cdots+2n_s}_{\SRL}\; , \quad X^{2n=2_1+\cdots+2_n}_{\SRL}\quad \text{and} \quad
X^{2n_R\times 2n_L}_{\SRL}\]
be respectively the partially reducible, orthogonal completely reducible and
nonorthogonal completely reducible lattice bisemispaces as introduced in lemma 4.1.5 and proposition 4.1.6. Then, we define the projective emergent isomorphisms of toroidal compactifications:
\begin{align*}
&\gamma^{c(2n)}_{R\times L} : \quad X^{2n=2n_1+\cdots+2n_s}_{\SRL}\\
&\quad \To\o X^{2n=2n_1+\cdots+2n_s}_{\SRL}
=\GL_{n=n_1+\cdots+n_s}({F^T_{R }}\times {F^T_L })\Big/ \GL_{n=n_1+\cdots+n_s}( (\ZZ\big/N\ \ZZ)^2 )\;, \\[11pt]
%
&\gamma^{c(2n)}_{R\times L} : \quad  X^{2n=2_1+\cdots+2_n}_{\SRL}\\
&\quad \To\o X^{2n=2_1+\cdots+2_n}_{\SRL}
=\GL_{2n=2_1+\cdots+2_n}({F^T_{R }}\times {F^T_L })\Big/ \GL_{2n=2_1+\cdots+2_n}( (\ZZ\big/N\ \ZZ)^2)\;, \\[11pt]
&\gamma^{c(2n_R\times 2n_L)}_{R\times L} : \quad  X^{2n_R\times 2n_L}_{\SRL}\\
&\quad \To\o X^{2n_R\times 2n_L}_{\SRL}
=\GL_{2n_R\times 2n_L}({F^T_{R }}\times 
{F^T_L })\Big/ \GL_{2n_R\times 2n_L}( (\ZZ\big/N\ \ZZ)^2)\;, \end{align*}
mapping these reducible lattice bisemispaces into their toroidal compactified equivalents in the sense of proposition 3.3.3 with the evident notations:
\[\o X^{n=n_1+\cdots+n_s}_{\SRL}=\mathop{\boxplus}\limits^{n_s}_{n_\ell=n_1}\o X^{n_\ell}_{\SRL}\;.\]\vskip 11pt

\subsubsection{Definition} 

The double coset decomposition of the reducible bilinear general semigroups are given by:
\begin{align*}
&\o S^{P_{n=n_1+\cdots+n_s}}_{K_{n=n_1+\cdots+n_s}} \\
& \quad =P_{n=n_1+\cdots+n_s}( {F^{T}_{\o\omega^1 }}\times {F^{T}_{\omega^1} } )\setminus \GL_{n=n_1+\cdots+n_s}
({F^T_{R }}\times {F^T_L })\Big/ \GL_{n=n_1+\cdots+n_s}((\ZZ\big/N\ \ZZ)^2)\;, \\[11pt]
&\o S^{P_{2n=2_1+\cdots+2_n}}_{K_{2n=2_1+\cdots+2_n}} \\
& \quad =P_{2n=2_1+\cdots+2_n}( {F^{T}_{\o\omega^1 }}\times {F^{T}_{\omega^1} } )\setminus \GL_{2n=2_1+\cdots+2_n}
({F^T_R  }\times {F^T_L})\Big/ \GL _{2n=2_1+\cdots+2_n}( (\ZZ\big/N\ \ZZ)^2)\;, \\[11pt]
&\o S^{P_{2n_R\times 2n_L}}_{K_{2n_R\times 2n_L}} \\
& \quad =P_{2n_R\times 2n_L}( {F^{T}_{\o\omega^1 }}\times {F^{T}_{\omega ^1}})\setminus \GL_{2n_R\times 2n_L}
({F^T_{R  }}\times {F^T_L   })\Big/ \GL_{2n_R\times 2n_L}( (\ZZ\big/N\ \ZZ)^2 )\;, \end{align*}
with the notations of definition 3.3.4.
\vskip 11pt

\subsubsection{Definition: induction} 

If the double coset decompositions of the reducible bilinear general semigroups are considered, there are evident inductions respectively from the representations of the complex bilinear reducible parabolic subgroups:
\[
P_{n=n_1+\cdots+n_s}( {F^{T}_{\o\omega^1 }}\times {F^{T}_{\omega^1} } )\;, \quad 
P_{2n=2_1+\cdots+2_n}( {F^{T}_{\o\omega^1 }}\times {F^{T}_{\omega^1} })\quad \text{and}\quad 
P_{2n_R\times 2n_L}( {F^{T}_{\o\omega^1 }}\times {F^{T}_{\omega^1} })\]
and from the respective representations of the bilinear reducible arithmetic subgroups:
\[\GL_{n=n_1+\cdots+n_s}( (\ZZ\big/N\ \ZZ)^2 )\;, \quad 
\GL_{2n=2_1+\cdots+2_n}( (\ZZ\big/N\ \ZZ)^2 )\quad \text{and}\quad 
\GL_{2n_R\times 2n_L}((\ZZ\big/N\ \ZZ)^2 )\]
to the respective representations of the reducible bilinear general semigroups: 
\[\GL_{n=n_1+\cdots+n_s} \;, \quad
\GL_{2n=2_1+\cdots+2_n}\quad \text{and} \quad
\GL_{2n_R\times 2n_L}\;.\]
\vskip 11pt

\subsubsection{Proposition} 

{\em\parindent=0pt
Let
\begin{align*}
\gamma^{c(2n)}_{R\times L} : \quad & \CY^{2n=2n_1+\cdots+2n_s}(X_{R})\times \CY^{2n=2n_1+\cdots+2n_s}
(X_{L}) =
\txt\bigoplus\limits^{2n_s}_{2n_\ell=2n_1}\CY^{2n_\ell}(X_{R})\times \CY^{2n_\ell}(X_{L}) 
 \\
&\quad \To
\CY^{2n}_T(\o X_{R})\times \CY_T^{2n}(\o X_{L})
= \txt\bigoplus\limits_{n_\ell}\CY^{2n_\ell}_T(\o X_{R})\times \CY^{2n_\ell}_T(\o X_{L})
\end{align*}
be the projective emergent isomorphism of compactifications applied to 
a product, right by left, of $2n$-dimensional partially reducible semicycles such that $\CY_T^{2n_\ell} (\o X_{L})$ (resp.\linebreak $\CY_T^{2n_\ell} (\o X_{R})$~) be a $n_\ell$-dimensional complex semitorus $T^{2n_\ell}_L$ (resp. $T^{2n_\ell}_R$~). 

Similarly, let
\begin{align*}
\gamma^{c(2n_{R\times L})}_{R\times L} : \quad & 
CY^{2n=2_1+\cdots+2_n}(X_{R})\times \CY^{2n=2_1+\cdots+2_n}(X_{L}) 
 \\
&\quad \To
\CY^{2n=2_1+\cdots+2_n}_T(\o X_{R})\times \CY_T^{2n=2_1+\cdots+2_n}(\o X_{L})
\end{align*}
and
\begin{align*}
\gamma^{c(2n_R\times 2n_L)}_{R\times L} : \quad & 
CY^{2n_R}(X_{R})\times \CY^{2n_L}(X_{L}) 
 \\
&\quad \To
\CY^{2n_R}_T(\o X_{R})\times \CY_T^{2n_L}(\o X_{L})
\end{align*}
be the two other considered projective emergent isomorphisms applied to products of completely reducible semicycles. 

Then, we have the evident decompositions of the cohomologies:
\Bi
\item $H^{*}(  \o S 
^{P_{n=n_1+\cdots+n_s}} _{K_{n=n_1+\cdots+n_s}} , \widehat M^{2n}_{T_{R_\oplus }}\otimes \widehat M^{2n}_{T_{L_\oplus }}) = \bigoplus\limits_{2n_\ell} H^{2n_\ell}(  \o S 
^{P_{n}} _{K_{n}} ,\widehat M^{2n_\ell}_{T_{R_\oplus }}\otimes \widehat M^{2n_\ell}_{T_{L_\oplus }})$

\hspace{1cm} $\simeq 
\bigoplus\limits_{n_\ell} \bigoplus\limits_{j}\bigoplus\limits_{m_j}
(\CY^{2n_\ell}_T(\o X_{R}[j,m_j])\times \CY_T^{2n_\ell}(\o X_{L}[j,m_j]))$\vskip 11pt

\item $H^{2n}(  \o S 
^{P_{2n=2_1+\cdots+2_n}} _{K_{2n=2_1+\cdots+2_n}} ,\widehat M^{2n}_{T_{R_\oplus }}\otimes \widehat M^{2n}_{T_{L_\oplus }})  = \bigoplus\limits_{2\ell_R=2\ell_L} H^{2}
(  \o S 
^{P_{2{n}}}_{K_{2{n}}} ,\widehat M^{2_{\ell_R}}_{T_{R_\oplus }}\otimes 
\widehat M^{2_{\ell_L}}_{T_{L_\oplus }})$

\hspace{1cm} $\simeq 
\bigoplus\limits_{2\ell_R=2\ell_L} \bigoplus\limits_{j}\bigoplus\limits_{m_j}
(\CY^{2_{\ell_R}}_T(\o X_{R}[j,m_j])\times \CY_T^{2_{\ell_L}}(\o X_{L}[j,m_j]))$\vskip 11pt

\item $H^{2n}(  \o S 
^{P_{2n_R\times 2n_L}} _{K_{2n_R\times 2n_L}} ,\widehat M^{2n_R}_{T_{R_\oplus }}\otimes \widehat M^{2n_L}_{T_{L_\oplus }})  = \bigoplus\limits_{2\ell_R=2\ell_L} H^{2}
(  \o S 
^{P_{2n_R\times 2n_L}}_{K_{2n_R\times 2n_L}} ,\widehat M^{2_{\ell_R}}_{T_{R_\oplus }}\otimes \widehat M^{2_{\ell_L}}_{T_{L_\oplus }})$

\hspace{2cm} $ \bigoplus\limits_{k_R\neq\ell_L} H^{2}
(  \o S 
^{P_{2n_R\times 2\ell _L}}_{K_{2n_R\times 2\ell _L}} ,\widehat M^{2_{k_R}}_{T_{R_\oplus }}\otimes \widehat M^{2_{\ell_L}}_{T_{L_\oplus }})$

\hspace{1cm} $\simeq 
\bigoplus\limits_{\ell_R=\ell_L} \bigoplus\limits_{j}\bigoplus\limits_{m_j}
(\CY^{2_{\ell_R}}_T(\o X_{R}[j,m_j])\times \CY_T^{2_{\ell_L}}(\o X_{L}[j,m_j]))$

\hspace{2cm} $ 
\bigoplus\limits_{k_R\neq\ell_L} \bigoplus\limits_{j_R}\bigoplus\limits_{m_{j_R}}
\bigoplus\limits_{j_L}\bigoplus\limits_{m_{j_L}}
(\CY^{2_{k_R}}_T(\o X_{R}[j_R,m_{j_R}])\times \CY_T^{2_{\ell_L}}(\o X_{L}[j_L,m_{j_L}]))$
\Ei
being equal to the sums of products, right by left, of the equivalence classes, counted with
their multiplicities, of the irreducible semicycles.
}

\vskip 11pt

\subsubsection{Proposition} 

 {\em 
\Bi
\item If \quad $z_{n_\ell}=\sum\limits^{2n_\ell}_{\alpha =1} z_{n_\ell}(\alpha )\ |\vec e_\alpha|\in \CC^{n_\ell}$~, then every left (resp. right) $n_\ell$-dimensional
complex semitorus has the analytic development:
\begin{align*}
T^{2n_\ell}_L[j_L,m_{j_L}] &\simeq
\lambda^{\half} (2n_\ell,j_L,m_{j_L}) \ e^{2\pi ij_Lz_{n_\ell}}\\
(\text{resp.} \quad 
T^{2n_\ell}_R[j_R,m_{j_R}] &\simeq 
\lambda^{\half}(2n_\ell,j_R,m_{j_R})\ e^{-2\pi ij_Rz_{n_\ell}}\ ).\end{align*}

\item On the other hand, let
\begin{align*}
\Eis_L (2n_\ell,j_L,m_{j_L}) &=
\txt\sum\limits^r_{j_L=1} \sum\limits_{m_{j_L}} \lambda^{\half} (2n_\ell,j_L,m_{j_L}) \ e^{2\pi ij_Lz_{n_\ell}}\\
(\text{resp.} \quad 
\Eis_R (2n_\ell,j_R,m_{j_R}) &=
\txt\sum\limits^r_{j_R=1} \sum\limits_{m_{j_R}} \lambda^{\half} (2n_\ell,j_R,m_{j_R}) \ e^{-2\pi ij_Lz_{n_\ell}}\ ).\end{align*}
be the (truncated) Fourier development of a $2n_\ell$-dimensional left (resp. right) cusp form of weight $k = 2$ restricted to the upper (resp. lower) half space.
\Ei

Then, we have the following analytical developments of the bilinear cohomologies:
 \Bi
\item $H^{*}(  \o S 
^{P_{n=n_1+\cdots+n_s}} _{K_{n=n_1+\cdots+n_s}} ,\widehat M^{2n}_{T_{R_\oplus }}\otimes \widehat M^{2n}_{T_{L_\oplus }}) =
 \bigoplus\limits_{n_\ell} (\Eis_R (2n_\ell,j_R,m_{j_R})\otimes \Eis_L (2n_\ell,j_L,m_{j_L})) $~;\vskip 11pt

\item $H^{2n}(  \o S 
^{P_{2n=2_1+\cdots+2_n}} _{K_{2n=2_1+\cdots+2_n}} ,\widehat M^{2n}_{T_{R_\oplus }}\otimes \widehat M^{2n}_{T_{L_\oplus }})  
=\bigoplus\limits_{\ell_R=\ell_L} (\Eis_R (2_{\ell_R},j_R,m_{j_R})\otimes \Eis_L 
(2_{\ell_L},j_L,m_{j_L})) $~;\vskip 11pt

\item $H^{2n}(  \o S 
^{P_{2n_R\times 2n_L}} _{K_{2n_R\times 2n_L}} ,\widehat M^{2n_R}_{T_{R_\oplus }}\otimes \widehat M^{2n_L}_{T_{L_\oplus }})  
=
 \bigoplus\limits_{\ell_R=\ell_L} (\Eis_R (2_{\ell_R},j_R,m_{j_R})\otimes \Eis_L (2_{\ell_L},j_L,m_{j_L})) $

\hspace{2cm} $  \bigoplus\limits_{k_R\neq\ell_L} (\Eis_R (2_{k_R},j_R,m_{j_R})\otimes 
\Eis_L (2_{\ell_L},j_L,m_{j_L})) $
\Ei
where $\Eis_L (2_{\ell_L},j_L,m_{j_L})$ corresponds to the (truncated) Fourier development of a cusp form of weight $k = 2$~.}
\vskip 11pt


\subsubsection{Reducible cusp biforms}

\Bi
\item Let $L^{1-1}_{\rm cusp}(G^{(2n=2n_1+\cdots+2n_s)}( {F^T_{\o \omega }}\times {F^T_{\o \omega }}))$ denote the bialgebra of cusp biforms on the partially reducible pseudo-ramified toroidal bilinear  semigroup 
$G^{(2n=2n_1+\cdots+2n_s)}( {F^T_{\o \omega }}\times {F^T_{\o \omega }})$~.\pagebreak

Referring to section 3.4.10, a partially reducible cusp biform on this ``diagonal'' bialgebra is given by:
\begin{align*}
& \phi ^{(2n=2n_1+\cdots+2n_s)}_{G^T_R} (x_{g^T_R}) \otimes_D \phi ^{(2n=2n_1+\cdots+2n_s)}_{G^T_L}(x_{g^T_L})\\[11pt]
& \quad = \bigoplus^{n_s}_{n_\ell =n_1} \phi ^{(2n_\ell )}_{G^T_R}\otimes_D \phi ^{(2n_\ell )}_{G^T_L}\\[11pt]
& \quad = \bigoplus^{n_s}_{n_\ell =n_1} \bigoplus^r_{j=1}\bigoplus_{m_j}
\bigl(\phi ^{(2n_\ell )}_{G^T_{j_R}} (x_{g^T_{j_R}})\otimes_D \phi ^{(2n_\ell )}_{G^T_{j_L}} (x_{g^T_{j_L}})\bigr)\\[11pt]
& \quad = \bigoplus_{n_\ell } ( \Eis_R(2n_\ell ,j_R,m_{j_R}) \otimes _{(D)}\Eis_L(2n_\ell ,j_L,m_{j_L}) )\;.\end{align*}

\item Similarly, let
$L^{1-1}_{\rm cusp}(G^{(2n=2_1+\cdots+2_n)}( {F^T_{\o \omega }}\times {F^T_{\o \omega }}))$ denote the bialgebra of orthogonally completely reducible 
cusp biforms on $G^{(2n=2_1+\cdots+2_n)}( {F^T_{\o \omega }}\times {F^T_{\o \omega }})$
given by:
\begin{align*}
& \phi ^{(2n=2_1+\cdots+2_n)}_{G^T_R} (x_{g^T_R}) \otimes_D \phi ^{(2n=2_1+\cdots+2_n)}_{G^T_L}(x_{g^T_L})\\[11pt]
& \quad = \bigoplus_{\ell _R=\ell _L} \phi ^{2_{\ell _R}} _{G^T_R}\otimes_D \phi ^{2_{\ell _L}} _{G^T_L}  \\[11pt]
& \quad = \bigoplus_{\ell _R=\ell _L}( \Eis_R(2_{\ell _R},j_R,m_{j_R})\otimes_{(D)}  \Eis_L(2_{\ell _L},j_L,m_{j_L}))
\;.\end{align*}

\item And, let
$L^{1-1}_{\rm cusp}(G^{(2n_R\times 2n_L)}( {F^T_{\o \omega }}\times {F^T_{\omega }}))$ be the bialgebra of nonorthogonally completely reducible 
cusp biforms on $G^{(2n_R\times 2n_L)}( {F^T_{\o \omega }}\times {F^T_{\omega }})$
given diagonally by:
\begin{align*}
& \phi ^{(2n_R=2_{1_R}+\cdots+2_{n_R})}_{G^T_R} (x_{g^T_R}) \otimes_D \phi ^{(2n_L=2_{1_L}+\cdots+2_{n_L})}_{G^T_L}(x_{g^T_L})\\[11pt]
& \quad = \bigoplus_{\ell _R=\ell _L}( \phi ^{2_{\ell _R}} _{G^T_R}\otimes_D \phi ^{2_{\ell _L}} _{G^T_L} )
\bigoplus_{k_R\neq \ell _L}( \phi ^{2_{k_R}} _{G^T_R}\otimes_D \phi ^{2_{\ell _L}} _{G^T_L}  )
 \\[11pt]
& \quad = \bigoplus_{\ell _R=\ell _L}( \Eis_R(2_{\ell _R},j_R,m_{j_R})\otimes_{(D)}  \Eis_L(2_{\ell _L},j_L,m_{j_L}))\\[11pt]
& \qquad \qquad 
\bigoplus_{k_R\neq\ell _L}( \Eis_R(2_{k_R},j_R,m_{j_R})\otimes_{(D)}  \Eis_L(2_{\ell _L},j_L,m_{j_L}))
\;.\end{align*}
\Ei
\vskip 11pt 

\subsubsection{Proposition}

{\em 
\Bi
\item Let 
\begin{gather*}
R_{G^{(2n=2n_1+\cdots+2n_s)}({F^T_{\o\omega }}\times {F^T_\omega })}
( \phi ^{(2n=2n_1+\cdots+2n_s)}_{G^T_R} \otimes \phi ^{(2n=2n_1+\cdots+2n_s)}_{G^T_L} )\;, \\[11pt]
R_{G^{(2n=2_1+\cdots+2_n)}({F^T_{\o\omega }}\times {F^T_\omega })}
( \phi ^{(2n=2_1+\cdots+2_n)}_{G^T_R} \otimes \phi ^{(2n=2_1+\cdots+2_n)}_{G^T_L} )\;,\\[11pt]
R_{G^{(2n_R\times 2n_L)}({F^T_{\o\omega }}\times {F^T_\omega })}
( \phi ^{2n_R}_{G^T_R} \otimes \phi ^{2n_L}_{G^T_L} )\end{gather*}
be the integral operators on the respective cusp biforms introduced in sections 4.2.6 and 3.4.10.
\vskip 11pt 

\item Then, the trace formulas of these integral operators are respectively given by:
\begin{align*}
& \Tr( R_{G^{(2n=2n_1+\cdots+2n_s)}({F^T_{\o\omega }}\times {F^T_\omega })}
( \phi ^{(2n=2n_1+\cdots+2n_s)}_{G^T_R} \otimes \phi ^{(2n=2n_1+\cdots+2n_s)}_{G^T_L} )\\[11pt]
& \quad = {\rm vol}(G^{(2n=2n_1+\cdots+2n_s)}({F^T_{\o\omega }}\times {F^T_\omega }))
\bigoplus_{n_\ell } ( \Eis_R(2n_\ell ,j_R,m_{j_R}),  \Eis_L(2n_\ell ,j_L,m_{j_L}))\\[11pt]
& \quad = {\rm vol}(G^{(2n=2n_1+\cdots+2n_s)}({F^T_{\o\omega }}\times {F^T_\omega }))
\bigoplus_{n_\ell } \bigoplus_{j} \bigoplus_{m_j}\lambda (2n_\ell ,j,m_j)\end{align*}
where $(\Eis_R(2n_\ell ,j_R,m_{j_R}),  \Eis_L(2n_\ell ,j_L,m_{j_L}))$ is a diagonal bilinear form from\linebreak
$\Red\cusp (\GL_{(2n=2n_1+\cdots+2n_s)}({F^T_{\o\omega }}\times {F^T_\omega }))$ to $\CC$~.

\begin{align*}
& \Tr( R_{G^{(2n=2_1+\cdots+2_n)}({F^T_{\o\omega }}\times {F^T_\omega })}
( \phi ^{(2n=2_1+\cdots+2_n)}_{G^T_R} \otimes \phi ^{(2n=2_1+\cdots+2_n)}_{G^T_L} )\\[11pt]
& \quad = {\rm vol}(G^{(2n=2_1+\cdots+2_n)}({F^T_{\o\omega }}\times {F^T_\omega }))
\bigoplus_{\ell_R=\ell _L } ( \Eis_R(2_{\ell _R},j_R,m_{j_R}),  \Eis_L(2_{\ell _L},j_L,m_{j_L}))\\[11pt]
& \quad = {\rm vol}(G^{(2n=2_1+\cdots+2_n)}({F^T_{\o\omega }}\times {F^T_\omega }))
\bigoplus_{\ell_R=\ell _L } \bigoplus_{j} \bigoplus_{m_j}\lambda (2_\ell ,j,m_j)\end{align*}
and finally

\begin{align*}
& \Tr( R_{G^{(2n_R\times 2n_L)}({F^T_{\o\omega }}\times {F^T_\omega })}
( \phi ^{(2n_R)}_{G^T_R} \otimes \phi ^{(2n_L)}_{G^T_L} )\\[11pt]
& \quad = {\rm vol}(G^{(2n_R\times 2n_L)}({F^T_{\o\omega }}\times {F^T_\omega }))
\L[
\bigoplus_{\ell_R=\ell _L } ( \Eis_R(2_{\ell _R},j_R,m_{j_R}),  \Eis_L(2_{\ell _L},j_L,m_{j_L}))
\R.\\[11pt]
& \qquad \quad 
+\L.
\bigoplus_{k_R=\ell _L } ( \Eis_R(2_{k_R},j_R,m_{j_R}),  \Eis_L(2_{\ell _L},j_L,m_{j_L}))
\R]
\\[11pt]
& \quad = {\rm vol}(G^{(2n_R\times 2n_L)}({F^T_{\o\omega }}\times {F^T_\omega }))
\L(
\bigoplus_{\ell_R=\ell _L } \bigoplus_{j} \bigoplus_{m_j}\lambda (2_\ell ,j,m_j)\R.
\\[11pt]
& \qquad \quad + \L.\bigoplus_{k_R\neq\ell _L }
\bigoplus_{j_R} \bigoplus_{j_L} \bigoplus_{m_{j_R}} \bigoplus_{m_{j_L}} 
 \lambda^{\half} (2_{k_R},j_R,m_{j_R}) \cdot 
\lambda^{\half} (2_{\ell _L},j_L,m_{j_L}) \R)
\end{align*}
\Ei
}\vskip 11pt 

\bpr 
\Be
\item The proof is clear from the preceding developments and, especially from proposition 3.4.10.  

\item These trace formulas are complete in the sense that they refer to a decomposition of these traces according to the direct sums of the irreducible representations of the bilinear algebraic semigroups, as in the classical case for the linear  groups.
\epr\Ee
\vskip 11pt


\subsubsection{Proposition: Langlands  global reducible correspondences} 

 {\em
\Bi
\item Let
\begin{align*}
\Red\Rep W^{2n=2n_1+\cdots+2n_s}_{\FRL} : \quad
& W^{2n=2n_1+\cdots+2n_s}_{F_{\o\omega}} \times W^{2n=2n_1+\cdots+2n_s}_{F_{\omega}}\\
& \quad \To \GL_{n=n_1+\cdots+n_s}( {\wt F_{\o\omega_\oplus }}\times {\wt F_{\omega_\oplus }})\; ,\\
\Red\Rep W^{2n=2_1+\cdots+2_n}_{\FRL} : \quad
& W^{2n=2_1+\cdots+2_n}_{F_{\o\omega}} \times W^{2n=2_1+\cdots+2_n}_{F_{\omega}}\\
& \quad \To \GL_{2n=2_1+\cdots+2_n}( {\wt F_{\o\omega_\oplus }}\times {\wt F_{\omega_\oplus }})\; ,\\
\Red\Rep W^{2n_R\times 2n_L}_{\FRL} : \quad
& W^{2n_R}_{F_{\o\omega}} \times W^{2n_L}_{F_{\omega}}\\
& \quad \To \GL_{2n_R\times 2n_L}( {\wt F_{\o\omega_\oplus }}\times {\wt F_{\omega_\oplus }})\end{align*}
be respectively the partially reducible, orthogonal completely reducible and nonorthogonal completely reducible representations of the bilinear global Weil groups
\[W^{2n=2n_1+\cdots+2n_s}_{F_{\o\omega}} \times W^{2n=2n_1+\cdots+2n_s}_{F_{\omega}}\;, \quad
W^{2n=2_1+\cdots+2_n}_{F_{\o\omega}} \times W^{2n=2_1+\cdots+2_n}_{F_{\omega}}\]
\[ \text{and}
\quad W^{2n_R}_{F_{\o\omega}} \times W^{ 2n_L}_{F_{\omega}}\;.\]
\vskip 11pt


\item Let
\begin{align*}
&\Red  \cusp (\GL_{n=n_1+\cdots+n_s} ({F^{T}_{\o \omega }} \times {F^{T}_{\omega }}))\\
& \qquad \qquad \qquad =\txt\bigoplus\limits_{n_\ell }( \Eis_R(2n_\ell,j_R,m_{j_R}) \otimes \Eis_L(2n_\ell,j_L,m_{j_L}) )\;,\\
&\Red  \cusp (\GL_{2n=2_1+\cdots+2_n} ({F^{T}_{\o \omega }} \times {F^{T}_{\omega }})\\
& \qquad \qquad \qquad =\txt\bigoplus\limits_{\ell_R=\ell_L }( \Eis_R(2_{\ell_R},j _R,m_{j_R}) \otimes \Eis_L(2_{\ell_L},j_L,m_{j_L}) )\;,\\
 &\Red  \cusp (\GL_{2n_R\times 2n_L} ({F^{T}_{\o \omega }} \times {F^{T}_{\omega }})\\
& \qquad \qquad \qquad =\txt\bigoplus\limits_{\ell_R=\ell_L }( \Eis_R(2_{\ell_R},j _R,m_{j_R}) \otimes \Eis_L(2_{\ell_L},j_L,m_{j _L}) )\\
& \qquad \qquad \qquad \qquad \txt\bigoplus\limits_{k_R\neq\ell_L }( \Eis_R(2_{k_R},j _R,m_{j _R}) \otimes \Eis_L(2_{\ell_L},j_L,m_{j_L}) )\;,\end{align*}
be respectively the sums of the products, right by left, of the equivalence classes of the partially reducible, orthogonal completely reducible and nonorthogonal completely reducible cuspidal representations of the corresponding reducible bilinear semigroups.\vskip 11pt

\item Let
\[
\CY^{2n=2n_1+\cdots+2n_s}_T(\o X_{R}) \times \CY^{2n=2n_1+\cdots+2n_s}_T(\o X_{L}) \;,\]
\[\CY^{2n=2_1+\cdots+2_n}_T(\o X_{R}) \times \CY^{2n=2_1+\cdots+2_n}_T(\o X_{L}) \;, \]
\[ \text{and} \quad \CY^{2n_R}_T(\o X_{R}) \times \CY^{2n_L}_T(\o X_{L}) \;, \]
be respectively the sums of the products, right by left, of the equivalence classes of the partially reducible, orthogonal completely reducible and nonorthogonal completely reducible $2n$-dimensional toroidal compactified cycles.\vskip 11pt

\item Then, we have the following {\bf{Langlands global reducible correspondences\/}}:
\Ei
\[\begin{array}[t]{rcl}
\Red\Rep^{(2n)}_{W_{F_{R\times L}}}(W^{2n=2n_1+\cdots+2n_s}_{F_{\o \omega }}\times W^{2n=2n_1+\cdots+2n_s}_{F_{\omega }})
& \to&
\Red\cusp(\GL_{n=n_1+\cdots+n_s} ({F^{T}_{\o \omega }}\times {F^{T}_{\omega }}))\\
\searrow \qquad\qquad  && \qquad\qquad  \nearrow\end{array}\]
 
\[ \CY ^{2n=2n_1+\cdots+2n_s}_T(\o X_{R}) \times 
\CY ^{2n=2n_1+\cdots+2n_s}_T(\o X_{L}) \]

\[\begin{array}[t]{rcl}
\Red\Rep^{(2n)}_{W_{F_{R\times L}}}(W^{2n=2_1+\cdots+2_n}_{F_{\o \omega }}\times W^{2n=2_1+\cdots+2_n}_{F_{\omega }})
& \to&
\Red\cusp(\GL_{2n=2_1+\cdots+2_n} ({F^{T}_{\o \omega }}\times {F^{T}_{\omega }}))\\
\searrow \qquad\qquad  && \qquad\qquad  \nearrow\end{array}\]
\[ \CY ^{2n=2_1+\cdots+2_n}_T(\o X_{R}) \times 
\CY ^{2n=2_1+\cdots+2_n}_T(\o X_{L})  \]

\[\begin{array}[t]{rcl}
\Red\Rep^{(2n)}_{W_{F_{R\times L}}}(W^{2n_R}_{F_{\o \omega }}\times W^{2n_L}_{F_{\omega  }})
& \to&
\Red\cusp(\GL_{2n_R\times 2n_L} ({F^{T}_{\o \omega }}\times {F^{T}_{\omega }}))\\
\searrow \qquad\qquad  && \qquad\qquad  \nearrow\end{array}\]
\[ \CY ^{2n_R}_T(\o X_{R}) \times 
\CY ^{2n_L}_T( \o X_{L}) \]
}\vskip 11pt

\bpr this proposition is an adaptation of proposition 3.4.14 to the reducible case by noticing that:
\Bi
\item $\CY ^{2n=2n_1+\cdots+2n_s}_T(\o X_{R}) \times 
\CY ^{2n=2n_1+\cdots+2n_s}_T(\o X_{L}) =\widehat G^{2n=2n_1+\cdots+2n_s)}
(F^T_{\o\omega _\oplus }\times F^T_{\omega _\oplus })=\linebreak \bigoplus_{n_\ell}(\Eis_R(2n_\ell,j_R,m_{j_R})\otimes \Eis_L((2n_\ell,j_L,m_{j_L}))$ according to proposition 4.1.7.

\item $\Red\Rep^{(2n)}_{W_{F_{R\times L}}}(W^{2n=2n_1+\cdots+2n_s}_{F_{\o \omega }}\times W^{2n=2n_1+\cdots+2n_s}_{F_{\omega }})=\widehat G^{(n=n_1+\cdots+n_s)}(F_{\o\omega _\oplus }\times F_{\omega _\oplus })$ is isomorphic to its toroidal equivalent $\widehat G^{2n=2n_1+\cdots+2n_s)}(F^T_{\o\omega _\oplus }\times F^T_{\omega \oplus })= \bigoplus_{n_\ell}(\Eis_R(2n_\ell,j_R,m_{j_R})\times \Eis_L(2n_\ell,j_L,m_{j_L}))$~.
%
\epr
\Ei\vskip 11pt

\subsubsection{Definition} 

The real boundaries of the reducible compactified bisemispaces 
\linebreak $\o X^{2n=2n_1+\cdots+2n_s} _{\SRL} $~, \quad $\o X^{2n=2_1+\cdots+2_n}_{\SRL}$ \quad and \quad
$\o X^{2n_R\times 2n_L}_{\SRL}$
result from the surjective morphisms:
\begin{align*}
&\gamma^{\delta (2n)}_{\RL} : \quad \o X^{2n=2n_1+\cdots+2n_s}_{\SRL}\\
&\quad \To\partial \o X^{2n=2n_1+\cdots+2n_s}_{\SRL}=\GL_{n=n_1+\cdots+n_s}({F^{+,T}_{R}}\times {F^{+,T}_L})\Big/ \GL_{n=n_1+\cdots+n_s}( (\ZZ\big/N\ \ZZ)^2 )\;, \\
&\gamma^{\delta (2n)}_{R\times L} : \quad \o X^{2n=2_1+\cdots+2_n}_{\SRL}\\
&\quad \To \partial \o X^{2n=2_1+\cdots+2_n}_{\SRL}=\GL_{2n=2_1+\cdots+2_n}({F^{+,T}_{R}}\times {F^{+,T}_L})\Big/ \GL_{2n=2_1+\cdots+2_n}((\ZZ\big/N\ \ZZ)^2)\;, \\
&\gamma^{\delta (2n_R\times 2n_L)}_{R\times L} : \quad \o X^{2n_R\times 2n_L}_{\SRL}\\
&\quad \To\partial \o X^{2n_R\times 2n_L}_{\SRL}=\GL_{2n_R\times 2n_L}({F^{+,T}_{R}}\times 
{F^{+,T}_L})\Big/ \GL_{2n_R\times 2n_L}( (\ZZ\big/N\ \ZZ)^2)\;, \end{align*}
as developed in section 3.5.1.
\vskip 11pt

\subsubsection{Definition: the reducible Shimura bisemivarieties} 

The double coset decompositions of
the reducible bilinear general semigroups over $({F^{+,T}_{R}}\times 
{F^{+,T}_L})$ are given by:
\begin{align*}
&\partial \o S^{P_{n=n_1+\cdots+n_s}}_{K_{n=n_1+\cdots+n_s}} \\
& \qquad =P_{n=n_1+\cdots+n_s}( {F^{+,T}_{\o v^1}}\times {F^{+,T}_{v^1}} )\setminus \GL_{n=n_1+\cdots+n_s}
({F^{+,T}_{R}}\times {F^{+,T}_L})\\
& \hspace{8cm}\Big/ \GL_{n=n_1+\cdots+n_s}( (\ZZ\big/N\ \ZZ)^2)\;, \\ 
&\partial \o S^{P_{2n=2_1+\cdots+2_n}}_{K_{2n=2_1+\cdots+2_n}} \\
& \qquad =P_{2n=2_1+\cdots+2_n}({F^{+,T}_{\o v^1}}\times {F^{+,T}_{v^1}} )\setminus \GL_{2n=2_1+\cdots+2_n}
({F^{+,T}_{R}}\times {F^{+,T}_L})\\
& \hspace{8cm}\Big/ \GL_{2n=2_1+\cdots+2_n}( (\ZZ\big/N\ \ZZ)^2)\;, \\
&\partial \o S^{P_{2n_R\times 2n_L}}_{K_{2n_R\times 2n_L}} \\
& \qquad =P_{2n_R\times 2n_L}( {F^{+,T}_{\o v^1}}\times {F^{+,T}_{v^1}} )\setminus \GL_{2n_R\times 2n_L}
({F^{+,T}_{R}}\times {F^{+,T}_L }  )\\
&\hspace{8cm}\Big/ \GL_{2n_R\times 2n_L}((\ZZ\big/N\ \ZZ)^2)\;; \end{align*}
they correspond to reducible Shimura bisemivarieties according to definition 3.5.2.
\vskip 11pt

\subsubsection{Proposition}  

{\em \parindent=0pt
The following Eisenstein cohomologies develop according to:
\Bi
\item $H^{*}( \partial  \o S 
^{P_{n=n_1+\cdots+n_s}} _{K_{n=n_1+\cdots+n_s}} , \widehat M^{2n}_{T_{v_{R_\oplus }}}\otimes \widehat M^{2n}_{T_{v_{L_\oplus }}}) = \bigoplus\limits_{n_\ell} H^{2n_\ell}( \partial  \o S 
^{P_{n}} _{K_{n}} ,\widehat M^{2n_\ell}_{T_{v_{R_\oplus }}}\otimes \widehat M^{2n_\ell}_{T_{v_{L_\oplus }}})$

\hspace{1cm} $\simeq 
\bigoplus\limits_{n_\ell} \bigoplus\limits_{j_\delta}\bigoplus\limits_{m_{j_\delta}}
(\CY^{n_\ell}_T(\partial \o X_{R}[j_\delta,m_{j_\delta}])\times \CY_T^{n_\ell}(\partial \o X_{L}[j_\delta,m_{j_\delta}]))$\vskip 11pt

\item $H^{2n}( \partial  \o S 
^{P_{2n=2_1+\cdots+2_n}} _{K_{2n=2_1+\cdots+2_n}} ,\widehat M^{2n}_{T_{v_{R_\oplus }}}\otimes \widehat M^{2n}_{T_{v_{L_\oplus }}})  = \bigoplus\limits_{\ell_R=\ell_L} H^{2_{\ell}}
( \partial  \o S 
^{P_{2{n}}}_{K_{2{n}}} ,\widehat M^{2_{\ell_R}}_{T_{v_{R_\oplus }}}\otimes 
\widehat M^{2_{\ell_L}}_{T_{v_{L_\oplus }}})$

\hspace{1cm} $\simeq 
\bigoplus\limits_{\ell_R=\ell_L} \bigoplus\limits_{j_\delta}\bigoplus\limits_{m_{j_\delta}}
(\CY^{1_{\ell_R}}_T(\partial \o X_{R}[j_\delta,m_{j_\delta}])\times \CY_T^{1_{\ell_L}}(\partial \o X_{L}[j_\delta,m_{j_\delta}]))$\vskip 11pt

\item $H^{2n}( \partial  \o S 
^{P_{2n_R\times 2n_L}} _{K_{2n_R\times 2n_L}} ,\widehat M^{2n_R}_{T_{v_{R_\oplus }}}\otimes \widehat M^{2n_L}_{T_{v_{L_\oplus }}})  = \bigoplus\limits_{\ell_R=\ell_L} H^{2_{\ell}}
( \partial  \o S 
^{P_{2{n_R},2 {n_L}}}_{K_{2{n_R},2{n_L}}} ,\widehat M^{2_{\ell_R}}_{T_{v_{R_\oplus }}}\otimes \widehat M^{2_{\ell_L}}_{T_{v_{L_\oplus }}})$

\hspace{2cm} $ \bigoplus\limits_{k_R\neq\ell_L} H^{2_{k_R},2_{\ell_L}}
( \partial  \o S 
^{P_{2n_R\times 2n_L}}_{K_{2n_R\times 2n_L}} ,\widehat M^{2_{k_R}}_{T_{v_{R_\oplus }}}\otimes \widehat M^{2_{\ell_L}}_{T_{v_{L_\oplus }}})$

\hspace{1cm} $\simeq 
\bigoplus\limits_{\ell_R=\ell_L} \bigoplus\limits_{j_\delta}\bigoplus\limits_{m_{j_\delta}}
(\CY^{1_{\ell_R}}_T(\partial \o X_{R}[j_\delta,m_{j _\delta }])\times \CY_T^{1_{\ell_L}}(\partial \o X_{L}[j_\delta,m_{j_\delta}]))$

\hspace{2cm} $ 
\bigoplus\limits_{k_R\neq\ell_L} \bigoplus\limits_{j^\delta _R}\bigoplus\limits_{m_{j^\delta _R}}
\bigoplus\limits_{j^\delta _L}\bigoplus\limits_{m_{j^\delta _L}}
(\CY^{1_{k_R}}_T(\partial \o X_{R}[j^\delta _R,m_{j^\delta _R}])\times \CY_T^{1_{\ell_L}}(\partial \o X_{L}[j^\delta _L,m_{j^\delta_L }]))$
\Ei
they are equal to the sums of products, right by left, of the equivalence classes, counted with\linebreak their multiplicities, of the irreducible cycles.}
\vskip 11pt

\bpr these decompositions of cohomologies are the real equivalents of those considered in proposition 4.2.4.\epr\vskip 11pt

\subsubsection{Proposition} 

{\em  
\Bean
\item If $ x_{n_\ell}=\sum\limits^{n_\ell}_{\beta =1} x_{n_\ell}(\beta )\ |\vec e_\beta| \in\RR^{n_\ell}$~,

 then, every left (resp. right) $n_\ell$-dimensional
real semitorus has the analytic development:
\begin{align*}
T^{n_\ell}_L[j^\delta_L,m_{j^\delta_L}] &\simeq
\lambda^{\half} (n_\ell,j^\delta_L,m_{j^\delta_L}) \ e^{2\pi ij^\delta_Lx_{n_\ell}}\\
(\text{resp.} \quad 
T^{n_\ell}_R[j^\delta_R,m_{j^\delta_R}] &\simeq
\lambda^{\half}(n_\ell,j^\delta_R,m_{j^\delta_R})\ e^{-2\pi ij^\delta_Rx_{n_\ell}}\ ).\end{align*}

\item On the other hand, let
\begin{align*}
\ELLIP_L (2n_\ell,j^\delta_L,m_{j^\delta_L}) &=
\txt\sum\limits^t_{j^\delta_L=1} \sum\limits_{m_{j^\delta_L}} \lambda^{\half} (n_\ell,j^\delta_L,m_{j^\delta_L}) \ e^{2\pi ij^\delta_Lx_{n_\ell}}\\
(\text{resp.} \quad 
\ELLIP_R (2n_\ell,j^\delta_R,m_{j^\delta_R}) &=
\txt\sum\limits^t_{j^\delta_R=1} \sum\limits_{m_{j^\delta_R}} \lambda^{\half} (n_\ell,j^\delta_R,m_{j^\delta_R}) \ e^{-2\pi ij^\delta_Rx_{n_\ell}})\ ), \quad t\le \infty\;,\end{align*}
be the (truncated) Fourier development of a left (resp. right) $n_\ell$-dimensional solvable global elliptic semimodule where
\Bi
\item[\textbullet] $\lambda (n_\ell,j_\delta,m_{j_\delta}) =\prod\limits^{n_\ell}_{c=1}
\lambda_c (n_\ell,j_\delta,m_{j_\delta}) \approx j^{n_\ell}\cdot N^{n_{\ell}}$~;
\vskip 11pt

\item[\textbullet] $e^{2\pi ix_{n_\ell}}\to e^{2\pi ij^\delta_Lx_{n_\ell}}$ (resp. $e^{-2\pi ix_{n_\ell}}\to e^{-2\pi ij^\delta_Rx_{n_\ell}}$~)
 is a left (resp. right) global Frobenius substitution.
\Ei

Then, we have the following analytic developments of the Eisenstein bilinear cohomologies of the equivalents of the reducible Shimura bisemivarieties:
\Bi
\item[\textbullet] $H^{*}( \partial  \o S 
^{P_{n=n_1+\cdots+n_s}} _{K_{n=n_1+\cdots+n_s}} ,
\widehat M^{2n}_{T_{v_{R_\oplus }}}\otimes \widehat M^{2n}_{T_{v_{L_\oplus }}}) 
$

 \hspace{1cm} $=\bigoplus\limits_{2n_\ell} (\ELLIP_R (2n_\ell,j^\delta_R,m_{j^\delta_R})\otimes \ELLIP_L 
(2n_\ell,j^\delta_L,m_{j^\delta_L}))$\vskip 11pt

\item[\textbullet] $H^{2n}( \partial  \o S 
^{P_{2n=2_1+\cdots+2_n}} _{K_{2n=2_1+\cdots+2_n}} ,\widehat M^{2n}_{T_{v_{R_\oplus }}}\otimes \widehat M^{2n}_{T_{v_{L_\oplus }}}) $

 \hspace{1cm} $=\bigoplus\limits_{\ell_R=\ell_L}   
(\ELLIP_R (2_{\ell_R},j^\delta_R,m_{j^\delta_R})\otimes \ELLIP_L (2_{\ell_L},j^\delta_L,m_{j^\delta_L}))$\vskip 11pt

\item[\textbullet] \mbox{$H^{2n}( \partial  \o S 
^{P_{2n_R\times 2n_L}} _{K_{2n_R\times 2n_L}} ,\widehat M^{2n_R}_{T_{v_{R_\oplus }}}\otimes \widehat M^{2n_L}_{T_{v_{L_\oplus }}}) $}

 \hspace{1cm} $=\bigoplus\limits_{\ell_R=\ell_L}   
(\ELLIP_R (2_{\ell_R},j^\delta_R,m_{j^\delta_R})\otimes \ELLIP_L (2_{\ell_L},j^\delta_L,m_{j^\delta_L}))$

\hspace{4cm} $\bigoplus\limits_{k_R\neq\ell_L}   
(\ELLIP_R (2_{k_R},j^\delta_R,m_{j^\delta_R})\otimes \ELLIP_L (2_{\ell_L},j^\delta_L,m_{j^\delta_L}))$
\Ei
\Ee}
\vskip 11pt

\bpr referring to proposition 3.5.4 and 3.5.5, each equivalence class representative\linebreak $\CY^{n_\ell}_T (\partial \o X_{R}[j_\delta ,m_{j_\delta }])$ (resp. $\CY^{n_\ell}_T (\partial \o X_{L}[j_\delta ,m_{j_\delta }])$~) of a semicycle of codimension $n_\ell$ corresponds to the analytic representation of a $n_\ell$-dimensional real semitorus
$T^{n_\ell}_R  [j_\delta ,m_{j_\delta }]$ (resp. $T^{n_\ell}_L  [j_\delta ,m_{j_\delta }]$~).\epr
\vskip 11pt

\subsubsection{Remark}

Trace formulas for the real reducible cases considered here can be developed similarly as it was done for the complex cases in sections
4.2.6 and 4.2.7.

\subsubsection{Proposition: Langlands global correspondences on the equivalents of the reducible Shimura bisemivarieties}

{\em \parindent=0pt
\Bi

\item Let
\begin{align*}
&\Red  \Rep^{(n)}_{W_{F^+_{R\times L}}} 
(W^{n=n_1+\cdots+n_s}_{F^{+}_{\o v}} \times W^{n=n_1+\cdots+n_s}_{F^{+}_{v}})\\
& \qquad \qquad \qquad \To H^*(\partial \o S^{P_{n=n_1+\cdots+n_s}}_{K_{n=n_1+\cdots+n_s}},
 \widehat M^{2n}_{T_{v_{R_\oplus }}}\otimes \widehat M^{2n}_{T_{v_{L_\oplus }}}) \\
& \hspace*{5cm}
=\FReps(\GL_{n=n_1+\cdots+n_s}({F^{+,T}_{\o v_\oplus }} \times {F^{+,T}_{v_\oplus }}))\;, \\[11pt]
&\Red  \Rep^{(n)}_{W_{F^+_{R\times L}}} 
(W^{2n=2_1+\cdots+2_n}_{F^{+}_{\o v}} \times W^{2n=2_1+\cdots+2_n}_{F^{+}_{v}})\\
& \qquad \qquad \qquad \To H^{2n}(\partial \o S^{P_{2n=2_1+\cdots+2_n}}_{K_{2n=2_1+\cdots+2_n}},\widehat M^{2n^*}_{T_{v_{R_\oplus }}}\otimes \widehat M^{2n^*}_{T_{v_{L_\oplus }}})  \\
& \hspace*{5cm}
=\FReps(\GL_{2n=2_1+\cdots+2_n}({F^{+,T}_{\o v_\oplus }} \times {F^{+,T}_{v_\oplus }}))\;,   \\[11pt]
%
&\Red  \Rep^{(n)}_{W_{F^{+}_{R\times L}}} 
(W^{2n_R}_{F^{+}_{\o v}} \times W^{2n_L}_{F^{+}_{v}})\\
& \qquad \qquad \qquad \To H^{2n}(\partial \o S^{P_{2n_R\times 2n_L}}_{K_{2n_R\times 2n_L}},\widehat M^{2n_R}_{T_{v_{R_\oplus }}}
\otimes \widehat M^{2n_L}_{T_{v_{L_\oplus }}}) \\
& \hspace*{5cm}
=\FReps(\GL_{2n_R\times 2n_L}({F^{+,T}_{\o v_\oplus }} \times {F^{+,T}_{v_\oplus }}))\;,   
\end{align*}
be respectively the sums of the products, right by left, of the equivalence classes of the partially reducible, orthogonal completely reducible and nonorthogonal completely reducible  representations of the bilinear global Weil groups \quad 
$(W^{n=n_1+\cdots+n_s}_{F^{+}_{\o v}} \times W^{n=n_1+\cdots+n_s}_{F^{+}_{v}})$~, \quad 
$(W^{2n=2_1+\cdots+2_n}_{F^{+}_{\o v}} \times W^{2n=2_1+\cdots+2_n}_{F^{+}_{v}})$ \quad 
and \quad 
$(W^{2n_R\times 2n_L}_{F^+_{\o v}} \times W^{2n_R\times 2n_L}_{F^+_{v}})$~.

\item Let
\begin{align*}
&\Red  \ELLIP (\GL_{n=n_1+\cdots+n_s} ({F^{+,T}_{\o v}} \times {F^{+,T}_{v}})\\
& \qquad \qquad =\txt\bigoplus\limits_{2n_\ell }( \ELLIP_R(2n_\ell,j^\delta _R,m_{j^\delta _R}) \otimes \ELLIP_L(n_\ell,j^\delta _L,m_{j^\delta _L}) )\;,\\\noalign{\vskip 11pt}
&\Red  \ELLIP (\GL_{2n=2_1+\cdots+2_n} ({F^{+,T}_{\o v}} \times {F^{+,T}_{v}})\\
& \qquad \qquad =\txt\bigoplus\limits_{\ell_R=\ell_L }( \ELLIP_R(2_{\ell_R},j^\delta _R,m_{j_R}) \otimes \ELLIP_L(2_{\ell_L},j^\delta _L,m_{j^\delta _L}) )\;,\\\noalign{\vskip 11pt}
& \Red  \ELLIP (\GL_{2n_R\times 2n_L} ({F^{+,T}_{\o v}} \times {F^{+,T}_{v}})\\
& \qquad \qquad =\txt\bigoplus\limits_{\ell_R=\ell_L }( \ELLIP_R(2_{\ell_R},j^\delta _R,m_{j^\delta _R}) \otimes \ELLIP_L(2_{\ell_L},j^\delta _L,m_{j^\delta _L}) )\\
& \qquad \qquad \qquad \txt\bigoplus\limits_{k_R\neq\ell_L }( \ELLIP_R(2_{k_R},j^\delta _R,m_{j^\delta _R}) \otimes \ELLIP_L(2_{\ell_L},j^\delta _L,m_{j^\delta _L}) )\;,\end{align*}
be respectively the sums of the products, right by left, of the equivalence classes of the partially reducible, orthogonal completely reducible and nonorthogonal completely reducible global elliptic representations of the corresponding reducible bilinear semigroups.

\item Let finally
\begin{gather*}
\CY^{n=n_1+\cdots+n_s}_T(\partial \o X_{R}) \times \CY^{n=n_1+\cdots+n_s}_T(\partial \o X_{L}) \;, \\[11pt]
\CY^{n=1_1+\cdots+1_n}_T(\partial \o X_{R}) \times \CY^{n=1_1+\cdots+1_n}_T(\partial \o X_{L}) \;, \\[11pt]
\text{and} \quad\CY^{n_R}_T(\partial \o X_{R}) \times \CY^{n_L}_T(\partial \o X_{L}) \;, \end{gather*}
be respectively the sums of the products, right by left, of the equivalence classes of the partially reducible, orthogonal completely reducible and nonorthogonal completely reducible $n$-dimensional toroidal compactified cycles over $\partial \o X_{R}$ and $\partial \o X_{L}$~.

\item Then, we have the following {\bfseries{Langlands global reducible correspondences\/}} on the equivalents of the reducible Shimura bisemivarieties:
\Ei
\[\begin{array}[t]{rcl}
\Red\Rep^{(n)}_{W_{F^{+}_{R\times L}}}
(W^{n=n_1+\cdots+n_s}_{F^{+}_{\o v}}\times W^{n=n_1+\cdots+n_s}_{F^{+}_{v}})
& \to&
\Red\ELLIP(\GL_{n=n_1+\cdots+n_s} ({F^{+,T}_{\o v}}\times {F^{+,T}_{v}}))\\
\searrow \qquad\qquad  && \qquad\qquad  \nearrow\end{array}\]
\[ \CY ^{n=n_1+\cdots+n_s}_T(\partial  \o X_{R}) \times 
\CY ^{n=n_1+\cdots+n_s}_T(\partial \o X_{L}) \]
\vskip 11pt

\[\begin{array}[t]{rcl}
\Red\Rep^{(n)}_{W_{F^{+}_{R\times L}}}
(W^{2n=2_1+\cdots+2_n}_{F^{+}_{\o v}}\times W^{2n=2_1+\cdots+2_n}_{F^{+}_{v}})
& \to&
\Red\ELLIP(\GL_{2n=2_1+\cdots+2_n} ({F^{+,T}_{\o v}}\times {F^{+,T}_{v}}))\\
\searrow \qquad\qquad  && \qquad\qquad  \nearrow\end{array}\]
\[ \CY ^{n=1_1+\cdots+1_n}_T(\partial  \o X_{R}) \times 
\CY ^{n=1_1+\cdots+1_n}_T(\partial  \o X_{L})  \]
\vskip 11pt

\[\begin{array}[t]{rcl}
\Red\Rep^{(n)}_{W_{F^{+}_{R\times L}}}
(W^{2n_R}_{F^{+}_{\o v}}\times W^{2n_L}_{F^{+}_{v  }})
& \to&
\Red\ELLIP(\GL_{2n_R\times 2n_L} ({F^{+,T}_{\o v}}\times {F^{+,T}_{v}}))\\
\searrow \qquad\qquad  && \qquad\qquad  \nearrow\end{array}\]
\[ \CY ^{n_R}_T(\partial  \o X_{R}) \times 
\CY ^{n_L}_T( \partial \o X_{L}) \]
}\vskip 11pt

\bpr this is an adaptation to the real case of proposition 4.2.8 (see also proposition 3.4.10).\epr

\section*{Appendix}
\addcontentsline{toc}{section}{Appendix}
\setcounter{subsection}{0}
{\def\thesubsection{\arabic{subsection}}

\subsection{Philosophy of bisemiobjects}

Instead of working with classical mathematical objects, the developments of this paper deal with the products of semiobjects.  In this respect, a mathematical symmetric object ``~$O$~'' will be cut into two semiobjects, labeled left and right, in such a way that the left semiobject $O_L$ be localized in (or refers to) the upper half space while the right symmetric semiobject $O_R$ is localized in (or refers to) the lower half space.  Then, informations concerning the internal mathematical structure of our object ``~$O$~'' can be more advantageously extracted from the product $O_R\times O_L$ of the two semiobjects $O_R$ and $O_L$~.  Indeed, the endomorphism of our object ``~$O$~'' can be decomposed into the product $E_R\times E_L$ of the endomorphisms $E_L:O_L\to O_L$ acting on the semiobject $O_L$ by the opposite endomorphism
$E_R:O_R\to O_R$ acting on $O_R$ such that $E_R=E_L^{-1}$~.  By this way, the action of $E_L$ is neutralized by the coaction of $E_R$~.

A classical analogue of this endormophism $E_R\times E_L$ on $O_R\times O_L$ is, for example, the normal endomorphism $E_N$ of a group $G$ given by:
\[a\ E_N(b)\ a^{-1}=E_N(a\ b\ a^{-1})\quad \text{for\ } a,b\in G\;.\]
(In the philosophy of bisemiobjects, the element $b$ would have to be decomposed into the product $b_R\times b_L$ of two semielements $b_R$ and $b_L$~, if mathematically feasible.)
\vskip 11pt 

\subsection{Semistructures}

So, semistructures will be recalled or introduced in the following.  The notation $L,R$ means  ``left'' or ``right''.

\Bi
\item We are concerned with a {\bf semigroup\/} $G\Lr$~, called left of right, which is a nonempty set of \lr elements, localized or referring to the upper (resp. lower) half space, together with a binary operation on $G\Lr$~, i.e. a function $G\Lr\times G\Lr\to G\Lr$~, or $G\Lr\times G_{L,L}\to G\Lr$~.

\item A \lr {\bf monoid\/} is a \lr semigroup $G\Lr$ which contains an identity element $e\Lr\in G\Lr$ such that:
\[ e_L\cdot a_L=a_L \quad \text{(resp.\ } a_R\cdot e_R=a_R\ ), \qquad \forall\ a\Lr\in G\Lr\;.\]

\item A \lr {\bf semiring\/} is a nonempty set $R\Lr$ together with two binary operations (addition and multiplication) such that:
\Bean
\item $(R\Lr,+)$ is an abelian \lr semigroup;
\item $(a\Lr\ b\Lr)\ c\Lr=a\Lr\ (b\Lr\ c\Lr) \quad \forall\ a\Lr,b\Lr,c\Lr\in R\Lr$  (associative multiplication);
\item $a\Lr\ (b\Lr+c\Lr)=a\Lr\ b\Lr+a\Lr\ c\Lr$ and $(a\Lr+b\Lr)\ c\Lr=a\Lr\ c\Lr + b\Lr\ c\Lr$ (left and right distribution).
\Ee

\item If $R\Lr$ is a commutative semiring with identity $\ung_{R\Lr}$ and no zero divisors, it will be called a \lr integral semidomain. 

 Furthermore, if every element of $R\Lr$ is a unit (i.e. left and right invertible), $R\Lr$ is a {\bf division semiring\/}.

And, a \lr {\bf semifield\/} is a commutative division semiring.

\item A {\bf \lr adele semiring\/} is the product of the completions of the \lr considered semifield.

\item Let $R\Lr$ be a \lr semiring.  A \lr {\bf\boldmath $R\Lr$-semimodule\/} is an additive abelian \lr semigroup $M\Lr$ together with a function $R\Lr\times M_L\to M_L$ (resp. $M_R\times R\Rl\to M_R$~) such that:
\Bean
\item $\begin{array}[t]{lll}
 & r_L\ (a_L+b_L)=r_L\ a_L+r_L\ b_L  \;, & \forall\ r_L\in R_L\ , \; a_L,b_L\in M_L\\
\text{(resp.}
 & (a_R+b_R)\ r_R= a_R\ r_R+b_R\ r_R)  \;, & \forall\ r_R\in R_R\ , \; a_R,b_R\in M_R\ );\end{array}$
\item $\begin{array}[t]{lll}
 & (r_L+s_L)\ a_L=r_L\ a_L+s_L\ a_L \;, & \forall\ s_L\in R_L\\
\text{(resp.}
 & a_R\ (r_R+s_R)=a_R\ r_R+a_R\ s_R \;, & \forall\ s_R\in R_R\ );\end{array}$
\item $\begin{array}[t]{ll}
 & r_L\ (s_L\ a_L)=(r_L\ s_L) \ a_L \\
\text{(resp.}
 & (a_R\ s_R)\ r_R= a_R\ (s_R\ r_R) \ ). \end{array}$
\Ee

If $R\Lr$ has an identity element $\ung\Lr$ such that $\ung_L\cdot a_L=a_L$ (resp. $a_R\cdot\ung_R=a_R$~), $M\Lr$ is a \lr unitary $R\Lr$-semimodule.

If $R\Lr$ is a \lr division semiring, then the unitary \lr unitary $R\Lr$-semimodule, is a \lr {\bf vector semispace\/}.

\item If $R\Lr$ is a commutative semiring with identity, a {\bf\boldmath $R\Lr$-semialgebra $A\Lr$} is a semiring $A\Lr$ such that:
\Bean
\item $(A\Lr,+)$ is a unitary \lr $R\Lr$-semimodule;
\item $\begin{array}[t]{lll}
 & r_L\ (a_L\ b_L)=(r_L\ a_L) \ b_L =a_L\ (r_L\ b_L)  \;, & \forall\ r_L\in R_L\ , \; a_L,b_L\in A_L\\
\text{(resp.}
 & (a_R\ b_R)\ r_R=a_R\ (b_R\ r_R)=b_R\ (a_R\ r_R)  \;, & \forall\ r_R\in R_R\ , \; a_R,b_R\in A_R
\ ). \end{array}$
\Ee

If $A\Lr$ is a division semiring, then $A\Lr$ is called a division semialgebra.
\Ei
\vskip 11pt 

\subsection{Bisemistructures \cite{Pie5}}

Before recalling the salient features of the algebraic bilinear semigroups, it is useful to give a precise definition of a $R_R-R_L$-bisemimodule $M_{R-L}$ and  of a $R_R\times R_L$-bisemimodule $M_{R}\otimes M_L$.

\Bi

\item Let $R_L$ (resp. $R_R$~) be a \lr semiring.

Then, a {\bbf $R_R-R_L$-bisemimodule\/} is an abelian (semi)group $M_{R-L}$ (under the addition) with the structure of both a left $R_L$-semimodule $M_L$ and a right $R_R$-semimodule $M_R$ and together with a bifunction $R_{L,R}\times M_{L-R}\times R\Rl\to M_{R-L}$ such that:
\Bean
\item $M_{L-R}\equiv M_{R-L}=M_R\oplus M_L$~.

\item $(r_L\ a_L)\ (a_R\ r_R)=r_L\ (a_L\ a_R\ r_R)=(r_L\ a_L\ a_R)\ r_R$~, $\forall\ a_L\in M_L$~, $a_R\in M_R$~, $r_L\in r_L$~, $r_R\in R_R$~.
\item the conditions, mentioned above for a \lr $R\Lr$-semimodule, are applied on the left and on the right of $M_{L-R}$~.
\Ee

\item {\bbf A $R_R\times R_L$-bisemimodule $M_R\otimes M_L$\/} is a  bilinear semigroup given by the map : $f:M_R \times M_L\to M_R\otimes M_L$~, where $f$ sends pairs $(a_{R_i},a_{L_i})\in M_R\times M_L$ of symmetric elements to their products $(a_{R_i}\times a_{L_i})$ such that:
\Bean
\item $f(a_{R_1}+a_{R_2},a_{L_1})=f(a_{R_1},a_{L_1})+f(a_{R_2},a_{L_1})$~;
\item $f(a_{R_1},a_{L_1}+a_{L_2})=f(a_{R_1},a_{L_1})+f(a_{L_1},a_{L_2})$~;
\item $f(a_R\ r_R,r_L\ a_L) \neq f(r_L\ a_R,a_L\ r_L)$~, 
$\forall\ a_{R_i}\in M_R$~, $i=1,2$~, $a_{L_i}\in M_L$~.
\Ee

 
 \item Taking into account the section 1 of the appendix and the introduction, we have to specify {\bbf the product $X_R\times X_L$ of a right affine semigroup scheme $X_R$ by its left symmetric equivalent $X_L$\/} as occurring from:
 \[ (X_R  \circ \phi _R^{'-1}) \times (X_L \circ \phi _L^{'-1})](\widetilde F_R\times \widetilde F_L)
 = \{(v_R\times v_L)\mid v_R\in V_R\ ,\; v_L\in V_L\}\]
 where $(v_R\times v_L)$ is a bipoint of an affine bisemispace $(V_R\otimes V_L)$ resulting from the bihomomorphism
 \[ X_R\times X_L:\quad Q_R\times Q_L\To V_R\otimes V_L\]
 which sends the product of the quotient $k$-algebras $Q_R$ and $Q_L$ into $(V_R\otimes V_L)$~.


\item It was seen in chapter 2 that $\GL_n(\wt F_{\o\omega }\times \wt F_\omega )$ is an {\bf algebraic bilinear semigroup\/} of matrices in the sense that its algebraic representation space is an $\GL_n(\wt F_{\o\omega }\times \wt F_\omega )$-bisemimodule $\wt M_R\otimes \wt M_L$ defined as follows:

Let $\wt M_L$ (resp. $\wt M_R$~) be an $n$-dimensional $T_n(\wt F_\omega )$- (resp. $T^t_n(\wt F_{\o\omega })$~)-semimodule.  Then, the canonical bilinear map $\wt M_R\times \wt M_L\to \wt M_R\otimes \wt M_L$ defines the bisemimodule $\wt M_R\otimes \wt M_L$~, provided that $\wt M_R\otimes \wt M_L$ is the bilinear tensor product of a left $T_n(\wt F_\omega )$-semimodule by a right $T^t_n(\wt F_{\o\omega })$-semimodule.

\item The {\bf complete algebraic bilinear  semigroup\/} $G^{(n)}(\wt F_{\o\omega }\times \wt F_\omega )$   is generated by the bilinear algebraic 
semigroup of matrices $\GL_n (\wt F_{\o\omega }\times \wt F_\omega )   \equiv T^t_n(\wt F_{\o\omega })\times \wt T_n(F_\omega )$ where $T_n(\wt F_\omega )$ (resp. $T^t_n(\wt F_{\o\omega })$~) is the group of upper (resp. lower) triangular matrices.

\item The complete algebraic bilinear  semigroup $G^{(n)}(F_{\o\omega }\times F_\omega )$ can be viewed as {\bf a bilinear Lie semigroup\/} defined as follows:
\Be
\item {\bbf $G^{(n)}(F_{\o\omega }\times F_\omega )$ is a topological bilinear semigroup\/}, that is to say that:
\Be
\item it is a {\bf bilinear semigroup\/} $G\RL$~:
\Be
\item[$\alpha$)] given by a function:
\begin{align*}
G_R\times G_L &\To G\RL\\
(g_{r_i},g_{L_i}) &\To (g_{r_i}\times g_{L_i})\end{align*}
sending pairs $(g_{r_i},g_{L_i})$ of symmetric elements, localized in (or referring to) the upper and lower half spaces, to their products $(g_{r_i}\times g_{L_i})$~, $i\in\NN$~.

\item[$\beta )$] submitted to a binary operation $\u\times$
\[ G\RL\u\times G\RL\To G\RL\]
defined by
\[ (g_{r_i}\times g_{L_i})\u\times(g_{r_j}\times g_{L_j})
\To (g_{r_i}+g_{r_j})\times (g_{L_i}+g_{L_j})\]
in such a way that cross products $(g_{r_i}\times g_{L_j})$ could be generated.
\Ee

\item {\bbf $G^{(n)}(F_{\o\omega }\times F_\omega )$ is a topological bisemispace\/}, product of a right topological semispace, restricted to the lower half space by its left symmetric equivalent.
\Ee

\item {\bbf $G^{(n)}(F_{\o\omega }\times F_\omega )$ is an abstract bisemivariety\/}, i.e. the product of a left semivariety restricted to the upper half space by the symmetric right semivariety restricted to the lower half space.
\Ee


\item Let $ P^{(n)}  (F_{\o\omega^1 }\times F_{\omega^1} )  =P^{(n)}  (F_{\o\omega^1 })\times P^{(n)}(F_{\omega^1} )$    be a bilinear subgroup of  $G^{(n)} (F_{\o\omega }\times F_\omega )$   such that $F_{\omega^1 }$  (resp. $F_{\o\omega^1 }$~) denotes the set of \lr  irreducible completions (see section 2.4.1).  

Let $a\Rl,b\Rl\in G^{(n)} (F_{\o\omega }\times F_\omega ) $~.  

$a_R$ is  right congruent to $b_R$ modulo $P^{(n)} (F_{\o\omega^1 })$~, denoted $a_R\equiv b_R\mod(P^{(n)} (F_{\o\omega^1 })$ 
if $a_R\ b_R^{-1}\in P^{(n)} (F_{\o\omega^1 })$~.

$a_L$ is  left congruent to $b_L$ modulo $P^{(n)} (F_{\omega^1 })$~, denoted $a_L\equiv b_L\mod(P^{(n)} (F_{\omega^1 })$ if $a_L\ b_L^{-1}\in P^{(n)} (F_{\omega^1 })$~.

Then, $ P^{(n)} (F_{\o\omega^1 }  \times F_{\omega^1 })$ is a {\bf bilinear normal subgroup\/} of 
$ G^{(n)} (F_{\o\omega }  \times F_{\omega })$ if left congruence modulo $ P^{(n)} (F_{\omega^1 }) $ correspond to right congruence modulo $ P^{(n)} (F_{\o\omega^1 })$~.

\item The {\bf endomorphism\/} of an algebraic bilinear semigroup $ G^{(n)} (\wt F_{\o\omega }  \times \wt F_{\omega })$ is a homomorphism
\[ E\RL: \quad G^{(n)} (\wt F_{\o\omega }  \times \wt F_{\omega })\To G^{(n)}  (\wt F_{\o\omega }  \times \wt F_{\omega })\]
which can be decomposed into:
\begin{align*}
E_R : \quad G^{(n)} ( \wt F_{\o\omega } ) &\To G^{(n)}( \wt F_{\o\omega } )\; ,\\
E_L : \quad G^{(n)}( \wt F_{\omega } ) &\To G^{(n)}( \wt F_{\omega } )\; ,\end{align*}
for $G^{(n)}( \wt F_{\o\omega } \times \wt F_{\omega } )= G^{(n)}( \wt F_{\o\omega } ) \times G^{(n)}( \wt F_{\omega } )$ such that:
\Be
\item $E_R(a_R\ b_R)=E_R(a_R)\cdot E_R(b_R)$ for all $a_R,b_R\in G^{(n)}(\wt  F_{\o\omega } )$~;\vskip 11pt 

\item $E_L(a_L\ b_L)=E_L(a_L)\cdot E_L(b_L)$ for all $a_L,b_L\in G^{(n)}(\wt  F_{\omega } )$~;\vskip 11pt 

\item $E\RL(a_R\ b_R\ a_L\ b_L)=E\RL(a_R\cdot a_L)\cdot E\RL(b_R\cdot b_L)$~.
\Ee

If $E_R=E_L^{-1}$~, then $E_R(a_R\ b_R)=E_L^{-1}(a_L^{-1}\ b_L^{-1})$ where $a_L^{-1},b_L^{-1}$ denote the opposites with respect to the addition or the complex conjugates of $a_L,b_L$~.

So, \begin{align*}
E\RL(a_R\ b_R\ a_L\ a_L)&=E_L^{-1}\times E_L(a_L^{-1}\ b_L^{-1}\ a_L\ b_L)\\
\text{or}\quad
E\RL(b_R\ a_R\ a_L\ b_L)&= E_L^{-1}\times E_L(b_L^{-1}\ a_L^{-1}\ a_L\ b_L)\;.\end{align*}
\Ei
}

       }
\end{document}